\documentclass[reqno,10pt]{amsart}
\usepackage{amssymb,amsmath,amsthm,}
\usepackage{amscd}
\usepackage{amsfonts}
\usepackage{amssymb}
\usepackage{latexsym}
\usepackage{color}
\usepackage{esint}

\newtheorem{theorem}{Theorem}

\newtheorem{corollary}[theorem]{Corollary}
\newtheorem{definition}[theorem]{Definition}
\newtheorem{lemma}[theorem]{Lemma}
\newtheorem{proposition}[theorem]{Proposition}
\newtheorem{remark}[theorem]{Remark}
\let\a=\alpha

\let\p=\partial

\let\O=\Omega

\let\o=\omega

\newcommand{\be}{\begin{equation}}
\newcommand{\bm}{\begin{multline}}
\newcommand{\ee}{\end{equation}}
\newcommand{\dd}{\mathrm{d}}

\numberwithin{equation}{section}
\numberwithin{theorem}{section}

\def\p{\partial}

\def\O{\Omega}

\def\B{\begin{equation}}
\def\E{\end{equation}}
\def\BN{\begin{eqnarray*}}
\def\EN{\end{eqnarray*}}

\begin{document}


\title{Uniform estimate of viscous free-boundary magnetohydrodynamics with zero vacuum magnetic field}
\author{Donghyun Lee}
\maketitle

\begin{abstract}
We consider viscous free-boundary magnetohydrodynamics(MHD) under vacuum in $\mathbb{R}^3$, especially when vacuum magnetic field is identically zero. It is a central problem in mathematics to perform vanishing viscosity limit to get a solution of hyperbolic inviscid system. However, boundary layer behavior happens near the free-boundary, so existence time $T^{\varepsilon} \rightarrow 0$ as kinematic viscosity $\varepsilon\rightarrow 0$ in standard sobolev space. Inspired by \cite{NMFR1}, we use sobolev conormal space to derive uniform regularity in viscosity $\varepsilon$. Finally, we get a solution of inviscid free-boundary magnetohydrodynamics when initial magnetic field is zero on the free-boundary and in vacuum. 

\end{abstract}


\section{Introduction}
In 1981, Beale \cite{TB} proved local existence result for free-boundary problem of Navier-Stokes equation. Also he proved global regularity for surface tension case in his later work \cite{TB2}. Similar work was done by many authors. We refer Allain \cite{GA} and Tani \cite{AT}, \cite{TT}. In \cite{TT}, using his local result \cite{AT}, Tani claimed global existence for both with and without surface tension. These works depend on Stokes regularity using diffusive effect of Navier-Stokes. We also refer some works by M. Padula and V.A.Solonnikov, for example, \cite{PS1}, \cite{PS2}.

Without kinematic viscosity, free-boundary Euler problem is much harder to solve. For irrotational case, many researches were performed. Using curl free and divergence free properties of velocity field, we can introduce scalar potential of velocity which solves Laplace's equation. Therefore, we can change the problem into the problem on the free-boundary. We refer some works by S.Wu \cite{SW1}, \cite{SW2}, and \cite{SW3}, and recent work by Germain, Masmoudi, and Shatah \cite{GMS}, in which they used space-time resonance method. For general rotational case, only local in time results are known. See Lindblad \cite{HL}, for example. He used apriori result in \cite{CL} and Nash-Moser technique to prove local existence. Also we refer Coutand and Shkoller \cite{CS}, Shatah and Zeng \cite{SZ}, and Masmoudi and Rousset \cite{NMFR1}. Especially Masmoudi and Rousset \cite{NMFR1} solved the problem by inviscid limit. 

If we consider conducting fluid, such as plasma, we gain magnetohydrodynamics(MHD). Movement of fluid itself generate electromagnetic force, called Lorentz force. Moreover, for magnetic field, we have Faraday's law which is part of Maxwell's equations of electromagnetism. Hence, we should consider two coupled PDEs with several boundary conditions and divergence free conditions. For viscous MHD, M.Padula and V.A.Solonnikov got some results in their works \cite{PS1}, \cite{PS2}, and \cite{VS}, where they used Stokes regularity which comes from smoothing effect of velocity $\triangle v$ and magnetic field $\triangle H$. However, there have been only few results for the free-boundary inviscid MHD. Because of nonlinear couplings between velocity and magnetic field, linearization by standard lagrangian map is not good way to approach our problem. Linearized compressible plasma-vacuum problem was studied by Trakhinin \cite{TRA2} and current-vortex sheet problem was studied in \cite{TRA1} and also in \cite{CHEN}. Recently, C.Hao and T.Luo \cite{HAO} got apriori estimate for inviscid case in the spirit of \cite{CL}. \\

\subsection{Free-boundary MHD with zero magnetic field on the free-boundary.} 
\indent Let us formulate inviscid free-boundary MHD problem in whole $\mathbb{R}^{3}$ with infinite depth. We write velocity field $u=(u_1,u_2,u_3)$ and magnetic field $H=(H_1,H_2,H_3)$. We use $\Omega(t)$ to denote domain of plasma at time $t\geq 0$ and write initial domain as $\Omega(0):= \Omega$. We also write free-surface as $S_F(t)$, and initial surface as $S_F(0):=S_F$. We use $h(t,y_{1}, y_{2})$ for profile of the free-surface, where $y:=(y_{1},y_{2})$ is two dimensional horizontal variables. Note that velocity $u$ is defined in plasma region $\O(t)$, whereas magnetic field $H$ is defined in whole space $\mathbb{R}^{3}$. In the following (\ref{lim mhd}), the first equation is well-known Euler equation with Lorentz force where $P$ is pressure including constant downward gravitational force. The second one is Faraday's law in Maxwell's equations. Both two vector fields $u$ and $H$ are divergence free, which mean incompressibility of fluid and non-existence of magnetic monopole from electromagnetic theory, which is part of Maxwell's equations. On the free-boundary $S_{F}(t)$, kinematic boundary condition and continuity of stress tensor condition are considered in sixth and fifth equation. In (\ref{lim mhd data}), we give compatible initial data $u_{0}$ and $H_{0}$, especially $H_{0}$ has zero value on the initial free-boundary $S_{F}$ in vacuum region $\mathbb{R}^{3}\backslash \O$.

\begin{equation} \label{lim mhd}
\begin{cases}
\partial_t u + (u\cdot\nabla)u + \nabla P = (H\cdot\nabla)H - \frac{1}{2}\nabla|H|^2,\quad\text{in}\quad\Omega(t),  \\
\partial_t H + (u\cdot\nabla)H = (H\cdot\nabla)u  ,\quad\text{in}\quad\Omega(t),  \\
\nabla\cdot u = 0,\quad\text{in}\quad\Omega(t), \\
\nabla\cdot H = 0,\quad\text{in}\quad\Omega(t), \\
P\mathbf{n} = gh\mathbf{n} + (H \otimes H - \frac{1}{2}I|H|^2)\mathbf{n},\quad\text{on}\quad S_{F}(t), \\
\partial_t h = u \cdot\mathbf{N},\quad\text{on}\quad S_{F}(t), \\
\end{cases}
\end{equation}
with initial data
\begin{equation} \label{lim mhd data}
\begin{cases}
u(0)=u_0,\quad H(0)=H_0,\quad\text{in}\quad \O, \\
\nabla\cdot u_{0} =0,\quad \nabla\cdot H_{0} = 0,\quad\text{in}\quad \O, \\
H_{0} = 0 \quad\text{on}\quad S_{F}\cup \{\mathbb{R}^{3}\backslash \O \},
\end{cases}
\end{equation}
where outward normal vector $\mathbf{N} := (-\p_{1}h, -\p_{2}h, 1)$ and $\mathbf{n} := \frac{\mathbf{N}}{|\mathbf{N}|}$. $g$ is gravitational constant. Note that downward gravitational force was combined with pressure $P$ and therefore we see constant $g$ in the fifth equation in (\ref{lim mhd}), instead of in the first equation. \\

Let us study special structure given by initial data $H_{0}\vert_{\p\O}=0$. If this holds, we have
\begin{equation*}
	\p_{t} H + (u\cdot\nabla)H = 0, \quad\text{on}\quad S_{F}:=\p\O.
\end{equation*}
Therefore, magnetic field is zero along particles which were on the initial free-boundary. Meanwhile, sixth equation in (\ref{lim mhd}) is equivalent to 
$$
\frac{DF}{Dt} = 0\quad\text{on}\quad S_F(t),
$$
where $S_F(t)$ is given by equation $F(t,y_{1},y_{2},z) := z - h(t,y_{1},y_{2})$ and $\frac{D}{Dt}$ is material derivative. So, a particle on initial boundary stays on the free-boundary as far as we have smooth solution of (\ref{lim mhd}). Finally, we get boundary condition 
\[
	H = 0 \quad\text{on}\quad S_{F}(t),
\] 
for (\ref{lim mhd}). In general, it is not natural to impose zero boundary value in (\ref{lim mhd}), because second equation is hyperbolic. Instead, this hidden boundary condition should be understood in the sense of propagation from initial data. 	\\

Meanwhile, there is no magnetic viscous dissipation crossing the free-boundary (or equivalently energy conservation). So initial condition 
\[
H_{vac} = 0 \quad\text{in}\quad \{\mathbb{R}^{3}\backslash \O \} 
\]
gives zero vacuum magnetic field
\begin{equation} \label{lim vac zero}
H_{vac} = 0,\quad \text{in}\quad \{\mathbb{R}^{3}\backslash \O(t) \}.
\end{equation}

\subsection{Viscous free-boundary MHD with zero magnetic field in vacuum.}
To solve (\ref{lim mhd}), we construct parabolic approximation system. Let $\varepsilon > 0$ be kinematic viscosity of Navier-Stokes equation and $\lambda > 0$ be magnetic diffusivity of Faraday's law. $\lambda$ is in fact $\lambda = \frac{1}{\mu\sigma}$, where $\mu$ is constant vacuum permeability and $\sigma$ is electric conductivity of material. So, magnetic diffusivity limit means $\sigma \rightarrow \infty$, which implies perfect conductivity limit of plasma. Considering viscous effects for fluid and magnetic field, we construct the following system \\
\begin{equation} \label{1.1}
\begin{cases}
\partial_t u + (u\cdot\nabla)u + \nabla P = \varepsilon\triangle u + (H\cdot\nabla)H - \frac{1}{2}\nabla|H|^2,\quad\text{in}\quad\Omega(t),\quad \varepsilon>0, \\
\partial_t H + (u\cdot\nabla)H = \lambda\triangle H + (H\cdot\nabla)u  ,\quad\text{in}\quad\Omega(t), \quad \lambda>0, \\
\nabla\cdot u = 0,\quad\text{in}\quad\Omega(t), \\
\nabla\cdot H = 0,\quad\text{in}\quad\Omega(t), \\
P\mathbf{n} - 2\varepsilon\mathbf{S}(u)\mathbf{n}= gh\mathbf{n} + (H \otimes H - \frac{1}{2}I|H|^2)\mathbf{n},\quad\text{on}\quad S_{F}(t), \\
\partial_t h = u \cdot\mathbf{N},\quad\text{on}\quad S_{F}(t), \\
H = 0,\quad\text{on}\quad S_{F}(t) \cup \{\mathbb{R}^{3}\backslash\O(t)\} , \\
\end{cases}
\end{equation}
with initial compatibility conditions,
\begin{equation} \label{1.1 data}
\begin{cases}
u(0)=u_0,\quad H(0)=H_0,\quad\text{in}\quad \O, \\
\nabla\cdot u_0 = \nabla\cdot H_0 = 0,\quad\text{in}\quad \O, \\
H_0 = 0,\quad\text{on}\quad S_{F} \cup \{\mathbb{R}^{3}\backslash\O \} , \\
\Pi\mathbf{S}(u_0)\mathbf{n} = 0,\quad\text{on}\quad S_F ,\\
\end{cases}
\end{equation}
where $\mathbf{S}(u)$ denotes symmetric part of $\nabla u$, 
\begin{equation} \label{def S}
\mathbf{S}(u) := \frac{\nabla u + (\nabla u)^T}{2},
\end{equation}
and $\Pi$ means tangential projection operator, $\Pi := \mathbf{I} - \mathbf{n}\otimes\mathbf{n}$.

	Let us explain boundary condition $H=0$ on the free-surface and in the vacuum. Let us use $H_{vac}$ and $H_{plasma}$ to denote magnetic fields in vacuum and plasma regions, respectively. When displacement current is not assumed, magnetic field in vacuum solves
	\[
		\nabla\cdot H_{vac} = 0,\quad \nabla\times H_{vac} = 0,
	\]
	by Maxwell's equation. Unlike to inviscid case (\ref{lim mhd}), zero boundary value of $H$ does not propagate from initial data. Therefore, zero magnetic value on the free-boundary and in vacuum should be understood in the sense of imposed constraint. Physically we may need some complicate equipments to realize such condition but this system makes sense mathematically. For hyperbolic system (\ref{lim mhd}) we cannot impose such constraint unless we have propagation from initial data but for parabolic system (\ref{1.1}), it is not overdetermined. Note that linearized equation in Lagrangian coordinates for Navier-Stokes and Faraday's law are just Stokes equation and heat equation  with zero boundary value, respectively. See (\ref{app stokes system}) and (\ref{app heat system}) for example.  

\begin{remark} \label{rmk boundary}	
	Since we are assuming nonzero magnetic diffusivity, plasma is not perfect conductor, which implies zero surface current. Therefore, magnetic field is curl free on the free-boundary. Of course, divergence is also divergence-free on the boundary by Maxwell's equations. Applying divergence-free condition and divergence theorem to the cylinder with infinitesimal height crossing the free-boundary, we easily get normal continuity
	\[
		H_{vac}\cdot\mathbf{n} = H_{{plasma}}\cdot\mathbf{n},\quad\text{on}\quad S_{F}(t).
	\]  
	Similarly, applying curl free condition and Stokes' theorem to closed circuit near the boundary, we get tangential continuity
	\[
		\big(\mathbf{I-\mathbf{n}\otimes\mathbf{n}} \big)H_{vac} = \big(\mathbf{I-\mathbf{n}\otimes\mathbf{n}} \big) H_{{plasma}},\quad\text{on}\quad S_{F}(t).
	\] 
	Therefore, magentic field is continuous on the boundary and special condition $H_{vac}=0$ gives $H=0$ on the free-boundary. Note that if we have $\lambda=0$, it corresponds to $\sigma=\infty$ which implies perfect conductor. In that case, existence of surface current gives normal continuity only. In mathematical aspect, when plasma is perfect conductor, $\lambda=0$ yield hyperbolic PDE for $H$ and we cannot give boundary value to the second PDE of (\ref{1.1}), in general. Instead, if initial data $H_0$ satisfies
	\[
		H_0\cdot\mathbf{n}=0\quad \text{on}\quad S_{F}(0),
	\]
	then we have
	\[
		H \cdot\mathbf{n}=0\quad \text{on}\quad S_{F}(t),
	\]
	as long as smooth solution exists. See \cite{HAO} for the detail, for example. Note that propagation of zero boundary value $H$ is special case of propagation of $H\cdot\mathbf{n}=0$ on the free-boundary. 
\end{remark}	
	
	In (\ref{1.1}), it is convenient to define total pressure $p$ as sum of pressure and magnetic pressure,
	$$
		p := P + \frac{1}{2}|H|^2.
	$$
Considering all of these, (\ref{1.1}) system becomes,
\begin{equation} \label{1.2}
\begin{cases}
\partial_t u + (u\cdot\nabla)u - (H\cdot\nabla)H + \nabla p = \varepsilon\triangle u,\,\,\,\,\text{in}\,\,\,\,\Omega(t), \\
\partial_t H + (u\cdot\nabla)H - (H\cdot\nabla)u  = \lambda\triangle H,\,\,\,\,\text{in}\,\,\,\,\Omega(t), \\
\nabla\cdot u = 0,\,\,\,\,\text{in}\,\,\,\,\Omega(t), \\
\nabla\cdot H = 0,\,\,\,\,\text{in}\,\,\,\,\Omega(t), \\
p\mathbf{n} - 2\varepsilon\mathbf{S}(u)\mathbf{n}= gh\mathbf{n},\,\,\,\,\text{on}\,\,\,\,S_{F}(t), \\
H = 0,\quad\text{on}\quad S_{F}(t) \cup \{\mathbb{R}^{3}\backslash\O(t)\} , \\
\partial_t h = u\cdot\mathbf{N},\,\,\,\,\,\,\,\text{on}\,\,\,\,S_{F}(t) ,\\
\end{cases}
\end{equation}
with initial compatibility conditions,
\begin{equation} \label{1.3}
\begin{cases}
u(0)=u_0,\,\,\,\,H(0)=H_0,\quad\text{in}\quad\Omega, 	\\
\nabla\cdot u_0 = \nabla\cdot H_0 = 0,\quad\text{in}\quad \O, \\
H_0 = 0,\quad\text{on}\quad S_{F} \cup \{\mathbb{R}^{3}\backslash\O \} , \\
\Pi\mathbf{S}(u_0)\mathbf{n} = 0,\quad\text{on}\quad S_F . 	\\
\end{cases}
\end{equation}



\subsection{Parametrization into a fixed domain}
We rewrite the system (\ref{1.2}) and (\ref{1.3}) in the fixed domain $S:=\{(y_{1},y_{2},z)|z<0\}$, lower half space in $\mathbb{R}^3$. We use ${x}:=(y_{1},y_{2},z)$ to denote a point in $\mathbb{R}^{3}$ and ${y}:=(y_{1},y_{2})$ to denote horizontal coordinate. Parametrization into the fixed domain $S$ is gained by  $\Phi(t,\cdot)$,
\begin{equation} \label{Phi}
\begin{split}
\Phi(t,\cdot): \ S \rightarrow \Omega(t), \quad \Phi(t, {y},z) = ( {y},\varphi(t, {y},z)),\quad z\leq 0,
\end{split}
\end{equation}
where $\varphi$ is to be defined in below (\ref{phi}). When $z>0$, we define 
\begin{equation*}  
\begin{split}
\Phi(t,\cdot): \ \mathbb{R}^{3}\backslash S \rightarrow \mathbb{R}^{3}\backslash\Omega(t), \quad \Phi(t, {y},z) = ( {y}, z + h(t,y)),\quad z>0.
\end{split}
\end{equation*}

 We use function $v$, $B$, and $q$ for velocity, magnetic field, and pressure in the fixed domain $S$,
\begin{equation*}
\begin{split}
v(t,{x}) &:=u(t,\Phi(t,{x})),\quad q(t,{x}):=p(t,\Phi(t,{x})),\quad {x}\in S \\
B(t,{x}) &:=H(t,\Phi(t,{x})),\quad  {x}\in \mathbb{R}^{3}. 
\end{split}
\end{equation*}
We have to decide $\varphi(t,\cdot)$ in (\ref{Phi}) so that $\Phi(t,\cdot)$ becomes a diffeomorphism between $S$ and $\O(t)$. From determinant of $\nabla\Phi$, we should have $\partial_z\varphi > 0$ for diffeomorphism. There are many ways to take $\varphi$. One easy option is to set $\varphi(t,y,z)=z+h(t,y)$. However, there is no gain of regularity, so this fits in the Euler equations case. Instead, we take a smoothing diffeomorphism similar as \cite{NMFR1}, \cite{DL}, \cite{BS}, and \cite{DHL}. We define 
\begin{equation} \label{phi}
\varphi(t,{y},z) := Az+\eta(t,{y},z).
\end{equation} 
To ensure that $\Phi(0,\cdot)$ is a diffeomorphism, $A$ should be picked so that 
\begin{equation} \label{diffeo}
\partial_z\varphi(t=0,{y},z)\geq 1,\quad\forall ({y},z)\in S,
\end{equation}
and $\eta$ is given by extension of $h$ to the inside of domain $S$, defined by 
\begin{equation} \label{eta def}
\hat{\eta}(\xi,z)=\chi(z\xi)\hat{h}(\xi),
\end{equation}
where $\hat{\cdot}$ is horizontal Fourier transform and $\xi$ is corresponding two dimensional frequency variable. $\chi$ is a smooth, compactly supported function which is 1 on the unit ball $B(0,1)$. This smoothing diffeomorphism was used in \cite{DL}, \cite{BS}, and also in \cite{NMFR1}. In Proposition \ref{proposition 3.7}, we will see that $\varphi$ has $\frac{1}{2}$ better regularity than $h$. \\ 

\indent We also define new derivatives of $v,B$ in $S$, to measure $\partial_i u,\partial_i B$ in the fixed domain $S$. Then we could rewrite the systems (\ref{1.2}) and (\ref{1.3}) in a fixed domain $S$. Using change of variable, we get,
\begin{equation*}
\begin{split}
(\partial_i u)(t,{y},\varphi) &= (\partial_i v - \frac{\partial_i\varphi}{\partial_z\varphi}\partial_z v)(t,{y},z),\quad i=t,1,2,  \\
(\partial_i H)(t,{y},\varphi) &= (\partial_i B - \frac{\partial_i\varphi}{\partial_z\varphi}\partial_z B)(t,{y},z),\quad i=t,1,2,  \\
(\partial_3 u)(t,{y},\varphi) &= (\frac{1}{\partial_z\varphi}\partial_z v)(t,{y},z),  \\
(\partial_3 H)(t,{y},\varphi) &= (\frac{1}{\partial_z\varphi}\partial_z B)(t,{y},z).
\end{split}
\end{equation*}
So it is convenient to define the following operator,
\begin{equation} \label{new deriva}
\partial_i^\varphi := \partial_i - \frac{\partial_i\varphi}{\partial_z\varphi}\partial_z,\quad\text{for}\quad i=t,1,2, \quad\text{and}\quad \partial_{z}^\varphi := \frac{1}{\partial_z\varphi}\partial_z.
\end{equation}
This definition implies that $\partial_i f \circ \Phi = \partial_i^\varphi (f\circ\Phi),\quad i=t,1,2,z$ for smooth $f$ defined in $S$. Hence, (\ref{1.2}) and (\ref{1.3}) are written in $S$ as  following.
\begin{equation} \label{1.4}
\begin{cases}
\partial_t^\varphi v + v\cdot\nabla^\varphi v - B\cdot\nabla^\varphi B + \nabla^\varphi q = 2\varepsilon\nabla^\varphi\cdot (\mathbf{S}^\varphi v),\quad \text{in}\quad S, \\
\partial_t^\varphi B + v\cdot\nabla^\varphi B - B\cdot\nabla^\varphi v = 2\lambda\nabla^\varphi\cdot(\mathbf{S}^\varphi B),\quad \text{in}\quad S, \\
\nabla^\varphi\cdot v = 0,\quad \text{in}\quad S, \\ 
\nabla^\varphi\cdot B = 0,\quad \text{in}\quad S, \\ 
q\mathbf{n} - 2\varepsilon(\mathbf{S}^\varphi v)\mathbf{n} = gh\mathbf{n},\quad \text{on}\quad \partial S, \\
\partial_t h = v\cdot \mathbf{N},\quad \text{on}\quad \partial S, \\
B = 0,\quad \text{on}\quad \partial S \cup \{\mathbb{R}^{3}\backslash S\}, \\
\end{cases}
\end{equation}
with initial compatibility condition
\begin{equation} \label{1.5}
\begin{cases}
v(0) = v_0,\quad B(0) = B_0,\quad\text{in}\quad S, 	\\
\nabla^\varphi\cdot v_0 = \nabla^\varphi\cdot B_0 = 0,\quad\text{in}\quad S, \\ 
B_0 = 0\quad\text{on}\quad \partial S \cup \{\mathbb{R}^{3}\backslash S\}, \\ 
\Pi\mathbf{S}^\varphi (v_0)\mathbf{n} = 0\quad \text{on}\quad \partial S. \\
\end{cases}
\end{equation}
Note that $\mathbf{S}^{\varphi}$ is defined by changing $\nabla$ into $\nabla^{\varphi}$ in the definition of (\ref{def S}). 

\begin{remark}
In this paper, we use $\nabla$ for full gradient $\nabla f := (\p_{1}, \p_{2}, \p_{3})f$ and $\nabla_{y}$ for horizontal gradient $(\p_{1}, \p_{2})$. And, since boundary profile $h(t,y)$ is extended into inside the domain, we define
\begin{equation}
	\mathbf{N} := (-\nabla_{y}\varphi,1),\quad \mathbf{n}:= \frac{\mathbf{N}}{|\mathbf{N}|} \quad \text{in}\quad S. 
\end{equation}
On the free-boundary $\p S$, $\mathbf{N}\vert_{z=0} = (-\nabla_y h, 1)$ and $\mathbf{n}\vert_{z=0}$ is just outward unit normal vector on the free-boundary as we used in the systems above. Sometimes we use upper index $\cdot^{b}$ to stress $f^{b} := f\big\vert_{z=0}$. However, we skip upper $\cdot^{b}$ notation when we have no confusion. 
\end{remark}


\subsection{Functional Framework and Notations}
We briefly introduce conormal space and other function spaces those are proper to our analysis. See \cite{NMFR1} for more explanation. First we define sobolev conormal derivatives.
\begin{definition} \label{definition 1.1}
We define conormal derivatives in $S$.
\begin{equation*}
Z_1 := \partial_{y_{1}},\,\,\,\,\,Z_2 := \partial_{y_{2}},\,\,\,\,\,Z_3 := \frac{z}{1-z}\partial_z,\,\,\,\,\,Z^\alpha := Z^{(\alpha_1,\alpha_2,\alpha_3)},
\end{equation*}
\begin{equation*}
\left|f\right|_{H_{co}^s(S)}^2 := \sum_{|\alpha|\leq s}\left|Z^\alpha f\right|_{L^2(S)}^2,\,\,\,\,\left|f\right|_{W_{co}^{s,\infty}(S)} := \sum_{|\alpha|\leq s}\left|Z^\alpha f\right|_{L^\infty(S)}.
\end{equation*} 
\end{definition}
In this paper we abbreviate the notation as $|\cdot|_s=|\cdot|_{H^s}$, $\|\cdot\|_s=|\cdot|_{H_{co}^s}$ and $\|\cdot\|=|\cdot|_{L^2} $. Similarly, $|\cdot|_{s,\infty}=|\cdot|_{W^{s,\infty}}$ and $\|\cdot\|_{s,\infty}=|\cdot|_{W_{co}^{s,\infty}}$. For horizontal component, we use $v_y := (v_1,v_2)$ and $\nabla_y := (\partial_1,\partial_2)$. Sometimes we may use notation $Z^m$. This means some $Z^\alpha=Z^{(\alpha_1,\alpha_2,\alpha_3)}$, where $|\alpha|:=|\alpha_1| + |\alpha_2| + |\alpha_3|=m$. We will add for all such $\alpha$, so we do not need to specify $\alpha$. 

\begin{definition} \label{definition 1.2}
For $m \geq 1$,
$$
E^m := \{f\in H^m_{co},\,\,\partial_z f \in H^{m-1}_{co} \}\quad \text{and}\quad E^{m,\infty} := \{f\in W^{m,\infty}_{co},\,\,\partial_z f \in W^{m-1,\infty}_{co} \},
$$
with norms,
$$
\|f\|^2_{E^m} := \|f\|^2_m + \|\partial_z f\|^2_{m-1},\,\,\,\|f\|_{E^{m,\infty}} := \|f\|_{m,\infty} + \|\partial_z f\|_{m-1,\infty}.
$$
\end{definition}

We also define tangential sobolev spaces which is weaker than conormal space.
\begin{definition} \label{definition 1.3}
$$
H^s_{tan}(S) := \{ f\in L^2(S),\,\,\,\Lambda^s f \in L^2(S)\},
$$
where $\Lambda^s$ is tangential Fourier multiplier by $(1+|\xi|^2)^{s/2}$, with norm
$$
\|f\|_{H^s_{tan}} := \|\Lambda^s f\|_{L^2}.
$$
\end{definition}
\noindent $\textbf{Notation}$ In this paper,  $\Lambda_0$ means $\Lambda(\frac{1}{c_0})$ where $\Lambda(\cdot)$ is a monotone increasing function and $\Lambda(\cdot,\cdot)$ is a monotone increasing function with all its variables. They may vary line to line.

\subsection{Main results}
We state two main results for this paper. First theorem states uniform energy estimate in Sobolev conormal space when $\varepsilon=\lambda$. We refer \cite{DHL} or \cite{Solonni} (general case with bounded domain) for local well-posedness of (\ref{1.1}). Or see Theorem \ref{app main theorem} in appendix for scheme of well-posedness.

\begin{theorem} \label{theorem 1.4}
Assume that (\ref{1.4}) is well-posed and consider $\varepsilon=\lambda \in (0,1]$ case. For fixed sufficiently large $m\geq 6$, let initial data $(v^{\varepsilon}_{0}, B^{\varepsilon}_{0}, h^{\varepsilon}_{0})$ be given so that
\begin{equation} \label{1.6}
\begin{split}
&\sup_{\varepsilon\in(0,1]} ( |h_0^\varepsilon|_{m} + \sqrt{\varepsilon}|h_0^\varepsilon|_{m+\frac{1}{2}} + \|v_0^\varepsilon\|_m + \|B_0^\varepsilon\|_m + \|\partial_z v_0^\varepsilon\|_{m-1} + \|\partial_z B_0^\varepsilon\|_{m-1}   \\
&\quad\quad\quad + \|\partial_z v_0^\varepsilon\|_{1,\infty} + \|\partial_z B_0^\varepsilon\|_{1,\infty} + \sqrt{\varepsilon}\|\partial_{zz} v_0^\varepsilon\|_{L^\infty} + \sqrt{\varepsilon}\|\partial_{zz} B_0^\varepsilon\|_{L^\infty} ) \leq R,
\end{split}
\end{equation}
and satisfy Taylor sign condition $-\partial_n q + g \geq c_0 > 0$ uniform in $\varepsilon$. Also, they should satisfy initial compatibility conditions
\begin{equation} \label{1.7}
\begin{split}
&\nabla^\varphi\cdot v_0^\varepsilon = 0,\,\,\,\,\nabla^\varphi\cdot B_0^\varepsilon = 0,\,\,\,\,\text{in}\,\,\,S,  \\
&\Pi \left( \mathbf{S}^{\varphi} v^\varepsilon_0\right)\mathbf{n} = 0,\quad \text{on}\quad \partial S, \\
& B_0^\varepsilon = 0 \quad\text{on}\quad \p S \cup \{\mathbb{R}^{3}\backslash S\}.
\end{split}
\end{equation}


Then there exist $T>0$ (uniform in $\varepsilon$) and some $C>0$ such that there exist a solution $(v^\varepsilon,B^\varepsilon,h^\varepsilon)$ for (\ref{1.4}), (\ref{1.7}), (\ref{phi}), and (\ref{eta def}) on $[0,T]$. Moreover the following energy estimates hold. For non-dissipation type, 	\\

\begin{equation} \label{1.8}
\begin{split}
&\sup_{t\in [0,T]}\left( |h^\varepsilon|^2_m + \|v^\varepsilon\|^2_m +  \|B^\varepsilon\|^2_m + \|\partial_z v^\varepsilon\|^2_m + \|\partial_z B^\varepsilon\|^2_m + \|\partial_z v^\varepsilon\|^2_{1,\infty} + \|\partial_z B^\varepsilon\|^2_{1,\infty} \right)  \\
&\quad \quad\quad + \|\partial_z v^\varepsilon\|^2_{L^4([0,T],H^{m-1}_{co})} + \|\partial_z B^\varepsilon\|^2_{L^4([0,T],H^{m-1}_{co})} \leq C.
\end{split}
\end{equation}
For dissipation type,
\begin{equation} \label{1.9}
\begin{split}
&\sup_{t\in [0,T]}\left( \varepsilon|h^\varepsilon|_{m+\frac{1}{2}}^2 + \varepsilon\|\partial_{zz}v^\varepsilon\|^2_{L^\infty} + \varepsilon\|\partial_{zz}B^\varepsilon\|^2_{L^\infty} \right)   \\
&\quad\quad\quad+ \varepsilon\int_0^T \left( \|\nabla v^\varepsilon\|^2_m + \|\nabla B^\varepsilon\|^2_m + \|\nabla \partial_z v^\varepsilon\|^2_{m-2} +  \|\nabla \partial_z B^\varepsilon\|^2_{m-2} \right) \leq C.
\end{split}
\end{equation}
\end{theorem}

Using Theorem \ref{theorem 1.4}, we get zero kinematic viscosity-magnetic diffusivity limit. 
\begin{theorem} \label{theorem 1.5}
Let us assume that we have a unique local solution $(v^\varepsilon,B^\varepsilon,h^\varepsilon) \in$ for (\ref{1.4}) with $\varepsilon=\lambda$ on $[0,T^{\varepsilon})$ for initial data $(v_0^\varepsilon, B_0^\varepsilon, h_0^\varepsilon)$ given in Theorem \ref{theorem 1.4}. We also assume
\begin{equation} \label{1.10}
\lim_{\varepsilon\rightarrow 0}\left( \left\|v_0^\varepsilon - v_0\right\|_{L^2(S)} + \left\|B_0^\varepsilon - B_0\right\|_{L^2(S)} + \left\|h_0^\varepsilon - h_0\right\|_{L^2(\partial S)} \right) = 0,
\end{equation}
where $(v_0, B_0, h_0)$ also satisfy assumptions of Theorem \ref{theorem 1.4}. Then there exist $(v,B,h)$ such that those are in 
\begin{equation} \label{1.11}
v,B \in L^\infty([0,T],H_{co}^m(S)),\,\,\partial_z v,\partial_z B \in L^\infty([0,T],H_{co}^{m-1}(S)),\,\,h \in L^\infty([0,T],H_{co}^{m}(\mathbb{R}^2))
\end{equation}
and
\begin{equation} \label{1.12}
\begin{split}
&\lim_{\varepsilon\rightarrow 0}\sup_{[0,T]} ( \left\|v^\varepsilon - v\right\|_{L^2(S)} + \left\|v^\varepsilon - v\right\|_{L^\infty(S)} + \left\|B^\varepsilon - B\right\|_{L^2(S)} \\
&\quad\quad\quad+ \left\|B^\varepsilon - B\right\|_{L^\infty(S)} + \left\|h^\varepsilon - h\right\|_{L^2(\partial S)} + \left\|h^\varepsilon - h\right\|_{W^{1,\infty}(\partial S)} ) = 0.
\end{split}
\end{equation}
Moreover, $(v,B,q)$ solves 
\begin{equation} \label{1.13}
\begin{cases}
\partial_t^\varphi v + v\cdot\nabla^\varphi v + \nabla^\varphi q = B\cdot\nabla^\varphi B,\quad \text{in}\quad S, \\
\partial_t^\varphi B + v\cdot\nabla^\varphi B = B\cdot\nabla^\varphi v,\quad\text{in}\quad S, \\
\nabla^\varphi\cdot v = 0,\quad\text{in}\quad S, \\ 
\nabla^\varphi\cdot B = 0,\quad\text{in}\quad S, \\ 
q \mathbf{n} = gh\mathbf{n},\quad\text{on}\quad\partial S, \\
\partial_t h = v\cdot \mathbf{N},\quad\text{on}\quad\partial S, \\
\end{cases}
\end{equation}
with initial data
\begin{equation} \label{1.13 data}
\begin{cases}
v(0) = v_0,\,\,\,B(0) = B_0,\quad \text{in}\quad S, 	\\
\nabla^{\varphi}\cdot v_0 = \nabla^{\varphi}\cdot B_0 = 0,\quad\text{in}\quad S, \\
B_0 = 0,\quad\text{on}\quad \p S \cup \{\mathbb{R}^{3}\backslash S \} .	\\
\end{cases}
\end{equation}
\end{theorem}

\subsection{Scheme of proof}
\indent The aim of this paper is to gain uniform estimate of viscous free-boundary magnetohydrodynamics when magentic field is zero in vacuum and on the free-boundary. From bondary layer behavior, we expect 
\[
u \sim u_{E} + \sqrt{\varepsilon} U(t,y,{z\over {\sqrt{\varepsilon}}}),
\]
when $u_{E}$ is solution of limit system $\varepsilon=0$ and $U$ is a some profile. Therefore, boundary layer behavior makes it impossible to get uniform estimate in standard sobolev space. Note that sobolev conormal space kills normal derivatives on the boundary by multiplying factor ${z\over 1-z}$, and is equivalent to standard sobolev space away from the boundary because $\p_{z}^{k}{z \over 1-z}$ is uniformly bounded for all order $k\in\mathbb{N}$ away from $z=0$. In \cite{NMFR1}, Masmoudi and Rousset used this functional framework to solve vanishing viscosity limit problem of Navier boundary problem in \cite{NMFR2} and free-boundary problem without surface tension in \cite{NMFR1}. For the free-boundary problem with surface tension, we refer Tarek and Donghyun \cite{TD}. \\
\\
\noindent\textit{1. Alinhac's unknown} We use fixed domain $S$ as defined in (\ref{1.4}) and apply conormal derivatives $Z^{m}$. Since we do not take surface tension into account, we get $\|v\|_m$ and $|h|_m$ as high order energy. However, we have some bad commutators which requires $|h|_{m+\frac{1}{2}}$ to be controlled, which is $\frac{1}{2}$-order higher than energy $|h|_{m}$. This problem can be fixed by introducing Alinhac's new unknown for transformed velocity $v=u\circ\Phi$ and magnetic field $B=H\circ\Phi$. We study structure of commutators and we will see that bad commutators have transport structure. We will write the system (\ref{1.4}) in terms of Alinhac's new variables,
\begin{equation} \label{new Alinhac}
(\mathcal{V}^\alpha,\mathcal{Q}^\alpha,\mathcal{B}^\alpha) := ( Z^\alpha v - \partial_z^\varphi v Z^\alpha \eta, Z^\alpha q - \partial_z^\varphi q Z^\alpha \eta,Z^\alpha B - \partial_z^\varphi B Z^\alpha \eta ),
\end{equation}
and then all commutators becomes low order in terms of $h$. Meanwhile, non-transport type nonlinear terms,
\[
B\cdot\nabla^{\varphi} B,\quad B\cdot\nabla^{\varphi} v,
\]
in (\ref{1.4}) generate some terms with higher order than energy. From (\ref{1.4}), we have two energy estimate from first two PDEs. Divergence free property and boundary condition of $B$ play critical roles to cancel these high order terms when we combine two energy estimates. This property hold for another variables whenever we use $L^{2}$-type energy estimates.  \\

\noindent\textit{2. Normal derivative} Second problem arises from the fact that conormal space is weaker than standard sobolev space. To control commutators we should control normal derivative terms, $\|\p_z v\|_{m-1}$ and $\|\p_z B\|_{m-1}$. These terms contain $m$ derivatives in total. From definition \ref{definition 1.1}, however, these cannot be controlled by $\|v\|_m$ or $\|B\|_m$, because conormal space is weaker than standard sobolev space. To estimate $\p_z v$, we introduce new variables,
\[
S^{v}_n := (\mathbf{I}-\mathbf{n}\otimes\mathbf{n})(\mathbf{S}^{\varphi}v)\mathbf{n} ,\quad S^{B}_n := (\mathbf{I}-\mathbf{n}\otimes\mathbf{n})(\mathbf{S}^{\varphi}B)\mathbf{n}.
\]
These variables are equivalent to $\p_z v$ and $\p_z B$ in conormal space. Moreover, their boundary values vanish by fifth equation and seventh equation in (\ref{1.4}). \textit{ However, equation of $S_n^{v}$ produce $\nabla^{2} p$ which require $|h|_{m+{1\over 2}}$ to be controlled. } Therefore, optimal regularity of $S_n^{v}$ and $S_n^{B}$ are $m-2$, not $m-1$. We cannot use these estimates to close energy estimate, but it can be used to control $L^\infty$-type terms in the next step. \\

\noindent\textit{3. $L^\infty$ type estimate} It is natural to expect finite low order $L^{\infty}$ terms in commutators. Unfortunately, we cannot use standard sobolev embedding in conormal space even with sufficient high order regularity. Instead, we use the high order energy of $ S^{v}_n$ and $S^{B}_n$ to control these finite low order $L^\infty$ type terms in commutators. Using conormal embedding with sufficiently large $m\geq 6$, $\|S_{n}^{v}\|_{m-2}$ and $\|S_{n}^{B}\|_{m-2}$ control these low order terms. Our basic approach is to use maximal principle of convection-diffusion equation, however, commutators between laplacian $\Delta^{\varphi}$ and conormal derivative $Z^{m}$ are not easy to control. Instead we use geometric re-parametrization which preserves second order normal derivative $\p_{zz}$ structure. \\

 Meanwhile, we have non-transport type nonlinear term,
\[
B\cdot\nabla^{\varphi} B,\quad B\cdot\nabla^{\varphi} v.
\]  
Moreover, we cannot expect cancellation between these two terms, because we use maximal principle, instead of standard $L^{2}$ type energy estimate. So we assume $\varepsilon=\lambda$ and combine two PDE's to get 
\begin{equation} \label{add and sub}
\begin{split}
\p_t^{\varphi} (v+B) + (v-B)\cdot\nabla^{\varphi}(v+B) + \nabla^{\varphi} q = \varepsilon\Delta^{\varphi} (v+B) ,\\
\p_t^{\varphi} (v-B) + (v+B)\cdot\nabla^{\varphi}(v-B) + \nabla^{\varphi} q = \varepsilon\Delta^{\varphi} (v-B) .\\
\end{split}
\end{equation}  
Some analysis with boundary condition $B \vert_{\p S}=0$ will give finite order estimates for 
\[
\|S_n^{v} + S_n^{B}\|_{k,\infty},\quad \|S_n^{v} - S_n^{B}\|_{k,\infty},\quad \text{for some small finite} \ \  k\in\mathbb{N},
\]
and therefore, we get estimates for each $\|S_n^{v}\|_{k,\infty}, \|S_n^{B}\|_{k,\infty}$.  \\

\noindent\textit{4. Vorticity estimate} Since optimal regularities of $S_n^{v}$ and $S_n^{B}$ were $m-2$, we should estimate $m-1$ order for $\p_z v$. The reason we could not reach to $m-1$ order was pressure, since we do not see gradient structure for $p$ in the equation of $S_{n}^{v}$. Equations for vorticity $\omega_{v}$ removes pressure term, so we get hope to get $m-1$ order conormal estimates. For general rotational fluid, however, vorticities 
\[
\o_v := \nabla^\varphi \times v,\quad \o_B := \nabla^\varphi \times B,
\]
do not vanish on the boundary $\p S$. Therefore $L^{2}$-type energy estimate generate boundary integral, which need $\frac{1}{2}$ more space regularity to control trace. To avoid these difficulty, in \cite{NMFR1}, the authors derived new $L^4$ in time estimate using Lagrangian map and microlocal symmetrizer, instead of strong $L^\infty$ type estimate. This is weaker norm in time, but is sufficient to control commutators, since commutator terms contains only $L^2$ in time. Also, we should use equations (\ref{add and sub}) again, because we cannot control each nonlinear terms otherwise. \\
\indent One critical idea of \cite{NMFR1} is to use Lagrangian transform which converts convection-diffusion equation into heat equation. However, we have two equations in (\ref{add and sub}), with different transport velocities, $v\pm B$. We consider two Lagrangian maps $Y_1$ and $Y_2$, so that
\[
\p_t Y_1 = u-H,\quad \p_t Y_2 = u+H,
\] 
then $\O(t)$ is transformed into two fixed initial domains $\O_1$ and $\O_2$. However, both $u$ and $H$ (or $v$ and $B$) must be defined in the same domain by definition of systems. One important remark is that domain $\O(t)$ is defined by boundary velocity only. And, propagation of boundary profile is determined by normal velocity on the free-boundary, \textit{see sixth equation in (\ref{1.4})}. Therefore, when we have $H\cdot\mathbf{n}=0$ (or $H=0$) on the free-boundary, boundary profile is determined only by velocity field $u$. So introducing two maps $Y_1$ and $Y_2$, equations in (\ref{add and sub}) are transformed into \textit{same domain} $\O_{1} = \O_{2} = \O$. Using these transforms and the new $L^4_T$ type result of \cite{NMFR1}, we get estimates for
\[
\|\o_v\|_{L^4_T H^{m-1}_{co}}\quad\text{and}\quad  \|\o_B\|_{L^4_T H^{m-1}_{co}}.
\]
We can check $m-1$ order conormal norm of vorticity is equivalent to $\|\p_z v\|_{m-1}$, and also for $B$ obviously. Finally we can close energy estimate and get uniform regularities for $u$, $B$, and $h$. \\

\noindent\textit{5. Uniqueness and vanihsing viscosity limit} Uniqueness is gained by $L^2$ energy estimate with higher order energy bound, which was obtained from previous steps. For vanishing-viscosity limit $\varepsilon=\lambda \rightarrow 0$, we use compactness argument to get weak limit. And for the limit system $\varepsilon=\lambda=0$, we have $L^{2}$ energy conservation which gives norm convergence in the limit process $\varepsilon\rightarrow 0$. Finally we get $L^2$ strong convergence.


\section{Formal differentiation and Alinhac's unknown}
Suppose that we have a smooth solution $(v,B,h)$. Let us study transport operator and commutators gained by $Z^{m}$. 
From definition of $\p_{i}^\varphi$,
$$
\partial_t^\varphi + v\cdot\nabla^\varphi = \partial_t + v_y\nabla_y + \frac{1}{\partial_z\varphi}(v\cdot\mathbf{N} - \partial_t\eta)\partial_z,\quad\text{where}\quad\mathbf{N} = (-\nabla_y\varphi,1).
$$ 
If we apply conormal derivatives $Z^{\alpha}$ with $|\alpha|=m$, then we see some commutators look like $\|Z^\alpha \mathbf{N}\| \sim |\nabla\varphi|_m \sim |\varphi|_{m+1} \sim |h|_{m+\frac{1}{2}}$ since $\varphi$ is $\frac{1}{2}$ smoother than $h$ by diffeomorphism (\ref{phi}). See \cite{NMFR1}, \cite{DHL}, or Proposition \ref{proposition 3.7} in section 4 for this estimate. Hence, we rewrite the system in terms of Alinhac's new unknowns, because this new unknown kills all these bad commutators. 
In this section, we  construct Alinhac's unknown systematically. First, let us define the following symbols.
\begin{equation} \label{2.1}
\begin{split}
\mathcal{N}(v,B,q,\varphi) &:= \partial_t^\varphi v + (v\cdot\nabla^\varphi)v + \nabla^\varphi q - 2\varepsilon\nabla^\varphi\cdot(\mathbf{S}^\varphi v) -  (B\cdot\nabla^\varphi)B, \\
\mathcal{F}(v,B,\varphi) &:= \partial_t^\varphi B + (v\cdot\nabla^\varphi)B - (B\cdot\nabla^\varphi)v - 2\lambda\nabla^\varphi\cdot(\mathbf{S}^\varphi B), \\
d_v (v,\varphi) &:= \nabla^\varphi\cdot v, \\
d_B (B,\varphi) &:= \nabla^\varphi\cdot B, \\
\mathcal{B}(v,B,q,\varphi) &:= (q-gh)\mathbf{N} - 2\varepsilon(\mathbf{S}^\varphi v)\mathbf{N}.
\end{split}
\end{equation}
In the following two propositions, $\sim$ means first order expansion in terms of $\dot{f},\dot{g}$.
\begin{proposition} \label{proposition 2.1}
We have the following first order expansion.
\begin{equation*}
(f+\dot{f})\cdot\nabla^{\varphi+\dot{\varphi}}(g+\dot{g}) \sim (f\cdot\nabla^\varphi)g + (f\cdot\nabla^\varphi)\dot{g} + (\dot{f}\cdot\nabla^\varphi)g - \partial_z^\varphi g(f\cdot(\nabla^\varphi \dot{\varphi})).
\end{equation*}
\end{proposition}
\begin{proof}
We abbreviate $\sum_{i=1,2}$ for terms with index $i$.
\begin{equation*}
\begin{split}
(f+\dot{f})\cdot\nabla^{\varphi+\dot{\varphi}}(g+\dot{g}) & = (f_i + \dot{f}_i)\left(\partial_i g + \partial_i\dot{g} - \frac{\partial_i\varphi + \partial_i\dot{\varphi}}{\partial_z\varphi + \partial_z\dot{\varphi}}(\partial_z g + \partial_z \dot{g}) \right) + (f_3 + \dot{f}_3)\frac{\partial_z g + \partial_z \dot{g}}{\partial_z\varphi + \partial_z\dot{\varphi}} \\
& \sim f_i \left( \partial_i g + \partial_i\dot{g} - \frac{(\partial_i\varphi + \partial_i\dot{\varphi})(\partial_z g + \partial_z\dot{g})}{\partial_z\varphi} \left( 1 - \frac{\partial_z\dot{\varphi}}{\partial_z\varphi} \right) \right) \\
&\quad + \dot{f}_i \left( \partial_i g + \partial_i\dot{g} - \frac{(\partial_i\varphi + \partial_i\dot{\varphi})(\partial_z g + \partial_z\dot{g})}{\partial_z\varphi} \right) + (f_3 + \dot{f}_3)\frac{\partial_z g + \partial_z \dot{g}}{\partial_z\varphi + \partial_z\dot{\varphi}} \\
& \sim f_i\partial_i g + f_i \partial_i \dot{g} - f_i \frac{(\partial_i\varphi\partial_z g + \partial_i\varphi\partial_z\dot{g} + \partial_i\dot{\varphi}\partial_z g)}{\partial_z\varphi}\left(1-\frac{\partial_z\dot{\varphi}}{\partial_z\varphi}\right) \\
&\quad + \dot{f}_i\partial_i g - \frac{\dot{f}_i}{\partial_z\varphi}\partial_i\varphi\partial_z g + \frac{f_3}{\partial_z\varphi}\left(\partial_z g + \partial_z\dot{g} - \partial_z g\frac{\partial_z\dot{\varphi}}{\partial_z\varphi}\right) + \frac{\dot{f}_3}{\partial_z\varphi}\partial_z g \\
& \sim f_i\partial_i^\varphi g + f_i\partial_i^\varphi \dot{g} + \dot{f}_i\partial_i^\varphi g - g_i\partial_i\dot{\varphi} + \dot{f}_i\frac{\partial_i\varphi\partial_z g}{\partial_z\varphi}\frac{\partial_z\dot{\varphi}}{\partial_z\varphi} \\
&\quad + f_3\partial_3^\varphi g + f_3\partial_3^\varphi \dot{g} + \dot{f}_3\partial_3^\varphi g - f_3\frac{\partial_z g \partial_z\dot{\varphi}}{\partial_z\varphi\partial_z\varphi} \\
& \sim f\cdot\nabla^\varphi g + f\cdot\nabla^\varphi\dot{g} + \dot{f}\cdot\nabla^\varphi g - \frac{\partial_z g}{\partial_z\varphi}(f_i\partial_i^\varphi\dot{\varphi} + f_3\partial_3^\varphi\dot{\varphi}) \\
& = (f\cdot\nabla^\varphi)g + (f\cdot\nabla^\varphi)\dot{g} + (\dot{f}\cdot\nabla^\varphi)g - \partial_z^\varphi g(f\cdot(\nabla^\varphi \dot{\varphi})).
\end{split}
\end{equation*}
\end{proof}

\begin{proposition} \label{proposition 2.2}
We have the following first order expansion.
\begin{equation*}
\partial_i^{\varphi+\dot{\varphi}}|f + \dot{f}|^2 \sim \partial_i^\varphi |f|^2 + 2\dot{f}\cdot(\partial_i^\varphi f) + 2f\cdot(\partial_i^\varphi \dot{f}) - 2\partial_i^\varphi \dot{\varphi}(f\cdot\partial_z^\varphi f)
\end{equation*}
\end{proposition}
\begin{proof}
We abbreviate $\sum_{j=1}^{3}$ for $j$ terms. For $i=1,2$,
\begin{equation*}
\begin{split}
\partial_i^{\varphi+\dot{\varphi}}|f + \dot{f}|^2 & \sim \left( \partial_i - \frac{\partial_i\varphi + \partial_i\dot{\varphi}}{\partial_z\varphi + \partial_z\dot{\varphi}}\partial_z \right)(f_j^2 + 2f_j\dot{f}_j) \\
& \sim 2f_j\partial_i f_j + 2\partial_if_j\dot{f}_j + 2f_j\partial_i \dot{f}_j - 2\frac{(\partial_i\varphi + \partial_i\dot{\varphi})(f_j\partial_z f_j + \partial_z f_j \dot{f}_j + f_j\partial_z \dot{f}_j)}{\partial_z\varphi}\left(1-\frac{\partial_z\dot{\varphi}}{\partial_z\varphi}\right) \\
& \sim \partial_i|f_j|^2 + 2\dot{f}_j\partial_i f_j + 2 f_j\partial_i\dot{f}_j - \frac{2}{\partial_z\varphi}(\partial_i\varphi f_j\partial_z f_j + \partial_i\varphi \partial_z f_j\dot{f}_j + \partial_i\varphi f_j\partial_z f_j + \partial_i\dot{\varphi}f_j\partial_z f_j) \\
&\quad + \frac{2}{\partial_z\varphi}\frac{\partial_z\dot{\varphi}}{\partial_z\varphi}\partial_i\varphi f_j\partial_z f_j \\
& \sim \left(\partial_i |f_j|^2 - \frac{\partial_i\varphi}{\partial_z\varphi}\partial_z |f_j|^2\right) + 2\dot{f}_j\left( \partial_i f_j - \frac{\partial_i\varphi}{\partial_z\varphi}\partial_z f_j \right) + 2f_j\left(\partial_i \dot{f}_j - \frac{\partial_i\varphi}{\partial_z\varphi}\partial_z \dot{f}_j\right) \\
&\quad - 2f_j\frac{\partial_z f_j}{\partial_z\varphi}\left(\partial_i\dot{\varphi} - \frac{\partial_i\varphi}{\partial_z\varphi}\partial_z\dot{\varphi}\right) \\
& = \partial_i^\varphi |f|^2 + 2\dot{f}\cdot(\partial_i^\varphi f) + 2f\cdot(\partial_i^\varphi \dot{f}) - 2\partial_i^\varphi \dot{\varphi}(f\cdot\partial_z^\varphi f).
\end{split}
\end{equation*}
For $i=3$,
\begin{equation*}
\begin{split}
\partial_3^{\varphi+\dot{\varphi}}|f + \dot{f}|^2 & \sim \frac{1}{\partial_z\varphi}\left( 1 - \frac{\partial_z\dot{\varphi}}{\partial_z\varphi}\right)\partial_z(f_j^2 + 2f_j\dot{f}_j) \\
& \sim \frac{2}{\partial_z\varphi}(f_j\partial_z f_j + \partial_z f_j\dot{f}_j + f_j\partial_z\dot{f}_j)\left( 1 - \frac{\partial_z\dot{\varphi}}{\partial_z\varphi}\right) \\
& \sim \frac{2}{\partial_z\varphi}\left(  f_j\partial_z f_j + \partial_z f_j\dot{f}_j + f_j\partial_z\dot{f}_j - \frac{\partial_z\dot{\varphi}}{\partial_z\varphi}f_j\partial_z f_j\right) \\
& = \partial_3^\varphi |f|^2 + 2\dot{f}\cdot(\partial_3^\varphi f) + 2f\cdot(\partial_3^\varphi \dot{f}) - 2\partial_3^\varphi \dot{\varphi}(f\cdot\partial_z^\varphi f).
\end{split}
\end{equation*}
\end{proof}

Using Proposition \ref{proposition 2.1} and \ref{proposition 2.2}, we get linearization for (\ref{2.1}),
\begin{equation} \label{linearization of symbols}
\begin{split}
D\mathcal{N}(v,B,q,\varphi)\cdot(\dot{v},\dot{B},\dot{q},\dot{\varphi}) & = \partial_t^\varphi\dot{v} + (v\cdot\nabla^\varphi)\dot{v} + \nabla^\varphi\dot{q} - 2\varepsilon\nabla^\varphi\cdot(\mathbf{S}^\varphi\dot{v}) - (B\cdot\nabla^\varphi)\dot{B} \\ 
&\quad + (\dot{v}\cdot\nabla^\varphi)v - \partial_z^\varphi v(\partial_t^\varphi\dot{\varphi} + v\cdot\nabla^\varphi\dot{\varphi}) - \partial_z^\varphi q\nabla^\varphi \dot{\varphi} \\
&\quad + 2\varepsilon\nabla^\varphi(\partial_z^\varphi v \otimes \nabla^\varphi\dot{\varphi} + \nabla^\varphi\dot{\varphi}\otimes\partial_z^\varphi v) + 2\varepsilon\partial_z^\varphi(\mathbf{S}^\varphi v)\nabla^\varphi \dot{\varphi} \\
&\quad - \left( \dot{B}\cdot\nabla^\varphi B - \partial_z^\varphi B(B\cdot\nabla^\varphi\dot{\varphi}) \right),   \\
Dd_v (v,\varphi)\cdot(\dot{v},\dot{\varphi}) &= \nabla^\varphi\cdot\dot{v} - \nabla^\varphi\dot{\varphi}\cdot\partial_z^\varphi v,  \\
Dd_B (B,\varphi)\cdot(\dot{B},\dot{\varphi}) &= \nabla^\varphi\cdot\dot{B} - \nabla^\varphi\dot{\varphi}\cdot\partial_z^\varphi B,  \\
D\mathcal{F}(v,B,\varphi)\cdot(\dot{v},\dot{B},\dot{\varphi}) & = \partial_t^\varphi\dot{B} + (v\cdot\nabla^\varphi)\dot{B} - (B\cdot\nabla^\varphi)\dot{v} - 2\lambda\nabla^\varphi\cdot(\mathbf{S}^\varphi\dot{B}) \\
&\quad + (\dot{v}\cdot\nabla^\varphi)B - (\dot{B}\cdot\nabla^\varphi)v - \partial_z^\varphi B(\partial_t^\varphi\dot{\varphi} + v\cdot\nabla^\varphi\dot{\varphi}) + \partial_z^\varphi v(B\cdot\nabla^\varphi\dot{\varphi}) \\
&\quad + 2\lambda\nabla^\varphi(\partial_z^\varphi B \otimes \nabla^\varphi\dot{\varphi} + \nabla^\varphi\dot{\varphi}\otimes\partial_z^\varphi B) + 2\lambda\partial_z^\varphi(\mathbf{S}^\varphi B)\nabla^\varphi \dot{\varphi},  \\
D\mathcal{B}(v,B,q,\varphi)\cdot(\dot{v},\dot{B},\dot{q},\dot{\varphi}) & = 2\varepsilon \mathbf{S}^\varphi\dot{v}\mathbf{N} - \partial_z^\varphi v \otimes \nabla^\varphi\dot{\varphi}\mathbf{N} - \nabla^\varphi\dot{\varphi} \otimes \partial_z^\varphi v\mathbf{N} - (\dot{q}-g\dot{h})\mathbf{N} \\
&\quad + \left(2\varepsilon \mathbf{S}^\varphi v - (q-gh)\right)\dot{\mathbf{N}}.
\end{split}
\end{equation}

On the right hand sides of above linearizations, we see $\nabla^\varphi \varphi$ which behaves like $|\nabla\varphi|_{m} \sim |h|_{m+\frac{1}{2}}$ for high order estimate. Now, we define Alinhac's new unknowns to remove $\nabla^\varphi \varphi$'s on the right hand side. For example, on the right hand side of $\mathcal{N}$, we see that $-(\partial_z^\varphi v) v\cdot\nabla^\varphi \dot{\varphi}$ is one of bad terms. But this term has $v\cdot\nabla^\varphi$ so this can be combined  with $(v\cdot\nabla^\varphi)\dot{v}$. Then $-(\partial_z^\varphi v) v\cdot\nabla^\varphi \dot{\varphi}$ gives 1 derivative $\nabla^\varphi$ to nonlinear structure $(v\cdot\nabla^\varphi)$ and remained $-(\partial_z^\varphi v) \dot{\varphi}$ is combined with $\dot{v}$ to generate a new variable.
\begin{equation} \label{2.2}
\mathcal{V} := \dot{v} - \partial_z^\varphi v\dot{\varphi},\,\,\,\,\mathcal{Q} := \dot{q} - \partial_z^\varphi q\dot{\varphi},\,\,\,\,\mathcal{B} := \dot{B} - \partial_z^\varphi B\dot{\varphi}
\end{equation}

\begin{lemma} \label{lemma 2.3}
Let us define
\begin{equation*}
\mathcal{A}_i(v,\varphi) := \partial_i^\varphi v,\,\,\,\,\mathcal{F}_{ij}(v,\varphi) := \partial_i^\varphi\partial_j^\varphi v,
\end{equation*}
then linearizations of $\mathcal{A}$ and $\mathcal{F}$ can be expressed by 
\begin{equation*}
D\mathcal{A}_i(v,\varphi)\cdot(\dot{v},\dot{\varphi}) = \partial_i^\varphi(\dot{v}-\partial_z^\varphi v\dot{\varphi}) + \dot{\varphi}\partial_z^\varphi(\mathcal{A}_i(v,\varphi)),
\end{equation*}
\begin{equation*}
D\mathcal{F}_{ij}(v,\varphi)\cdot(\dot{v},\dot{\varphi}) = \partial_{ij}^\varphi(\dot{v}-\partial_z^\varphi v\dot{\varphi}) + \dot{\varphi}\partial_z^\varphi(\mathcal{F}_{ij}(v,\varphi)).
\end{equation*}
\end{lemma}
\begin{proof}
This is simple calculations which use commutativity property of $\partial_i^\varphi$. See proof of Lemma 2.7 in \cite{NMFR1} for detail.
\end{proof}

Using Lemma \ref{lemma 2.3}, we have the following proposition. Note that on the right hand side, all bad terms (which behaves like $|h|_{m+\frac{1}{2}}$) are removed.
\begin{proposition} \label{proposition 2.4}
Linearization of (\ref{2.1}) can be expressed as the following, using new unknowns $\mathcal{V},\mathcal{Q}$, and $\mathcal{B}$ in (\ref{2.2}),
\begin{equation*}
\begin{split}
D\mathcal{N}(v,B,q,\varphi)\cdot(\dot{v},\dot{B},\dot{q},\dot{\varphi}) & = (\partial_t^\varphi + (v\cdot\nabla^\varphi) - 2\varepsilon\nabla\cdot(\mathbf{S}^\varphi\cdot))\mathcal{V} + \nabla^\varphi\mathcal{Q} - (B\cdot\nabla^\varphi)\mathcal{B} \\
& + (\dot{v}\cdot\nabla^\varphi)v - (\dot{B}\cdot\nabla^\varphi)B + \dot{\varphi}\{ \partial_z^\varphi\mathcal{N}(v,B,P,\varphi) - (\partial_z^\varphi v\cdot\nabla^\varphi)v + (\partial_z^\varphi B\cdot\nabla^\varphi)B \},  \\
Dd_v (v,\varphi)\cdot(\dot{v},\dot{\varphi}) &= \nabla^\varphi\cdot \mathcal{V} - \dot{\varphi}\partial_z^\varphi(d_v(v,\varphi)),  \\
Dd_B (B,\varphi)\cdot(\dot{B},\dot{\varphi}) &= \nabla^\varphi\cdot \mathcal{B} - \dot{\varphi}\partial_z^\varphi(d_B(B,\varphi)),  \\
D\mathcal{F}(v,B,\varphi)\cdot(\dot{v},\dot{B},\dot{\varphi}) & = (\partial_t^\varphi + (v\cdot\nabla^\varphi) - 2\lambda\nabla\cdot(\mathbf{S}^\varphi\cdot))\mathcal{B} - (B\cdot\nabla^\varphi)\mathcal{V} \\
& + (\dot{v}\cdot\nabla^\varphi)B - (\dot{B}\cdot\nabla^\varphi)v + \dot{\varphi}\left\{\partial_z^\varphi(\mathcal{F}(v,B,\varphi)) - (\partial_z^\varphi v\cdot\nabla^\varphi)B + (\partial_z^\varphi B\cdot\nabla^\varphi)v \right\},  \\
D\mathcal{B}(v,B,q,\varphi)\cdot(\dot{v},\dot{B},\dot{q},\dot{\varphi}) & = 2\varepsilon \mathbf{S}^\varphi\mathcal{V}\mathbf{N} + 2\varepsilon\dot{\varphi}\partial_z(\mathbf{S}^\varphi v)\mathbf{N} - (\dot{q}-g\dot{h})\mathbf{N} + \left(2\varepsilon \mathbf{S}^\varphi - (q-gh)\right)\dot{\mathbf{N}}. \\
\end{split}
\end{equation*}
\end{proposition}
\begin{proof}
We use linearization (\ref{linearization of symbols}) and Lemma \ref{lemma 2.3}.
\end{proof}

\section{Preliminaries estimates}
In this section we collect some necessary propositions and preliminary estimates from \cite{NMFR1}. Every functions are defined in the fixed domain $S$.
\begin{proposition} \label{proposition 3.1}
We have the following products, and commutator estimates. \\
$\bullet$ For $u,v \in L^\infty \cap H^k_{co}$, $k\geq 0$,
\begin{equation*}  \label{3.1}
\|Z^{\alpha_1}u Z^{\alpha_2}v\| \lesssim \|u\|_{L^\infty}\|v\|_k + \|v\|_{L^\infty}\|u\|_k,\,\,\,|\alpha_1|+|\alpha_2|=k.
\end{equation*}
$\bullet$ For $1\leq |\alpha| \leq k$, $g\in H^{k-1}_{co} \cap L^\infty$, $f\in H^k_{co}$ such that $Zf\in L^\infty$, we have
\begin{equation*} \label{3.2}
\|[Z^\alpha,f]g\|\lesssim \|Zf\|_{k-1}\|g\|_{L^\infty} + \|Zf\|_{L^\infty}\|g\|_{k-1}.
\end{equation*}
$\bullet$ For $|\alpha|=k\geq 2$, we define the symmetric commutator $[Z^\alpha,f,g] = Z^\alpha(fg) - (Z^\alpha f)g - fZ^\alpha g$. Then we have the estimate
\begin{equation*} \label{3.3}
\|[Z^\alpha,f,g]\| \lesssim \|Zf\|_{L^\infty}\|Zg\|_{k-2} + \|Zg\|_{L^\infty}\|Zf\|_{k-2}.
\end{equation*}
\end{proposition}

The following proposition states embedding and trace estimate.
\begin{proposition} \label{proposition 3.2}
$\bullet$ For $s_1 \geq 0, s_2 \geq 0$ such that $s_1+s_2 > 2$ and $f$ such that $f\in H^{s_1}_{tan}$, $\partial_z f \in H^{s_2}_{tan}$, we have the anisotropic sobolev embedding. 
\begin{equation*} \label{3.4}
\|f\|^2_{L^\infty} \lesssim \|\partial_z f\|_{H^{s_2}_{tan}} \|f\|_{H^{s_1}_{tan}}.
\end{equation*}
$\bullet$ For $f\in H^1(S)$, we have the trace estimates,
\begin{equation*} \label{3.5}
|f(\cdot,0)|_{H^s(\mathbb{R}^2)} \leq C\|\partial_z f\|^{1/2}_{H^{s_2}_{tan}} \|f\|^{1/2}_{H^{s_1}_{tan}},
\end{equation*}
with $s_1+s_2 = 2s \geq 0$.
\end{proposition}

We have similar estimates on boundary $\p S=\mathbb{R}^2$.
\begin{proposition} \label{proposition 3.3}
	When $f,g$ are defined in $\mathbb{R}^2$, we have the following commutator estimates.
\begin{equation*} \label{3.6}
\begin{split}
|\Lambda^s(fg)|_{L^2(\mathbb{R}^2)} &\leq C_s\left( |f|_{L^\infty(\mathbb{R}^2)}|g|_{H^s(\mathbb{R}^2)} + |g|_{L^\infty(\mathbb{R}^2)}|f|_{H^s(\mathbb{R}^2)} \right),  \\
\|[\Lambda^s,f]\nabla g\|_{L^2(\mathbb{R}^2)} &\leq C_s\left( |\nabla f|_{L^\infty(\mathbb{R}^2)}|g|_{H^s(\mathbb{R}^2)} + |\nabla g|_{L^\infty(\mathbb{R}^2)}|f|_{H^s(\mathbb{R}^2)} \right),  \\
|uv|_{\frac{1}{2}} &\lesssim |u|_{1,\infty}|v|_{\frac{1}{2}}.
\end{split}
\end{equation*}
\end{proposition}

From (\ref{Phi}), Jacobian of change of variable $\Phi$ is $\partial_z \varphi$. Let us define volume element $dV_t$ by
$$
\int_{\Omega(t)} F dydz = \int_S f (\partial_z\varphi) dydz := \int_S f dV_t,\quad\text{where}\quad F(t,\Phi(t,y,z)) = f(t,y,z).
$$
Now, we state integration by part for $\int_S \partial_i^\varphi f g  dV_t$.
\begin{proposition} \label{proposition 3.4}
In $S$, we have the following integration by parts rules.
\begin{equation*} \label{3.9}
\begin{split}
\int_S \partial_i^\varphi f g dV_t &= -\int_S f\partial_i^\varphi g dV_t + \int_{\partial S} fg \mathbf{N}_i dy,\,\,\,i=1,2,3,  \\
\int_S \partial_t^\varphi f g dV_t &= \partial_t\int_S fg dV_t -\int_S f\partial_t^\varphi g dV_t - \int_{\partial S} fg \partial_t h.
\end{split}
\end{equation*}
\end{proposition}

Using Proposition \ref{proposition 3.4}, we gain following Corollary which is useful to get $L^{2}$-type energy estimate for two main PDE's in (\ref{1.4}).
\begin{corollary} \label{corollary 3.5}
Let $v$ be a vector field such that $\nabla^\varphi\cdot v = 0$. For every smooth function $f,g$ and smooth vector field $u,w$, we have
\begin{equation} \label{3.11}
\begin{split}
\int_S (\partial_t f + v\cdot\nabla^\varphi f)f dV_t &= \frac{1}{2}\partial_t \int_S |f|^2 dV_t - \frac{1}{2}\int_{\partial S} |f|^2(\partial_t h - v\cdot\mathbf{N})dy,  \\
\int_S \triangle^\varphi f g dV_t &= -\int_S \nabla^\varphi f\cdot \nabla^\varphi g dV_t + \int_{\partial S} \nabla^\varphi f\cdot\mathbf{N} g dy,  \\
\int_S \nabla^\varphi\cdot(\mathbf{S}^\varphi u)\cdot w dV_t &= -\int_S\mathbf{S}^\varphi u \cdot \mathbf{S}^\varphi w dV_t + \int_{\partial S}(\mathbf{S}^\varphi u \mathbf{N})\cdot w dy.
\end{split}
\end{equation}
\end{corollary}

The following propositions is about Adapted Korn's inequality in $S$.
\begin{proposition} \label{proposition 3.6}
Let $\partial_z\varphi \geq c_0$, $\|\nabla\varphi\|_{L^\infty} + \|\nabla^2\varphi\|_{L^\infty} \leq \frac{1}{c_0}$ for some $c_0 > 0$, then there exists $\Lambda_0 = \Lambda(1/c_0)>0$, such that for every $v\in H^{1}(S)$, we have
\begin{equation} \label{3.14}
\|\nabla v\|^2_{L^2(S)} \leq \Lambda_0 \left( \int_S |\mathbf{S}^\varphi v|^2 dV_t + \|v\|^2_{L^2(S)} \right).
\end{equation}
\end{proposition}

As explained before, gain of regularity for $\eta$ is very important in this paper. We give estimates for $\eta$.
\begin{proposition}  \label{proposition 3.7}
We have the following estimates for $\eta$.
\begin{equation} \label{3.15}
\begin{split}
\|\nabla\eta(t)\|_{H^s(S)} &\leq C_s |h(t)|_{s+\frac{1}{2}},\quad\forall s \geq 0,  \\
\|\nabla\partial_t\eta(t)\|_{H^s(S)} &\leq C_s \left( 1 + \|v\|_{L^\infty} + |\nabla_{y} h|_{L^\infty} \right)\left( \|v\|_{E^{s+1}} + |\nabla_{y} h|_{s+\frac{1}{2}} \right),\quad \forall s \in \mathbb{N},  \\
\|\eta\|_{W^{s,\infty}} &\leq C_s |h|_{s,\infty},\quad \forall s \in \mathbb{N},  \\
\|\partial_t \eta\|_{W^{s,\infty}} &\leq C_s \left( 1 + |h|_{s,\infty} \right) \|v\|_{s,\infty}, \forall s \in \mathbb{N}.
\end{split}
\end{equation}
\end{proposition}

To treat fraction terms, the following proposition is very useful.
\begin{proposition} \label{proposition 3.8}
For every $m\in\mathbb{N}$, we have,
\begin{equation} \label{3.19}
\Big\| \frac{f}{\partial_z\varphi} \Big\|_m \leq \Lambda(\frac{1}{c_0},|h|_{1,\infty} + \|f\|_{L^\infty})\left(|h|_{m+\frac{1}{2}} + \|f\|_m \right).
\end{equation}
\end{proposition}

With vanishing factor $\varepsilon$, $h$ gets one half more regularity. 
\begin{proposition} \label{proposition 3.9}
For every $m\in\mathbb{N}$, $\varepsilon \in (0,1)$, we have the estimate,
\begin{equation}
\varepsilon|h(t)|^2_{m+\frac{1}{2}} \leq \varepsilon|h_0|^2_{m+\frac{1}{2}} + \varepsilon\int_0^t|v^b|^2_{m+\frac{1}{2}} + \int_0^t \Lambda(|\nabla_{y} h|_{L^\infty(\mathbb{R}^2)} + \|v\|_{1,\infty})\left( \|v\|^2_m + \varepsilon|h|^2_{m+\frac{1}{2}} \right) d\tau,
\end{equation}
where $v^{b} := v\vert_{z=0}$.
\end{proposition}

\section{High order equations}
In this section, we apply high order conormal derivatives, $Z^\alpha$ with $|\alpha|=m$, to (\ref{1.4}) and (\ref{1.5}) and rewrite the high order system in terms of new unknowns
\begin{equation} \label{high new unknown}
(\mathcal{V}^\alpha,\mathcal{Q}^\alpha,\mathcal{B}^\alpha) := ( Z^\alpha v - \partial_z^\varphi v Z^\alpha \eta, Z^\alpha q - \partial_z^\varphi q Z^\alpha \eta,Z^\alpha B - \partial_z^\varphi B Z^\alpha \eta ),
\end{equation}
using linearization results. 

\subsection{A Commutator estimate}
For $i=1,2$ we write,
\begin{equation*}
Z^\alpha \partial_i^\varphi f = \partial_i^\varphi Z^\alpha f - \partial_z^\varphi f \partial_i^\varphi Z^\alpha\eta + \mathcal{C}_{i}^\alpha(f),
\end{equation*}
\begin{equation*}
\mathcal{C}_{i}^\alpha(f) := \mathcal{C}_{i,1}^\alpha(f) + \mathcal{C}_{i,2}^\alpha(f) + \mathcal{C}_{i,3}^\alpha(f),
\end{equation*}
where
\begin{equation} \label{comm symbols}
\begin{cases}
\mathcal{C}_{i,1}^\alpha := -[Z^\alpha,\frac{\partial_i\varphi}{\partial_z \varphi},\partial_z f], \\
\mathcal{C}_{i,2}^\alpha := -\partial_z f[Z^\alpha,\partial_i\varphi,\frac{1}{\partial_z\varphi}] - \partial_i\varphi\left( Z^\alpha\left(\frac{1}{\partial_z\varphi}\right) + \frac{Z^\alpha\partial_z\eta}{(\partial_z\varphi)^2}\right)\partial_z f, \\
\mathcal{C}_{i,3}^\alpha := -\frac{\partial_i\varphi}{\partial_z\varphi}[Z^\alpha,\partial_z]f + \frac{\partial_i\varphi}{(\partial_z\varphi)^2}\partial_z f[Z^\alpha,\partial_z]\eta. \\
\end{cases}
\end{equation}
For $i=3$, result is very similar and we suffice to replace $\partial_i\varphi$ by $1$ in above terms. We need to estimate commutators.
\begin{lemma} \label{lemma 4.1}
For $1\leq |\alpha| \leq m$, $i=1,2,3$, we have
$$
\|\mathcal{C}_{i}^\alpha(f)\| \leq \Lambda \left(\frac{1}{c_0},|h|_{2,\infty} + \|\nabla f\|_{1,\infty} \right)\left( \|\nabla f\|_{m-1} + |h|_{m-\frac{1}{2}} \right).
$$
\end{lemma}
\begin{proof}
See Lemma 5.1 in \cite{NMFR1}
\end{proof}

\subsection{Divergence free condition for $\mathbf(v,B)$}
By applying $Z^{\alpha}$, we have,
$$
Z^\alpha (\nabla^\varphi\cdot v) = 0.
$$
Using notations in (\ref{comm symbols}),
\begin{equation*}
\begin{split}
&\nabla^\varphi\cdot(Z^\alpha v) - \partial_z^\varphi v\cdot\nabla^\varphi(Z^\alpha \varphi) + \sum_{i=1}^3 \mathcal{C}_i^\alpha(v_i) = 0,  \\
&\nabla^\varphi\cdot(Z^\alpha v - \partial_z^\varphi v Z^\alpha\varphi) - (\nabla^\varphi\cdot\partial_z^\varphi v)Z^\alpha\varphi + \sum_{i=1}^3 \mathcal{C}_i^\alpha(v_i) = 0.
\end{split}
\end{equation*}
Second term is zero since $\partial_i^\varphi$'s commute and therefore, $\nabla^\varphi\cdot v = 0$. Hence we get
\begin{equation} \label{4.1}
\begin{split}
&\nabla^\varphi\cdot\mathcal{V}^\alpha + \mathcal{C}^\alpha(d_v) = 0,\,\,\,\,\mathcal{C}^\alpha(d_v) := \sum_{i=1}^3 \mathcal{C}_i^\alpha(v_i),  \\
&\nabla^\varphi\cdot\mathcal{B}^\alpha + \mathcal{C}^\alpha(d_B) = 0,\,\,\,\,\mathcal{C}^\alpha(d_B) := \sum_{i=1}^3 \mathcal{C}_i^\alpha(B_i),
\end{split}
\end{equation}
where commutators $\mathcal{C}^\alpha(d_v)$ and $\mathcal{C}^\alpha(d_B)$ satisfy the following estimates by Lemma \ref{lemma 4.1}.
\begin{equation} \label{4.3}
\begin{split}
&\|\mathcal{C}^\alpha(d_v)\| \leq \Lambda \left(\frac{1}{c_0},|h|_{2,\infty} + \|\nabla v\|_{1,\infty} \right)\left( \|\nabla v\|_{m-1} + |h|_{m-\frac{1}{2}} \right),  \\
&\|\mathcal{C}^\alpha(d_B)\| \leq \Lambda \left(\frac{1}{c_0},|h|_{2,\infty} + \|\nabla B\|_{1,\infty} \right)\left( \|\nabla B\|_{m-1} + |h|_{m-\frac{1}{2}} \right).
\end{split}
\end{equation}

\subsection{Navier-Stokes Equation with Lorentz force}
Applying high order conormal derivatives,
$$
Z^\alpha\left\{ \partial_z^\varphi v + (v\cdot\nabla^\varphi)v + \nabla^\varphi q - 2\varepsilon\nabla^\varphi\cdot(\mathbf{S}^\varphi v) - (B\cdot\nabla^\varphi)B \right\} = 0.
$$
$\textit{Transport}$ \ 
Let us use notation $V_z = \frac{1}{\partial_z\varphi}(v\cdot\mathbf{N}^\varphi -  \partial_t\eta)$ and $\mathbf{N}^\varphi := (-\nabla_y\eta,1)$. Then,
\begin{equation} \label{4.5}
\begin{split}
Z^\alpha\left(\partial_t^\varphi + (v\cdot\nabla^\varphi) \right)v &= Z^\alpha(\partial_t + v_y\cdot\nabla_y v + V_z\partial_z) v,  \\
&= (\partial_t + v_y\cdot\nabla_y + V_z\partial_z)Z^\alpha v + (v\cdot Z^\alpha\mathbf{N}^\varphi - \partial_t Z^\alpha\eta)\partial_z^\varphi v - \partial_z^\varphi Z^\alpha\eta(v\cdot\mathbf{N}^\varphi - \partial_t\eta)\partial_z^\varphi v + \mathcal{C}^\alpha(\mathcal{T}_v)  \\
&= (\partial_t^\varphi + v\cdot\nabla^\varphi)(Z^\alpha v) - \partial_z^\varphi v(\partial_t^\varphi + v\cdot\nabla^\varphi)(Z^\alpha\varphi) + \mathcal{C}^\alpha(\mathcal{T}_v),
\end{split}
\end{equation}
where $\mathcal{C}^\alpha(\mathcal{T}_v)$ is defined by
\begin{equation*}
\mathcal{C}^\alpha(\mathcal{T}_v) := \sum_{i=1}^6 \mathcal{T}_i^\alpha.
\end{equation*}
Each terms are given by,
\begin{equation*}
\begin{split}
\mathcal{T}_1^\alpha &= [Z^\alpha, v_y,\partial_y v],\,\,\,\,\mathcal{T}_2^\alpha = [Z^\alpha, V_z,\partial_z v],\,\,\,\,\mathcal{T}_3^\alpha = \frac{1}{\partial_z\varphi}[Z^\alpha,v_z]\partial_z v,  \\
\mathcal{T}_4^\alpha &= \left(Z^\alpha\left(\frac{1}{\partial_z\varphi}\right) + \frac{\partial_z Z^\alpha\eta}{(\partial_z\varphi)^2}\right)v_z\partial_z v,\,\,\,\,\mathcal{T}_5^\alpha = v_z\partial_z v\frac{[Z^\alpha,\partial_z]\eta}{(\partial_z\varphi)^2} + V_z[Z^\alpha,\partial_z]v,  \\
\mathcal{T}_6^\alpha &= [Z^\alpha,v_z,\frac{1}{\partial_z\varphi}]\partial_z v.
\end{split}
\end{equation*}
Estimate for $\mathcal{C}^\alpha(\mathcal{T}_v)$ is given as following using Propositions in section 3. 
\begin{equation} \label{4.6}
\|\mathcal{C}^\alpha(\mathcal{T}_v)\| \leq \Lambda\left(\frac{1}{c_0},|h|_{2,\infty} + \|v\|_{E^{2,\infty}}\right)\left(\|v\|_{E^m} + |h|_{m-\frac{1}{2}}\right).
\end{equation}
$\textit{Pressure}$ \ For pressure term,
\begin{equation} \label{4.7}
Z^\alpha\nabla^\varphi q = \nabla^\varphi(Z^\alpha q) - \partial_z^\varphi q \nabla^\varphi(Z^\alpha \varphi) + \mathcal{C}^\alpha(q),
\end{equation}
with estimate
\begin{equation} \label{4.8}
\|\mathcal{C}^\alpha(q)\| \leq \Lambda \left(\frac{1}{c_0},|h|_{2,\infty} + \|\nabla q\|_{1,\infty} \right)\left( \|\nabla q\|_{m-1} + |h|_{m-\frac{1}{2}} \right).
\end{equation}
$\textit{Diffusion}$ \ For diffusion term,
\begin{equation} \label{4.9}
\begin{split}
Z^\alpha\left(-2\varepsilon\nabla^\varphi\cdot(\mathbf{S}^\varphi v)\right) &= -2\varepsilon\nabla^\varphi\cdot(\mathbf{S}^\varphi Z^\alpha v) + 2\varepsilon\nabla^\varphi\cdot\left( \partial_z^\varphi v \otimes \nabla^\varphi Z^\alpha\varphi + \nabla^\varphi Z^\alpha\varphi \otimes \partial_z^\varphi v \right)  \\
&\quad + 2\varepsilon\partial_z^\varphi(\mathbf{S}^\varphi v)\nabla^\varphi(Z^\alpha\varphi) - \varepsilon\mathcal{D}^\alpha(\mathbf{S}^\varphi v) - \varepsilon\nabla^\varphi\cdot(\mathcal{E}^\alpha v),
\end{split}
\end{equation}
where $\mathcal{D}^\alpha(\mathbf{S}^\varphi v)$ and $(\mathcal{E}^\alpha v)$ are defined by
\[
(\mathcal{E}^\alpha v)_{ij} := \mathcal{C}_i^\alpha(v_j) + \mathcal{C}_j^\alpha(v_i),\,\,\,\,\mathcal{D}^\alpha(\mathbf{S}^\varphi v)_i := 2\mathcal{C}_j^\alpha (\mathbf{S}^\varphi v)_{ij},
\]
with estimate for $\mathcal{E}^\alpha (v)$,
\begin{equation} \label{4.10}
\|\mathcal{E}^\alpha (v)\| \leq \Lambda\left(\frac{1}{c_0},|h|_{2,\infty} + \|\nabla v\|_{1,\infty}\right)\left(\|v\|_{m} + \|\partial_z v\|_{m-1} + |h|_{m-\frac{1}{2}}\right).
\end{equation}
$\textit{Lorentz\,\,force}$ \ For Lorentz force term,
\begin{equation} \label{4.11}
\begin{split}
-Z^\alpha(B\cdot\nabla^\varphi B) &= -Z^\alpha\left(\sum_{i=1}^3 B_i\partial_i^\varphi B\right) = -\sum_{i=1}^3\left(Z^\alpha B_i\partial_i^\varphi B + B_i Z^\alpha\partial_i^\varphi B + [Z^\alpha,B_i,\partial_i^\varphi B] \right)  \\
&= -(Z^\alpha B)\cdot\nabla^\varphi B - \sum_{i=1}^3 B_i\left( \partial_i^\varphi (Z^\alpha B) - \partial_z^\varphi B\partial_i^\varphi(Z^\alpha\varphi) + \mathcal{C}_i^\alpha(B) \right)
- \sum_{i=1}^3[Z^\alpha,B_i,\partial_i^\varphi B]   \\
&= -(Z^\alpha B)\cdot\nabla^\varphi B - (B\cdot\nabla^\varphi)(Z^\alpha B) + (\partial_z^\varphi B)(B\cdot\nabla^\varphi(Z^\alpha\varphi)) - \sum_{i=1}^3 B_i\mathcal{C}_i^\alpha(B) - \sum_{i=1}^3[Z^\alpha,B_i,\partial_i^\varphi B].
\end{split}
\end{equation}
Let us define,
\begin{equation*}
\begin{split}
\mathcal{C}^\alpha(\mathcal{T}_B) & := \sum_{i=1}^3[Z^\alpha,B_i,\partial_i^\varphi B] = \sum_{i=1}^3[Z^\alpha,B_i,\partial_i B] - \sum_{i=1}^3[Z^\alpha,B_i,\frac{\partial_i\varphi}{\partial_z\varphi}\partial_z B].
\end{split}
\end{equation*}
Using Lemma \ref{lemma 4.1}, we have an estimate for $\mathcal{C}^\alpha(\mathcal{T}_{B})$,
\begin{equation} \label{4.12}
\|\mathcal{C}^\alpha(\mathcal{T}_B)\| \leq \Lambda\left(\frac{1}{c_0},|h|_{2,\infty} + \|\nabla B\|_{1,\infty} \right)\left( \|\nabla B\|_{m-1} + |h|_{m-\frac{1}{2}} \right).
\end{equation}
Now putting (\ref{4.5}), (\ref{4.7}), (\ref{4.9}), and (\ref{4.11}) together, and using linearization of $\mathcal{N}$, to get
\begin{equation*}
\begin{split}
0 &= D\mathcal{N}(v,B,q,\varphi)\cdot(Z^\alpha v,Z^\alpha B, Z^\alpha q, Z^\alpha\varphi) - (Z^\alpha v\cdot\nabla^\varphi) v  \\
&+ \mathcal{C}^\alpha(\mathcal{T}_v) + \mathcal{C}^\alpha(q) - \varepsilon\mathcal{D}^\alpha(\mathbf{S}^\varphi v) - \varepsilon\nabla^\varphi\cdot(\mathcal{E}^\alpha v) - \sum_{i=1}^3 B_i\mathcal{C}_i^\alpha(B) - \mathcal{C}^\alpha(\mathcal{T}_B).
\end{split}
\end{equation*}
By Proposition \ref{proposition 2.4} for $D\mathcal{N}(v,B,q,\varphi)\cdot(Z^\alpha v,Z^\alpha B, Z^\alpha q, Z^\alpha\varphi)$, we get the following.
\begin{equation} \label{4.13}
\begin{split}
&(\partial_t^\varphi + v\cdot\nabla^\varphi - 2\varepsilon\nabla^\varphi\cdot(\mathbf{S}^\varphi\cdot))\mathcal{V}^\alpha + \nabla^\varphi\mathcal{Q}^\alpha - (B\cdot\nabla^\varphi)\mathcal{B}^\alpha  \\
&= (Z^\alpha B\cdot\nabla^\varphi)B + Z^\alpha\varphi\left((\partial_z^\varphi v\cdot\nabla^\varphi)v - (\partial_z^\varphi B\cdot\nabla^\varphi)B \right) + \sum_{i=1}^3 B_i\mathcal{C}_i^\alpha(B) + \mathcal{C}^\alpha(\mathcal{T}_B)  \\
&+ \varepsilon\mathcal{D}^\alpha(\mathbf{S}^\varphi v) + \varepsilon\nabla^\varphi\cdot(\mathcal{E}^\alpha v) - \mathcal{C}^\alpha(\mathcal{T}_v) - \mathcal{C}^\alpha(q) .
\end{split}
\end{equation}

\subsection{Faraday law} Similar as above, we apply $Z^{\alpha}$ and perform computation for each terms.
$$
Z^\alpha\left(\partial_t^\varphi B + (v\cdot\nabla^\varphi)B - (B\cdot\nabla^\varphi)v - 2\lambda\nabla^\varphi\cdot(\mathbf{S}^\varphi B) \right) = 0
$$
$\textit{Transport}$ \ When $V_z := \frac{v\cdot \mathbf{N}^\varphi - \partial_t\varphi}{\partial_z\varphi}$,
\begin{equation*}
\begin{split}
&\quad  Z^\alpha(\partial_t^\varphi + v\cdot\nabla^\varphi)B \\
&= Z^\alpha(\partial_t + v_y\cdot\nabla_y + V_z\partial_z)B,  \\
&= \partial_t(Z^\alpha B) + \sum_{i=1}^2 \left(Z^\varphi v_i\partial_i B + v_i\partial_i Z^\alpha B + [Z^\alpha,v_i,\partial_i B]\right) + Z^\alpha V_z\cdot\partial_z B + V_z Z^\alpha\partial_z B + [Z^\alpha,V_z,\partial_z B]  \\
&= (\partial_t + v_y\cdot\nabla_y)(Z^\alpha B) + Z^\alpha v_y\cdot\nabla_y B + \sum_{i=1}^2[Z^\alpha,v_i,\partial_i B] + Z^\alpha\left(\frac{v\cdot \mathbf{N}^\varphi - \partial_t\varphi}{\partial_z\varphi}\right)\cdot \partial_z B + V_z\partial_z Z^\alpha B \\
&\quad + V_z[Z^\alpha,\partial_z]B + [Z^\alpha,V_z,\partial_z B]  \\
&= (\partial_t^\varphi + v_y\cdot\nabla_y + V_z\partial_z)(Z^\alpha B) + Z^\alpha\left(\frac{v\cdot \mathbf{N}^\varphi - \partial_t\varphi}{\partial_z\varphi}\right)\cdot\partial_z B + \mathcal{R}_1	\\
&= (\partial_t^\varphi + v\cdot\nabla^\varphi)(Z^\alpha B) + \partial_z B\left(\frac{1}{\partial_z\varphi}v\cdot Z^\alpha \mathbf{N}^\varphi + Z^\alpha\left(\frac{1}{\partial_z\varphi}\right)v\cdot \mathbf{N}^\varphi - \frac{1}{\partial_z\varphi}\partial_t Z^\alpha\varphi - \partial_t\varphi Z^\alpha\left(\frac{1}{\partial_z\varphi}\right) \right)  \\
&= (\partial_t^\varphi + v\cdot\nabla^\varphi)(Z^\alpha B) + \partial_z B\left( \frac{1}{\partial_z\varphi}Z^\alpha v\cdot \mathbf{N}^\varphi +  \sum_{i=1}^3[Z^\alpha,\frac{1}{\partial_z\varphi},v_i,\mathbf{N}^\varphi_i] - [Z^\alpha,\partial_t\varphi,\frac{1}{\partial_z\varphi}] \right) + \mathcal{R}_1,
\end{split}
\end{equation*}
where
\begin{equation} \label{mathcal R1}
\mathcal{R}_1 := V_z[Z^\alpha,\partial_z]B + Z^\alpha v_y\cdot\nabla_y B + \sum_{i=1}^2[Z^\alpha,v_i,\partial_i B] + [Z^\alpha,V_z,\partial_z B].
\end{equation}
Note that
\begin{equation*}
\begin{split}
Z^\alpha\left(\frac{1}{\partial_z\varphi}\right) &= Z^{\alpha-1}\left(-\frac{Z\partial_z\varphi}{(\partial_z\varphi)^2}\right)   \\
&= -\frac{1}{(\partial_z\varphi)^2}\left(\partial_z Z^\alpha\varphi + [Z^\alpha,\partial_z]\varphi\right) - Z\partial_z\varphi Z^{\alpha-1}\left(\frac{1}{\partial_z\varphi}\right)^2 - [Z^{\alpha-1},\frac{1}{(\partial_z\varphi)^2}, Z\partial_z\varphi].
\end{split}
\end{equation*}
Hence, transport part becomes
\begin{equation*}
\begin{split}
& Z^\alpha(\partial_t^\varphi + v\cdot\nabla^\varphi)B \\
&= (\partial_t^\varphi + v\cdot\nabla^\varphi)(Z^\alpha B) + (v\cdot Z^\alpha \mathbf{N}^\varphi - \partial_t Z^\alpha\varphi)\partial_z^\varphi B + \partial_z B Z^\alpha\left(\frac{1}{\partial_z\varphi}\right) (v\cdot \mathbf{N}^\varphi - \partial_t\varphi)  \\
&\quad + \partial_z B\left(\frac{1}{\partial_z\varphi}Z^\alpha v\cdot \mathbf{N}^\varphi + \sum_{i=1}^3[Z^\alpha,\frac{1}{\partial_z\varphi},v_i,\mathbf{N}^\varphi_i] - [Z^\alpha,\partial_t\varphi,\frac{1}{\partial_z\varphi}] \right) + \mathcal{R}_1  \\
&= (\partial_t^\varphi + v\cdot\nabla^\varphi)(Z^\alpha B) + (v\cdot Z^\alpha N - \partial_t Z^\alpha\varphi)\partial_z^\varphi B   \\
&\quad + \partial_z B(v\cdot \mathbf{N}^\varphi - \partial_t\varphi)\left( -\frac{1}{(\partial_z\varphi)^2}\left(\partial_z Z^\alpha\varphi + [Z^\alpha,\partial_z]\varphi\right) - Z\partial_z\varphi Z^{\alpha-1}\left(\frac{1}{\partial_z\varphi}\right)^2 - [Z^{\alpha-1},\frac{1}{(\partial_z\varphi)^2}, Z\partial_z\varphi] \right) + \tilde{\mathcal{R}}_1,
\end{split}
\end{equation*}
where
$$
\tilde{\mathcal{R}}_1 := \partial_z B\left(\frac{1}{\partial_z\varphi}Z^\alpha v\cdot \mathbf{N}^\varphi + \sum_{i=1}^3[Z^\alpha,\frac{1}{\partial_z\varphi},v_i,\mathbf{N}^\varphi_i] - [Z^\alpha,\partial_t\varphi,\frac{1}{\partial_z\varphi}] \right) + \mathcal{R}_1.
$$
So,
$$
Z^\alpha(\partial_t^\varphi + v\cdot\nabla^\varphi)B = (\partial_t^\varphi + v\cdot\nabla^\varphi)(Z^\alpha B) + (v\cdot Z^\alpha \mathbf{N}^\varphi - \partial_t Z^\alpha\varphi)\partial_z^\varphi B  - \partial_z^\varphi Z^\alpha\varphi(v\cdot \mathbf{N}^\varphi - \partial_t\varphi)\partial_z^\varphi B + \mathcal{C}^\alpha(\mathcal{T}_F).
$$
and at result
\begin{equation} \label{4.14}
Z^\alpha(\partial^\varphi + v\cdot\nabla^\varphi)B = (\partial_t^\varphi + v\cdot\nabla^\varphi)Z^\alpha B - \partial_z^\varphi B(\partial^\varphi + v\cdot\nabla^\varphi)Z^\alpha\varphi + \mathcal{C}^\alpha(\mathcal{T}_F),
\end{equation}
where
\begin{equation*}
\begin{split}
\mathcal{C}^\alpha(\mathcal{T}_F) &= \partial_z B(v\cdot \mathbf{N}^\varphi - \partial_t\varphi)\left( -\frac{[Z^\alpha,\partial_z]\varphi}{(\partial_z\varphi)^2} - Z\partial_z\varphi Z^{\alpha-1}\left(\frac{1}{\partial_z\varphi}\right)^2 - [Z^{\alpha-1},\frac{1}{(\partial_z\varphi)^2}, Z\partial_z\varphi] \right)  \\
&\quad + \partial_z B\left(\frac{1}{\partial_z\varphi}Z^\alpha v\cdot \mathbf{N}^\varphi + \sum_{i=1}^3[Z^\alpha,\frac{1}{\partial_z\varphi},v_i,\mathbf{N}^\varphi_i] - [Z^\alpha,\partial_t\varphi,\frac{1}{\partial_z\varphi}] \right)   \\
&\quad + V_z[Z^\alpha,\partial_z]B + Z^\alpha v_y\cdot\nabla_y B + \sum_{i=1}^2[Z^\alpha,v_i,\partial_i B] + [Z^\alpha,V_z,\partial_z B].
\end{split}
\end{equation*}

\noindent Using Propositions in section 4 and Lemma \ref{lemma 4.1} we get,
\begin{equation} \label{4.15}
\left\| \mathcal{C}^\alpha(\mathcal{T}_F) \right\| \leq \Lambda\left(\frac{1}{c_0},|h|_{2,\infty} + \|v\|_{1,\infty} + \|B\|_{E^{2,\infty}}\right)\left( |h|_{m-\frac{1}{2}} + \|B\|_{E^m} + \|v\|_{m} \right).
\end{equation}
$\textit{Forcing term} \ \mathit{(B\cdot\nabla^\varphi)v}$  
\begin{equation} \label{4.16}
\begin{split}
-Z^\alpha(B\cdot\nabla^\varphi v) &= -Z^\alpha\left(\sum_{i=1}^3 B_i\partial_i^\varphi v\right) = -\sum_{i=1}^3\left(Z^\alpha B_i\partial_i^\varphi v + B_i Z^\alpha\partial_i^\varphi v + [Z^\alpha,B_i,\partial_i^\varphi v] \right)   \\
&= -(Z^\alpha B)\cdot\nabla^\varphi v - \sum_{i=1}^3 B_i\left( \partial_i^\varphi (Z^\alpha v) - \partial_z^\varphi v\partial_i^\varphi(Z^\alpha\varphi) + \mathcal{C}_i^\alpha(v) \right)
- \sum_{i=1}^3[Z^\alpha, B_i,\partial_i^\varphi v]  \\
&= -(Z^\alpha B)\cdot\nabla^\varphi v - (B\cdot\nabla^\varphi)(Z^\alpha v) + (\partial_z^\varphi v)(B\cdot\nabla^\varphi(Z^\alpha\varphi)) - \sum_{i=1}^3 B_i\mathcal{C}_i^\alpha(v) - \mathcal{C}^\alpha(\mathcal{T}_I),  \\
\end{split}
\end{equation}
where
\begin{equation*}
\begin{split}
\mathcal{C}^\alpha(\mathcal{T}_I) & := \sum_{i=1}^3[Z^\alpha,B_i,\partial_i^\varphi v] = \sum_{i=1}^3[Z^\alpha,B_i,\partial_i v] - \sum_{i=1}^3[Z^\alpha,B_i,\frac{\partial_i\varphi}{\partial_z\varphi}\partial_z v].
\end{split}
\end{equation*}
Then using Lemma \ref{lemma 4.1}, we have an estimate for $\mathcal{C}^\alpha(\mathcal{T}_{I})$.
\begin{equation} \label{4.17}
\|\mathcal{C}^\alpha(\mathcal{T}_I)\| \leq \Lambda\left(\frac{1}{c_0},|h|_{2,\infty} + \|\nabla v\|_{1,\infty} + \|\nabla B\|_{1,\infty} \right)\left( \|\nabla v\|_{m-1} + \|\nabla B\|_{m-1} + |h|_{m-\frac{1}{2}} \right).
\end{equation}
$\textit{Diffusion}$ \ Diffusion part is same as Navier-Stokes part.
\begin{equation} \label{4.18}
\begin{split}
Z^\alpha\left(-2\lambda\nabla^\varphi\cdot(\mathbf{S}^\varphi B)\right) &= -2\lambda\nabla^\varphi\cdot(\mathbf{S}^\varphi Z^\alpha B) + 2\lambda\nabla^\varphi\cdot\left( \partial_z^\varphi B \otimes \nabla^\varphi Z^\alpha\varphi + \nabla^\varphi Z^\alpha\varphi \otimes \partial_z^\varphi B \right)  \\
&\quad + 2\lambda\partial_z^\varphi(\mathbf{S}^\varphi B)\nabla^\varphi(Z^\alpha\varphi) - \varepsilon\mathcal{D}^\alpha(\mathbf{S}^\varphi B) - \varepsilon\nabla^\varphi\cdot(\mathcal{E}^\alpha B),
\end{split}
\end{equation}
where definitions and estimates of $\mathcal{D}^\alpha(\mathbf{S}^\varphi B)$ and $(\mathcal{E}^\alpha B)$ are sams as before.
$$
(\mathcal{E}^\alpha B)_{ij} := \mathcal{C}_i^\alpha(B_j) + \mathcal{C}_j^\alpha(B_i),\,\,\,\,\mathcal{D}^\alpha(\mathbf{S}^\varphi B)_i := 2\mathcal{C}_j^\alpha (\mathbf{S}^\varphi B)_{ij},
$$
with estimate for $\mathcal{E}^\alpha (B)$,
\begin{equation} \label{4.19}
\|\mathcal{E}^\alpha (B)\| \leq \Lambda\left(\frac{1}{c_0},|h|_{2,\infty} + \|\nabla B\|_{1,\infty}\right)\left(\|v\|_{m} + \|\partial_z B\|_{m-1} + |h|_{m-\frac{1}{2}}\right).
\end{equation}
Combining (\ref{4.14}), (\ref{4.16}), and (\ref{4.18}) together, we get
\begin{equation*}
\begin{split}
& D\mathcal{F}(v,B,\varphi)\cdot(Z^\alpha v, Z^\alpha B, Z^\alpha\varphi) - (Z^\alpha v\cdot\nabla^\varphi)B \\
&= -\mathcal{C}^\alpha(\mathcal{T}_F) + \sum_{i=1}^3 B_i\mathcal{C}_i^\alpha(v) + \mathcal{C}^\alpha(\mathcal{T}_I) + \lambda\mathcal{D}^\alpha(\mathbf{S}^\varphi B) + \lambda\nabla^\varphi\cdot(\mathcal{E}^\alpha B).
\end{split}
\end{equation*}

\noindent Using Proposition \ref{proposition 2.4},
\begin{equation} \label{4.20}
\begin{split}
&(\partial_t^\varphi + v\cdot\nabla^\varphi - 2\lambda\nabla^\varphi\cdot(\mathbf{S}^\varphi\cdot))\mathcal{B}^\alpha - (B\cdot\nabla^\varphi)\mathcal{V}^\alpha  \\
&= (Z^\alpha B\cdot\nabla^\varphi)v + Z^\alpha\varphi\left((\partial_z^\varphi v\cdot\nabla^\varphi)B + (\partial_z^\varphi B\cdot\nabla^\varphi)v\right) -\mathcal{C}^\alpha(\mathcal{T}_F)   \\
&\quad + \sum_{i=1}^3 B_i\mathcal{C}_i^\alpha(v) + \mathcal{C}^\alpha(\mathcal{T}_I) + \lambda\mathcal{D}^\alpha(\mathbf{S}^\varphi B) + \lambda\nabla^\varphi\cdot(\mathcal{E}^\alpha B).
\end{split}
\end{equation}

\subsection{Kinematic Boundary}
For the boundary condition, $Z^\alpha = D_{y}^\alpha$, i.e. $\alpha_3 = 0$. It is easy to check the following as we did before. In fact, this is just same as Lemma 5.7 in \cite{NMFR1},
\begin{equation} \label{4.21}
\partial_z Z^\alpha h - v^b\cdot Z^\alpha \mathbf{N} - \mathcal{V}^\alpha\cdot \mathbf{N} = \mathcal{C}^\alpha(h),
\end{equation}
where
$$
\mathcal{C}^\alpha(h) := -[Z^\alpha,v_y^b,\nabla_y h] - \frac{(\partial_z v)^b}{\partial_z\varphi}\cdot\mathbf{N}Z^\alpha h,
$$
with estimate for $ \mathcal{C}^\alpha(h)$
\begin{equation} \label{4.22}
|\mathcal{C}^\alpha(h)|_{L^2} \leq \Lambda\left(\frac{1}{c_0},\|v\|_{E^{1,\infty}} + |h|_{2,\infty}\right)\left( \|v\|_{E^m} + |h|_m \right).
\end{equation}

\subsection{Continuity of Stress-Tensor}
This is also same as Lemma 5.6 in \cite{NMFR1},
$$
Z^\alpha\left((q-gh)\mathbf{N} - 2\varepsilon \mathbf{S}^\varphi v\mathbf{N}\right) = 0.
$$
Using new variable $\mathcal{V}^{\alpha}$, 
\begin{equation} \label{4.23}
2\varepsilon \mathbf{S}^\varphi\mathcal{V}^\alpha\mathbf{N} - (Z^\alpha q - g Z^\alpha h)\mathbf{N} + (2\varepsilon \mathbf{S}^\varphi v - (q-gh))Z^\alpha\mathbf{N} = \mathcal{C}^\alpha(\partial) - 2\varepsilon Z^\alpha h\partial_z^\varphi(\mathbf{S}^\varphi v)\mathbf{N},
\end{equation}
where 
$$
\mathcal{C}^\alpha(\partial) := -\varepsilon\mathcal{E}^\alpha(v) - \sum_{\substack{\beta+\gamma=\alpha,\\0<|\beta|<|\alpha|}} \varepsilon Z^\beta(\mathbf{S}^\varphi v)Z^\gamma\mathbf{N} + \sum_{\substack{\beta+\gamma=\alpha,\\0<|\beta|<|\alpha|}} \varepsilon Z^\beta(q-gh)Z^\gamma\mathbf{N},
$$
with estimate for $ \mathcal{C}^\alpha(\partial)$
\begin{equation} \label{4.24}
|\mathcal{C}^\alpha(\partial)|_{L^2} \leq \varepsilon\Lambda\left(\frac{1}{c_0},\|v\|_{E^{2,\infty}} + |h|_{2,\infty}\right)\left( |v^b|_{m} + |h|_m \right).
\end{equation}

\section{Pressure estimates}
In this section, we get estimate for total pressure $p$, for any smooth solution $(v,B,q,h)$. We decompose $q$ into $q = q^E + q^{NS}$, where
\begin{equation*}
\begin{split}
\triangle^\varphi q^E &= -\nabla^\varphi\cdot(v\cdot\nabla^\varphi v) + \nabla^\varphi\cdot(B\cdot\nabla^\varphi B),\quad q^E|_{z=0} = gh,   \\
\triangle^\varphi q^{NS} &= 0,\quad q^{NS}|_{z=0} = 2\varepsilon \mathbf{S}^\varphi v\mathbf{n} \cdot \mathbf{n}.
\end{split}
\end{equation*}

We express $\triangle^\varphi$ as elliptic operator.
$$
\triangle^\varphi f = \frac{1}{\partial_z\varphi}\nabla\cdot(E\nabla f),
$$
where
$$
E = \begin{pmatrix}
\partial_z\varphi & 0 & -\partial_1\varphi \\
0 & \partial_z\varphi & -\partial_2\varphi \\
-\partial_1\varphi & -\partial_2\varphi & \frac{1+(\partial_1\varphi)^2+(\partial_2\varphi)^2}{\partial_z\varphi}
\end{pmatrix} = \frac{1}{\partial_z\varphi}PP*,
$$
and
$$
P = \begin{pmatrix}
\partial_z\varphi & 0 & 0 \\
0 & \partial_z\varphi & 0 \\
-\partial_1\varphi & -\partial_2\varphi & 1
\end{pmatrix}.
$$
Matrix $E$ is positive symmetric and there exists $\delta(c_0) > 0$ such that
$$
EX\cdot X \geq \delta|X|^2,\,\,\,\,\forall\in \mathbb{R}^3,
$$
if $\|\nabla_y\varphi\|_{L^\infty} \leq \frac{1}{c_0}$, and $\partial_z\varphi \leq c_0 > 0$. We have an estimate
\begin{equation} \label{5.1}
\|E\|_{W^{k,\infty}} \leq \Lambda(\frac{1}{c_0},|h|_{k+1,\infty}).
\end{equation}
Also, using the following decomposition, 
$$
E = \mathbf{I}_A + \tilde{E},\,\,\,\,\tilde{E} = \begin{pmatrix}
\partial_z\eta & 0 & -\partial_1\eta \\
0 & \partial_z\eta & -\partial_2\eta \\
-\partial_1\eta & -\partial_2\eta & \frac{A((\partial_1\eta)^2+(\partial_1\eta)^2)-\partial_z\eta}{A\partial_z\varphi}
\end{pmatrix},\,\,\,\,\mathbf{I}_A = \text{diag}(A,A,1/A),
$$
we also get an estimate,
\begin{equation} \label{5.2}
\|\tilde{E}\|_{H^s} \leq \Lambda(\frac{1}{c_0},|h|_{1,\infty})|h|_{s+\frac{1}{2}}.
\end{equation}

We employ the following lemmas about elliptic problem, from \cite{NMFR1}. First Lemma is for Euler part, $q^E$.
\begin{lemma} \label{lemma 5.1}
For elliptic equation in $S$,
$$
-\nabla\cdot(E\nabla\rho) = \nabla\cdot F,\,\,\,\,\rho|_{z=0} = 0,
$$
we have the estimates :
\begin{equation} \label{5.3}
\begin{split}
\|\nabla\rho\| &\leq \Lambda(\frac{1}{c_0},|h|_{1,\infty})\|F\|_{L^2},\,\,\,\,\|\nabla^2\rho\| \leq \Lambda(\frac{1}{c_0},|h|_{2,\infty})(\|\nabla\cdot F\| + \|F\|_1),  \\
\|\nabla\rho\|_k &\leq \Lambda(\frac{1}{c_0},|h|_{2,\infty} + |h|_3 + \|F\|_{H_{tan}^2} + \|\nabla\cdot F\|_{H_{tan}^1} )(|h|_{k+\frac{1}{2}} + \|F\|_k),\,\,\,\,k\geq 1,  \\
\|\partial_{zz}\rho\|_{k-1} &\leq \Lambda(\frac{1}{c_0},|h|_{2,\infty}+ |h|_3 + \|F\|_{H_{tan}^2} + \|\nabla\cdot F\|_{H_{tan}^1})(|h|_{k+\frac{1}{2}} + \|F\|_k + \|\nabla\cdot F\|_{k-1}),\,\,\,\,k\geq 2.  
\end{split}
\end{equation}
\end{lemma}
\begin{proof}
See Lemma 6.1 in \cite{NMFR1}.
\end{proof}

Second lemma is about Navier-Stokes part, $q^{NS}$.
\begin{lemma} \label{lemma 5.2}
For elliptic equation in $S$,
$$
-\nabla\cdot(E\nabla\rho) = 0,\,\,\,\,\rho|_{z=0} = f^b,
$$
we have the estimates :
\begin{equation} \label{5.6}
\|\nabla\rho\|_{H^k} \leq \Lambda(\frac{1}{c_0},|h|_{2,\infty} + |h|_3 + |f^b|_{1,\infty} + |f^b|_{5/2})(|h|_{k+1/2} + |f^b|_{k+1/2}).
\end{equation}
\end{lemma}
\begin{proof}
See Lemma 6.2 in \cite{NMFR1}.
\end{proof}

Using above two lemmas, we can get estimates for $q^E,q^{NS}$.
\begin{proposition} \label{proposition 5.3}
For $q^E$, we have the estimates : 
\begin{equation} \label{5.7}
\begin{split}
\|\nabla q^E\|_{m-1} + \|\partial_{zz}q^E\|_{m-2} &\leq \Lambda(\frac{1}{c_0},|h|_{2,\infty} + |h|_3 + \|v\|_{E^{1,\infty}} + \|B\|_{E^{1,\infty}} + \|v\|_{E^3} + \|B\|_{E^3} )  \\
&\times ( \|v\|_{E^m} + \|B\|_{E^m} + |h|_{m-\frac{1}{2}}),  \\
\|\nabla q^E\|_{1,\infty} + \|\partial_{zz}q^E\|_{L^\infty} &\leq \Lambda(\frac{1}{c_0},|h|_{2,\infty} + |h|_4 + \|v\|_{E^{1,\infty}} + \|B\|_{E^{1,\infty}} + \|v\|_{E^4} + \|B\|_{E^4} ),  \\
\|\nabla q^E\|_{2,\infty} &\leq \Lambda(\frac{1}{c_0},|h|_{2,\infty} + |h|_5 + \|v\|_{E^{1,\infty}} + \|B\|_{E^{1,\infty}} + \|v\|_{E^5} + \|B\|_{E^5} ).
\end{split}
\end{equation}
\end{proposition}
\begin{proof}
See Proposition 6.4 in \cite{NMFR1}. We just suffice to add same types of norms for $B$.
\end{proof}

\begin{proposition} \label{proposition 5.4}
For $q^{NS}$, we have the estimates for $m\leq 1$ : 
\begin{equation} \label{5.10}
\begin{split}
|\nabla q^{NS}|_{H^{m-1}} &\leq \varepsilon\Lambda(\frac{1}{c_0},|h|_{2,\infty} + |h|_4 + \|v\|_{E^{2,\infty}} + \|v\|_{E^4})( |v^b|_{m+\frac{1}{2}} + |h|_{m+\frac{1}{2}} ),  \\
\|\nabla q^{NS}\|_{L^\infty} &\leq \varepsilon\Lambda(\frac{1}{c_0},|h|_{2,\infty} + |h|_4 + \|v\|_{E^{2,\infty}} + \|v\|_{E^4} ).
\end{split}
\end{equation}
\end{proposition}
\begin{remark}
Note that $q^{NS}$ can be estimated in standard sobolev space, not necessarily in conormal one.
\end{remark}
\begin{proof}
See Proposition 6.3 in \cite{NMFR1}.
\end{proof}

The following proposition will be used for Taylor sign condition.
\begin{proposition} \label{proposition 5.6}
For $T\in [0,T^\varepsilon)$, we have the following estimate.
\begin{equation} \label{5.12}
\begin{split}
\int_0^T |(\partial_z\partial_t q^E)^b|_{L^\infty} & \leq \int_0^T \Lambda(\frac{1}{c_0}, |h|_6 + |h|_{3,\infty} + \|v\|_6 + \|B\|_6 + \|\partial_z v\|_4 + \|\partial_z B\|_4 + \|v\|_{E^{2,\infty}} + \|B\|_{E^{2,\infty}} ) \\
& \quad \times ( 1 + \varepsilon\|\partial_{zz}v\|_{L^\infty} + \lambda\|\partial_{zz}B\|_{L^\infty} + \varepsilon\|\partial_{zz}v\|_3 + \lambda\|\partial_{zz}B\|_3 ) d\tau.
\end{split}
\end{equation}
\end{proposition}
\begin{proof}
This is a version of Proposition 6.5 in \cite{NMFR1}. 
First,
$$
\triangle^\varphi q^E = -\nabla^\varphi\cdot(v\cdot\nabla^\varphi v) + \nabla^\varphi\cdot(B\cdot\nabla^\varphi B),\,\,\,q^E|_{z=0} = gh.
$$
Taking $\partial_t$,
$$
\nabla\cdot(E\nabla \partial_t q^E) = -\nabla\cdot(\partial_t(P(v\cdot\nabla^\varphi v))) + \nabla\cdot(\partial_t(P(B\cdot\nabla^\varphi B))) - \nabla\cdot(\partial_t E \nabla P^E),\,\,\,\partial_t q^E |_{z=0} = g\partial_t h.
$$
We divide into $q^E := q^i + q^B$ so that,
\begin{equation*}
\begin{split}
\nabla\cdot(E\nabla \partial_t P^i) &= -\nabla\cdot(\partial_t(P(v\cdot\nabla^\varphi v))) + \nabla\cdot(\partial_t(P(B\cdot\nabla^\varphi B))) - \nabla\cdot(\partial_t E \nabla P^E),\,\,\,\partial_t P^i |_{z=0} = 0, \\
\nabla\cdot(E\nabla \partial_t P^B) &= 0,\,\,\,\partial_t P^B |_{z=0} = g\partial_t h.
\end{split}
\end{equation*}
Estimate for $q^B$ is exactly same as \cite{NMFR1}, and for $q^i$, we use Lemma 6.6 in \cite{NMFR1}, where $F$ also includes similar structure for $B$. We get
\begin{equation*}
\begin{split}
\int_0^T |(\partial_z\partial_t q^i)^b|_{L^\infty} &\leq \int_0^T \Lambda\left(\frac{1}{c_0},\|v\|_{E^5} + \|v\|_{E^{2,\infty}} + \|v\|_6 + \|B\|_{E^5} + \|B\|_{E^{2,\infty}} + \|B\|_6 + |h|_6 + |h|_{3,\infty} \right)  \\
&\times\left(1 + \|\partial_t v\|_{L^\infty} + \|\partial_t v\|_3 + \|\partial_t B\|_{L^\infty} + \|\partial_t B\|_3 \right) ds.
\end{split}
\end{equation*}
Estimate for $\|\partial_z v\|_3 + \|\partial_t v\|_{L^\infty} + \|\partial_z B\|_3 + \|\partial_t B\|_{L^\infty}$ is gained from the two main equations of (\ref{1.1}), Proposition \ref{proposition 5.3}, Proposition \ref{proposition 5.4}, and Proposition \ref{proposition 3.2}.
\begin{equation}
\begin{split}
\|\partial_z v\|_3 + \|\partial_t v\|_{L^\infty} + \|\partial_z B\|_3 + \|\partial_t B\|_{L^\infty} &\leq \Lambda\left( \frac{1}{c_0}, \|B\|_{E^5} + \|B\|_{E^{2,\infty}} + \|v\|_{E^5} + \|v\|_{E^{2,\infty}} + |h|_5 + |h|_{2,\infty} \right)  \\
&\quad\quad\quad \times \left( 1 + \varepsilon\|\partial_{zz}v\|_{L^\infty} + \varepsilon\|\partial_{zz}v\|_3 + \lambda\|\partial_{zz}B\|_{L^\infty} + \lambda\|\partial_{zz}B\|_3 \right).
\end{split}
\end{equation}
Euler part $q^B$ is given by Lemma \ref{proposition 5.3}.
Putting altogether, we get the result.
\begin{equation} \label{5.13}
\begin{split}
\int_0^T |(\partial_z\partial_t q^E)^b|_{L^\infty} &\leq \int_0^T \Lambda\left(\frac{1}{c_0},\|v\|_{6} + \|\partial_z v\|_4 + \|v\|_{E^{2,\infty}} + \|B\|_{6} + \|\partial_z B\|_4 + \|B\|_{E^{2,\infty}} + |h|_6 + |h|_{3,\infty} \right)  \\
&\quad\quad\quad\times\left( 1 + \varepsilon\|\partial_{zz}v\|_{L^\infty} + \varepsilon\|\partial_{zz}v\|_3 + \lambda\|\partial_{zz}B\|_{L^\infty} + \lambda\|\partial_{zz}B\|_3 \right) ds.
\end{split}
\end{equation}
\end{proof}

\section{$L^2$ energy estimate}
\begin{proposition} \label{L2 est}
For any smooth solution of (\ref{1.4}), we have the following zero-order energy estimate.
$$
\frac{d}{dt}\left( \int_S |v|^2 dV_t + g\int_{\partial S} |h|^2 dy + \int_S |B|^2 dV_t \right) + 4\varepsilon\int_S |\mathbf{S}^\varphi v|^2 dV_t + 4\lambda\int_S |\mathbf{S}^\varphi B|^2 dV_t = 0.
$$
\end{proposition}
\begin{proof}
Multiplying $v$ and integrating on $S$ for Navier-Stokes, and also using boundary condition, we get
$$
\frac{d}{dt}\int_S |v|^2 dV_t + 4\varepsilon\int_S |\mathbf{S}^\varphi v|^2 dV_t - \int_{\partial S}|v|^2(h_t - v^b\cdot\mathbf{N})dy = 2\int_{\partial S}(2\varepsilon \mathbf{S}^\varphi v - q\mathbf{I})\mathbf{N}\cdot v dy + 2\int_S v\cdot(B\cdot\nabla^\varphi)B dV_t.
$$
Then using kinematic boundary condition and Continuity of stress tensor condition,
$$
\frac{d}{dt}\int_S |v|^2 dV_t + 4\varepsilon\int_S |\mathbf{S}^\varphi v|^2 dV_t = 2\int_S v\cdot(B\cdot\nabla^\varphi)B dV_t - 2\int_{\partial S}) gh(\mathbf{N}\cdot v) dy,
$$
\begin{equation} \label{6.1}
\frac{d}{dt}\left( \int_S |v|^2 dV_t + g\int_{\partial S}|h|^2 dy \right)+ 4\varepsilon\int_S |\mathbf{S}^\varphi v|^2 dV_t = 2\int_S v\cdot(B\cdot\nabla^\varphi)B dV_t.
\end{equation}
Multiplying $B$ and integrating on $S$ for Faraday's Law, we get
$$
\frac{d}{dt}\int_S |B|^2 dV_t + 4\lambda\int_S |\mathbf{S}^\varphi B|^2 dV_t = 2\int_S B\cdot(B\cdot\nabla^\varphi)v dV_t + 4\lambda\int_{\partial S} \mathbf{S}^\varphi B (\mathbf{N}\cdot B) dV_t.
$$
Using divergence free condition, we know that
\begin{equation} \label{6.2}
\begin{split}
\int_{S} B\cdot(B\cdot\nabla^\varphi) v dV_t &= \int_{\partial S} (B\cdot v)(B\cdot\mathbf{n}) dy - \int_{S} v\cdot(B\cdot\nabla^\varphi)B dV_t   \\
&= - \int_{S} v\cdot(B\cdot\nabla^\varphi)B dV_t.
\end{split}
\end{equation}
Therefore, Faraday's law gives
\begin{equation} \label{6.3}
\frac{d}{dt}\int_S |B|^2 dV_t + 4\lambda\int_S |\mathbf{S}^\varphi B|^2 dV_t = - 2\int_S v\cdot(B\cdot\nabla^\varphi)B dV_t.
\end{equation}
We add (\ref{6.1}) and (\ref{6.3}), and cancel both right hand sides to finish the proof.
\end{proof}

\section{Higher order Energy estimate}
Using the results of section 3, 4, and 5, we can make high order energy estimates. We define the $\Lambda_{\infty}(t)$ which contains all low order terms.
\begin{equation} \label{lambda infty}
\Lambda_{\infty}(t) := \Lambda\left(\frac{1}{c_0},\|v
(t)\|_{E^{2,\infty}} + \sqrt{\varepsilon}\|\partial_{zz}v(t)\|_{L^\infty} + \|B
(t)\|_{E^{2,\infty}} + \sqrt{\lambda}\|\partial_{zz}B(t)\|_{L^\infty} + \|v(t)\|_{E^4} + |h|_4 \right).
\end{equation}
\subsection{Navier-Stokes Equation with Lorentz force}
We apply $L^2$ energy estimate in section 6 to (\ref{4.1}), (\ref{4.13}), (\ref{4.21}), and (\ref{4.23}), to get 
\begin{equation} \label{7.1}
\begin{split}
&\frac{d}{dt}\int_S |\mathcal{V}^\alpha|^2 dV_t + 4\varepsilon\int_S |\mathbf{S}^\varphi\mathcal{V}^\alpha|^2 dV_t - 2\int_S (B\cdot\nabla^\varphi)\mathcal{B}^\alpha\cdot\mathcal{V}^\alpha dV_t  \\
&= \mathcal{R}_S + \mathcal{R}_C + \mathcal{R}_P + 2\int_{z=0}(2\varepsilon \mathbf{S}^\varphi\mathcal{V}^\alpha - \mathcal{Q}^\alpha Id)\mathbf{n}\cdot\mathcal{V}^\alpha dy,
\end{split}
\end{equation}
where $\mathcal{R}_S$ and $\mathcal{R}_C$ are defined as
\begin{equation} \label{7.2}
\begin{split}
\mathcal{R}_S &:= 2\int_S \left( \varepsilon\mathbf{D}^\alpha(\mathbf{S}^\varphi v) + \varepsilon\nabla^\varphi\cdot(\mathcal{E}^\alpha(v)) \right)\cdot\mathcal{V}^\alpha dV_t,  \\
\mathcal{R}_C &:= -2\int_S\left( (\mathcal{C}^\alpha(\mathcal{T}_v) + \mathcal{C}^\alpha(q))\cdot\mathcal{V}^\alpha - \mathcal{C}^\alpha(d_v)\mathcal{Q}^\alpha \right) dV_t,
\end{split}
\end{equation}
and $\mathcal{R_P}$ is defined as
\begin{equation} \label{7.4}
\mathcal{R}_P := 2\int_S \left\{ (Z^\alpha B\cdot\nabla^\varphi)B + Z^\alpha\varphi\left( (\partial_z^\varphi v\cdot\nabla^\varphi)v + (\partial_z^\varphi B\cdot\nabla^\varphi)B \right) + \sum_{i=1}^3 B_i\mathcal{C}^\alpha_i(B) + \mathcal{C}^\alpha(\mathcal{T}_B) \right\}\cdot\mathcal{V}^\alpha dV_t.
\end{equation}

\subsection{Faraday law}
We perform similar work and get the following energy estimate using (\ref{4.1}) and (\ref{4.20}).
\begin{equation} \label{7.5}
\begin{split}
&\frac{d}{dt}\int_S |\mathcal{B}^\alpha|^2 dV_t + 4\lambda\int_S |\mathbf{S}^\varphi\mathcal{B}^\alpha|^2 dV_t + 2\int_S (B\cdot\nabla^\varphi)\mathcal{B}^\alpha\cdot\mathcal{V}^\alpha dV_t   \\
&= 2\int_S \{ (Z^\alpha B\cdot\nabla^\varphi)v + Z^\alpha\varphi\left((\partial_z^\varphi v\cdot\nabla^\varphi)B + (\partial_z^\varphi B\cdot\nabla^\varphi)v\right)  \\
&\quad - \mathcal{C}^\alpha(\mathcal{T}_F) + \sum_{i=1}^3 v_i\mathcal{C}_i^\alpha(B) + \mathcal{C}^\alpha(\mathcal{T}_I) + \varepsilon\mathbf{D}^\alpha(\mathbf{S}^\varphi B) + \varepsilon\nabla^\varphi\cdot(\mathcal{E}^\alpha (B)) \} \cdot\mathcal{B}^\alpha dV_t  \\
&= 4\lambda\int_{\p S} (\mathbf{S}^\varphi\mathcal{B}^\alpha)\mathbf{n}\cdot\mathcal{B}^\alpha dy + \mathcal{R}_{S_B} + \mathcal{R}_{C_B} + \mathcal{R}_{P_B} = \mathcal{R}_{S_B} + \mathcal{R}_{C_B} + \mathcal{R}_{P_B},  \\
\end{split}
\end{equation}
where
\begin{equation} \label{Rs for B}
\begin{split}
\mathcal{R}_{S_B} &:= 2\int_S \{\varepsilon\mathbf{D}^\alpha(\mathbf{S}^\varphi B) + \varepsilon\nabla^\varphi\cdot(\mathcal{E}^\alpha (B)) \} \cdot\mathcal{B}^\alpha dV_t  \\
\mathcal{R}_{C_B} &:= -2\int_S \mathcal{C}^\alpha(\mathcal{T}_F)\cdot\mathcal{B}^\alpha dV_t  \\
\mathcal{R}_{P_B} &:= 2\int_S \left\{ (Z^\alpha B\cdot\nabla^\varphi)v + Z^\alpha\varphi\left((\partial_z^\varphi v\cdot\nabla^\varphi)B + (\partial_z^\varphi B\cdot\nabla^\varphi)v\right) + \sum_{i=1}^3 v_i\mathcal{C}_i^\alpha(B) + \mathcal{C}^\alpha(\mathcal{T}_I) \right\}\cdot\mathcal{B}^\alpha dV_t.
\end{split}
\end{equation}


\subsection{More commutator estimates}
We need to know the estimates of $\mathcal{R}_{S}$, $\mathcal{R}_{C}$, $\mathcal{R}_P$, $\mathcal{R}_{S_B}$, $\mathcal{R}_{C_B}$, and $\mathcal{R}_{P_B}$, those are defined in (\ref{7.2}), (\ref{7.4}), and (\ref{Rs for B}). All the estimates come from Propositions in section 3 and Lemma \ref{lemma 4.1}. \\
\noindent 1) Estimate of $\mathcal{R}_S$. For $\mathcal{R}_S$, using (\ref{7.2}) and (\ref{4.10}),
\begin{equation} \label{7.6}
\begin{split}
\|\mathcal{R}_S\| &\leq \varepsilon \ \Lambda_{\infty}(t) \big\{ \|\nabla\mathcal{V}^\alpha\| ( \|\mathcal{V}^\alpha\| + \|\partial_z v\|_{m-1} + |h|_{m+\frac{1}{2}} )  \\
&+ ( \|v\|^2_{E^m} + |h|^2_{m+\frac{1}{2}} ) + \|\partial_{zz}v\|_{L^\infty}( |h|^2_m + \|\mathcal{V}^\alpha\|^2 ) \big\} ,
\end{split}
\end{equation}
where $\Lambda_{\infty}(t)$ is defined in (\ref{lambda infty}). \\

\noindent 2) Estimate of $\mathcal{R}_C$. We use (\ref{4.6}), (\ref{4.8}), (\ref{4.3}), and Proposition \ref{proposition 5.3} to get,
$$
\|\mathcal{C}^\alpha(q^E)\|\|\mathcal{V}^m\| \leq \Lambda_{\infty}(t) \left( \|v\|_{E^{m}} + \|B\|_{E^{m}} + |h|_m \right) \|\mathcal{V}^m\|.
$$ 
So, we get 
\begin{equation} \label{7.7}
\begin{split}
\|\mathcal{R}_C\| &\leq \ \Lambda_{\infty}(t) \left(\|v\|_{E^{m}} + \|B\|_{E^m} + |h|_m + \varepsilon |h|_{m+\frac{1}{2}} + \varepsilon |v^b|_{m+\frac{1}{2}} \right)\|\mathcal{V}^m\| \\
&\leq \Lambda_{\infty}(t) \left( \varepsilon\|\nabla \mathcal{V}^m\|\|\mathcal{V}^m\| + \|v\|_{E^m}^2 + \|B\|_{E^m}^2 + \|\mathcal{V}^m\|^2 + |h|_m^2 + \varepsilon |h|_{m+\frac{1}{2}}^2 \right). 
\end{split}
\end{equation}
\noindent 3) Estimate of $\mathcal{R}_P$.
We know that
$$
\|\mathcal{C}_i^\alpha(B)\| \leq \Lambda\left(\frac{1}{c_0},|h|_{2,\infty} + \|\nabla B\|_{1,\infty}\right)\left( \|\nabla B\|_{m-1} + |h|_{m-\frac{1}{2}} \right),
$$
from Lemma \ref{lemma 4.1}. For $\mathcal{C}^\alpha(\mathcal{T}_B)$, using Proposition \ref{proposition 3.1},
\begin{equation*}
\begin{split}
\mathcal{C}^\alpha(\mathcal{T}_B) &= \sum_{i=1}^3[Z^\alpha,B_i,\partial_i^\varphi B] \leq \Lambda\left( \frac{1}{c_0}, \|B\|_{1,\infty} + \|\partial_i^\varphi B\|_{1,\infty} \right) \left( \|B\|_{m-1} + \|\partial_i^\varphi B\|_{m-1} \right)  \\
&\leq \Lambda\left( \frac{1}{c_0}, \|B\|_{E^{2,\infty}} + |h|_{2,\infty} \right) \left( \|B\|_{E^{m}} + |h|_{m} \right).
\end{split}
\end{equation*}
Hence, we get
\begin{equation} \label{7.8}
\left\|\mathcal{R}_P\right\| \leq \Lambda\left( \frac{1}{c_0}, \|v\|_{E^{1,\infty}} + \|B\|_{E^{2,\infty}} + |h|_{2,\infty} \right) \left( \|\mathcal{V}^m\| + \|B\|_{E^{m}} + |h|_{m} \right).
\end{equation}

\noindent 4) Estimate of $\int_{z=0}(2\varepsilon \mathbf{S}^\varphi\mathcal{V}^\alpha - \mathcal{Q}^\alpha Id)\mathbf{n}\cdot\mathcal{V}^\alpha dy$.
\begin{equation} \label{7.111}
\begin{split}
\int_{z=0}(2\varepsilon \mathbf{S}^\varphi\mathcal{V}^\alpha - \mathcal{Q}^\alpha Id)\mathbf{n}\cdot\mathcal{V}^\alpha dy &= \underbrace{ \int_{\p S} (-g Z^{\alpha}h + \p_{z}^{\varphi}q   Z^{\alpha} h ) \mathbf{N}\cdot\mathcal{V}^{\alpha} }_{(B)}  \\
&- \underbrace{ 2\int_{\p S} ( 2\varepsilon \mathbf{S}^{\varphi} v - (q-gh)I) Z^{\alpha}\mathbf{N}\cdot\mathcal{V}^{\alpha} }_{(A)} + \mathcal{R}_{B},  \\	
\end{split}
\end{equation}
where 
\begin{equation} \label{RB}
R_{B} := \int_{\p S} (\mathcal{C}^{\alpha}(\p)-2\varepsilon Z^{\alpha}h \p_{z}^{\varphi}(\mathbf{S}^{\varphi} v)\mathcal{N} )\cdot\mathcal{V}^{\alpha},
\end{equation}
and is estimated by 
\begin{equation} \label{RB est}
\begin{split}
|R_{B}| &\leq \varepsilon\|\nabla\mathcal{V}^{m}\|\|\mathcal{V}^{m}\| + \Lambda_{\infty}(t)\|\p_z v\|_{m-1} |h|_{m}  \\
&\quad + \Lambda_{\infty}(t) (1 + |(\p_z \p_t q^{E})^b|_{L^\infty} ) (|h|^{2}_{m} + \varepsilon|h|^{2}_{m+\frac{1}{2}} + \|\mathcal{V}^{m}\| ),
\end{split}
\end{equation}
by (\ref{4.23}) and (\ref{4.24}). $(A)$ can be estimated by
\begin{equation} \label{(A) est}
\begin{split}
|(A)| &\leq \Big|  \int_{z= 0} \big(2 \varepsilon \mathbf{S}^{\varphi} v -  q^{NS}{I} \big) Z^\alpha \mathbf
N \cdot \mathcal{V}^\alpha \, dy \Big|  \\
& \leq |  Z^\alpha \nabla_{y} h |_{-\frac{1}{2}} \, | \big(2 \varepsilon \mathbf{S}^{\varphi} v - q^{NS}I \big)(\mathcal{V}^\alpha)^b|_{\frac{1}{2}}  \\
&\leq \varepsilon\Lambda_{\infty}(t) |h|_{m+{\frac{1}{2}}} |(\mathcal{V}^\alpha)^b |_{\frac{1}{2}}.
\end{split}
\end{equation}
Meanwhile, for $(B)$,
\begin{equation} \label{(B) est}
\begin{split}
(B) &:= \int_{\p S} \big( -g Z^\alpha h + \p_{z}^{\varphi} q\, Z^\alpha h \big)  \mathbf{N}\cdot
\mathcal{V}^\alpha \\
&=   \int_{z=0}\big( -g Z^\alpha h + \p_{z}^{\varphi} q^E\, Z^\alpha h \big)  \mathbf{N}\cdot
\mathcal{V}^\alpha \, dy + \int_{z=0} \p_{z}^{\varphi} q^{NS} Z^\alpha h  \mathcal{V}^\alpha \cdot \mathbf{N}  \\
&\leq \int_{\p S}\big( -g Z^\alpha h + \p_{z}^{\varphi} q^E\, Z^\alpha h \big)\partial_{t} Z^\alpha h \\
&\quad  - \int_{\p S}\big( -g Z^\alpha h + \p_{z}^{\varphi} q^E\, Z^\alpha h \big)v^b \cdot (Z^\alpha \nabla_{y}h - \mathcal{C}^\alpha(h))\, dy +  \Lambda_{\infty}(t) \varepsilon |h|_{m} |(\mathcal{V}^{\alpha})^{b}|  \\
&\leq  \int_{\p S}\big( -g Z^\alpha h + \p_{z}^{\varphi} q^E\, Z^\alpha h \big)\partial_{t} Z^\alpha h +  \Lambda_{\infty}(t)\big( |h|_{m} + \|v\|_{E^m})|h|_{m} +  \Lambda_{\infty}(t) \varepsilon |h|_{m} |(\mathcal{V}^{\alpha})^{b}| \\
&= - {\frac{1}{2}} {\frac{\dd}{\dd t}}  \int_{\p S} ( g- \p_{z}^{\varphi} q^E) |Z^\alpha h |^2 - \int_{\p S} \partial_{t}\big( \p_{z}^{\varphi} q^E \big) |Z^\alpha h|^2  +  \Lambda_{\infty}(t)\big( |h|_{m} + \|v\|_{E^m})|h|_{m} +  \Lambda_{\infty}(t) \varepsilon |h|_{m} |(\mathcal{V}^{\alpha})^{b}|, \\
\end{split}
\end{equation}
where we used Proposition \ref{proposition 5.4}, (\ref{4.21}), and (\ref{4.22}). Therefore, using (\ref{7.111}), (\ref{7.14}), (\ref{(A) est}), and (\ref{(B) est}), we can write this term  as
\begin{equation} \label{7.112}
\begin{split}
2\int_{z=0}(2\varepsilon \mathbf{S}^\varphi\mathcal{V}^\alpha - \mathcal{Q}^\alpha Id)\mathbf{n}\cdot\mathcal{V}^\alpha dy &= - \frac{1}{2}\frac{\dd}{\dd t}\int_{\p S} (g-\p_{z}^{\varphi}q^{E}) |Z^{\alpha} h|^2 + \tilde{\mathcal{R}}_{B} ,
\end{split}
\end{equation}
with estimate,
\begin{equation} \label{tilde RB}
\begin{split}
| \tilde{\mathcal{R}}_{B}|  &\leq  \Lambda_{\infty}(t) \Big\{  \varepsilon \big(1+   \|\partial_{zz} v \|_{L^\infty} \big)|h|_{m} + \varepsilon |v^b |_{m}\big) |(V^\alpha)^b| + \varepsilon |h|_{m+{\frac{1}{2}}} |(\mathcal{V}^\alpha)^b|_{{\frac{1}{2}}} \\
&\quad + (1 + |(\partial_{z} \partial_{t} q^E)^b |_{L^\infty})|h|_{m}^2 + \|v\|_{E^m} |h|_{m} \Big\},
\end{split}
\end{equation}
where we used Proposition \ref{proposition 5.4}.  \\

\noindent 5) Estimate of $\mathcal{R}_{S_B}$. We use (\ref{4.10}). We suffice to replace $v$ in (\ref{7.6}) into $B$.
\begin{equation} \label{7.9}
\begin{split}
\|\mathcal{R}_{S_B}\| &\leq \varepsilon\Lambda_\infty(t) \Big\{ \|\nabla\mathcal{B}^\alpha\| ( \|\mathcal{B}^\alpha\| + \|\partial_z B\|_{m-1} + |h|_{m+\frac{1}{2}} )  \\
&\quad + ( \|B\|^2_{E^m} + |h|^2_{m+\frac{1}{2}} ) + \|\partial_{zz}B\|_{L^\infty}( |h|^2_m + \|\mathcal{B}^\alpha\|^2 ) \Big\}.
\end{split}
\end{equation}

\noindent 6) Estimate of $\mathcal{R}_{C_B}$.  Simply we get,
\begin{equation} \label{7.10}
\|\mathcal{R}_{C_B}\| \leq \Lambda_{\infty}(t)(\|B\|_m + \|\partial_z B\|_{m-1} + |h|_m)\|\mathcal{B}^\alpha\|.
\end{equation}
\noindent 7) Estimate of $\mathcal{R}_{P_B}$. This can be estimated similar as $\mathcal{R}_{P}$,
\begin{equation} \label{7.11}
\left\|\mathcal{R}_{P_B}\right\| \leq \Lambda\left( \frac{1}{c_0}, \|v\|_{E^{1,\infty}} + \|B\|_{E^{2,\infty}} + |h|_{2,\infty} \right) \left( \|\mathcal{B}^m\| + \|B\|_{E^{m}} + |h|_{m} \right).
\end{equation}

\subsection{Energy estimate for Navier-Stokes Equation with Lorentz Force}
Using (\ref{7.1}) and (\ref{7.112}),
\begin{equation} \label{7.12}
\begin{split}
&\frac{d}{dt}\frac{1}{2}\left(\int_S |\mathcal{V}^\alpha|^2 dV_t + \int_{\partial S}(g - \partial_z^\varphi q^E)|Z^\alpha h|^2 dy \right) + 4\varepsilon\int_S |\mathbf{S}^\varphi\mathcal{V}^\alpha|^2 dV_t - 2\int_S (B\cdot\nabla^\varphi)\mathcal{B}^\alpha\cdot\mathcal{V}^\alpha dV_t  \\
&= \mathcal{R}_S + \mathcal{R}_C + \mathcal{R}_P + \tilde{\mathcal{R}}_B,
\end{split}
\end{equation}
where estimates for four terms on the RHS are given by (\ref{7.6}), (\ref{7.7}), (\ref{7.8}), and (\ref{tilde RB}). High order regularity of $h$ requires positivity of $g - \partial_z^\varphi q^E$, which is known as Rayleigh-Taylor sign condition. Therefore, if we assume,
\begin{equation} \label{7.14}
\partial_z\varphi \geq c_0,\,\,\,|h|_{2,\infty} \leq \frac{1}{c_0},\,\,\,g-(\partial_z^\varphi q^E)|_{z=0} \geq \frac{c_0}{2},\,\,\,\forall t\in [0,T^\varepsilon],
\end{equation}
we get
\begin{equation*}
\begin{split}
&\|\mathcal{V}^m(t)\|^2 + |h(t)|_m^2 + 4\varepsilon \int_S \|\mathbf{S}^\varphi \mathcal{V}^m\|^2 dV_t - \sum_{\forall \alpha} 2 \int_0^t \int_S (B\cdot\nabla^\varphi)\mathcal{B}^\alpha\cdot\mathcal{V}^\alpha dV_t ds   \\
&\leq \Lambda_0 \left( \|\mathcal{V}^m(0)\|^2 + |h(0)|_m^2 \right) + \int_0^t \varepsilon\Lambda_{\infty}(s) \|\nabla\mathcal{V}^m\| \left( \|\mathcal{V}^m\| + \|v\|_{E^m} + |h|_m + |h|_{m+\frac{1}{2}} \right) \dd s \\
&\quad + \int_0^t \Lambda_{\infty}(s) \left( 1 + |(\partial_z\partial_t q^E)^b|_{L^\infty} \right) \left( \|\mathcal{V}^m\|^2 + \|v\|_{E^m}^2 + \|B\|_{E^m}^2 + |h|_m^2 + \varepsilon |h|_{m+\frac{1}{2}}^2 \right) \dd s ,
\end{split}
\end{equation*}
where
\begin{equation*}
\|\mathcal{V}^m(t)\|^2 := \sum_{|\alpha|\leq m}\|\mathcal{V}^\alpha (t)\|^2,\quad \|\mathbf{S}^\varphi\mathcal{V}^m(t)\|^2 := \sum_{|\alpha|\leq m}\|\mathbf{S}^\varphi\mathcal{V}^\alpha (t)\|^2.
\end{equation*}
Using Young's inequality, Proposition \ref{proposition 3.6}, and Proposition \ref{proposition 3.9}, we get the following proposition.

\begin{proposition} \label{proposition 7.2}
Under the assumption of (\ref{7.14}),
we have the following estimate.
\begin{equation} \label{7.15}
\begin{split}
&\|\mathcal{V}^m(t)\|^2 + |h(t)|_m^2 + \varepsilon |h(t)|_{m+\frac{1}{2}}^2 + \varepsilon \int_{0}^{t} \|\nabla\mathcal{V}^m\|^2 ds - \sum_{\forall \alpha} 2\int_0^t \int_S (B\cdot\nabla^\varphi)\mathcal{B}^\alpha\cdot\mathcal{V}^\alpha dV_t ds  \\
&\leq \Lambda_0 \left( \|\mathcal{V}^m(0)\|^2 + |h(0)|_m^2 + \varepsilon |h(0)|_{m+\frac{1}{2}}^2 \right) + \int_0^t \Lambda_{\infty}(s) \left( 1 + |(\partial_z\partial_t q^E)^b|_{L^\infty} \right) \left( \|\mathcal{V}^m\|^2 + \|B\|_{m}^2 + |h|_m^2 + \varepsilon |h|_{m+\frac{1}{2}}^2 \right) \dd s \\
&\quad + \int_0^t \Lambda_{\infty}(s)\left( \|\partial_z v\|_{m-1}^2 + \|\partial_z B\|_{m-1}^2 \right) \dd s.
\end{split}
\end{equation}
Note that $\|v\|_{E^m}$ and $\|B\|_{E^m}$ are absorbed into $\|\mathcal{V}^m\|$, $\|\partial_z v\|_{m-1}$, $\|\mathcal{B}^m\|$, and $\|\partial_z B\|_{m-1}$ by definition.  \\
\end{proposition}

\subsection{Energy estimate for Faraday Law}
Using (\ref{7.5}), 
$$
\frac{d}{dt}\int_S |\mathcal{B}^\alpha|^2 dV_t + 4\lambda\int_S |\mathbf{S}^\varphi\mathcal{B}^\alpha|^2 dV_t + 2\int_S (B\cdot\nabla^\varphi)\mathcal{B}^\alpha\cdot\mathcal{V}^\alpha dV_t = \mathcal{R}_{S_B} + \mathcal{R}_{C_B} + \mathcal{R}_{P_B}.
$$
Using estimates (\ref{7.9}), (\ref{7.10}), and (\ref{7.11}), we get,
\begin{equation*} 
\begin{split}
&\|\mathcal{B}^m\|^2 + 4\lambda\int_0^t \int_S |\mathbf{S}^\varphi\mathcal{B}^{m}|^2 dV_t ds + \sum_{\forall \alpha} 2\int_0^t \int_S (B\cdot\nabla^\varphi)\mathcal{B}^\alpha\cdot\mathcal{V}^\alpha dV_t ds   \\
&\leq \Lambda_0 \left(\|\mathcal{B}^m\|^2\right) + \int_0^t \varepsilon\Lambda_{\infty}(s) \|\nabla\mathcal{B}^m\| \left( \|\mathcal{B}^m\| + \|B\|_{E^m} + |h|_m + |h|_{m+\frac{1}{2}} \right) \dd s  \\
&\quad + \int_0^t \Lambda_{\infty}(s) \left( \|\mathcal{B}^m\|^2 + \|v\|_{m}^2 + \|B\|_{E^m}^2 + |h|_m^2 + \varepsilon |h|_{m+\frac{1}{2}}^2 \right) \dd s.
\end{split}
\end{equation*}
There is no $\partial_z\partial_t q^E$ part here since it come from $\tilde{\mathcal{R}}_B$ which come from (\ref{7.9}). Again, using Young's inequality, Proposition \ref{proposition 3.6}, and Proposition \ref{proposition 3.9}, we get the following proposition.

\begin{proposition} \label{proposition 7.3}
Under the assumption of 
\begin{equation*} 
\partial_z\varphi \geq c_0,\,\,\,|h|_{2,\infty} \leq \frac{1}{c_0},\quad \forall t\in [0,T^\varepsilon],
\end{equation*}
we have
\begin{equation} \label{7.16}
\begin{split}
&\|\mathcal{B}^m\|^2 + \lambda\int_0^t \|\nabla\mathcal{B}^{m}\|^2 ds + \sum_{\forall \alpha} 2\int_0^t \int_S (B\cdot\nabla^\varphi)\mathcal{B}^\alpha\cdot\mathcal{V}^\alpha dV_t ds   \\
&\leq \Lambda_0 \left(\|\mathcal{B}^m(0)\|^2\right) + \int_0^t \Lambda_{\infty}(s) \left( \|\mathcal{B}^m\|^2 + \|v\|_{m}^2 + |h|_m^2 + \varepsilon |h|_{m+\frac{1}{2}}^2 \right) \dd s  + \int_0^t \Lambda_{\infty}(s)( \|\partial_z B\|_{m-1}^2 ) \dd s.
\end{split}
\end{equation}
\end{proposition}
From two main estimates (\ref{7.15}) and (\ref{7.16}), then we see that we should estimate some terms in $\Lambda_\infty$ (such as $\|\partial_z v\|_{k,\infty}$) and $\|\partial_z v\|_{m-1}^2 + \|\partial_z B\|_{m-1}^2$. Note that we should use Proposition \ref{proposition 5.6} to estimate $|(\partial_z\partial_t q^E)^b|_{L^\infty}$ on the right hand side.

\section{Normal estimate}
From Propositions \ref{proposition 7.2} and \ref{proposition 7.3}, we should estimate $\|\partial_z v\|_{m-1}$ and $\|\partial_z B\|_{m-1}$, because they are not controlled by $\|v\|_m$ and $\|B\|_m$. First, the following Lemma {\ref{lemma 8.1}} is true for both $v$ and $B$. 
\begin{lemma} \label{lemma 8.1}
For every integer $m \geq 1$, normal part of $\partial_z v,\partial_z B$ can be estimated as follow.
$$
\|\partial_z v\cdot\mathbf{n}\|_{m-1} \leq \Lambda\left(\frac{1}{c_0},\|\nabla v\|_{L^\infty}\right)\left( \|\mathcal{V}^m\| + |h|_{m-\frac{1}{2}} \right),
$$
$$
\|\partial_z B\cdot\mathbf{n}\|_{m-1} \leq \Lambda\left(\frac{1}{c_0},\|\nabla B\|_{L^\infty}\right)\left( \|\mathcal{B}^m\| + |h|_{m-\frac{1}{2}} \right).
$$
\end{lemma}
\begin{proof}
This is derived easily from divergence free condition of $v$ and $B$. See \cite{NMFR1}.
\end{proof}

\begin{lemma} \label{lemma 8.2}
For every integer $k \geq 0$, when we define,
$$
S^v_n := \Pi \mathbf{S}^\varphi v\mathbf{n},\quad S^B_n := \Pi \mathbf{S}^\varphi B\mathbf{n},\quad \Pi := I - \mathbf{n}\otimes\mathbf{n}.
$$
Then we get
\begin{equation} \label{8.1}
\begin{split}
\|\partial_z v\|_k &\leq \Lambda\left(\frac{1}{c_0},\|\nabla v\|_{L^\infty}\right)\left( \|S^v_n\|_k + |h|_{k+\frac{1}{2}} + \|v\|_{k+1} \right),  \\
\|\partial_{zz}v\|_k &\leq \Lambda\left(\frac{1}{c_0},\|v\|_{E^{2,\infty}}\right)\left(\|\nabla S^v_n\|_k + |h|_{k+\frac{3}{2}} + \|v\|_{k+2}\right).
\end{split}
\end{equation}
This is exactly same as $B$, because these come from definition and divergence free condition.
\begin{equation} \label{8.2}
\begin{split}
\|\partial_z B\|_k &\leq \Lambda\left(\frac{1}{c_0},\|\nabla B\|_{L^\infty}\right)\left( \|S^B_n\|_k + |h|_{k+\frac{1}{2}} + \|B\|_{k+1} \right),  \\
\|\partial_{zz} B\|_k &\leq \Lambda\left(\frac{1}{c_0},\|B\|_{E^{2,\infty}}\right)\left(\|\nabla S^B_n\|_k + |h|_{k+\frac{3}{2}} + \|B\|_{k+2}\right).
\end{split}
\end{equation}
\end{lemma}
\begin{proof}
	We suffice to use definition (\ref{new deriva}), divergence free condition of $v$ and $B$, and Lemma \ref{lemma 8.1}. See \cite{NMFR1}.
\end{proof}
From this lemma, we estimate $S_{n}^{v}$, $\nabla S_{n}^{v}$, $S_{n}^{B}$, and $\nabla S_{n}^{B}$, instead of $\p_z v$, $\p_{zz} v$, $\p_z B$, $\p_{zz} B$. So we make equations of $S_{n}^{v}$ and $S_{n}^{B}$, and estimate them.
\begin{proposition} \label{proposition 8.3}
We have the following estimate.
\begin{equation}
\begin{split}
&\|S^v_n\|_{m-2}^2 + 2\varepsilon\int_0^T \|\nabla^\varphi S^v_n\|_{m-2}^2 ds - 2\int_0^T\int_S Z^\alpha S^v_n\cdot (B\cdot\nabla^\varphi)Z^\alpha S^B_n dV_t ds   \\
&\leq \Lambda_0\|S^v_n(0)\|_{m-2}^2 + \int_0^T \Lambda_\infty (s) \left(\|v\|_{E^m} + \|B\|_{E^m} + |h|_{m-\frac{1}{2}} + \sqrt{\varepsilon}|h|_{m+\frac{1}{2}}\right) \left(\|S^v_n\|_{m-2} + \|S^B_n\|_{m-2} + |h|_{m-\frac{1}{2}}\right) ds  \\
&+ 2\varepsilon\int_0^T |v^b|_{m+\frac{1}{2}}\|S^v_n\|_{m-2} ds + \Lambda_0\varepsilon\int_0^T \|\nabla S^v_n\|_{m-3}^2 ds.
\end{split}
\end{equation}
Note that from $|h|_{m-\frac{1}{2}}$, we cannot gain $m-1$ order estimate. $m-2$ is optimal regularity for $S_n^{v}$.
\end{proposition}
\begin{proof}
We apply $\nabla^\varphi$ to the system (\ref{1.4}) to get,
$$
\partial_t^\varphi\nabla^\varphi v + (v\cdot\nabla^\varphi)\nabla^\varphi v + (\nabla^\varphi v)^2 - (B\cdot\nabla^\varphi)\nabla^\varphi B - (\nabla^\varphi B)^2 + (D^\varphi)^2 q - \varepsilon\triangle^\varphi\nabla^\varphi v = 0.
$$
We take symmetric part, then we get ($\mathbf{D}^\varphi$ is symmetric part.)
\begin{equation*}
\begin{split}
&\partial_t^\varphi \mathbf{S}^\varphi v + (v\cdot\nabla^\varphi)\mathbf{S}^\varphi v + \frac{1}{2}\left( (\nabla^\varphi v)^2 + ((\nabla^\varphi v)^T)^2 \right)  \\
&\quad - (B\cdot\nabla^\varphi)\mathbf{S}^\varphi B - \frac{1}{2}\left( (\nabla^\varphi B)^2 + ((\nabla^\varphi B)^T)^2 \right) + (D^\varphi)^2 q - \varepsilon\triangle^\varphi(\mathbf{S}^\varphi v)= 0.
\end{split}
\end{equation*}
Applying $\cdot\mathbf{n}$ and tangential projection operator $\Pi$, we get
\begin{equation} \label{8.4}
\partial_t^\varphi S^v_n + (v\cdot\nabla^\varphi)S^v_n - (B\cdot\nabla^\varphi)S^B_n - \varepsilon\triangle^\varphi(S^v_n) = F_S,
\end{equation}
where $F_S := F_S^1 + F_S^2 + F_S^3$,
\begin{equation*}
\begin{split}
F_S^1 &:= -\frac{1}{2}\Pi\left( (\nabla^\varphi v)^2 + ((\nabla^\varphi v)^T)^2 \right)\mathbf{n} + (\partial_t\Pi + v\cdot\nabla^\varphi\Pi)\mathbf{S}^\varphi v\mathbf{n} + \Pi \mathbf{S}^\varphi v (\partial_t\mathbf{n} + v\cdot\nabla^\varphi\mathbf{n}),  \\
F_S^2 &:= -2\varepsilon\partial_i^\varphi\Pi\partial_i^\varphi(\mathbf{S}^\varphi v\mathbf{n}) - 2\varepsilon\Pi(\partial_i^\varphi(\mathbf{S}^\varphi v)\partial_i^\varphi\mathbf{n}) -\varepsilon(\triangle^\varphi\Pi)\mathbf{S}^\varphi v\mathbf{n} - \varepsilon\Pi \mathbf{S}^\varphi v \triangle^\varphi\mathbf{n} - \Pi((D^\varphi)^2 q)\mathbf{n},  \\
F_S^3 &:= \frac{1}{2}\Pi\left( (\nabla^\varphi B)^2 + ((\nabla^\varphi B)^T)^2 \right)\mathbf{n} - (B\cdot\nabla^\varphi\Pi)\mathbf{S}^\varphi B\mathbf{n} - \Pi \mathbf{S}^\varphi B ( B\cdot\nabla^\varphi\mathbf{n}).
\end{split}
\end{equation*}
Note that,
\begin{equation*}
\begin{split}
(\partial_t^\varphi \mathbf{S}^\varphi v)\mathbf{n} &:= \partial_t^\varphi(\mathbf{S}^\varphi v\mathbf{n}) - \mathbf{S}^\varphi v \partial_t\mathbf{n}, \\
\Pi(\partial_t^\varphi \mathbf{S}^\varphi v)\mathbf{n} &:= \Pi(\partial_t^\varphi(\mathbf{S}^\varphi v\mathbf{n}) - \mathbf{S}^\varphi v \partial_t\mathbf{n}) \\
&= \Pi\partial_t^\varphi(\mathbf{S}^\varphi v\mathbf{n}) - \Pi(\mathbf{S}^\varphi v \partial_t\mathbf{n}) = \partial_t^\varphi(\Pi \mathbf{S}^\varphi v \mathbf{n}) - (\partial_t\Pi)\mathbf{S}^\varphi\mathbf{n} - \Pi(\mathbf{S}^\varphi v\partial_t\mathbf{n}) \\
&= \partial_t^\varphi(S^v_n) - (\partial_t\Pi)\mathbf{S}^\varphi\mathbf{n} - \Pi \mathbf{S}^\varphi v\partial_t\mathbf{n}.
\end{split}
\end{equation*}
We can easily estimate $F_S^1$, $F_S^2$, and $F_S^3$.
\begin{equation} \label{8.5}
\begin{split}
\|F_S^1\|_{m-2} &\leq \Lambda_{\infty}(t) \left( \|S^v_n\|_{m-2} + |h|_{m-\frac{1}{2}} + \|v\|_{m-1} \right),  \\
\|F_S^2\|_{m-2} &\leq \Lambda_{\infty}(t) \varepsilon\left( \|\nabla^\varphi S^v_n\|_{m-2} + |h|_{m+\frac{1}{2}} + |v^b|_{m+\frac{1}{2}} \right) + \Lambda_{\infty}(t) \left( \|v\|_{E^m} + |h|_{m-\frac{1}{2}} \right),  \\
\|F_S^3\|_{m-2} &\leq \Lambda_{\infty}(t) \left( \|S^B_n\|_{m-2} + |h|_{m-\frac{1}{2}} + \|B\|_{m-1} \right).
\end{split}
\end{equation}

\begin{remark}
In above three estimates, the order of $v$ and $h$ give critical optimal criteria. For $h$,
\begin{equation*}
\begin{split}
F_S^1 &\sim \nabla\mathbf{n} \sim \nabla\nabla\varphi \sim |h|_{m-2+2-\frac{1}{2}} \sim |h|_{m-\frac{1}{2}},		\\
F_S^2 &\sim \varepsilon\triangle\Pi \sim \varepsilon\triangle\nabla\varphi \sim \varepsilon|\varphi|_{m-2+3} \sim \varepsilon|h|_{m+\frac{1}{2}}.
\end{split}
\end{equation*}
We already got full regularity of $h$, so cannot raise its order. $F_S^3$ is similar as $F_S^1$.
For $v$,
\begin{equation*}
\begin{split}
	F_S^1 &\sim \|v\|_{m-1},\|S^v_n\|_{m-2}, 	\\
	F_S^2 &\sim \Lambda_{\infty}(t) \|\nabla q^E\|_{E^{m-1}} \sim \|v\|_{E^m}.	
\end{split}
\end{equation*}
Regularity of $v$ in $F_S^2$ is also maximal, although we cannot try ($m-1$) order. $F_S^3$ is similar.
\end{remark}

Note that from boundary compatibility condition, we have,
\begin{equation} \label{Snv bdry cond}
	S^v_n|_{z=0} = 0.
\end{equation}
Therefore, from revised basic $L^2$ energy estimate, Proposition \ref{L2 est}, 
\begin{equation} \label{8.6}
\frac{1}{2}\frac{d}{dt}\int_S |S^v_n|^2 dV_t + \varepsilon\int_S |\nabla^\varphi S^v_n|^2 dV_t = \int_S F_S\cdot S^v_n dV_t + \int_S S^v_n\cdot (B\cdot\nabla^\varphi)S^B_n dV_t.
\end{equation}
Applying $Z^\alpha$ to (\ref{8.4}), we get
$$
\partial_t^\varphi Z^\alpha S^v_n + (v\cdot\nabla^\varphi)Z^\alpha S^v_n - (B\cdot\nabla^\varphi)Z^\alpha S^B_n - \varepsilon\triangle^\varphi Z^\alpha(S^v_n) = Z^\alpha(F_S) + \mathcal{C}_S,
$$
where 
$$
\mathcal{C}_S := \mathcal{C}_S^1 + \mathcal{C}_S^2 + \mathcal{C}_S^3,
$$
with
\begin{equation*}
\begin{split}
\mathcal{C}_S^1 &:= [Z^\alpha v_y]\cdot\nabla_y S^v_n + [Z^\alpha,V_z]\partial_z S^v_n := \mathcal{C}^1_{S_y} + \mathcal{C}^1_{S_z},  \\
\mathcal{C}_S^2 &:= \varepsilon[Z^\alpha,\triangle^\varphi]S^v_n,  \\
\mathcal{C}_S^3 &:= -[Z^\alpha B_y]\cdot\nabla_y S^B_n + [Z^\alpha,\frac{B\cdot N}{\partial_z\varphi}]\partial_z S^B_n := \mathcal{C}^3_{S_y} + \mathcal{C}^3_{S_z}.
\end{split}
\end{equation*}
High order estimate becomes,
\begin{equation*}
\frac{1}{2}\frac{d}{dt}\int_S |Z^\alpha S^v_n|^2 dV_t + \varepsilon\int_S |Z^\alpha\nabla^\varphi S^v_n|^2 dV_t = \int_S (Z^\alpha(F_S) + \mathcal{C}_S)\cdot Z^\alpha S^v_n dV_t + \int_S Z^\alpha S^v_n\cdot (B\cdot\nabla^\varphi)Z^\alpha S^B_n dV_t.
\end{equation*}
Note that $Z^{\a} S_{n}^{v} \vert_{z=0} = 0$ by (\ref{Snv bdry cond}) and definition $Z_{3} := {z\over 1-z} \p_{z}$ and hence, boundary integration vanishes. The last term will be canceled with similar term from faraday's Law. Estimates of $\mathcal{C}_S^1$ and $\mathcal{C}_S^2$ are given in \cite{NMFR1}, using some variants of Hardy's inequality (Lemma 8.4 in \cite{NMFR1}), which is valid only when function is zero at $z=0$. It is important that we have such condition of $S_{n}^{v}$, from continuity of stress tensor boundary condition, i.e. fifth equation in (\ref{1.4}).
\begin{equation*}
\begin{split}
\|\mathcal{C}_S^1\| &\leq \Lambda_\infty (t) (\|S_n^v\|_{m-2} + \|v\|_{E^{m-1}} + |h|_{m-\frac{1}{2}} ),  \\
\left| \int_S \mathcal{C}_S^2\cdot Z^\alpha S_n^v dV_t \right| &\leq \Lambda_0 (\varepsilon^{1/2}\|\nabla Z^\alpha S_n^v\| + \|S_n\|_{m-2} )  \\
&\quad \times \left( \varepsilon^{1/2}\|\nabla S_n^v\|_{m-3} + \|S_n\|_{m-2} + \Lambda_\infty (t)(|h|_{m-\frac{3}{2}} + \varepsilon^{1/2}|h|_{m-\frac{1}{2}})\right).
\end{split}
\end{equation*}
Differences to \cite{NMFR1} are
$$
\int_S Z^\alpha(F_S^3)\cdot Z^\alpha S^v_n dV_t,\quad\text{and}\quad\int_S \mathcal{C}^3_S \cdot Z^\alpha S^v_n dV_t.
$$
Since structures of $F_S^3,\,\mathcal{C}^3_S$ are nearly similar to $F_S^1,\,\mathcal{C}^1_S$, we get similar results, by replacing $v$ into $B$. Hence, when $|\alpha| = m-2$,
\begin{equation*}
\begin{split}
&\frac{1}{2}\frac{d}{dt}\int_S |Z^\alpha S^v_n|^2 dV_t + \varepsilon\int_S |Z^\alpha\nabla^\varphi S^v_n|^2 dV_t \leq \int_S Z^\alpha S^v_n\cdot (B\cdot\nabla^\varphi)Z^\alpha S^B_n dV_t  \\
&\quad + \Lambda_\infty (t)\left(\|v\|_{E^m} + \|B\|_{E^m} + |h|_{m-\frac{1}{2}} + \sqrt{\varepsilon}|h|_{m+\frac{1}{2}}\right) \left(\|S^v_n\|_{m-2} + \|S^B_n\|_{m-2} + |h|_{m-\frac{1}{2}}\right)   \\
&\quad + \varepsilon |v^b|_{m+\frac{1}{2}}\|S^v_n\|_{m-2} + \Lambda_0\varepsilon\|\nabla S^v_n\|_{m-3}^2.
\end{split}
\end{equation*}
This implies,
\begin{equation*}
\begin{split}
&\|S^v_n\|_{m-2}^2 + 2\varepsilon\int_0^T \|\nabla^\varphi S^v_n\|_{m-2}^2 ds - 2\int_0^T\int_S Z^\alpha S^v_n\cdot (B\cdot\nabla^\varphi)Z^\alpha S^B_n dV_t ds   \\
&\leq \Lambda_0\|S^v_n(0)\|_{m-2}^2 + \int_0^T \Lambda_\infty (s) \left(\|v\|_{E^m} + \|B\|_{E^m} + |h|_{m-\frac{1}{2}} + \sqrt{\varepsilon}|h|_{m+\frac{1}{2}}\right) \left(\|S^v_n\|_{m-2} + \|S^B_n\|_{m-2} + |h|_{m-\frac{1}{2}}\right) ds  \\
&\quad + 2\varepsilon\int_0^T |v^b|_{m+\frac{1}{2}}\|S^v_n\|_{m-2} ds + \Lambda_0\varepsilon\int_0^T \|\nabla S^v_n\|_{m-3}^2 ds.
\end{split}
\end{equation*}
\end{proof}

We perform similar estimate for Faraday's law.
\begin{proposition} \label{proposition 8.5}
We have the estimate. 
\begin{equation} \label{8.7}
\begin{split}
&\|S^B_n\|_{m-2}^2 + 2\lambda\int_0^T \|\nabla^\varphi S^B_n\|_{m-2}^2 ds + 2\int_0^T\int_S Z^\alpha S^v_n\cdot (B\cdot\nabla^\varphi)Z^\alpha S^B_n dV_t ds   \\ 
&\leq \Lambda_0\|S^B_n(0)\|_{m-2}^2 + \int_0^T \Lambda_\infty (s)\left(\|v\|_{E^m} + \|B\|_{E^m} + |h|_{m-\frac{1}{2}} + \sqrt{\lambda}|h|_{m+\frac{1}{2}}\right) \left(\|S^v_n\|_{m-2} + \|S^B_n\|_{m-2} + |h|_{m-\frac{1}{2}}\right) ds  \\
&\quad + \Lambda_0\lambda\int_0^T \|\nabla S^v_n\|_{m-3}^2 ds.
\end{split}
\end{equation}
\end{proposition}
\begin{proof}
Applying $\nabla^\varphi$ to the equation, we get
$$
\partial_t^\varphi\nabla^\varphi B + (v\cdot\nabla^\varphi)\nabla^\varphi B + (\nabla^\varphi B)(\nabla^\varphi v) - (B\cdot\nabla^\varphi)\nabla^\varphi v - (\nabla^\varphi v)(\nabla^\varphi B) - \lambda\triangle^\varphi\nabla^\varphi B = 0.
$$
We take transpose to get 
$$
\partial_t^\varphi \mathbf{S}^\varphi B + (v\cdot\nabla^\varphi)\mathbf{S}^\varphi B - (B\cdot\nabla^\varphi)\mathbf{S}^\varphi v - \lambda\triangle^\varphi(\mathbf{S}^\varphi B) 
$$
$$
+ \frac{1}{2}\left( (\nabla^\varphi B)(\nabla^\varphi v) - (\nabla^\varphi v)(\nabla^\varphi B) + (\nabla^\varphi v)^T(\nabla^\varphi B)^T - (\nabla^\varphi B)^T(\nabla^\varphi v)^T \right) = 0.
$$
Applying $\mathbf{n}$ and $\Pi$, we get
$$
\partial_t^\varphi S^B_n + (v\cdot\nabla^\varphi)S^v_B - (B\cdot\nabla^\varphi)S^v_n - \lambda\triangle^\varphi(S^B_n) = E_S,
$$
where $E_S = E_S^1 + E_S^2 + E_S^3$ and
\begin{equation} \label{8.8}
\begin{split}
E_S^1 &:= (\partial_t\Pi + v\cdot\nabla^\varphi\Pi)\mathbf{S}^\varphi B\mathbf{n} + \Pi \mathbf{S}^\varphi B (\partial_t\mathbf{n} + v\cdot\nabla^\varphi\mathbf{n}) \\
&\quad - \frac{1}{2}\left( (\nabla^\varphi B)(\nabla^\varphi v) - (\nabla^\varphi v)(\nabla^\varphi B) + (\nabla^\varphi v)^T(\nabla^\varphi B)^T - (\nabla^\varphi B)^T(\nabla^\varphi v)^T \right) = 0,  \\
E_S^2 &:= -2\lambda\partial_i^\varphi\Pi\partial_i^\varphi(\mathbf{S}^\varphi B\mathbf{n}) - 2\lambda\Pi(\partial_i^\varphi(\mathbf{S}^\varphi B)\partial_i^\varphi\mathbf{n}) -\lambda(\triangle^\varphi\Pi)\mathbf{S}^\varphi B\mathbf{n} - \lambda\Pi \mathbf{S}^\varphi B \triangle^\varphi\mathbf{n},  \\
E_S^3 &:= - (B\cdot\nabla^\varphi\Pi)\mathbf{S}^\varphi v\mathbf{n} - \Pi \mathbf{S}^\varphi v ( B\cdot\nabla^\varphi\mathbf{n}).
\end{split}
\end{equation}
We estimate $E_S^1$, $E_S^2$, and $E_S^3$. 
\begin{equation} \label{8.9}
\begin{split}
\|E_S^1\|_{m-2} &\leq \Lambda_{\infty}(t) \left( |h|_{m-\frac{1}{2}} + \|S^v_n\|_{m-2} + \|S^B_n\|_{m-2} + \|B\|_{m-1} + \|v\|_{m-2} \right),  \\
\|E_S^2\|_{m-2} &\leq \Lambda_{\infty}(t) \lambda\left( \|\nabla^\varphi S^B_n\|_{m-2} + |h|_{m+\frac{1}{2}} \right),  \\
\|E_S^3\|_{m-2} &\leq \Lambda_{\infty}(t) \left( |h|_{m-\frac{1}{2}} + \|S^v_n\|_{m-2} + \|B\|_{m-2} + \|v\|_{m-1} \right).
\end{split}
\end{equation}
Since $B$ is uniformly zero in $\mathbb{R}^{3}\backslash S$, we have boundary condition 
\begin{equation} \label{SnB bdry cond}
	S^B_n|_{z=0} = 0.
\end{equation}
And therefore $Z^{\alpha} S_{n}^{B} = 0$ holds similar as $S_{n}^{v}$ case. It is easy to get $L^2$ estimate
\begin{equation} \label{8.10}
\frac{1}{2}\frac{d}{dt}\int_S |S^B_n|^2 dV_t + \lambda\int_S |\nabla^\varphi S^B_n|^2 dV_t = \int_S E_S\cdot S^v_N dV_t + \int_S S^B_n\cdot (B\cdot\nabla^\varphi)S^v_n dV_t.
\end{equation}
Applying $Z^\alpha$, we get
$$
\partial_t^\varphi Z^\alpha S^B_n + (v\cdot\nabla^\varphi)Z^\alpha S^B_n - (B\cdot\nabla^\varphi)Z^\alpha S^v_n - \lambda\triangle^\varphi Z^\alpha(S^B_n) = Z^\alpha(E_S) + \bar{\mathcal{C}}_S,
$$
where 
$$
\bar{\mathcal{C}}_S := \bar{\mathcal{C}}_S^1 + \bar{\mathcal{C}}_S^2 + \bar{\mathcal{C}}_S^3,
$$
with
\begin{equation}
\begin{split}
\bar{\mathcal{C}}_S^1 &:= [Z^\alpha v_y]\cdot\nabla_y S^B_n + [Z^\alpha,V_z]\partial_z S^B_n := \bar{\mathcal{C}}^1_{S_y} + \bar{\mathcal{C}}^1_{S_z},  \\
\bar{\mathcal{C}}_S^2 &:= \lambda[Z^\alpha,\triangle^\varphi]S^B_n,  \\
\bar{\mathcal{C}}_S^3 &:= -[Z^\alpha B_y]\cdot\nabla_y S^v_n + [Z^\alpha,\frac{B\cdot N}{\partial_z\varphi}]\partial_z S^v_n := \bar{\mathcal{C}}^3_{S_y} + \bar{\mathcal{C}}^3_{S_z}.
\end{split}
\end{equation}
High order estimate becomes,
\begin{equation*}
\frac{1}{2}\frac{d}{dt}\int_S |Z^\alpha S^B_n|^2 dV_t + \lambda\int_S |Z^\alpha\nabla^\varphi S^B_n|^2 dV_t = \int_S (Z^\alpha(E_S) + \bar{\mathcal{C}}_S)\cdot Z^\alpha S^B_n dV_t + \int_S Z^\alpha S^B_n\cdot (B\cdot\nabla^\varphi)Z^\alpha S^v_n dV_t.
\end{equation*}
Estimates of these terms are similar, so when $\alpha = m-2$, we get
\begin{equation*}
\begin{split}
&\frac{1}{2}\frac{d}{dt}\int_S |Z^\alpha S^B_n|^2 dV_t + \lambda\int_S |Z^\alpha\nabla^\varphi S^B_n|^2 dV_t \leq \int_S Z^\alpha S^B_n\cdot (B\cdot\nabla^\varphi)Z^\alpha S^v_n dV_t 	\\
&\quad + \Lambda_\infty (t)\left(\|v\|_{E^m} + \|B\|_{E^m} + |h|_{m-\frac{1}{2}} + \sqrt{\lambda}|h|_{m+\frac{1}{2}}\right) \left(\|S^v_n\|_{m-2} + \|S^B_n\|_{m-2} + |h|_{m-\frac{1}{2}}\right) 	\\
&\quad + \Lambda_0\lambda\|\nabla S^B_n\|_{m-3}^2.
\end{split}
\end{equation*}
(Comparing with previous Proposition \ref{proposition 8.3}, $\lambda |v^b|_{m+\frac{1}{2}}\|S^B_n\|_{m-2}$ does not appear, since it comes from pressure estimate $q^{NS}$. However, there is no pressure term in Faraday's law.)
\begin{equation*}
\begin{split}
&\|S^B_n\|_{m-2}^2 + 2\lambda\int_0^T \|\nabla^\varphi S^B_n\|_{m-2}^2 ds - 2\int_0^T\int_S Z^\alpha S^B_n\cdot (B\cdot\nabla^\varphi)Z^\alpha S^v_n dV_t ds   \\
&\leq \Lambda_0\|S^B_n(0)\|_{m-2}^2 + \int_0^T \Lambda_\infty(s)\left(\|v\|_{E^m} + \|B\|_{E^m} + |h|_{m-\frac{1}{2}} + \sqrt{\lambda}|h|_{m+\frac{1}{2}}\right) \left(\|S^v_n\|_{m-2} + \|S^B_n\|_{m-2} + |h|_{m-\frac{1}{2}}\right) ds  \\
&\quad + \Lambda_0\lambda\int_0^T \|\nabla S^v_n\|_{m-3}^2 ds. 
\end{split}
\end{equation*}
To make cancellation with previous proposition, we calculate
$$
\int_S f\cdot (B\cdot\nabla^\varphi)g dV_t =\int_{\partial S} f\cdot g (B\cdot \mathbf{n}) - \int_S \nabla^\varphi\cdot B (f\cdot g) - \int_S g\cdot (B\cdot\nabla^\varphi)f = - \int_S g\cdot (B\cdot\nabla^\varphi)f.
$$
Hence,
\begin{equation*}
\begin{split}
&\|S^B_n\|_{m-2}^2 + 2\lambda\int_0^T \|\nabla^\varphi S^B_n\|_{m-2}^2 ds + 2\int_0^T\int_S Z^\alpha S^v_n\cdot (B\cdot\nabla^\varphi)Z^\alpha S^B_n dV_t ds   \\
&\leq \Lambda_0\|S^B_n(0)\|_{m-2}^2 + \int_0^T \Lambda_{\infty}(s)\left(\|v\|_{E^m} + \|B\|_{E^m} + |h|_{m-\frac{1}{2}} + \sqrt{\lambda}|h|_{m+\frac{1}{2}}\right) \left(\|S^v_n\|_{m-2} + \|S^B_n\|_{m-2} + |h|_{m-\frac{1}{2}}\right) ds  \\
&\quad + \Lambda_0\lambda\int_0^T \|\nabla S^v_n\|_{m-3}^2 ds.
\end{split}
\end{equation*}
\end{proof}
Now we add two Proposition \ref{proposition 8.3} and \ref{proposition 8.5} to  cancel the last terms on the LHS. And then we use induction for $\varepsilon\int_0^T \|\nabla S^v_n\|_{m-3}^2 ds,\lambda\int_0^T \|\nabla S^v_n\|_{m-3}^2 ds$.

\begin{proposition} \label{8.6}
We have the following estimate.
\begin{equation} \label{8.11}
\begin{split}
&\|S^v_n\|_{m-2}^2 + \|S^B_n\|_{m-2}^2 + 2\varepsilon\int_0^T \|\nabla^\varphi S^v_n\|_{m-2}^2 ds + 2\lambda\int_0^T \|\nabla^\varphi S^B_n\|_{m-2}^2 ds   \\
&\leq \Lambda_0 \left( \|S^v_n(0)\|_{m-2}^2 + \|S^B_n(0)\|_{m-2}^2 \right)  \\
&\quad + \int_0^T \Lambda_{\infty}(s) \left(\|\mathcal{V}^m\|^2 + \|\mathcal{B}^m\|^2 + \|S^v_n\|_{m-2}^2 + \|S^B_n\|_{m-2}^2 + |h|_m^2 + \varepsilon |h|_{m+\frac{1}{2}}^2 \right)  \\
&\quad + \int_0^T \Lambda_{\infty}(s) \left( \|\partial_z v\|_{m-1}^2 + \|\partial_z B\|_{m-1}^2 \right) + \varepsilon\int_0^T \|\nabla\mathcal{V}^m\|^2. 
\end{split}
\end{equation}
\end{proposition}
\begin{proof}
We sum Proposition \ref{proposition 8.3} and \ref{proposition 8.5} to get
\begin{equation*}
\begin{split}
&\|S^v_n\|_{m-2}^2 + \|S^B_n\|_{m-2}^2 + 2\varepsilon\int_0^T \|\nabla^\varphi S^v_n\|_{m-2}^2 ds + 2\lambda\int_0^T \|\nabla^\varphi S^B_n\|_{m-2}^2 ds   \\
&\leq \Lambda_0 \left( \|S^v_n(0)\|_{m-2}^2 + \|S^B_n(0)\|_{m-2}^2 \right)  \\
&\quad + \int_0^T \Lambda_{\infty}(s) \left(\|v\|_{E^m} + \|B\|_{E^m} + |h|_{m-\frac{1}{2}} + \sqrt{\varepsilon}|h|_{m+\frac{1}{2}}\right) \left(\|S^v_n\|_{m-2} + \|S^B_n\|_{m-2} + |h|_{m-\frac{1}{2}}\right) ds  \\
&\quad + 2\varepsilon\int_0^T |v^b|_{m+\frac{1}{2}}\|S^v_n\|_{m-2} ds + \Lambda_0\left( \varepsilon\int_0^T \|\nabla S^v_n\|_{m-3}^2 ds + \lambda\int_0^T \|\nabla S^B_n\|_{m-3}^2 ds \right).
\end{split}
\end{equation*}
Using trace estimate Proposition \ref{proposition 3.2}, we have,
$$
|v^b|_{m+\frac{1}{2}} \leq \|\nabla\mathcal{V}^m\| + \|\mathcal{V}\|_m + \Lambda_{\infty}(t) |h|_{m+\frac{1}{2}},
$$
and therefore,
$$
\varepsilon\int_0^T |v^b|_{m+\frac{1}{2}}\|S^v_n\|_{m-2} \leq \varepsilon\int_0^T \|\nabla\mathcal{V}^m\|^2 + \int_0^T \Lambda_{\infty}(s) \left( \|\mathcal{V}^m\|^2 + \|S^v_n\|_{m-2}^2 + \varepsilon |h|_{m+\frac{1}{2}}^2 \right).
$$
We can replace $\|\nabla^\varphi S^v_n\|_k$ by $\|\nabla S^v_n\|_k$, and use induction for $\left( \varepsilon\int_0^T \|\nabla S^v_n\|_{m-3}^2 ds + \lambda\int_0^T \|\nabla S^B_n\|_{m-3}^2 ds \right)$ to finish proof.
\end{proof}

\section{$L^\infty$ type estimate}
From this section, we set $\varepsilon=\lambda$. We estimate all $L^\infty$ type terms in $\Lambda_{\infty}(t)$. First we state basic properties of Proposition 9.1 in \cite{NMFR1}. 
\begin{proposition}
	We have the following estimates.
	\begin{equation}
	\begin{split}
		|h|_{k,\infty} + \sqrt{\varepsilon}|h|_{k+1,\infty} &\lesssim |h|_{2+k} + \sqrt{\varepsilon}|h|_{2+k+\frac{1}{2}},\quad k\in\mathbb{N},	\\
		\|v(t)\|_{2,\infty} &\leq \Lambda\big(\frac{1}{c_0}, |h|_{4,\infty} + \|\mathcal{V}^{4}\| + \|S_n^{v}\|_{3} \big),	\\
		\|\p_{z} v\|_{1,\infty} &\leq \Lambda_{0} \big( \|S_{n}^{v}\|_{1,\infty} + \|v\|_{2,\infty} \big),	\\
		\sqrt{\varepsilon}\|\p_{zz}v\|_{\infty} &\leq \Lambda_{0}\big( \sqrt{\varepsilon}\|\p_{z}S_n^v\|_{\infty} + \|S_n^v\|_{1,\infty} + \|v\|_{2,\infty} \big),
	\end{split}
	\end{equation}
and estimate 
Note that last three inequalities hold for $B$ version because they come from Proposition \ref{proposition 3.2}, structure of $\mathcal{V}^{\a}, \mathcal{B}^{\a}$, interpolation, and Young's inequality.
\end{proposition}

Inspired by above proposition, we define the following quantity.
\begin{equation*}
\begin{split}
\mathcal{Q}_m (t) &:= |h|_{m}^2 + \varepsilon |h|_{m+\frac{1}{2}}^2 + \|\mathcal{V}^m\|^2 + \|\mathcal{B}^m\|^2 + \|S^v_n\|_{m-2}^2 + \|S^v_n\|_{1,\infty}^2 + \varepsilon\|\partial_z S^v_n\|_{L^\infty}^2  \\
&\quad + \|S^B_n\|_{m-2}^2 + \|S^B_n\|_{1,\infty}^2 + \varepsilon\|\partial_z S^B_n\|_{L^\infty}^2.
\end{split}
\end{equation*}
 
From Proposition \ref{proposition 7.3}, we know that we should estimate $\|\partial_z v\|_{m-1}, \|\partial_z B\|_{m-1}$. But we have only $m-2$ order estimate in section 8. We will get $\|\partial_z v\|_{m-1}, \|\partial_z B\|_{m-1}$ at section 10. In this section, we control some $L^\infty$ type terms here and $\|S^v_n\|_{m-2}, \|S^B_n\|_{m-2}$ would be sufficient to estimate them, for sufficiently large $m$. First we state a corollary which resembles corollary 9.3 in \cite{NMFR1}. These terms come from $L^\infty$ type terms in $\Lambda_{\infty}(t)$.
\begin{corollary} \label{corollary 9.1}
When $m \geq 6$, for each time $t$,
\begin{equation} \label{9.1}
\|v\|_{2,\infty} + \|\partial_z v\|_{1,\infty} + \sqrt{\varepsilon}\|\partial_{zz}v\|_{L^\infty} +  \|B\|_{2,\infty} + \|\partial_z B\|_{1,\infty} + \sqrt{\varepsilon}\|\partial_{zz}B\|_{L^\infty} + |h|_{4,\infty} \leq \Lambda\left(\frac{1}{c_0},\mathcal{Q}_m\right).
\end{equation}
\end{corollary}
\noindent Therefore, by above Corollay \ref{corollary 9.1}, we get
\begin{equation} \label{infty equiv}
\begin{split}
	\Lambda_{m,\infty}(t) &:= \Lambda (\frac{1}{c_0},\|v\|_{m} + \|\partial_z v\|_{m-2} + \|B\|_{m} + \|\partial_z B\|_{m-2} + |h|_m + \sqrt{\varepsilon}|h|_{m+\frac{1}{2}} + |h|_{4,\infty}
	\\
	& \quad + \|v\|_{E^{2,\infty}} + \sqrt{\varepsilon}\|\partial_{zz}v\|_{L^\infty} + \|B\|_{E^{2,\infty}} + \sqrt{\varepsilon}\|\partial_{zz}B\|_{L^\infty} ) \leq \Lambda\left(\frac{1}{c_0},\mathcal{Q}_m\right),\,\,\,\,m \geq 6.
\end{split}
\end{equation}

\noindent Hence we control $\mathcal{Q}_m$ instead of $\Lambda_{m,\infty}$. Now we start with estimates of $\|S^v_n\|_{1,\infty}$ and $\|S^B_n\|_{1,\infty}$.

\subsection{Zero order estimate}
First we note that zero-order estimate is just come from maximal principle of transport equation, estimate (\ref{9.12}), but we should be careful that we treat the term from Lorentz force. From (\ref{8.4}),
$$
\partial_t^\varphi S^v_n + (v\cdot\nabla^\varphi)S^v_n - (B\cdot\nabla^\varphi)S^B_n - \varepsilon\triangle^\varphi(S^v_n) = F_S,
$$
and therefore we get,
$$
\|S^v_n(t)\|_{L^\infty} \leq \|S^v_n(0)\|_{L^\infty} + \int_0^t \left( \|F_S\|_{L^\infty} + \|(B\cdot\nabla)S^B_n\|_{L^\infty} \right).
$$
Main problem is that $\|\nabla S^B_n\|_{L^\infty}$ is not controlled by $\Lambda_{m,\infty}$. In fact, we need $\|\partial_{zz}B\|_{L^\infty}$ to control this. This is same for Faraday law. From the equation of $S_n^B$,
$$
\partial_t^\varphi S^B_n + (v\cdot\nabla^\varphi)S^v_B - (B\cdot\nabla^\varphi)S^v_n - \varepsilon\triangle^\varphi(S^B_n) = E_S
$$
and we get
$$
\|S^B_n(t)\|_{L^\infty} \leq \|S^B_n(0)\|_{L^\infty} + \int_0^t \left( \|E_S\|_{L^\infty} + \|(B\cdot\nabla)S^v_n\|_{L^\infty} \right).
$$
We also need $\|\partial_{zz}v\|_{L^\infty}$ to control $\|(B\cdot\nabla)S^v_n\|_{L^\infty}$. Instead of attacking $S^v_n$ and $S^B_n$ separately, we treat these terms by adding and subtracting two equations, \begin{equation} \label{9.2}
\begin{split}
&\partial_t^\varphi(S^v_n + S^B_n) + ((v-B)\cdot\nabla^\varphi)(S^v_n + S^B_n) - \varepsilon\triangle^\varphi(S^v_n + S^B_n) = F_S + E_S  \\
&\partial_t^\varphi(S^v_n - S^B_n) + ((v+B)\cdot\nabla^\varphi)(S^v_n - S^B_n) - \varepsilon\triangle^\varphi(S^v_n - S^B_n) = F_S - E_S.
\end{split}
\end{equation}
Using maximal principle for (\ref{9.2}), we get
\begin{equation} \label{9.4}
\|(S^v_n + S^B_n)(t)\|_{L^\infty} \leq \|(S^v_n + S^B_n)(0)\|_{L^\infty} + \int_0^t \|F_S + E_S\|_{L^\infty},
\end{equation}
where
$$
\int_0^t \|F_S\|_{L^\infty} \leq \varepsilon\int_0^t \|\mathbf{S}^\varphi\mathcal{V}^m\|^2 + (1+\varepsilon)\int_0^t\Lambda_{m,\infty}(s),\quad \int_0^t \|E_S\|_{L^\infty} \leq \int_0^t \Lambda_{m,\infty}(s).
$$
This result is nearly same for $(S^v_n - S^B_n)$ case. Therefore,
\begin{equation} \label{9.5}
\|(S^v_n + S^B_n)(t)\|_{L^\infty} \lesssim \|S^v_n(0)\|_{L^\infty} + \|S^B_n(0)\|_{L^\infty} + \varepsilon\int_0^t \|\mathbf{S}^\varphi\mathcal{V}^m\|^2 + (1+\varepsilon)\int_0^t\Lambda_{m,\infty},
\end{equation}
\begin{equation} \label{9.6}
\|(S^v_n - S^B_n)(t)\|_{L^\infty} \leq \|S^v_n(0)\|_{L^\infty} + \|S^B_n(0)\|_{L^\infty} + \varepsilon\int_0^t \|\mathbf{S}^\varphi\mathcal{V}^m\|^2 + (1+\varepsilon)\int_0^t\Lambda_{m,\infty}.
\end{equation}
By (\ref{9.5}) and (\ref{9.6}),
\begin{equation} \label{9.7}
\|S^v_n(t)\|_{L^\infty}, \ \|S^B_n(t)\|_{L^\infty} \leq \|S^v_n(0)\|_{L^\infty} + \|S^B_n(0)\|_{L^\infty} + \varepsilon\int_0^t \|\mathbf{S}^\varphi\mathcal{V}^m\|^2 + (1+\varepsilon)\int_0^t\Lambda_{m,\infty}.
\end{equation}

\subsection{First order estimate}
To estimate first order terms, we divide thin layer near the boundary of S, then we can apply sobolev embedding to lower part, since it loose essential information for conormal norms. Main part is $L^\infty$ estimate for $\|\chi Z S^v_n\|_{L^\infty}$, where $\chi$ is zero away from thin layer near the boundary. Here, we know that direct maximal principle of transport equation is not good way because of commutators between $Z$ and $\triangle^\varphi$. \\

Let us define transformation $\Psi$,
\begin{equation} \label{9.8}
\Psi(t,\cdot):S = \mathbb{R}^2 \times (-\infty,0) \rightarrow \Omega(t),
\end{equation}
$$
x=(y,z) \mapsto \begin{pmatrix} y\\ h(t,y) \end{pmatrix} + z\mathbf{n}^b(t,y),
$$
where $\mathbf{n}^b$ is unit normal at the boundary, i.e. $(-\nabla_{y} h,1)/|\mathbf{N}|$. To ensure that this is diffeomorphism near the boundary, we check
$$
D\Psi(t,\cdot) = \begin{pmatrix} 1 & 0 & -\partial_1 h \\ 
                                 0 & 1 & -\partial_2 h \\
                                 \partial_1 h & \partial_2 h & 1 \end{pmatrix} + \begin{pmatrix} -z\partial_{11} h & -z\partial_{12} h & 0\\ 
                                -z\partial_{21} h & -z\partial_{22} h & 0\\
                                 0 & 0 & 1 \end{pmatrix}.
$$
This is diffeomorphism near the boundary since norm of second matrix is controlled by $|h|_{2,\infty}$. So, we restrict $\Psi(t,\cdot)$ on $\mathbb{R}^2 \times (-\delta,0)$ so that it is diffeomorphism. (Note that $\delta$ depends on $c_0$. Function $\chi$ is gained by $\chi(z) = \kappa(\frac{z}{\delta(c_0)}) \in [0,1]$, where $\kappa$ is smooth compactly supported function which is $1$ near $\p S$. Next, we write laplacian $\triangle^\varphi$ with respect to Riemannian metric of above parametrization. Riemannian metric becomes,
\begin{equation} \label{9.9}
g(y,z) = \begin{pmatrix}
\tilde{g}(y,z) & 0 \\
0 & 1 
\end{pmatrix},
\end{equation}
where $\tilde{g}$ is $2\times 2 $ block matrix. And with this metric, laplacian becomes,
\begin{equation} \label{9.10}
\triangle_g f = \partial_{zz} f + \frac{1}{2}\partial_z(\text{ln}|g|)\partial_z f + \triangle_{\tilde{g}} f,
\end{equation}
where
\begin{equation*}
\triangle_{\tilde{g}} f = \frac{1}{|\tilde{g}|^{\frac{1}{2}}}\sum_{1\leq i,j \leq 2} \partial_{y^i}(\tilde{g}^{ij}|\tilde{g}|^{\frac{1}{2}}\partial_{y^j} f),
\end{equation*}
where $\tilde{g}^{ij}$ is inverse matrix element of $\tilde{g}$. Notice that this map is invertible near the boundary, thin layer of thickness $\delta$ which depends on $c_0$. And we localize $\mathbf{S}^\varphi v$ by multiplying $\chi(z) = \kappa(\frac{z}{\delta(c_0)})$, that means this is 1 at thin layer and then smoothly decay to zero. We define 
\begin{equation} \label{9.11}
S_v^\chi := \chi(z)\mathbf{S}^\varphi v \quad\text{and}\quad S_B^\chi := \chi(z)\mathbf{S}^\varphi B.
\end{equation}
We find the equation for the $S_v^\chi$ and $S_B^\chi$,
\begin{equation} \label{9.12}
\partial_t^\varphi S_v^\chi + (v\cdot\nabla^\varphi)S_v^\chi - (B\cdot\nabla^\varphi)S_B^\chi - \varepsilon\triangle^\varphi(S_v^\chi) = F_{S^\chi} := F^\chi + F_v ,
\end{equation}
where
\begin{equation*}
\begin{split}
F^\chi &:= (V_z\partial_z\chi)\mathbf{S}^\varphi v - \left(\frac{B\cdot\mathbf{N}}{\partial_z\varphi}\partial_z\chi\right)\mathbf{S}^\varphi B - \varepsilon\nabla^\varphi\chi\cdot\nabla^\varphi \mathbf{S}^\varphi v - \varepsilon\triangle^\varphi\chi \mathbf{S}^\varphi v,  \\
F_v &:= -\chi(\mathbf{D}^\varphi)^2 q - \frac{\chi}{2}\left((\nabla^\varphi v)^2 + ((\nabla^\varphi v)^T)^2 \right) + \frac{\chi}{2}\left((\nabla^\varphi B)^2 + ((\nabla^\varphi B)^T)^2 \right).
\end{split}
\end{equation*}
Note that $F^\chi$ has $\mathbf{S}^\varphi v$ and $\mathbf{S}^\varphi B$ where as $F_v$ has only non symmetric parts. Note that $F^\chi$ is supported away from the boundary, because of $\nabla\chi$ and 
$$
\|F^\chi\|_{1,\infty} \leq \Lambda_{m,\infty} (t).
$$
For Faraday's law equation,
\begin{equation} \label{9.13}
\partial_t^\varphi S_B^\chi + (v\cdot\nabla^\varphi)S_B^\chi - (B\cdot\nabla^\varphi)S_v^\chi - \varepsilon\triangle^\varphi(S_B^\chi) = E_{S^\chi},
\end{equation}
where
$$
E_{S^\chi} := (V_z\partial_z\chi)\mathbf{S}^\varphi B - \left(\frac{B\cdot\mathbf{N}}{\partial_z\varphi}\partial_z\chi\right)\mathbf{S}^\varphi v - \varepsilon\nabla^\varphi\chi\cdot\nabla^\varphi \mathbf{S}^\varphi B - \varepsilon\triangle^\varphi\chi \mathbf{S}^\varphi B.
$$
Note that for Faraday's law, equation of $S^B_n$ is much simpler than $S^v_n$ and we do not have $F_v$ type commutators. Note that $E_{S^\chi}$ is supported away from the boundary because of $\nabla\chi$ and 
$$
\|E_{S^\chi}\|_{1,\infty} \leq \Lambda_{m,\infty} (t).
$$
Now, we define $\tilde{S}_{v}$, $\tilde{S}_{B}$ in $\Omega(t)$ and $S^\chi_{v}$, $S^\chi_{B}$, $S^\Psi_{v}$, and $S^\Psi_{B}$ on $S$ which are localized $\mathbf{S}^\varphi (v)$ and $\mathbf{S}^\varphi (B)$ near the boundary. $S^\Psi_{v}$ and $S^\Psi_{B}$ are main terms to estimate and this measure $S^\chi_{v} = \mathbf{S}_{loc}^\varphi (v)$ and $S^\chi_{v} = \mathbf{S}_{loc}^\varphi (B)$ in corresponding point,
\begin{equation} \label{9.14}
S^\Psi_{v,B} (t,y,z) = \tilde{S}_{v,B}(t,\Psi(t,y,z)) = S_{v,B}^\chi(t,(\Phi^{-1}\circ\Psi)(t,y,z)).
\end{equation}
So $\tilde{S}_{v}$ and $\tilde{S}_{B}$ solve (note that $\varphi$ come from $\Phi$ so $\tilde{S}$ solves similar equation in original domain.)
\begin{equation} \label{9.15}
\partial_t \tilde{S}_v + (u\cdot\nabla)\tilde{S}_v - (H\cdot\nabla)\tilde{S}_B - \varepsilon\triangle\tilde{S}_v = F_{S^{\chi}}(t,\Phi^{-1}(t,\cdot)),
\end{equation}
\begin{equation} \label{9.16}
\partial_t \tilde{S}_B + (u\cdot\nabla)\tilde{S}_B - (H\cdot\nabla)\tilde{S}_v - \varepsilon\triangle\tilde{S}_v = E_{S^{\chi}}(t,\Phi^{-1}(t,\cdot)).
\end{equation}
We use Laplacian (\ref{9.9}) to transform above equation via $\Psi$. Then  $S^{\Psi}_{v}$ and $S^{\Psi}_{B}$ solve
\begin{equation} \label{9.17}
\partial_t S^\Psi_v + ({w}_v\cdot\nabla)S^\Psi_v - ({w}_B\cdot\nabla)S^\Psi_B - \varepsilon\left( \partial_{zz}S^\Psi_v + \frac{1}{2}\partial_z(\ln |g|)\partial_z S^\Psi_v + \triangle_{\tilde{g}}S^\Psi_v \right) = F_{S^{\chi}}(t,(\Phi^{-1}\circ\Psi)(t,\cdot)),
\end{equation}
\begin{equation} \label{9.18}
\partial_t S^\Psi_B + ({w}_v\cdot\nabla)S^\Psi_B - ({w}_B\cdot\nabla)S^\Psi_v - \varepsilon\left( \partial_{zz}S^\Psi_B + \frac{1}{2}\partial_z(\ln |g|)\partial_z S^\Psi_B + \triangle_{\tilde{g}}S^\Psi_B \right) = E_{S^{\chi}}(t,(\Phi^{-1}\circ\Psi)(t,\cdot)),
\end{equation}
where
\begin{equation} \label{temp omega}
\begin{split}
{w}_v &:= \bar{\chi}(\mathbf{D}\Psi)^{-1}\left( v(t,\Phi^{-1}\circ\Psi) - \partial_t\Psi \right),  \\
{w}_B &:= \bar{\chi}(\mathbf{D}\Psi)^{-1} B(t,\Phi^{-1}\circ\Psi).
\end{split}
\end{equation}
$S^\Psi$ is compactly supported near the boundary and function $\bar{\chi}$ is slightly larger support in z such that $\bar{\chi}S^\Psi = S^\Psi$. Note that this function allows us to have ${w}$ which is also supported near the boundary. Now we set the alternatives for $S^{v}_n$ and $S^{B}_n$, which are $S^\Psi_{v,n}$ and $S^\Psi_{B,n}$,
\begin{equation} \label{9.19}
\begin{split}
S^\Psi_{v,n} &:= \Pi^b(t,y)S^\Psi_v\mathbf{n}^b(t,y),  \\
S^\Psi_{B,n} &:= \Pi^b(t,y)S^\Psi_B\mathbf{n}^b(t,y),
\end{split}
\end{equation}
where $\Pi^b = \mathbf{I} - \mathbf{n}^b\otimes\mathbf{n}^b$, which means tangential projection to the boundary.\
We get the equations for $S^\Psi_{v,n}$ and $S^\Psi_{B,n}$:
\begin{equation} \label{9.20}
\partial_t S^\Psi_{v,n} + ({w}_v\cdot\nabla)S^\Psi_{v,n} - ({w}_B\cdot\nabla)S^\Psi_{B,n} -\varepsilon\left(\partial_{zz} + \frac{1}{2}\partial_z(\ln |g|)\partial_z\right)S^\Psi_{v,n} = F^\Psi_n,
\end{equation}
where
$$
F^\Psi_n := \left( \Pi^b F_{S^\chi}\mathbf{n}^b + F^{\Psi,1}_n + F^{\Psi,2}_n + F^{\Psi,3}_n \right),
$$
with
\begin{equation*}
\begin{split}
F^{\Psi,1}_n &:= \left( (\partial_t + {w}_{v,y}\cdot\nabla_y)\Pi^b\right)S^\Psi_v\mathbf{n}^b + \Pi^b S^\Psi_v(\partial_t + {w}_{v,y}\cdot\nabla_y)\mathbf{n}^b,  \\
F^{\Psi,2}_n &:= -\varepsilon\Pi^b(\triangle_{\tilde{g}}S^\Psi_v)\mathbf{n}^b,  \\
F^{\Psi,3}_n &:= \left( ({w}_{B,y}\cdot\nabla_y)\Pi^b\right)S^\Psi_B\mathbf{n}^b + \Pi^b S^\Psi_B({w}_{B,y}\cdot\nabla_y)\mathbf{n}^b,
\end{split}
\end{equation*}
\begin{equation} \label{9.21}
\partial_t S^\Psi_{B,n} + ({w}_v\cdot\nabla)S^\Psi_{B,n} - ({w}_B\cdot\nabla)S^\Psi_{v,n} -\varepsilon\left(\partial_{zz}S^\Psi_{B,n} + \frac{1}{2}\partial_z(\ln |g|)\partial_z\right)S^\Psi_{B,n} = E^\Psi_n,
\end{equation}
where
$$
E^\Psi_n := \left( \Pi^b E_{S^\chi}\mathbf{n}^b + E^{\Psi,1}_n + E^{\Psi,2}_n + E^{\Psi,3}_n \right)
$$
with
\begin{equation*}
\begin{split}
E^{\Psi,1}_n &:= \left( (\partial_t + {w}_{v,y}\cdot\nabla_y)\Pi^b\right)S^\Psi_B\mathbf{n}^b + \Pi^b S^\Psi_B(\partial_t + {w}_{v,y}\cdot\nabla_y)\mathbf{n}^b,  \\
E^{\Psi,2}_n &:= -\varepsilon\Pi^b(\triangle_{\tilde{g}}S^\Psi_B)\mathbf{n}^b,  \\
E^{\Psi,3}_n &:= \left( ({w}_{B,y}\cdot\nabla_y)\Pi^b\right)S^\Psi_v\mathbf{n}^b + \Pi^b S^\Psi_v({w}_{B,y}\cdot\nabla_y)\mathbf{n}^b,
\end{split}
\end{equation*}
and boundary condition
\begin{equation} \label{9.22}
S^\Psi_{v,n}|_{z=0} = 0,\,\,\,S^\Psi_{B,n}|_{z=0} = 0.
\end{equation}
This is because $S^\Psi_{v,n} = \Pi \mathbf{S}^\varphi v\mathbf{n} := S^v_n$ and $S^\Psi_{B,n} = \Pi \mathbf{S}^\varphi B\mathbf{n} := S^B
_n$ on the boundary. Now we state the Lemma \ref{lemma 9.5} in \cite{NMFR1}.
\begin{lemma} \label{lemma 9.2}
Consider $\mathcal{T} : S \rightarrow S$ such that $\mathcal{T}(y,0) = y,\,\,\forall y\in \mathbb{R}^2$ and let $g(x) = f(\mathcal{T}x)$. Then for every $k \geq 1$, we have the estimate
$$
\|g\|_{k,\infty} \leq \Lambda\left(\|\nabla\mathcal{T}\|_{k-1,\infty}\right)\|f\|_{k,\infty}.
$$
\end{lemma}
\begin{remark}
Meaning of this lemma : Sobolev conormal spaces are invariant by diffeomorphism which preserve the boundary. i.e If f is conormal $k,\infty$, then $f\circ\mathcal{T}$ is also conormal $k,\infty$ if $\mathcal{T}$ preserves the boundary. Similar holds for $\|\cdot\|_m$ type sobolev space.
\end{remark}

\noindent We use above lemma (and remark) to show that equivalency of $S^v_n$ and $S^\Psi_{v,n}$, and also for $B$.
\begin{equation} \label{9.23}
\begin{split}
\|S^\Psi_{v,n}\|_{1,\infty} &\leq \Lambda_0(\|S^v_n\|_{1,\infty} + \|v\|_{2,\infty}),\,\,\|S^v_n\|_{1,\infty} \leq \Lambda_0\left(\|S^\Psi_{v,n}\|_{1,\infty} + \Lambda(\frac{1}{c_0},|h|_m + \|\mathcal{V}^m\| + \|S^v_n\|_{m-2})\right),  \\
\|S^\Psi_{B,n}\|_{1,\infty} &\leq \Lambda_0(\|S^B_n\|_{1,\infty} + \|B\|_{2,\infty}),\,\,\|S^B_n\|_{1,\infty} \leq \Lambda_0\left(\|S^\Psi_{B,n}\|_{1,\infty} + \Lambda(\frac{1}{c_0},|h|_m + \|\mathcal{B}^m\| + \|S^B_n\|_{m-2})\right).
\end{split}
\end{equation}
We should get estimate for ${w}_{v,B}$. Using the same argument and definition (\ref{temp omega}),
\begin{equation} \label{9.24}
\begin{split}
\|{w}_v\|_{1,\infty} &\leq \Lambda\left(\frac{1}{c_0},|h|_{3,\infty} + \|v\|_{1,\infty} + \|\partial_t\Psi\|_{1,\infty} \right) \leq \Lambda\left(\frac{1}{c_0},|h|_{3,\infty} + \|v\|_{2,\infty} + |\partial_t h|_{2,\infty} \right) \leq \Lambda_{m,\infty}(t),  \\
\|{w}_B\|_{1,\infty} &\leq \Lambda\left(\frac{1}{c_0},|h|_{3,\infty} + \|B\|_{1,\infty} \right) \leq \Lambda\left(\frac{1}{c_0},|h|_{3,\infty} + \|v\|_{2,\infty} \right) \leq \Lambda_{m,\infty}(t).
\end{split}
\end{equation}
Now we make proposition about estimate of $\|Z_i S^\Psi_{v,n}\|_{L^\infty},\,\|Z_i S^\Psi_{B,n}\|_{L^\infty}$. Note that $\|S^\Psi_{v,n}\|_{L^\infty},\,\|S^\Psi_{B,n}\|_{L^\infty}$ is already given above by just maximal principle.
\begin{proposition} \label{proposition 9.4}
We have the folloiwng estimate.
\begin{equation}
\begin{split}
&\|Z_3 S^v_n (t)\|_{L^\infty},\,\|Z_3 S^B_n (t)\|_{L^\infty} \leq \Lambda_0\left( \|v(0)\|_{E^{2,\infty}} + \|B(0)\|_{E^{2,\infty}} \right)  \\
&+ \Lambda\left(\frac{1}{c_0},\|\mathcal{V}^m\| + \|\mathcal{B}^m\| + \|S^v_n\|_{m-2} + \|S^B_n\|_{m-2} + |h|_m \right)  \\
&+ \int_0^t \Lambda_{m,\infty}(s)\left( 1 + \|S^v_n\|_{1,\infty} + \|S^B_n\|_{1,\infty} + \varepsilon\|\nabla\mathcal{V}^m\| + \varepsilon\|S^v_n\|_{3,\infty} + \varepsilon\|v\|_{4,\infty} + \varepsilon\|S^B_n\|_{3,\infty} + \varepsilon\|B\|_{4,\infty} \right).
\end{split}
\end{equation}
\end{proposition}
\begin{proof}
1) When $i=1,2$, we apply $\partial_i$ to get these two.
\begin{equation} \label{9.26}
\begin{split}
&\partial_t\partial_i S^\Psi_{v,n} + ({w}_v\cdot\nabla)\partial_i S^\Psi_{v,n} - ({w}_B\cdot\nabla)\partial_i S^\Psi_{B,n} - \varepsilon\left(\partial_{zz} + \frac{1}{2}\partial_z(\ln |g|)\partial_z\right)\partial_i S^\Psi_{v,n} 	\\
&\quad = \partial_i F^\Psi_n - \partial_i{w}_v\cdot\nabla S^\Psi_{v,n} + \partial_i{w}_B\cdot\nabla S^\Psi_{B,n} - \frac{\varepsilon}{2}\partial_z S^\Psi_{v,n} \partial_{iz}^2 (\ln |g|),
\end{split}
\end{equation}
\begin{equation} \label{9.27}
\begin{split}
&\partial_t\partial_i S^\Psi_{B,n} + ({w}_v\cdot\nabla)\partial_i S^\Psi_{B,n} - ({w}_B\cdot\nabla)\partial_i S^\Psi_{v,n} - \varepsilon\left(\partial_{zz} + \frac{1}{2}\partial_z(\ln |g|)\partial_z\right)\partial_i S^\Psi_{B,n} 	\\
&\quad = \partial_i E^\Psi_n - \partial_i{w}_v\cdot\nabla S^\Psi_{B,n} + \partial_i{w}_B\cdot\nabla S^\Psi_{v,n} - \frac{\varepsilon}{2}\partial_z S^\Psi_{B,n} \partial_{iz}^2 (\ln |g|).
\end{split}
\end{equation}

We also add and subtract these two equations.
\begin{equation} \label{9.28}
\begin{split}
&\partial_t\partial_i(S^\Psi_{v,n}+S^\Psi_{B,n}) + (({w}_v-{w}_B)\cdot\nabla)\partial_i(S^\Psi_{v,n}+S^\Psi_{B,n})  - \varepsilon\left(\partial_{zz} + \frac{1}{2}\partial_z(\ln |g|)\partial_z\right)\partial_i(S^\Psi_{v,n}+S^\Psi_{B,n}) 	\\
&\quad = \partial_i(F^\Psi_n+E^\Psi_n) + (\partial_i{w}_B- \partial_i{w}_v)\cdot\nabla(S^\Psi_{v,n}+S^\Psi_{B,n}) - \frac{\varepsilon}{2}\partial_z(S^\Psi_{v,n}+S^\Psi_{B,n}) \partial_{iz}^2 (\ln |g|),
\end{split}
\end{equation}

\begin{equation} \label{9.29}
\begin{split}
&\partial_t\partial_i(S^\Psi_{v,n}-S^\Psi_{B,n}) + (({w}_v+{w}_B)\cdot\nabla)\partial_i(S^\Psi_{v,n}-S^\Psi_{B,n})  - \varepsilon\left(\partial_{zz} + \frac{1}{2}\partial_z(\ln |g|)\partial_z\right)\partial_i(S^\Psi_{v,n}-S^\Psi_{B,n}) 	\\
&\quad = \partial_i(F^\Psi_n-E^\Psi_n) - ( \partial_i{w}_B + \partial_i{w}_v)\cdot\nabla(S^\Psi_{v,n}-S^\Psi_{B,n}) - \frac{\varepsilon}{2}\partial_z(S^\Psi_{v,n}-S^\Psi_{B,n}) \partial_{iz}^2 (\ln |g|).
\end{split}
\end{equation}

Maximal principle yields
\begin{equation} \label{9.30}
\begin{split}
\|\partial_i(S^\Psi_{v,n}+S^\Psi_{B,n})(t)\|_{L^\infty} &\leq \|\partial_i(S^\Psi_{v,n}+S^\Psi_{B,n})(0)\|_{L^\infty} + \int_0^t \Big( \|\partial_i (F^\Psi_n+E^\Psi_n)\|_{L^\infty}  	\\
&\quad + \|(\partial_i{w}_B- \partial_i{w}_v)\cdot\nabla(S^\Psi_{v,n}+S^\Psi_{B,n})\|_{L^\infty} + \varepsilon\|\partial_z(S^\Psi_{v,n}+S^\Psi_{B,n})\partial_{iz}^2 (\ln |g|)\|_{L^\infty} \Big) 	\\
&\leq \|\partial_i S^\Psi_{v,n}(0)\|_{L^\infty} + \|\partial_i S^\Psi_{B,n}(0)\|_{L^\infty} + \int_0^t \Big( \|\partial_i F^\Psi_n\|_{L^\infty} + \|\partial_i E^\Psi_n\|_{L^\infty}  	\\
&\quad + \|(\partial_i{w}_B- \partial_i{w}_v)\cdot\nabla(S^\Psi_{v,n}+S^\Psi_{B,n})\|_{L^\infty} + \varepsilon\Lambda(\frac{1}{c_0},|h|_{3,\infty}) \big( \|\partial_z S^\Psi_{v,n}\|_{L^\infty} + \|\partial_z S^\Psi_{B,n}\|_{L^\infty} \big)\Big) 	\\
&\leq \|\partial_i S^\Psi_{v,n}(0)\|_{L^\infty} + \|\partial_i S^\Psi_{B,n}(0)\|_{L^\infty} + \int_0^t \Big( \|\partial_i F^\Psi_n\|_{L^\infty} + \|\partial_i E^\Psi_n\|_{L^\infty} 	\\
&\quad + \|(\partial_i{w}_B- \partial_i{w}_v)\cdot\nabla(S^\Psi_{v,n}+S^\Psi_{B,n})\|_{L^\infty} + \Lambda_{m,\infty}(s) \Big),
\end{split}
\end{equation}

\begin{equation} \label{9.31}
\begin{split}
\|\partial_i(S^\Psi_{v,n}-S^\Psi_{B,n})(t)\|_{L^\infty} &\leq \|\partial_i(S^\Psi_{v,n}-S^\Psi_{B,n})(0)\|_{L^\infty} + \int_0^t\Big( \|\partial_i (F^\Psi_n-E^\Psi_n)\|_{L^\infty} 	\\
&\quad + \|(\partial_i{w}_B+\partial_i{w}_v)\cdot\nabla(S^\Psi_{v,n}-S^\Psi_{B,n})\|_{L^\infty} + \varepsilon\|\partial_z(S^\Psi_{v,n}-S^\Psi_{B,n})\partial_{iz}^2 (\ln |g|)\|_{L^\infty} \Big) 	\\
&\leq \|\partial_i S^\Psi_{v,n}(0)\|_{L^\infty} + \|\partial_i S^\Psi_{B,n}(0)\|_{L^\infty} + \int_0^t\Big( \|\partial_i F^\Psi_n\|_{L^\infty} + \|\partial_i E^\Psi_n\|_{L^\infty} 	\\
&\quad + \|(\partial_i{w}_B+\partial_i{w}_v)\cdot\nabla(S^\Psi_{v,n}-S^\Psi_{B,n})\|_{L^\infty} + \varepsilon\Lambda(\frac{1}{c_0},|h|_{3,\infty}) \big(\|\partial_z S^\Psi_{v,n}\|_{L^\infty} + \|\partial_z S^\Psi_{B,n}\|_{L^\infty} \big)\Big) 	\\
&\leq \|\partial_i S^\Psi_{v,n}(0)\|_{L^\infty} + \|\partial_i S^\Psi_{B,n}(0)\|_{L^\infty} + \int_0^t\Big( \|\partial_i F^\Psi_n\|_{L^\infty} + \|\partial_i E^\Psi_n\|_{L^\infty} 	\\
&\quad + \|(\partial_i{w}_B+\partial_i{w}_v)\cdot\nabla(S^\Psi_{v,n}-S^\Psi_{B,n})\|_{L^\infty} + \Lambda_{m,\infty}(s) \Big).
\end{split}
\end{equation}
We estimate high order terms in the RHS. \\
I) $\|\partial_i{w}_v\cdot\nabla S^\Psi_{v,n}\|_{L^\infty}$ estimate.
$$
\|\partial_i{w}_v\cdot\nabla S^\Psi_{v,n}\|_{L^\infty} \leq \|{w}_v\|_{1,\infty}\|S^\Psi_{v,n}\|_{1,\infty} + \|\partial_i{w}_{v,3}\partial_z S^\Psi_{v,n}\|_{L^\infty} \leq \Lambda_{m,\infty}(t) + \|\partial_i{w}_{v,3}\partial_z S^\Psi_{v,n}\|_{L^\infty}.
$$
Note that,
$$
w^b = (\mathbf{D}\Phi(t,y,0))^{-1}(v^b - (0,\partial_t h)),
$$
and 
$$
{w}_{v,3}^b = \frac{1}{|\mathbf{N}|}(v^b\cdot\mathbf{N} - \partial_t h) = 0,
$$
by boundary condition. So $\partial_i{w}_{v,3}$ also vanishes on the boundary since $i=1,2$. Then we get the estimate of the last term of above,
\begin{equation} \label{9.32}
\begin{split}
\|\partial_i{w}_{v,3}\partial_z S^\Psi_{v,n}\|_{L^\infty} &\leq \|\frac{1-z}{z}\partial_i{w}_{v,3}\frac{z}{1-z}\partial_z S^\Psi_{v,n}\|_{L^\infty} \leq \|\frac{1-z}{z}\partial_i{w}_{v,3}\|_{L^\infty} \|S^\Psi_{v,n}\|_{1,\infty} 	\\
&\leq \|\frac{1-z}{z}(0 + |\partial_z\partial_i{w}_{v,3}|_{L^\infty_{z,loc}}z )\|_{L^\infty} \|S^\Psi_{v,n}\|_{1,\infty} \leq \|\partial_z\partial_i {w}_{v,3}\|_{L^\infty} \|S^\Psi_{v,n}\|_{1,\infty} 	\\ 
&\leq \Lambda\left(\frac{1}{c_0},|h|_{3,\infty} + \|v\|_{E^{2,\infty}} \right) \leq \Lambda_{m,\infty}(t).
\end{split}
\end{equation}
Hence
$$
\|\partial_i{w}_v\cdot\nabla S^\Psi_{v,n}\|_{L^\infty} \leq \Lambda_{m,\infty}(t).
$$
Similarly, we get
$$
\|\partial_i{w}_v\cdot\nabla S^\Psi_{B,n}\|_{L^\infty} \leq \Lambda_{m,\infty}(t).
$$
II) $\|\partial_i{w}_B\cdot\nabla S^\Psi_{B,n}\|_{L^\infty}$ estimate.
$$
\|\partial_i{w}_B\cdot\nabla S^\Psi_{B,n}\|_{L^\infty} \leq \|{w}_B\|_{1,\infty}\|S^\Psi_{B,n}\|_{1,\infty} + \|\partial_i{w}_{B,3}\partial_z S^\Psi_{B,n}\|_{L^\infty} \leq \Lambda_{m,\infty}(t) + \|\partial_i{w}_{B,3}\partial_z S^\Psi_{B,n}\|_{L^\infty}.
$$
Similar as I), ${w}_B^b = 0$, since $B$ vanish on the boundary, so 
$$
{w}_{B}^b = (\mathbf{D}\Psi(t,y,0))^{-1}(B^b) = 0.
$$
Using zero boundary value property of $\partial_i{w}_B^b$, we get
$$
\|\partial_i{w}_B\cdot\nabla S^\Psi_{B,n}\|_{L^\infty} \leq \Lambda_{m,\infty}(t).
$$
and similarly,
$$
\|\partial_i{w}_B\cdot\nabla S^\Psi_{v,n}\|_{L^\infty} \leq \Lambda_{m,\infty}(t).
$$
III) $\|\partial_i F^\Psi_n\|_{L^\infty}$ estimate.
$$
\partial_i F^\Psi_n = \partial_i\left(\Pi^b F_{S^\chi}\mathbf{n}^b \right) + \partial_i F^{\Psi,1}_n + \partial_i F^{\Psi,2}_n + \partial_i F^{\Psi,3}_n,
$$
$$
\|F^{\Psi,1}_n\|_{1,\infty},\,\,\|F^{\Psi,3}_n\|_{1,\infty} \leq \Lambda_{m,\infty}(t),
$$
$$
\|F^{\Psi,2}_n\|_{1,\infty} \leq \varepsilon\Lambda_{m,\infty}(t) \left( \|S^v_n\|_{3,\infty} + \|v\|_{4,\infty} \right).
$$
Considering $ \partial_i\left(\Pi^b F_{S^\chi}\mathbf{n}^b \right) $, we get 
\begin{equation*}
\begin{split}
	\|F^\Psi_n\|_{1,\infty} &\leq \Lambda_{m,\infty}(t) \left( 1 + \varepsilon\|S^v_n\|_{3,\infty} + \varepsilon\|v\|_{4,\infty} + \|\Pi^b\left((\mathbf{D}^\varphi)^2 q\right)\mathbf{n}^b\|_{1,\infty} \right), 	\\
	\|\Pi^b\left((\mathbf{D}^\varphi)^2 q\right)\mathbf{n}^b\|_{1,\infty} &\leq \Lambda_0 \left( \|\nabla q^E\|_{2,\infty} + \|\nabla q^{NS}\|_{2,\infty} \right) \leq \Lambda_0 \left( 1 + \varepsilon\|\mathbf{S}^\varphi\mathcal{V}^m\|\right).
\end{split}
\end{equation*}
Hence,
$$
\|F^\Psi_n\|_{1,\infty} \leq \Lambda_{m,\infty}(t) \left( 1 + \varepsilon\|S^v_n\|_{3,\infty} + \varepsilon\|v\|_{4,\infty} + \varepsilon\|\nabla\mathcal{V}^m\| \right).
$$
IV) $\|\partial_i E^\Psi_n\|_{L^\infty}$ estimate.
$$
\partial_i E^\Psi_n = \partial_i\left(\Pi^b E_{S^\chi}\mathbf{n}^b \right) + \partial_i E^{\Psi,1}_n + \partial_i E^{\Psi,2}_n + \partial_i E^{\Psi,3}_n,
$$
$$
\|E^{\Psi,1}_n\|_{1,\infty},\,\,\|E^{\Psi,3}_n\|_{1,\infty} \leq \Lambda_{m,\infty}(t),
$$
$$
\|E^{\Psi,2}_n\|_{1,\infty} \leq \lambda\Lambda_{m,\infty}(t) \left( \|S^B_n\|_{3,\infty} + \|B\|_{4,\infty} \right).
$$
Considering $ \partial_i\left(\Pi^b E_{S^\chi}\mathbf{n}^b \right) $ (These is no terms like $F_v$, so there is no terms of pressure), we get simply
$$
\|E^\Psi_n\|_{1,\infty} \leq \Lambda_{m,\infty}(t) \left( 1 + \varepsilon\|S^B_n\|_{3,\infty} + \varepsilon\|B\|_{4,\infty} \right).
$$
Combining above results of $I)\sim IV)$ above, we get
\begin{equation} \label{9.33}
\begin{split}
\|\partial_i(S^\Psi_{v,n}+S^\Psi_{B,n})(t)\|_{L^\infty} &\leq \|\partial_i S^\Psi_{v,n}(0)\|_{L^\infty} + \|\partial_i S^\Psi_{B,n}(0)\|_{L^\infty} 	\\
&\quad + \int_0^t \Lambda_{m,\infty}(s)\left( 1 + \varepsilon\|\nabla \mathcal{V}^m\| + \varepsilon\|S^v_n\|_{3,\infty} + \varepsilon\|S^B_n\|_{3,\infty} + \varepsilon\|v\|_{4,\infty} + \varepsilon\|B\|_{4,\infty} \right),
\end{split}
\end{equation}
\begin{equation} \label{9.34}
\begin{split}
\|\partial_i(S^\Psi_{v,n}-S^\Psi_{B,n})(t)\|_{L^\infty} &\leq \|\partial_i S^\Psi_{v,n}(0)\|_{L^\infty} + \|\partial_i S^\Psi_{B,n}(0)\|_{L^\infty} 	\\
&\quad + \int_0^t \Lambda_{m,\infty}(s)\left( 1 + \varepsilon\|\nabla \mathcal{V}^m\| + \varepsilon\|S^v_n\|_{3,\infty} + \varepsilon\|S^B_n\|_{3,\infty} + \varepsilon\|v\|_{4,\infty} + \varepsilon\|B\|_{4,\infty} \right).
\end{split}
\end{equation}
Therefore, we get
\begin{equation} \label{9.35}
\begin{split}
&\|\partial_i S^\Psi_{v,n}(t)\|_{L^\infty},\,\,\|\partial_i S^\Psi_{B,n}(t)\|_{L^\infty} \leq \|\partial_i S^\Psi_{v,n}(0)\|_{L^\infty} + \|\partial_i S^\Psi_{B,n}(0)\|_{L^\infty} 	\\
&\quad + \int_0^t \Lambda_{m,\infty}(s)\left( 1 + \varepsilon\|\nabla \mathcal{V}^m\| + \varepsilon\|S^v_n\|_{3,\infty} + \varepsilon\|S^B_n\|_{3,\infty} + \varepsilon\|v\|_{4,\infty} + \varepsilon\|B\|_{4,\infty} \right).
\end{split}
\end{equation}
2) When $i=3$, we apply $Z_3 = \frac{z}{1-z}\partial_z$. commutator between $Z_3,\varepsilon\partial_{zz}$ should be treated carefully. 
It is convenient to eliminate $\varepsilon \partial_z(\ln |g|)\partial_z$ in modified laplacian. This is done by defining,
\begin{equation} \label{9.36}
\begin{split}
\rho_v(t,y,z) &= |g|^{\frac{1}{4}}S^\Psi_{v,n} = |g|^{\frac{1}{4}}\Pi^b S^\Psi_v \mathbf{n}^b, 	\\
\rho_B(t,y,z) &= |g|^{\frac{1}{4}}S^\Psi_{B,n} = |g|^{\frac{1}{4}}\Pi^b S^\Psi_B \mathbf{n}^b.
\end{split}
\end{equation}
Since these solve,
\begin{equation} \label{9.37}
\partial_t\rho_v + {w}_v\cdot\nabla\rho_v - {w}_B\cdot\nabla\rho_B - \varepsilon\partial_{zz}\rho_v = |g|^{\frac{1}{4}}\left( F^\Psi_n + F_g \right) := \mathcal{H}_1,
\end{equation}
where 
$$
F_g := \frac{\rho_v}{|g|^{\frac{1}{2}}} \left({w}_v\cdot\nabla - \varepsilon\partial_{zz}\right)|g|^{\frac{1}{4}} + \frac{\rho_B}{|g|^{\frac{1}{2}}} \left({w}_B\cdot\nabla \right)|g|^{\frac{1}{4}},
$$
and 
\begin{equation} \label{9.38}
\partial_t\rho_B + {w}_v\cdot\nabla\rho_B - {w}_B\cdot\nabla\rho_v - \varepsilon\partial_{zz}\rho_B = |g|^{\frac{1}{4}}\left( E^\Psi_n + E_g \right) := \mathcal{H}_2,
\end{equation}
where
$$
E_g := \frac{\rho_B}{|g|^{\frac{1}{2}}} \left({w}_v\cdot\nabla - \varepsilon\partial_{zz}\right)|g|^{\frac{1}{4}} + \frac{\rho_v}{|g|^{\frac{1}{2}}} \left({w}_B\cdot\nabla \right)|g|^{\frac{1}{4}}.
$$
We treat $\rho_v$ instead of $Z_3 S^\Psi_{v,n}$ (same for $B$) since we have equivalency
\begin{equation} \label{9.39}
\begin{split}
\|Z_3 S^\Psi_{v,n}\|_{L^\infty} &\leq \Lambda_0\|\rho_v\|_{1,\infty},\quad \|\rho_v\|_{1,\infty} \leq \Lambda_0\|S^\Psi_{v,n}\|_{1,\infty}, 	\\
\|Z_3 S^\Psi_{B,n}\|_{L^\infty} &\leq \Lambda_0\|\rho_B\|_{1,\infty},\quad \|\rho_B\|_{1,\infty} \leq \Lambda_0\|S^\Psi_{B,n}\|_{1,\infty}.
\end{split}
\end{equation}
So we get equivalent relation between $\|S^\Psi_{v,n}\|_{1,\infty}$ and $\|\rho_v\|_{1,\infty}$. This is same for $B$. Note that $\rho_v,\rho_B$ are also zero on the boundary $z=0$. The following is Lemma 9.6 in \cite{NMFR1}.

\begin{lemma} \label{lemma 9.5}
For smooth function $\rho$,
\begin{equation} \label{9.40}
\partial_t\rho + {w}\cdot\nabla\rho = \varepsilon\partial_{zz}\rho + \mathcal{H},
\end{equation}
where ${w}_3$ vanishes on the boundary. Assume $\rho$ and $\mathcal{H}$ are compactly supported in $z$, then we have,
$$
\|Z_i\rho(t)\|_{L^\infty} \leq \|Z_i\rho_0\|_{L^\infty} + \|\rho_0\|_{L^\infty} + \int_0^t\left( (\|{w}\|_{E^{2,\infty}} + \|\partial_{zz}{w}_3\|_{L^\infty})(\|\rho\|_{1,\infty} + \|\rho\|_4) + \|\mathcal{H}\|_{1,\infty} \right),\,\,\,\,i=1,2,3.
$$
\end{lemma}

\noindent By adding and subtracting two $\rho$ equations, we have
\begin{equation} \label{9.41}
\partial_t(\rho_v+\rho_B) + ({w}_v-{w}_B)\cdot\nabla(\rho_v+\rho_B) - \varepsilon\partial_{zz}(\rho_v+\rho_B) = \mathcal{H}_1 + \mathcal{H}_2,
\end{equation}
\begin{equation} \label{9.42}
\partial_t(\rho_v-\rho_B) + ({w}_v+{w}_B)\cdot\nabla(\rho_v-\rho_B) - \varepsilon\partial_{zz}(\rho_v-\rho_B) = \mathcal{H}_1 - \mathcal{H}_2.
\end{equation}
Using Lemma \ref{lemma 9.5}, 
\begin{equation} \label{9.43}
\begin{split}
\|Z_3(\rho_v+\rho_B)(t)\|_{L^\infty} &\leq \|Z_3(\rho_v+\rho_B)(0)\|_{L^\infty} + \|(\rho_v+\rho_B)(0)\|_{L^\infty} 	\\
&\quad + \int_0^t \Big( (\|{w}_v-{w}_B\|_{E^{2,\infty}} + \|\partial_{zz}({w}_v-{w}_B)_3\|_{L^\infty})(\|\rho_v+\rho_B\|_{1,\infty} + \|\rho_v+\rho_B\|_4) 	\\
&\quad + \|\mathcal{H}_1\|_{1,\infty} + \|\mathcal{H}_2\|_{1,\infty} \Big),
\end{split}
\end{equation}
\begin{equation} \label{9.44}
\begin{split}
\|Z_3(\rho_v-\rho_B)(t)\|_{L^\infty} &\leq \|Z_3(\rho_v-\rho_B)(0)\|_{L^\infty} + \|(\rho_v-\rho_B)(0)\|_{L^\infty} 	\\
&\quad + \int_0^t \Big( (\|{w}_v+{w}_B\|_{E^{2,\infty}} + \|\partial_{zz}({w}_v+{w}_B)_3\|_{L^\infty})(\|\rho_v-\rho_B\|_{1,\infty} + \|\rho_v-\rho_B\|_4) 	\\
&\quad + \|\mathcal{H}_1\|_{1,\infty} + \|\mathcal{H}_2\|_{1,\infty} \Big).
\end{split}
\end{equation}
We estimate terms on the RHS of (\ref{9.43}) and (\ref{9.44}).  \\
\noindent I) $\|\mathcal{H}_1\|_{1,\infty},\,\|\mathcal{H}_1\|_{2,\infty}$ estimates.
$$
\|\mathcal{H}_1\|_{1,\infty} \leq \|F^\Psi_n\|_{1,\infty} + \|F_g\|_{1,\infty}\quad\text{and}\quad\|\mathcal{H}_2\|_{1,\infty} \leq \|E^\Psi_n\|_{1,\infty} + \|E_g\|_{1,\infty}.
$$
We already estimated $\|F^\Psi_n\|_{1,\infty},\,\|E^\Psi_n\|_{1,\infty}$ above. So,
\begin{equation*}
\begin{split}
\|\mathcal{H}_1\|_{1,\infty} &\leq \Lambda_{m,\infty}(t) \left( 1 + \varepsilon\|\mathbf{S}^\varphi\mathcal{V}^m\| + \varepsilon\|S^v_n\|_{3,\infty} + \varepsilon\|v\|_{4,\infty} \right),  \\
\|\mathcal{H}_2\|_{1,\infty} &\leq \Lambda_{m,\infty}(t) \left( 1 + \varepsilon\|S^B_n\|_{3,\infty} + \varepsilon\|B\|_{4,\infty} \right).
\end{split}
\end{equation*}
II) $\|\rho_v\|_4,\,\|\rho_B\|_4$ estimates.
\begin{equation*}
\begin{split}
\|\rho_v\|_4 \leq \Lambda\left(\frac{1}{c_0}, |h|_6 + \|S^\Psi_{v,n}\|_4 \right) &\leq \Lambda_{m,\infty}(t), 		\\
\|\rho_B\|_4 \leq \Lambda\left(\frac{1}{c_0}, |h|_6 + \|S^\Psi_{B,n}\|_4 \right) &\leq \Lambda_{m,\infty}(t).
\end{split}
\end{equation*}
III) $\|{w}_v\|_{E^{2,\infty}},\,\|{w}_B\|_{E^{2,\infty}}$ estimates. From definition,
$$
\|{w}_v\|_{E^{2,\infty}},\,\,\|{w}_B\|_{E^{2,\infty}} \leq \Lambda_{m,\infty}(t).
$$
IV) $\|\partial_{zz}{w}_{v,3}\|_{L^\infty},\,\|\partial_{zz}{w}_{B,3}\|_{L^\infty}$ estimates. It looks that this term has two normal derivatives. First,
$$
\|\partial_{zz}\left(\bar{\chi}(\mathbf{D}\Psi^{-1}\partial_t\Psi)\right)\|_{L^\infty} \leq \Lambda\left(\frac{1}{c_0}, |h|_{2,\infty} + |\partial_t h|_{2,\infty} \right) \leq \Lambda_{m,\infty}(t).
$$
Main part is third component,
$$
\|\partial_{zz}\left(\bar{\chi}D\Psi^{-1}v(t,\Phi^{-1}\circ\Psi) \right)_3\|_{L^\infty}.
$$
Key point is that this is bounded by term with $(D\Psi^{-1})^b$, 
\begin{equation*}
\begin{split}
\|\partial_{zz}\left(\bar{\chi}D\Psi^{-1}v(t,\Phi^{-1}\circ\Psi) \right)_3\|_{L^\infty} &\leq \|\bar{\chi}\partial_{zz}\left((D\Phi(t,y,0))^{-1}v(t,\Phi^{-1}\circ\Psi) \right)_3\|_{L^\infty} + \Lambda_{m,\infty}(t), 		\\
&\leq \|\bar{\chi}\partial_{zz}\left( v(t,\Phi^{-1}\circ\Psi)\cdot\mathbf{n}^b\right)\|_{L^\infty} + \Lambda_{m,\infty}(t),
\end{split}
\end{equation*}
where we used (\ref{9.35}). We write $v(t,\Phi^{-1}\circ\Psi) = u(t,\Psi) := u^{\Psi}(t,y,z)$ ,then using divergence free condition of $u$, we can change 1-normal derivative to tangential derivative so that $\partial_{zz} \rightarrow \partial_{iz}$ to be controlled by $\Lambda_{m,\infty}$.
At result, from \cite{NMFR1},
$$
\|\partial_{zz}{w}_{v,3}\|_{L^\infty}, \|\partial_{zz}{w}_{B,3}\|_{L^\infty} \leq \Lambda_{m,\infty}(t).
$$
So we get the estimates for $\|Z_3(\rho_v+\rho_B)\|_{L^\infty},\,\|Z_3(\rho_v-\rho_B)\|_{L^\infty}$ in the same form. 
\begin{equation} \label{9.45}
\begin{split}
\|Z_3(\rho_v+\rho_B)\|_{L^\infty},\,\|Z_3(\rho_v-\rho_B)\|_{L^\infty} &\leq \|Z_3(\rho_v+\rho_B)(0)\|_{L^\infty}  + \int_0^t \Lambda_{m,\infty}(s) \Big( 1 + \|\rho_v\|_{1,\infty} + \|\rho_B\|_{1,\infty} 	\\
&\quad + \varepsilon\|\mathbf{S}^\varphi\mathcal{V}^m\| + \varepsilon\|S^v_n\|_{3,\infty} + \varepsilon\|v\|_{4,\infty} + \varepsilon\|S^B_n\|_{3,\infty} + \varepsilon\|B\|_{4,\infty} \Big),
\end{split}
\end{equation}
with
$$
\|Z_3(\rho_v+\rho_B)(0)\|_{L^\infty} \leq \Lambda_0\left( \|v(0)\|_{E^{2,\infty}} + \|B(0)\|_{E^{2,\infty}} \right).
$$
Using same technique again, we get the same estimates for each $\|Z_3\rho_v\|_{L^\infty},\,\|Z_3\rho_B\|_{L^\infty}$.
\begin{equation} \label{9.46}
\begin{split}
\|Z_3\rho_v(t)\|_{L^\infty},\,\|Z_3\rho_B(t)\|_{L^\infty} &\leq \Lambda_0 \Big( \|v(0)\|_{E^{2,\infty}} + \|B(0)\|_{E^{2,\infty}} \Big) + \int_0^t \Lambda_{m,\infty}(s) \Big( 1 + \|\rho_v\|_{1,\infty} + \|\rho_B\|_{1,\infty} 	\\
&\quad + \varepsilon\|\mathbf{S}^\varphi\mathcal{V}^m\| + \varepsilon\|S^v_n\|_{3,\infty} + \varepsilon\|v\|_{4,\infty} + \varepsilon\|S^B_n\|_{3,\infty} + \varepsilon\|B\|_{4,\infty} \Big).
\end{split}
\end{equation}
3) We now know that $\|Z_i S^\Psi_{
(v,B),n} \|_{L^{\infty}}$ estimate, and $\|Z_3\rho_{v,B}\|_{L^\infty}$. Moreover, we have equivalent relation,
$$
\rho_{v,B} \sim S^\Psi_{(v,B),n} \sim S^{v,B}_n,
$$
with help of 
$$
\Lambda\left(\frac{1}{c_0},\|\mathcal{V}^m\| + \|\mathcal{B}^m\| + \|S^{v}_n\|_{m-2} + \|S^{B}_n\|_{m-2} + |h|_m \right).
$$
Finally we get
\begin{equation} \label{9.47}
\begin{split}
\|Z_3 S^v_n (t)\|_{L^\infty}, \ \|Z_3 S^B_n (t)\|_{L^\infty} &\leq \Lambda_0 \Big( \|v(0)\|_{E^{2,\infty}} + \|B(0)\|_{E^{2,\infty}} \Big) 	\\
&\quad + \Lambda\left(\frac{1}{c_0},\|\mathcal{V}^m\| + \|\mathcal{B}^m\| + \|S^v_n\|_{m-2} + \|S^B_n\|_{m-2} + |h|_m \right)  \\
&\quad + \int_0^t \Lambda_{m,\infty}(s) \Big( 1 + \|S^v_n\|_{1,\infty} + \|S^B_n\|_{1,\infty} 	\\
&\quad + \varepsilon\|\nabla\mathcal{V}^m\| + \varepsilon\|S^v_n\|_{3,\infty} + \varepsilon\|v\|_{4,\infty} + \varepsilon\|S^B_n\|_{3,\infty} + \varepsilon\|B\|_{4,\infty} \Big).
\end{split}
\end{equation}
\end{proof} 

From zero and first order estimates, we get an estimate for $\|S^{v,B}_n\|_{1,\infty}$. 
\begin{proposition} \label{proposition 9.6}
We have the following estimates for $\|S^{v}_n\|_{1,\infty}$ and $\|S^{B}_n\|_{1,\infty}$.
\begin{equation} \label{9.48}
\begin{split}
\|S^v_n (t)\|_{1,\infty}^2, \ \|S^B_n (t)\|_{1,\infty}^2 &\leq \Lambda_0\left( \|S^v_n(0)\|_{1,\infty}^2 + \|S^B_n(0)\|_{1,\infty}^2 \right)  \\
&\quad + \Lambda\left(\frac{1}{c_0},\|\mathcal{V}^m\| + \|\mathcal{B}^m\| + \|S^v_n\|_{m-2} + \|S^B_n\|_{m-2} + |h|_m \right)  \\
&\quad + ( 1 + \varepsilon ) \int_0^t \Lambda_{m,\infty}(s) + \int_0^t \Big( \|S^v_n\|_{1,\infty}^2 + \|S^B_n\|_{1,\infty}^2 	\\
&\quad + \varepsilon\|\nabla\mathcal{V}^m\|^2 + \varepsilon\|\nabla\mathcal{B}^m\|^2 + \varepsilon\|\nabla S^v_n\|_{m-2}^2 + \varepsilon\|\nabla S^B_n\|_{m-2}^2 \Big).
\end{split}
\end{equation}
\end{proposition}
\begin{proof}
By Proposition \ref{proposition 9.4} and (\ref{9.4}),
\begin{equation} \label{9.49}
\begin{split}
\|S^v_n (t)\|_{1,\infty},\,\|S^B_n (t)\|_{1,\infty} &\leq \Lambda_0\left( \|v(0)\|_{E^{2,\infty}} + \|B(0)\|_{E^{2,\infty}} \right)   \\
&\quad + \Lambda\left(\frac{1}{c_0},\|\mathcal{V}^m\| + \|\mathcal{B}^m\| + \|S^v_n\|_{m-2} + \|S^B_n\|_{m-2} + |h|_m \right)  \\
&\quad + ( 1 + \varepsilon ) \int_0^t \Lambda_{m,\infty}(s) \Big( 1 + \|S^v_n\|_{1,\infty} + \|S^B_n\|_{1,\infty} 	\\
&\quad + \varepsilon\|\nabla\mathcal{V}^m\| + \varepsilon\|S^v_n\|_{3,\infty} + \varepsilon\|v\|_{4,\infty} + \varepsilon\|S^B_n\|_{3,\infty} + \varepsilon\|B\|_{4,\infty} \Big).
\end{split}
\end{equation}
We estimate $\varepsilon\|S^v_n\|_{3,\infty}$, $\varepsilon\|v\|_{4,\infty}$, $\varepsilon\|S^B_n\|_{3,\infty}$, and $\varepsilon\|B\|_{4,\infty}$. For $\varepsilon\|S^v_n\|_{3,\infty}$ and $\varepsilon\|S^B_n\|_{3,\infty}$, by embedding,
\begin{equation*}
\begin{split}
\sqrt{\varepsilon}\|S^v_n\|_{3,\infty} &\leq \Lambda_{m,\infty}(t) +  \sqrt{\varepsilon}\|\nabla\mathcal{V}^m\| + \sqrt{\varepsilon}\|\nabla S^v_n\|_{m-2}, 	\\
\sqrt{\varepsilon}\|S^B_n\|_{3,\infty} &\leq \Lambda_{m,\infty}(t) +  \sqrt{\varepsilon}\|\nabla\mathcal{B}^m\| + \sqrt{\varepsilon}\|\nabla S^B_n\|_{m-2}.
\end{split}
\end{equation*}
For $\varepsilon\|v\|_{4,\infty}$ and $\varepsilon\|B\|_{4,\infty}$,
\begin{equation*}
\begin{split}
\sqrt{\varepsilon}\|v\|_{4,\infty} &\leq \Lambda_{m,\infty}(t) + \sqrt{\varepsilon}\|\nabla\mathcal{V}^m\|, 	\\
\sqrt{\varepsilon}\|B\|_{4,\infty} &\leq \Lambda_{m,\infty}(t) + \sqrt{\varepsilon}\|\nabla\mathcal{B}^m\|.
\end{split}
\end{equation*}
By using above 4 estimates and Young's inequality, we get the result.
\end{proof}

\subsection{$\sqrt{\varepsilon}\|\partial_z S^v_n\|_{L^\infty}$,$\sqrt{\varepsilon}\|\partial_z S^B_n\|_{L^\infty}$ estimates}
\begin{proposition} \label{proposition 9.7}
We have the estimate for $\sqrt{\varepsilon}\|\partial_z S^v_n\|_{L^\infty}$ and $\sqrt{\varepsilon}\|\partial_z S^B_n\|_{L^\infty}$,
\begin{equation} \label{9.51}
\begin{split}
\varepsilon\|\partial_z\rho_v\|_{L^\infty}^2, \ \varepsilon\|\partial_z\rho_B\|_{L^\infty}^2 &\leq \Lambda_{m,\infty}(0) + 2\int_0^t \left( \varepsilon\|\nabla\mathcal{V}^m\|^2 + \varepsilon\|\nabla\mathcal{B}^m\|^2 + \varepsilon\|\nabla S^v_n\|_{m-2}^2 + + \varepsilon\|\nabla S^B_n\|_{m-2}^2 \right)  	\\
&\quad + ( 1 + 16\sqrt{t})\int_0^t \frac{\Lambda_{m,\infty}(s)}{\sqrt{t-\tau}} d\tau.
\end{split}
\end{equation}
\end{proposition}
\begin{proof}
This estimate corresponds to estimate of $\sqrt{\varepsilon}\|\partial_{zz}v\|_{L^\infty}$. Our strategy is to derive the estimate for $\sqrt{\varepsilon}\|\partial_z\rho_{v,B}\|_{L^\infty}$, because, for both $v,B$,
$$
\partial_z S^\Psi_n = \Pi^b\frac{\partial}{\partial z}S^\chi(t,\Phi^{-1}\circ\Psi)\mathbf{n}^b.
$$
We can apply Lemma \ref{lemma 9.2}, so we get similar control,
$$
\|\partial_z S^\Psi_n \|_{1,\infty} \leq \Lambda_0 \|\partial_z\Pi^b S^\chi(t,\Phi^{-1}\circ\Psi)\mathbf{n}^b\|_{1,\infty}.
$$
Then using $|\Pi-\Pi^b| + |\mathbf{n}-\mathbf{n}^b| = \mathcal{O}(z)$, 
$$
\|\partial_z S^\Psi_n \|_{1,\infty} \leq \Lambda_0\left( \|\partial_z S_n\|_{1,\infty} + \|v\|_{2,\infty}\right).
$$
What we need is inverse argument. since the map $\mathcal
T$ in Lemma \ref{lemma 9.2} conserves boundary, 
$$
\|\partial_z S_n \|_{1,\infty} \leq \Lambda_0\left( \|\partial_z S^\Psi_n\|_{1,\infty} + \|v\|_{2,\infty}\right),
$$
and
$$
\|\partial_z\rho\|_{1,\infty} = \|\partial_z\left(|g|^{\frac{1}{2}}\Pi^b S^\Psi\mathbf{n}^b \right)\|_{1,\infty}.
$$
On the right hand side, $\partial_z$ hit $|g|$ and $S^\Psi$ and we have,
$$
\|\partial_z S^\Psi_n\|_{1,\infty} \leq \Lambda_0(\|\partial_z\rho\|_{1,\infty}).
$$
Hence, we get the control what we expected for both $v$ and $B$.
\begin{equation} \label{9.52}
\begin{split}
\sqrt{\varepsilon}\|\partial_z S_n\|_{1,\infty} &\leq \Lambda_0\left( \sqrt{\varepsilon}\|\partial_z\rho_v\|_{1,\infty} + \|v\|_{2,\infty}\right),	\\
\sqrt{\varepsilon}\|\partial_z S_n\|_{1,\infty} &\leq \Lambda_0\left( \sqrt{\varepsilon}\|\partial_z\rho_B\|_{1,\infty} + \|B\|_{2,\infty}\right).
\end{split}
\end{equation}
As we know, $\rho_v,\,\rho_B$ solve
\begin{equation*}
\begin{split}
&\partial_t(\rho_v+\rho_B) + ({w}_v-{w}_B)\cdot\nabla(\rho_v+\rho_B) - \varepsilon\partial_{zz}(\rho_v+\rho_B) = \mathcal{H}_1 + \mathcal{H}_2, 	\\
&\partial_t(\rho_v-\rho_B) + ({w}_v+{w}_B)\cdot\nabla(\rho_v-\rho_B) - \varepsilon\partial_{zz}(\rho_v-\rho_B) = \mathcal{H}_1 - \mathcal{H}_2.
\end{split}
\end{equation*}
This is heat equation with respect to z-direction, with zero boundary data on $z=0$. We use heat kernel,
$$
G(t,y,z) = \frac{1}{\sqrt{4\pi t}}\left( e^{-\frac{(z-z')^2}{4t}} - e^{-\frac{(z+z')^2}{4t}} \right).
$$
Using initial data $\rho_0$ and source $(\mathcal{H}_1 + \mathcal{H}_2) - ({w}_v-{w}_B)\cdot\nabla(\rho_v+\rho_B)$, we get
\begin{equation} \label{9.53}
\begin{split}
&\sqrt{\varepsilon}\partial_z(\rho_v + \rho_B)(t,y,z) = \int_{-\infty}^0 \sqrt{\varepsilon}\partial_z G(t,z,z')\rho_0(y,z')dz' 	\\
&\quad + \int_0^t\int_{\infty}^0 \sqrt{\varepsilon}\partial_z G(t-\tau,z,z')\left\{\mathcal{H}_1 + \mathcal{H}_2 - ({w}_v-{w}_B)\cdot\nabla(\rho_v+\rho_B)\right\}(\tau,y,z') dz'd\tau.
\end{split}
\end{equation}
Since $G$ has a gaussian form, we get
\begin{equation*}
\begin{split}
&\sqrt{\varepsilon}\|\partial_z(\rho_v + \rho_B)(t)\|_{L^\infty} \leq \sqrt{\varepsilon}\|\partial_z(\rho_v + \rho_B)(0)\|_{L^\infty} 	\\
&\quad + \frac{1}{\sqrt{4\pi}}\int_0^t\frac{1}{\sqrt{t-\tau}}\left( \|\mathcal{H}_1 + \mathcal{H}_2\|_{L^\infty} + \|({w}_v-{w}_B)\cdot\nabla(\rho_v+\rho_B)\|_{L^\infty} \right) d\tau.
\end{split}
\end{equation*}
I) Using previous estimates for $\|\mathcal{H}_1\|_{1,\infty}$ and $\|\mathcal{H}_2\|_{1,\infty}$,
\begin{equation*}
\begin{split}
\|\mathcal{H}_1\|_{1,\infty} &\leq \Lambda_{m,\infty}(t) \left( 1 + \varepsilon\|\mathbf{S}^\varphi \mathcal{V}^m\| + \varepsilon\|S^v_n\|_{2,\infty} + \varepsilon\|v\|_{3,\infty}\right), 	\\
\|\mathcal{H}_2\|_{1,\infty} &\leq \Lambda_{m,\infty}(t) \left( 1 + \varepsilon\|S^B_n\|_{2,\infty} + \varepsilon\|B\|_{3,\infty} \right).
\end{split}
\end{equation*}
We should control $\varepsilon\|S^v_n\|_{2,\infty}$, $\varepsilon\|v\|_{3,\infty}$, $\varepsilon\|S^B_n\|_{2,\infty}$, and $\varepsilon\|B\|_{3,\infty}$.
\begin{equation*}
\begin{split}
\varepsilon\|S^v_n\|_{2,\infty} &\leq \varepsilon\|\nabla S^v_n\|_3^{\frac{1}{2}}\|S^v_n\|_4 \leq \Lambda_{m,\infty}(t) \varepsilon\|\nabla S^v_n\|_{m-2}^{\frac{1}{2}},  \\
\varepsilon\|S^B_n\|_{2,\infty} &\leq \varepsilon\|\nabla S^B_n\|_3^{\frac{1}{2}}\|S^B_n\|_4 \leq \Lambda_{m,\infty}(t) \varepsilon\|\nabla S^B_n\|_{m-2}^{\frac{1}{2}},  \\
\varepsilon\|v\|_{3,\infty} &\leq \varepsilon\|\nabla v\|_4^{\frac{1}{2}}\|v\|_5^{\frac{1}{2}} \leq \Lambda_{m,\infty}(t) \left( 1 + \|\nabla\mathcal{V}^m\|\right),  \\
\varepsilon\|B\|_{3,\infty} &\leq \varepsilon\|\nabla B\|_4^{\frac{1}{2}}\|B\|_5^{\frac{1}{2}} \leq \Lambda_{m,\infty}(t) \left( 1 + \|\nabla\mathcal{B}^m\|\right).
\end{split}
\end{equation*}
II) We have,
$$
\sqrt{\varepsilon}\|\partial_z\rho_v(0)\|_{L^\infty} \leq \Lambda_{m,\infty}(0),\,\,\sqrt{\varepsilon}\|\partial_z\rho_B(0)\|_{L^\infty} \leq \Lambda_{m,\infty}(0).
$$
III) Using the fact that ${w}$ is zero on the boundary, we give $\partial_z$ to ${w}$ and tame the second term into conormal regularity not the $\partial_z$ regularity.
$$
\|({w}_v-{w}_B)\cdot\nabla(\rho_v+\rho_B)\|_{L^\infty} \leq \|{w}_v-{w}_B\|_{E^{1,\infty}}\|\rho_v+\rho_B\|_{1,\infty} \leq \Lambda_{m,\infty}(t).
$$
Hence, similar as coupled equations,
\begin{equation} \label{9.54}
\begin{split}
&\sqrt{\varepsilon}\|\partial_z(\rho_v + \rho_B)(t)\|_{L^\infty},\,\sqrt{\varepsilon}\|\partial_z(\rho_v - \rho_B)(t)\|_{L^\infty} \leq   \Lambda_{m,\infty}(0) 	\\
&\quad + \int_0^t \frac{\Lambda_{m,\infty}(s)}{\sqrt{t-\tau}} \left( 1 + \varepsilon\|\nabla\mathcal{V}^m\|^{\frac{1}{2}}  + \varepsilon\|\nabla\mathcal{B}^m\|^{\frac{1}{2}} + \varepsilon\|\nabla S^v_n\|_{m-2}^{\frac{1}{2}} + \varepsilon\|\nabla S^B_n\|_{m-2}^{\frac{1}{2}} \right)d\tau.
\end{split}
\end{equation}
Therefore, similar as before,
\begin{equation} \label{9.55}
\begin{split}
&\sqrt{\varepsilon}\|\partial_z\rho_v(t)\|_{L^\infty},\,\sqrt{\varepsilon}\|\partial_z\rho_B(t)\|_{L^\infty} \leq   \Lambda_{m,\infty}(0) 	\\
&\quad + \int_0^t \frac{\Lambda_{m,\infty}}{\sqrt{t-\tau}} \left( 1 + \varepsilon\|\nabla\mathcal{V}^m\|^{\frac{1}{2}}  + \varepsilon\|\nabla\mathcal{B}^m\|^{\frac{1}{2}} + \varepsilon\|\nabla S^v_n\|_{m-2}^{\frac{1}{2}} + \varepsilon\|\nabla S^B_n\|_{m-2}^{\frac{1}{2}} \right)d\tau.
\end{split}
\end{equation}
Now we square these inequalities. Main stuff is squaring the last term of these inequalities. There are two cases. First, when we product two different terms we can use young's inequality to get the terms like
$\varepsilon\int_0^t\|\nabla\mathcal{V}^m\|d\tau$, what we want. When we squre this terms, we should be careful, because if we use Holder's inequality for $L^2-L^2$ then we may get terms like $\int_0^t\frac{\Lambda_{m,\infty}}{(t-\tau)^2}d\tau$. This is bad, since it blows up near zero. So we use $L^4 L^4 L^2$ Holder inequality to get
\begin{equation*}
\begin{split}
\left( \int_0^t \frac{\Lambda_{m,\infty}(s)}{\sqrt{t-\tau}}\|\nabla S^v_n\|_{m-2}^{\frac{1}{2}} d\tau \right)^2 &= \left( \int_0^t \frac{\Lambda_{m,\infty}}{(t-\tau)^{\frac{1}{8}}} \|\nabla S^v_n\|_{m-2}^{\frac{1}{2}} \frac{1}{(t-\tau)^{\frac{3}{8}}} d\tau \right)^2 	\\
&\leq \left( \left(\int_0^t \frac{\Lambda^4_{\infty,m}}{\sqrt{t-\tau}} \right)^{\frac{1}{4}} \left(\int_0^t \|\nabla S^v_n\|_{m-2}^2 \right)^{\frac{1}{4}} \left(\int_0^t \frac{1}{(t-\tau)^{\frac{3}{4}}} \right)^{\frac{1}{2}} \right)^2 	\\
&\leq \left(\int_0^t \frac{\Lambda^4_{\infty,m}}{\sqrt{t-\tau}} \right)^{\frac{1}{2}} \left(\int_0^t \|\nabla S^v_n\|_{m-2}^2 \right)^{\frac{1}{2}} 4t^{\frac{1}{4}}.
\end{split}
\end{equation*}
Now we can use Young's inequality to get the terms what we want. We skip other terms, since they are nearly same. Finally we get, 
\begin{equation} \label{9.56}
\begin{split}
\varepsilon\|\partial_z\rho_v\|_{L^\infty}^2, \ \varepsilon\|\partial_z\rho_B\|_{L^\infty}^2 &\leq \Lambda_{m,\infty}(0) + 2\int_0^t \left( \varepsilon\|\nabla\mathcal{V}^m\|^2 + \varepsilon\|\nabla\mathcal{B}^m\|^2 + \varepsilon\|\nabla S^v_n\|_{m-2}^2 + + \varepsilon\|\nabla S^B_n\|_{m-2}^2 \right) 	\\
&\quad + ( 1 + 16\sqrt{t})\int_0^t \frac{\Lambda_{m,\infty}(s)}{\sqrt{t-\tau}} d\tau.
\end{split}
\end{equation}
\end{proof}

\subsection{$\int_0^t \sqrt{\varepsilon}\|\nabla^2 v\|_{1,\infty}$, $\int_0^t \sqrt{\varepsilon}\|\nabla^2 B\|_{1,\infty}$ estimates}
We need estimates of $\int_0^t \sqrt{\varepsilon}\|\nabla^2 v\|_{1,\infty}$, $\int_0^t \sqrt{\varepsilon}\|\nabla^2 B\|_{1,\infty}$ later.
\begin{lemma} \label{lemma 9.8}
Let $m \geq 6$ and $\sup_{[0,T]}\Lambda_{m,\infty}(t) \leq M$. Then 
\begin{equation} \label{9.57}
\sqrt{\varepsilon}\int_0^t \|\partial_{zz} v\|_{1,\infty} \leq \Lambda(M)(1+16\sqrt{t})\sqrt{t} \Big( 1 + \varepsilon\int_0^T( \|\nabla\mathcal{V}^m\|^2 + \|\nabla\mathcal{B}^m\|^2 + \|\nabla S_n^v\|^2_{m-2} + \|\nabla S_n^B\|^2_{m-2} )\Big).
\end{equation}
This is same for $\int_0^t \sqrt{\varepsilon}\|\partial_{zz} B\|_{1,\infty}$.
\end{lemma}
\begin{proof}
It is similar as Proposition \ref{proposition 9.7}. One more conormal derivative order changes nearly nothing, since our $m$ is sufficiently large. Integration for $t$ changes $\frac{1}{\sqrt{t-\tau}}$ into $\sqrt{t}$. 
\end{proof}

\begin{proposition} \label{proposition 9.9}
Under the same assumption, we get the following.
\begin{equation} \label{9.58}
\begin{split}
\sqrt{\varepsilon}\int_0^t \|\nabla^2 v\|_{1,\infty},\,\, \sqrt{\varepsilon}\int_0^t \|\nabla^2 B\|_{1,\infty} \leq \Lambda(M)(1+t)^2\left( 1 + \varepsilon\int_0^t( \|\nabla\mathcal{V}^m\|^2 + \|\nabla\mathcal{B}^m\|^2 + \|\nabla S_n^v\|^2_{m-2} + \|\nabla S_n^B\|^2_{m-2} ) \right).
\end{split}
\end{equation}
\end{proposition}
\begin{proof}
First, estimate of $\sqrt{\varepsilon}\int_0^t \|\nabla v\|_{2,\infty}$ is easy, because we already know,
$$
\|\nabla v\|_{2,\infty} \leq \Lambda_{m,\infty}(t) ( \|S_n^v\|_{2,\infty} + \|v\|_{3,\infty} ).
$$
Therefore, the following is easy.
$$
\sqrt{\varepsilon}\int_0^t \|\nabla v\|_{2,\infty} \leq \int_0^t \left( \Lambda_{m,\infty}(s) + \varepsilon\|\nabla\mathcal{V}^m\|^2 + \varepsilon\|\nabla S_n^v\|^2_{m-2} \right).
$$
Next, with help of Lemma \ref{lemma 9.8}, we have everything we need to finish the proof. We skip the detail.
\end{proof}

\section{Vorticity estimate}
We could not get estimate for $\|\partial_z v\|_{H^{m-1}_{co}}$ in previous sections. Instead $m-2$ was optimal. However, if we weaken $L^\infty_t$ to $L^{\alpha \geq 2}_t$  we may get the $m-1$ regularity. We define 
\begin{equation} \label{10.1}
\begin{split}
\omega_{v} &:= \nabla^\varphi\times v,\quad \omega_{v} = (\nabla\times u)(t,\Phi),	\\
\omega_{B} &:= \nabla^\varphi\times B,\quad \omega_{B} = (\nabla\times H)(t,\Phi).
\end{split}
\end{equation}
Since 
\begin{equation} \label{10.2}
\begin{split}
\omega_{v} \times \mathbf{n} &= {\frac{1}{2}} \big( D^\varphi \mathbf{n}  - (D^\varphi v)^t \mathbf{n} \big) = \mathbf{S}^{\varphi} v \mathbf{n} -  (D^\varphi v)^t \mathbf{n},
\end{split}
\end{equation}
we get
$$ \omega_{v} \times \mathbf{n} =   {1 \over 2} \partial_{n} u  - g^{ij}\big( \partial_{j}v \cdot \mathbf{n} \big) \partial_{y^i} .$$
By using 
\begin{equation} \label{10.2}
\begin{split}
	 \p_{N} u = {1 + |\partial_{1}\varphi|^2 + |\partial_{2} \varphi |^2 \over \partial_{z}\varphi} \partial_{z} v - \partial_{1}\varphi \partial_{1}v
- \partial_{2} \varphi \partial_{2} v,
\end{split}
\end{equation}
 
\noindent we get the following estimate. Estimate for $B$ is exactly same as $v$.
\begin{equation} \label{10.2}
\begin{split}
\|Z^{m-1} \partial_z v\| &\leq \Lambda_{m,\infty}(t)(\|v\|_{m} + |h|_{m-\frac{1}{2}} + \|\omega_v\|_{m-1} ),	\\
\|Z^{m-1} \partial_z B\| &\leq \Lambda_{m,\infty}(t)(\|B\|_{m} + |h|_{m-\frac{1}{2}} + \|\omega_B\|_{m-1} ).
\end{split}
\end{equation}
These imply that we suffice to control $\|\omega_{v}\|_{m-1}$ and $\|\omega_{B}\|_{m-1}$ to control $\|\partial_z v\|_{m-1}$ and $\|\partial_z B\|_{m-1}$, respectively. Applying $\nabla^\varphi\times$ kill pressure term in Navier-Stokes equation, therefore we get similar structure from two main PDEs.
\begin{equation} \label{10.2}
\begin{split}
& \partial_t^\varphi\omega_v + (v\cdot\nabla^\varphi)\omega_v - (B\cdot\nabla^\varphi)\omega_B - (\omega_v\cdot\nabla^\varphi)v + (\omega_B\cdot\nabla^\varphi)B = \varepsilon\triangle^\varphi\omega_v,	\\
& \partial_t^\varphi\omega_B + (v\cdot\nabla^\varphi)\omega_B - (B\cdot\nabla^\varphi)\omega_v - (\omega_v\cdot\nabla^\varphi)B + (\omega_B\cdot\nabla^\varphi)v = \varepsilon\triangle^\varphi\omega_B.
\end{split}
\end{equation}
Situation is quite different to $S_n$ case. The reason we used $S^{v}_n$ instead of $\partial_z v$ is that it is equivalent to $\partial_z v$ and moreover, it is zero boundary condition. For vorticity, we have
$$
\omega\times\mathbf{n} = \Pi(\omega\times\mathbf{n}) \neq 0,\quad\text{on}\quad\partial S,
$$
in general. Moreover equation of $\omega\times\mathbf{n}$ is more complicate than the equation of $\omega$. This means $\omega\times\mathbf{n}$ has no advantage than $\omega$. Thus we just use $\|\omega\|_{m-1}$ directly. Applying $Z^{\alpha}$ gives, ($|\alpha|\leq m-1$)
\begin{equation} \label{10.3}
\partial_t^\varphi Z^{\alpha}\omega_v + (v\cdot\nabla^\varphi)Z^{\alpha}\omega_v - (B\cdot\nabla^\varphi)Z^{\alpha}\omega_B - \varepsilon\triangle^\varphi Z^\alpha\omega_v =  F,
\end{equation}
where
\begin{equation*}
\begin{split}
F &:= Z^\alpha(\omega_v\cdot\nabla^\varphi v) - Z^\alpha(\omega_B\cdot\nabla^\varphi B) + \mathcal{C}_S,  \\
\mathcal{C}_S &:= \mathcal{C}_S^1 + \mathcal{C}_S^2 + \mathcal{C}_S^3,
\end{split}
\end{equation*}
with
\begin{equation*}
\begin{split}
\mathcal{C}_S^1 &:= [Z^\alpha v_y]\cdot\nabla_y\omega_v + [Z^\alpha,V_z]\partial_z\omega_v := \mathcal{C}^1_{S_y} + \mathcal{C}^1_{S_z},  \\
\mathcal{C}_S^2 &:= \varepsilon[Z^\alpha,\triangle^\varphi]\omega_v,  \\
\mathcal{C}_S^3 &:= -[Z^\alpha B_y]\cdot\nabla_y\omega_B + [Z^\alpha,\frac{B\cdot N}{\partial_z\varphi}]\partial_z\omega_B := \mathcal{C}^3_{S_y} + \mathcal{C}^3_{S_z}.
\end{split}
\end{equation*}
For $\o_{B}$,
\begin{equation} \label{10.4}
\partial_t^\varphi Z^{\alpha}\omega_B + (v\cdot\nabla^\varphi)Z^{\alpha}\omega_B - (B\cdot\nabla^\varphi)Z^{\alpha}\omega_v - \varepsilon\triangle^\varphi Z^{\alpha}\omega_B =  E,
\end{equation}
where
\begin{equation*}
\begin{split}
E &:= Z^\alpha(\omega_v\cdot\nabla^\varphi B) - Z^\alpha(\omega_B\cdot\nabla^\varphi v) + \bar{\mathcal{C}}_S,  \\
\bar{\mathcal{C}}_S &:= \bar{\mathcal{C}}_S^1 + \bar{\mathcal{C}}_S^2 + \bar{\mathcal{C}}_S^3,
\end{split}
\end{equation*}
with
\begin{equation*}
\begin{split}
\bar{\mathcal{C}}_S^1 &:= [Z^\alpha v_y]\cdot\nabla_y\omega_B + [Z^\alpha,V_z]\partial_z\omega_B := \bar{\mathcal{C}}^1_{S_y} + \bar{\mathcal{C}}^1_{S_z},  \\
\bar{\mathcal{C}}_S^2 &:= \varepsilon[Z^\alpha,\triangle^\varphi]\omega_B,  \\
\bar{\mathcal{C}}_S^3 &:= -[Z^\alpha B_y]\cdot\nabla_y\omega_v + [Z^\alpha,\frac{B\cdot \mathbf{N}}{\partial_z\varphi}]\partial_z \omega_v := \bar{\mathcal{C}}^3_{S_y} + \bar{\mathcal{C}}^3_{S_z}.
\end{split}
\end{equation*}
And for boundary data we use the following boundary estimate
\[
	|\nabla v(\cdot,0)|_{s} \leq \Lambda (\|v\|_{1,\infty} + |h|_{2,\infty}) ( |v(\cdot,0)|_{s+1} + |h|_{s+1} ),
\]
from Lemma 5.5 in \cite{NMFR1} to get
\begin{equation} \label{10.5}
|(Z^\alpha\omega_v)^b| \leq \Lambda_{6,\infty}(t)(|v|^b_m + |h|_m),
\end{equation}
\begin{equation} \label{10.6}
|(Z^\alpha\omega_B)^b| \leq \Lambda_{6,\infty}(t)(|B|^b_m + |h|_m).
\end{equation}
Using Proposition \ref{proposition 3.2},
\begin{equation*}
\begin{split}
|(Z^\alpha\omega_v)^b| &\leq \Lambda_{6,\infty}(t)\left( \|\nabla\mathcal{V}^m\|^{\frac{1}{2}}\|\mathcal{V}^m\|^{\frac{1}{2}} + \|\mathcal{V}^m\| + |h|_m \right), 	\\
|(Z^\alpha\omega_B)^b| &\leq \Lambda_{6,\infty}(t)\left( \|\nabla\mathcal{B}^m\|^{\frac{1}{2}}\|\mathcal{B}^m\|^{\frac{1}{2}} + \|\mathcal{B}^m\| + |h|_m \right),
\end{split}
\end{equation*}
and surely,
\begin{equation} \label{10.7}
\begin{split}
\sqrt{\varepsilon}\int_0^t |(Z^\alpha\omega_v)^b|^2 &\leq  \Lambda_{6,\infty}(t) \sqrt{\varepsilon} \int_0^t \left( \|\nabla\mathcal{V}^m\|\|\mathcal{V}^m\| + \|\mathcal{V}^m\|^2 + |h|_m^2 \right),	\\
\sqrt{\varepsilon}\int_0^t |(Z^\alpha\omega_B)^b|^2 &\leq   \Lambda_{6,\infty}(t) \sqrt{\varepsilon} \int_0^t \left( \|\nabla\mathcal{B}^m\|\|\mathcal{B}^m\| + \|\mathcal{B}^m\|^2 + |h|_m^2 \right).
\end{split}
\end{equation}
Using Young's inequality,
\begin{equation*}
\begin{split}
\sqrt{\varepsilon}\int_0^t |(Z^\alpha\omega_v)^b|^2 &\leq \varepsilon\int_0^t\|\nabla\mathcal{V}^m\|^2 + \int_0^t\Lambda_{6,\infty}(s)\left(\|\mathcal{V}^m\|^2 + |h|_m^2\right), 	\\
\sqrt{\varepsilon}\int_0^t |(Z^\alpha\omega_B)^b|^2 &\leq \varepsilon \int_0^t\|\nabla\mathcal{B}^m\|^2 + \int_0^t\Lambda_{6,\infty}(s)\left(\|\mathcal{B}^m\|^2 + |h|_m^2\right).
\end{split}
\end{equation*}
We split
$$
Z^\alpha\omega_v := \omega^\alpha_{v,h} + \omega^\alpha_{v,nh},\quad Z^\alpha\omega_B := \omega^\alpha_{B,h} + \omega^\alpha_{B,nh},
$$
and estimate each four terms on the right hand sides.  \\
I) $\omega^\alpha_{v,nh},\,\omega^\alpha_{B,nh}$ solve, non-homogeneous equation,
\begin{equation} \label{10.8}
\partial_t^\varphi \omega^\alpha_{v,nh} + (v\cdot\nabla^\varphi)\omega^\alpha_{v,nh} - (B\cdot\nabla^\varphi)\omega^\alpha_{B,nh} - \varepsilon\triangle^\varphi\omega^\alpha_{v,nh} =  F,
\end{equation}
with initial and zero-boundary condition,
$$
(\omega^\alpha_{v,nh})^b = 0,\,\,\,(\omega^\alpha_{v,nh})_{t=0} = \omega_v(0),
$$
and
\begin{equation} \label{10.9}
\partial_t^\varphi \omega^\alpha_{B,nh} + (v\cdot\nabla^\varphi)\omega^\alpha_{B,nh} - (B\cdot\nabla^\varphi)\omega^\alpha_{v,nh} - \varepsilon\triangle^\varphi\omega^\alpha_{B,nh} =  E,
\end{equation}
with initial and zero-boundary condition,
$$
(\omega^\alpha_{B,nh})^b = 0,\,\,\,(\omega^\alpha_{B,nh})_{t=0} = \omega_B(0).
$$
II) Meanwhile, $\omega^\alpha_{v,h},\,\omega^\alpha_{B,h}$ solve, homogeneous equation,
\begin{equation} \label{10.10}
\partial_t^\varphi \omega^\alpha_{v,h} + (v\cdot\nabla^\varphi)\omega^\alpha_{v,h} - (B\cdot\nabla^\varphi)\omega^\alpha_{B,h} - \varepsilon\triangle^\varphi\omega^\alpha_{v,h} =  0,
\end{equation}
with zero initial and general boundary condition,
$$
(\omega^\alpha_{v,h})^b = (Z^\alpha\omega_v)^b,\,\,\,(\omega^\alpha_{v,h})_{t=0} = 0,
$$
and
\begin{equation} \label{10.11}
\partial_t^\varphi \omega^\alpha_{B,h} + (v\cdot\nabla^\varphi)\omega^\alpha_{B,h} - (B\cdot\nabla^\varphi)\omega^\alpha_{v,h} - \varepsilon\triangle^\varphi\omega^\alpha_{B,h} =  0,
\end{equation}
with zero initial and general boundary condition,
$$
(\omega^\alpha_{B,h})^b = (Z^\alpha\omega_B)^b,\,\,\,(\omega^\alpha_{B,h})_{t=0} = 0.
$$
We state the energy estimates for these two vorticity terms. For non-homogeneous terms, we define
\begin{equation} \label{10.12}
\begin{split}
\|\omega_{v,nh}^{m-1}\|^2 &:= \sum_{|\alpha|\leq m-1}\|\omega_{v,nh}^\alpha\|^2,\quad \int_0^t\|\nabla\omega_{v,nh}^{m-1}\|^2 := \int_0^t\sum_{|\alpha|\leq m-1}\|\nabla\omega_{v,nh}^\alpha\|^2,  \\
\|\omega_{B,nh}^{m-1}\|^2 &:= \sum_{|\alpha|\leq m-1}\|\omega_{B,nh}^\alpha\|^2,\quad \int_0^t\|\nabla\omega_{B,nh}^{m-1}\|^2 := \int_0^t\sum_{|\alpha|\leq m-1}\|\nabla\omega_{B,nh}^\alpha\|^2.
\end{split}
\end{equation}
For homogeneous terms, we define similarly
\begin{equation} \label{10.13}
\begin{split}
\|\omega_{v,h}^{m-1}\|^2 &:= \sum_{|\alpha|\leq m-1}\|\omega_{v,h}^\alpha\|^2 \quad
\int_0^t\|\nabla\omega_{v,h}^{m-1}\|^2 := \int_0^t\sum_{|\alpha|\leq m-1}\|\nabla\omega_{v,h}^\alpha\|^2, 	\\
\|\omega_{B,h}^{m-1}\|^2 &:= \sum_{|\alpha|\leq m-1}\|\omega_{B,h}^\alpha\|^2 \quad
\int_0^t\|\nabla\omega_{B,h}^{m-1}\|^2 := \int_0^t\sum_{|\alpha|\leq m-1}\|\nabla\omega_{B,h}^\alpha\|^2.
\end{split}
\end{equation}

\noindent In following subsections, we estimate $\omega^{\alpha}_{v,nh}$, $\omega^{\alpha}_{B,nh}$, $\omega^{\alpha}_{v,h}$, and $\omega^{\alpha}_{B,h}$ since 
\begin{equation} \label{10.14}
\begin{split}
\|\omega\|_{m-1}^2 &= \sum_{|\alpha|\leq m-1} \|Z^\alpha\omega\|^2 = \sum_{|\alpha|\leq m-1} \|\omega^\alpha_{h} + \omega^\alpha_{nh}\|^2 \leq 2\|\omega^{m-1}_{nh}\|^2 + 2\|\omega^{m-1}_{h}\|^2.  \\
\end{split}
\end{equation}

\subsection{Non-homogeneous estimate}
We estimate $\omega^\alpha_{v,nh},\omega^\alpha_{B,nh}$.
\begin{proposition} \label{proposition 10.1}
We have the following vorticity estimate for $\omega^\alpha_{v,nh}$.
\begin{equation} \label{10.15}
\begin{split}
&\|\omega^{m-1}_{v,nh}\|^2 + 2\varepsilon\int_0^t \int_S |\nabla^\varphi\omega^{m-1}_{v,nh}|^2 dV_t ds - 2\int_0^t\int_S(B\cdot\nabla^\varphi) \omega^{m-1}_{B,nh} \cdot\omega^{m-1}_{v,nh} dV_t ds  \\
&\leq \Lambda_0\|\omega^{m-1}_{v,nh}(0)\|^2 + \int_0^t \Lambda_{m,\infty}(s) \left( \|v\|_{E^m}^2 + \|B\|_{E^m}^2 + \|\omega_v\|_{m-1}^2 + \|\omega_B\|_{m-1}^2 + |h|_m^2 + \varepsilon |h|_{m+\frac{1}{2}}^2 \right)  \\
&\quad + \varepsilon\int_0^t\Lambda_{m,\infty}(s) \left( \|\nabla S^v_n\|_{m-2}^2 + \|\nabla S^B_n\|_{m-2}^2 \right) ds.
\end{split}
\end{equation}
\end{proposition}
\begin{proof}
Using equation for $\omega^\alpha_{v,nh}$, (with dirichlet boundary condition) we get $L^2$ type energy estimate.
\begin{equation} \label{10.16}
\frac{1}{2}\frac{d}{dt}\int_S |\omega^\alpha_{v,nh}|^2 dV_t + \varepsilon\int_S |\nabla^\varphi\omega^\alpha_{v,nh}|^2 dV_t - \int_S(B\cdot\nabla^\varphi) \omega^\alpha_{B,nh} \cdot\omega^\alpha_{v,nh} dV_t = \int_S F \cdot \omega^\alpha_{v,nh} dV_t.
\end{equation}
I) $\|Z^\alpha(\omega_v\cdot\nabla^\varphi v)\|, \ \|Z^\alpha(\omega_B\cdot\nabla^\varphi B)\|$ estimates.\\
Simply we get
\begin{equation*} 
\begin{split}
\|Z^\alpha(\omega_v\cdot\nabla^\varphi v)\| &\leq \Lambda_{m,\infty}(t) \left( \|w_v\|_{m-1} + \|v\|_m + |h|_{m-\frac{1}{2}} \right), 	\\
\|Z^\alpha(\omega_B\cdot\nabla^\varphi B)\| &\leq \Lambda_{m,\infty}(t) \left( \|w_B\|_{m-1} + \|B\|_m + |h|_{m-\frac{1}{2}} \right).
\end{split}
\end{equation*}
II) $\left|\int_S\mathcal{C}_S^2 \cdot\omega^\alpha_{v,nh} dV_t\right|$ estimate.\\
As like in $S^{v,B}_n$ estimate,
$$
\left|\int_S\mathcal{C}_S^2 \cdot\omega^\alpha_{v,nh} dV_t\right| \leq \Lambda_0\left( \sqrt{\varepsilon}\|\nabla^\varphi\omega^\alpha_{v,nh}\| + \|\omega_v\|_{m-1}\right)\left(\sqrt{\varepsilon}\|\nabla\omega_v\|_{m-2} + \|\omega_v\|_{m-1} + \Lambda_{m,\infty}\left( |h|_{m-\frac{1}{2}} + \sqrt{\varepsilon}|h|_{m+\frac{1}{2}} \right) \right).
$$
III) $\|\mathcal{C}_{S_y}^1\|, \ \|\mathcal{C}_{S_y}^3\|$ estimates. \\
This is also similar as previous $S^{v,B}_n$ estimate,
\begin{equation*} 
\begin{split}
\|\mathcal{C}_{S_y}^1\| &\leq \Lambda_{m,\infty}(t) \left( \|\omega_v\|_{m-1} + \|v\|_m + |h|_m \right), 	\\
\|\mathcal{C}_{S_y}^3\| &\leq \Lambda_{m,\infty}(t) \left( \|\omega_B\|_{m-1} + \|B\|_m + |h|_m \right).
\end{split}
\end{equation*}
IV) $\|\mathcal{C}_{S_z}^1\|, \ \|\mathcal{C}_{S_z}^3\|$ estimates.\\
This is the main part of proof. In $S^{v,B}_n$ section, $|\alpha|=m-2$ was optimal, since $|h|_{m-\frac{1}{2}}$ and $\|v\|_{E^m}$. Note that in pressure estimate, if we use $\nabla^\varphi\cdot(v\cdot\nabla^\varphi)v = \nabla^\varphi v : (\nabla^\varphi v)^T$ then we could get 1 regularity for $v$. Nevertheless $m-1$ is not available because of the worst $h$ regularity term which come from $\|\nabla q^E\|$ in $F_S$ and $C_{S_z}$. But now, since we consider vorticity, $\|\nabla q^E\|$ does not appear, so only the problem is $|h|_{m-\frac{1}{2}}$ in $C_{S_z}$. What we show here is that we can get $\frac{1}{2}$ regularity of $h$ in fact. (We did not have to do same thing in $S_n$ section, since pressure generate $|h|_{m-\frac{1}{2}}$.)\\
To estimate $\mathcal{C}_{S_z}^1$, we should estimate terms like (with $|\gamma|+|\beta|\leq m-1,\,\,|\gamma|\leq m-2$)
$$
\|Z^\beta V_z\partial_z Z^\gamma \omega_v\|.
$$
We write this as
$$
c_{\tilde{\beta}} Z^{\tilde{\beta}}\left(\frac{1-z}{z} V_z\right) Z_3 Z^{\tilde{\gamma}} \omega_v,
$$
with $|\tilde{\gamma}|+|\tilde{\beta}|\leq m-1,\,\,|\tilde{\gamma}|\leq m-2$. Then using Lemma 8.4 in \cite{NMFR1}(variant of Hardy's inequality when function is zero at $z=0$),
$$
\leq \Lambda_{m,\infty}(t) \Big(\|\omega_v\|_{m-1} + |h|_{m-\frac{1}{2}} + \left\|\frac{1-z}{z}Z(v\cdot \mathbf{N} - \partial_z\varphi)\right\|_{m-2}\Big).
$$
Last term is main term.
$$
\left\|\frac{1-z}{z}Z(v\cdot \mathbf{N} - \partial_z\varphi)\right\|_{m-2} \leq \Lambda_{m,\infty}(t) \Big( \|v\|_{E^m} + \sum_{|\alpha|\leq m-1}\left\|v\cdot \partial_z Z^\alpha\mathbf{N} - \partial_z Z^\alpha\partial_t\varphi\right\|_{L^2} \Big).
$$
If we brutely estimate this, then we get
$$
\|\partial_z Z^\alpha\mathbf{N}\| \sim |\nabla\varphi|_{m-1+1} \sim |\varphi|_{m+1} \sim |h|_{m+\frac{1}{2}},
$$
so we loose $\frac{1}{2}$ regularity. We treat this carefully. First, we should see that $Z_3$ does not lose regularity of $h$. From definition of $\varphi(t,y,z) = Az + \eta(t,y,z)$, 
$$
|Z_3 \hat{\eta}| = \left|\frac{z}{1-z}\partial_z\hat{\eta}(\xi,z)\right| = \left|\frac{z}{1-z}\partial_z\left(\chi(\xi z)\hat{h}(\xi)\right)\right| \leq \left|\frac{z}{1-z}\xi\cdot\nabla\chi(\xi z)\right| \leq \left| \chi_2(\xi z)\hat{h}\right|,
$$
where $\chi_2$ is 1 on $B_2(0)$, which is bigger support than $\chi$.\\
This is because, $\nabla\chi$ has support on annular domain, $B_2(0)-B_1(0)$. Above regularity means applying $Z_3$ does not reduce 1 regularity. So, $\alpha_3 \neq 0$ case is not harmful, which means, 
$$
\left\|v\cdot \partial_z Z^\alpha\mathbf{N} - \partial_z Z^\alpha\partial_t\varphi\right\|_{L^2} \leq \Lambda_{m,\infty}(t) \left( |h|_{m-\frac{1}{2}} + |\partial_t h|_{m-\frac{3}{2}}\right) \leq \Lambda_{m,\infty}(t) \left( |h|_{m-\frac{1}{2}} + \|v\|_{E^m}\right),\,\,\alpha_3 \neq 0.
$$ 
When  $\alpha_3 = 0$,
$$
v\cdot\partial_z Z^\alpha\mathbf{N} - \partial_z Z^\alpha\partial_t\eta = -v_1\partial_z\left(\psi_z *_y \partial_1 Z^\alpha h\right) -v_2\partial_z\left(\psi_z *_y \partial_2 Z^\alpha h\right) - \partial_z\left(\psi_z *_y \partial_t Z^\alpha h\right) := \mathcal{T}_{\alpha},
$$
where
$$
\hat{\psi}_z = \chi(\xi z).
$$
So by inverse fourier transformation with respect to horizontal variable, 
$$
\psi_z(y) = \frac{1}{z^2}\check{\chi}\left(\frac{y}{z}\right).
$$
Note that $\chi$ is in Schwartz class, so $\check{\chi}$ is also in Schwartz class. Moreover, when $|z|\leq 1$,
\begin{equation} \label{10.17}
Z^\alpha\psi_z(z) = \mathbf{D}^\alpha\psi_z(z) = -\frac{2}{z^3}\check{\chi}\left(\frac{y}{z}\right) + \frac{1}{z^2}(\nabla\check{\chi})\left(\frac{y}{z}\right)\cdot y \leq \frac{2}{z^3}(1 + |y|)\zeta(\frac{y}{z}).
\end{equation}
For some function $\zeta$ with compact suuport neat origin. This is because $\nabla\check{\chi}$ is also in Schwartz class.\\
i) When $z \leq -1$,
$$
\left\|\mathcal{T}_\alpha\right\|_{L^2(S\cap{|z|\geq 1)}}  = \left\| 
-v_1\partial_z\left(\mathbf{D}^\alpha\psi_z *_y \partial_1 h\right) -v_2\partial_z\left(\mathbf{D}^\alpha\psi_z *_y \partial_2 h\right) - \partial_z\left(\mathbf{D}^\alpha\psi_z *_y \partial_t \right)
\right\|
$$
$$
\leq \Lambda_{m,\infty}(t) \left( \left\|\partial_z\left(\frac{1}{z^{m-1}}(\mathbf{D}^\alpha\psi_z) *_y \nabla_{y} h\right)\right\|_{L^2(S\cap{|z|\geq 1)}} + \left\|\partial_z\left(\frac{1}{z^{m-1}}(\mathbf{D}^\alpha\psi_z) *_y \partial_t h\right)\right\|_{L^2(S\cap{|z|\geq 1)}} \right).
$$
Since $D^\alpha\psi_z$ is also similar as $\psi_z$, we get
\begin{equation} \label{10.18}
\left\|\mathcal{T}_\alpha\right\|_{L^2(S\cap{|z|\geq 1)}} \leq \Lambda_{m,\infty}(t) \left( \|v\|_{E^1} + |h|_{1} \right).
\end{equation}
ii) When $z \geq -1$
\begin{equation*}
\begin{split}
\mathcal{T}_\alpha &= -v_1\partial_z\left(\psi_z *_y \partial_1 Z^\alpha h\right) -v_2\partial_z\left(\psi_z *_y \partial_2 Z^\alpha h\right) - \partial_z\left(\psi_z *_y \partial_t Z^\alpha h\right)  \\
&= -v_1^b\partial_z\left(\psi_z *_y \partial_1 Z^\alpha h\right) -v_2^b\partial_z\left(\psi_z *_y \partial_2 Z^\alpha h\right) - \partial_z\left(\psi_z *_y \partial_t Z^\alpha h\right) + \mathcal{R},
\end{split}
\end{equation*}
where
$$
\mathcal{R} := (v_1^b - v_1)\partial_z\left(\psi_z *_y \partial_1 Z^\alpha h\right) + (v_2^b - v_2)\partial_z\left(\psi_z *_y \partial_2 Z^\alpha h\right).
$$
Using taylor expansion with respect to z, we get
$$
|\mathcal{R}| \leq \Lambda_{m,\infty}(t) |z|\left|\partial_z(\psi_z *_y Z^\alpha\nabla_{y} h)\right| \leq \left|Z_3(\psi_z *_y Z^\alpha\nabla_{y} h)\right|.
$$
Note that $Z_3$ do nothing about regularity to $\psi_z$, so
\begin{equation} \label{10.19}
\|\mathcal{R}\| \leq \Lambda_{m,\infty}(t) |h|_{m-\frac{1}{2}}.
\end{equation}
For first part in $\mathcal{T}_\alpha$,
\begin{equation*}
\begin{split}
&-v_1^b\partial_z\left(\psi_z *_y \partial_1 Z^\alpha h\right) -v_2^b\partial_z\left(\psi_z *_y \partial_2 Z^\alpha h\right) - \partial_z\left(\psi_z *_y \partial_t Z^\alpha h\right)  \\
&\quad = \partial_z\left(\psi_z *_y \left( -v_1^b\partial_1 Z^\alpha h - v_2^b\partial_2 Z^\alpha h - \partial_t Z^\alpha h \right)\right)  \\
&\quad\quad + \partial_z\int_S\Big( (v_1(t,y',0)-v_1(t,y,0))\psi_z(y-y')\partial_1 Z^\alpha h(t,y') \\
&\quad\quad  + (v_2(t,y',0)-v_2(t,y,0))\psi_z(y-y')\partial_2 Z^\alpha h(t,y') \Big)dy.
\end{split}
\end{equation*}
By taylor expansion, 
$$
|(v_i^b(t,y',0)- v_i^b(t,y,0))\partial_z\psi_z(y-y')| \leq |\nabla_y v_i^b|_{L^\infty} \left|(y-y')\frac{2}{z^3}\partial_z\check{\psi}(\frac{y-y'}{z}) \right|.
$$
So, for $\forall z \in (0,1]$, 
\begin{equation*}
\begin{split}
&\sup_{z \in (0,1]} \left\|\partial_z\int_{\partial S} (v_i(t,y',0)-v_i(t,y,0))\psi_z(y-y')\partial_i Z^\alpha h(t,y')\right\|  \\
&\quad \leq \|\nabla_y v_i^b\|_{L^\infty} \sup_{z \in (0,1]} \int_{\partial S} \left|\frac{y-y'}{z^3}\right|^2 \zeta^2(\frac{y-y'}{z})z^2 d\left(\frac{y-y'}{z}\right)  \\
&\quad \leq \|\nabla_y v_i^b\|_{L^\infty} \sup_{z \in (0,1]} \int_{\partial S} y^2\zeta^2(y)dy \leq \Lambda_{m,\infty}(t).
\end{split}
\end{equation*}
For the first term of $\mathcal{T}_\alpha$,
$$
\|-v_1^b\partial_1 Z^\alpha h - v_2^b\partial_2 Z^\alpha h - \partial_t Z^\alpha h\|_{H^{\frac{1}{2}}} = \| \mathcal{C}^\alpha(h) - (\mathcal{V}^\alpha)^b - v_3^b \|_{H^{\frac{1}{2}}} \leq \Lambda_{m,\infty}(t) \left( \|v\|_{E^m} + |h|_{m-\frac{1}{2}} \right),
$$
thus 
$$
\left\| \partial_z\left(\psi_z *_y \left( -v_1^b\partial_1 Z^\alpha h - v_2^b\partial_2 Z^\alpha h - \partial_t Z^\alpha h \right)\right) \right\| \leq \Lambda_{m,\infty}(t) \left( \|v\|_{E^m} + |h|_{m-\frac{1}{2}} \right).
$$
Finally we get
$$
\left\|\mathcal{T}_\alpha\right\|_{L^2(S\cap{|z|\neq 1)}} \leq \Lambda_{m,\infty}(t) \left( \|v\|_{E^m} + |h|_m \right).
$$
Considering i) and ii), and $\alpha_3 \neq 0$,
\begin{equation} \label{10.20}
\|\mathcal{C}_{S_Z}^1\| \leq \Lambda_{m,\infty}(t) \left( \|\omega_v\|_{m-1} + \|v\|_{E^m} + |h|_m \right).
\end{equation}
For $\|\mathcal{C}_{S_z}^3\|$, we can do similar thing for $B$, since $B$ is zero on the boundary, so control is better. 
\begin{equation} \label{10.21}
\|\mathcal{C}_{S_Z}^3\| \leq \Lambda_{m,\infty}(t) \left( \|\omega_B\|_{m-1} + \|B\|_{E^m} + |h|_m \right).
\end{equation}
Considering I) $\sim$ IV),
\begin{equation*}
\begin{split}
&\|\omega^\alpha_{v,nh}\|^2 + 2\varepsilon\int_0^t \int_S |\nabla^\varphi\omega^\alpha_{v,nh}|^2 dV_t ds - 2\int_0^t\int_S(B\cdot\nabla^\varphi) \omega^\alpha_{B,nh} \cdot\omega^\alpha_{v,nh} dV_t ds  \\
&\leq \Lambda_0\|\omega^\alpha_{v,nh}(0)\|^2 + \int_0^t \Lambda_{m,\infty}(s)\left( \|v\|_{E^m}^2 + \|B\|_{E^m}^2 + \|\omega_v\|_{m-1}^2 + \|\omega_B\|_{m-1}^2 + |h|_m^2 + \varepsilon|h|_{m+\frac{1}{2}}^2 \right)  \\
&\quad + \Lambda_0\varepsilon\int_0^t \left( \|\nabla\omega_v\|_{m-2}^2 + \|\nabla\omega_B\|_{m-2}^2 \right) ds.
\end{split}
\end{equation*}
For the last term,
$$
\sqrt{\varepsilon}\|\nabla\omega_v\|_{m-2} \leq \Lambda_{m,\infty}\left( \sqrt{\varepsilon}\|\partial_{zz}v\|_{m-2} + \sqrt{\varepsilon}\|\partial_{z}v\|_{m-1}\right) + \Lambda_{m,\infty}|h|_{m-\frac{1}{2}}.
$$
This gives, (by summing for all indices)
\begin{equation*} \label{10.22}
\begin{split}
&\|\omega^{m-1}_{v,nh}\|^2 + 2\varepsilon\int_0^t \int_S |\nabla^\varphi\omega^{m-1}_{v,nh}|^2 dV_t ds - 2\int_0^t\int_S(B\cdot\nabla^\varphi) \omega^{m-1}_{B,nh} \cdot\omega^{m-1}_{v,nh} dV_t ds  \\
&\leq \Lambda_0\|\omega^{m-1}_{v,nh}(0)\|^2 + \int_0^t \Lambda_{m,\infty}(s) \left( \|v\|_{E^m}^2 + \|B\|_{E^m}^2 + \|\omega_v\|_{m-1}^2 + \|\omega_B\|_{m-1}^2 + |h|_m^2 + \varepsilon|h|_{m+\frac{1}{2}}^2 \right)  \\
&+ \varepsilon\int_0^t\Lambda_{m,\infty}(s) \left( \|\nabla S^v_n\|_{m-2}^2 + \|\nabla S^B_n\|_{m-2}^2 \right) ds.
\end{split}
\end{equation*}
\end{proof}
\noindent For $\omega^\alpha_{B,nh}$, similarly we get,
\begin{proposition} \label{proposition 10.2}
We have the following vorticity estimates for $\omega^\alpha_{B,nh}$.
\begin{equation} \label{10.23}
\begin{split}
&\|\omega^{m-1}_{B,nh}\|^2 + 2\varepsilon\int_0^t \int_S |\nabla^\varphi\omega^{m-1}_{B,nh}|^2 dV_t ds + 2\int_0^t\int_S(B\cdot\nabla^\varphi) \omega^{m-1}_{B,nh} \cdot\omega^{m-1}_{v,nh} dV_t ds  \\
&\leq \Lambda_0\|\omega^{m-1}_{B,nh}(0)\|^2 + \int_0^t \Lambda_{m,\infty}(s) \left( \|v\|_{E^m}^2 + \|B\|_{E^m}^2 + \|\omega_v\|_{m-1}^2 + \|\omega_B\|_{m-1}^2 + |h|_m^2 + \varepsilon|h|_{m+\frac{1}{2}}^2 \right)  \\
&\quad + \varepsilon\int_0^t\Lambda_{m,\infty}(s) \left( \|\nabla S^v_n\|_{m-2}^2 + \|\nabla S^B_n\|_{m-2}^2 \right) ds.
\end{split}
\end{equation}
\end{proposition}
\begin{proof}
Using equation for $\omega^\alpha_{B,nh}$, (with dirichlet boundary condition) we get $L^2$ type energy estimate.
\begin{equation} \label{10.24}
\frac{1}{2}\frac{d}{dt}\int_S |\omega^\alpha_{B,nh}|^2 dV_t + \varepsilon\int_S |\nabla^\varphi\omega^\alpha_{B,nh}|^2 dV_t + \int_S(B\cdot\nabla^\varphi) \omega^\alpha_{B,nh} \cdot\omega^\alpha_{v,nh} dV_t = \int_S E \cdot \omega^\alpha_{B,nh} dV_t.
\end{equation}
I) $\|Z^\alpha(\omega_v\cdot\nabla^\varphi B)\|, \ \|Z^\alpha(\omega_B\cdot\nabla^\varphi v)\|$ estimates.  \\
Simply we get
\begin{equation*}
\begin{split}
\|Z^\alpha(\omega_v\cdot\nabla^\varphi B)\| &\leq \Lambda_{m,\infty}(t) \left( \|w_v\|_{m-1} + \|B\|_m + |h|_{m-\frac{1}{2}} \right), 	\\
\|Z^\alpha(\omega_B\cdot\nabla^\varphi v)\| &\leq \Lambda_{m,\infty}(t) \left( \|w_B\|_{m-1} + \|v\|_m + |h|_{m-\frac{1}{2}} \right).
\end{split}
\end{equation*}
II) $\left|\int_S\bar{\mathcal{C}}_S^2 \cdot\omega^\alpha_{B,nh} dV_t\right|$ estimate.  \\
As like in $S^{v,B}_n$ estimate,
$$
\left|\int_S \bar{\mathcal{C}}_S^2 \cdot\omega^\alpha_{B,nh} dV_t\right| \leq \Lambda_0\left( \sqrt{\varepsilon}\|\nabla^\varphi\omega^\alpha_{B,nh}\| + \|\omega_B\|_{m-1}\right)\left(\sqrt{\varepsilon}\|\nabla\omega_B\|_{m-2} + \|\omega_B\|_{m-1} + \Lambda_{m,\infty}(t) \big( |h|_{m-\frac{1}{2}} + \sqrt{\varepsilon}|h|_{m+\frac{1}{2}} \big) \right).
$$
III) $\|\bar{\mathcal{C}}_{S_y}^1\|, \ \|\bar{\mathcal{C}}_{S_y}^3\|$ estimates.  \\
This is also similar to previous $S^{v,B}_n$ estimate,
\begin{equation*}
\begin{split}
\|\bar{\mathcal{C}}_{S_y}^1\| &\leq \Lambda_{m,\infty}\left( \|\omega_B\|_{m-1} + \|v\|_m + |h|_m \right),	\\
\|\bar{\mathcal{C}}_{S_y}^3\| &\leq \Lambda_{m,\infty}\left( \|\omega_v\|_{m-1} + \|B\|_m + |h|_m \right).
\end{split}
\end{equation*}
IV) $\|\bar{\mathcal{C}}_{S_z}^1\|, \ \|\bar{\mathcal{C}}_{S_z}^3\|$ estimates.  \\
We skip the detail, since it is nearly same as previous Proposition \ref{proposition 10.1}.
\begin{equation} \label{10.25}
\|\mathcal{C}_{S_Z}^1\| \leq \Lambda_{m,\infty}\left( \|\omega_v\|_{m-1} + \|B\|_{E^m} + |h|_m \right),
\end{equation}
\begin{equation} \label{10.26}
\|\mathcal{C}_{S_Z}^3\| \leq \Lambda_{m,\infty}\left( \|\omega_B\|_{m-1} + \|v\|_{E^m} + |h|_m \right).
\end{equation}
Hence as like previous Proposition, we get the result.
\end{proof}

Summing above two Propositions \ref{proposition 10.1} and \ref{proposition 10.2}, we get the following.
\begin{proposition} \label{proposition 10.3}
Non-homogeneous part estimate.
\begin{equation} \label{10.27}
\begin{split}
&\|\omega^{m-1}_{v,nh}(t)\|^2 + \|\omega^{m-1}_{B,nh}(t)\|^2 + \varepsilon\int_0^t \int_S |\nabla^\varphi\omega^{m-1}_{v,nh}|^2 dV_t ds + \varepsilon\int_0^t \int_S |\nabla^\varphi\omega^{m-1}_{B,nh}|^2 dV_t ds 	\\
&\leq \Lambda_0\left( \|\omega^{m-1}_{v,nh}(0)\|^2 + \|\omega^{m-1}_{B,nh}(0)\|^2 \right) + \varepsilon\int_0^t\Lambda_{m,\infty} \left( \|\nabla S^v_n\|_{m-2}^2 + \|\nabla S^B_n\|_{m-2}^2 \right) ds 	\\
&\quad + \int_0^t \Lambda_{m,\infty}\left( \|v\|_{E^m}^2 + \|B\|_{E^m}^2 + \|\omega_v\|_{m-1}^2 + \|\omega_B\|_{m-1}^2 + |h|_m^2 + \varepsilon|h|_{m+\frac{1}{2}}^2 \right) ds.
\end{split}
\end{equation}
\end{proposition}
Hence, $\omega_{nh}$ has zero boundary condition, we get the $L^\infty$ type energy estimate. Main difficulty of this section is how to estimate $\omega_{h}$ which has nonzero boundary condition. In this case, we get only $L^4_t$ type estimate.

\subsection{Homogeneous estimate}
In this section, we estimate $\omega^{\alpha}_{v,h},\omega^{\alpha}_{B,h}$ those have nonzero boundary condition as we see in (\ref{10.10}) and (\ref{10.11}). 
We define two maps $Y_i$ ($i=1,2$),
\begin{equation} \label{10.28}
\begin{split}
& Y_1 : \Omega \rightarrow \mathbb{R}^{3},\,\,\,\partial_t Y_1(t,x) = (u-H)(t,Y_1(t,x)),\,\,\,Y_1(0,x) = x,  \\
& Y_2 : \Omega \rightarrow \mathbb{R}^{3},\,\,\,\partial_t Y_2(t,x) = (u+H)(t,Y_2(t,x)),\,\,\,Y_2(0,x) = x.
\end{split}
\end{equation}
Let us compare image of these two maps. Images $Y_1(t,\O)$ and $Y_2(t,\O)$ are defined only by boundary values of vector fields, $(u\pm H)^b$. If we write boundary graphs as $h_1$ and $h_2$, then
\begin{equation} \label{same image}
h_{1,2}(t) := h_{1,2}(0) + \int_{0}^{t} (u\pm H)^b \cdot \mathbf{N} = h(0) + \int_{0}^{t} u^b \cdot \mathbf{N} = h(t),\quad H^b = 0,\quad\text{on} \quad \p\O.
\end{equation}
Now, we have equations for $\omega^\alpha_{v,h},\omega^\alpha_{B,h}$,
\begin{equation} \label{10.29}
\partial_t^\varphi (\omega^\alpha_{v,h}+\omega^\alpha_{B,h}) + (v-B)\cdot\nabla^\varphi (\omega^\alpha_{v,h}+\omega^\alpha_{B,h}) - \varepsilon\triangle^\varphi(\omega^\alpha_{v,h}+\omega^\alpha_{B,h}) =  0,
\end{equation}
with initial and boundary condition,
$$
(\omega^\alpha_{v,h}+\omega^\alpha_{B,h})^b = (Z^\alpha\omega_v)^b + (Z^\alpha\omega_B)^b,\,\,\,(\omega^\alpha_{v,h}+\omega^\alpha_{B,h})_{t=0} = 0.
$$
And,
\begin{equation} \label{10.30}
\partial_t^\varphi (\omega^\alpha_{v,h}-\omega^\alpha_{B,h}) + (v+B)\cdot\nabla^\varphi (\omega^\alpha_{v,h}-\omega^\alpha_{B,h}) - \varepsilon\triangle^\varphi(\omega^\alpha_{v,h}-\omega^\alpha_{B,h}) =  0,
\end{equation}
with initial and boundary condition,
$$
(\omega^\alpha_{v,h}-\omega^\alpha_{B,h})^b = (Z^\alpha\omega_v)^b - (Z^\alpha\omega_B)^b,\quad (\omega^\alpha_{v,h}-\omega^\alpha_{B,h})_{t=0} = 0.
$$
(\ref{10.29}) can be transformed into
\begin{equation} \label{10.31}
\begin{split}
&\partial_t \left((\omega^\alpha_{v,h}+\omega^\alpha_{B,h})\circ\Phi^{-1}\right) + ((v-B)\circ\Phi^{-1})\cdot\nabla \left((\omega^\alpha_{v,h}+\omega^\alpha_{B,h})\circ\Phi^{-1}\right) - \varepsilon\triangle\left((\omega^\alpha_{v,h}+\omega^\alpha_{B,h})\circ\Phi^{-1}\right) =  0, 	\\
&\partial_t \left((\omega^\alpha_{v,h}+\omega^\alpha_{B,h})\circ\Phi^{-1}\right) + (u-H)\cdot\nabla \left((\omega^\alpha_{v,h}+\omega^\alpha_{B,h})\circ\Phi^{-1}\right) - \varepsilon\triangle\left((\omega^\alpha_{v,h}+\omega^\alpha_{B,h})\circ\Phi^{-1}\right) =  0.
\end{split}
\end{equation}
$u+H$ case is also similar, so we get,
\begin{equation} \label{lar equation}
\begin{split}
& \partial_t \left((\omega^\alpha_{v,h}+\omega^\alpha_{B,h})\circ\Phi^{-1}\circ Y_{1}\right) - \varepsilon\frac{1}{a_0}\partial_i \left(a_{ij}\partial_j \left((\omega^\alpha_{v,h}+\omega^\alpha_{B,h})\circ\Phi^{-1}\circ Y_{1}\right)\right) = 0,\\
& \partial_t \left((\omega^\alpha_{v,h}-\omega^\alpha_{B,h})\circ\Phi^{-1}\circ Y_{2}\right) - \varepsilon\frac{1}{b_0}\partial_i \left(b_{ij}\partial_j \left((\omega^\alpha_{v,h}-\omega^\alpha_{B,h})\circ\Phi^{-1}\circ Y_{2}\right)\right) = 0,\\
\end{split}
\end{equation}
where $Y_{1,2}$ are defined in (\ref{10.28}) and,
\[
	a_{0}:= |J^{1}_{0}|^{\frac{1}{2}} := |\det\nabla Y_{1}(0,x)|^{\frac{1}{2}},\quad b_{0}:= |J^{2}_{0}|^{\frac{1}{2}} := |\det\nabla Y_{2}(0,x)|^{\frac{1}{2}}.
\] 
Matrix $a_{ij}$ and $b_{ij}$ are defined by 
\begin{equation} \label{matrix ab}
\begin{split}  
a_{ij} =  |J^{1}_{0}|^{\frac{1}{2}} P^{-1}, &\quad P_{ij}=  \partial_{i} Y_{1} \cdot \partial_{j} Y_{1}  \\
b_{ij} =  |J^{2}_{0}|^{\frac{1}{2}} P^{-1}, &\quad P_{ij}=  \partial_{i} Y_{2} \cdot \partial_{j} Y_{2}.  \\
\end{split}
\end{equation}
Note that in above equations (\ref{10.31}), we use two different transformation $\Phi\circ Y_{1}$ and $\Phi\circ Y_{2}$. However, in these two transforms, we can use same $\Phi$, since image of $Y_{1}$ and $Y_{2}$ are identical by (\ref{same image}).  \\

\noindent Now we define, 
$$
\left((\omega^\alpha_{v,h} \pm \omega^\alpha_{B,h})\circ\Phi^{-1}\circ X\right) := \mathcal{W}_{\pm},
$$
and these solve,
\begin{equation} 
\begin{split} 
&\partial_t \mathcal{W}_+ - \varepsilon\frac{1}{a_0}\partial_i \left(a_{ij}\partial_j \mathcal{W}_+ \right) = 0,  \\
&\partial_t \mathcal{W}_- - \varepsilon\frac{1}{b_0}\partial_i \left(b_{ij}\partial_j \mathcal{W}_- \right) = 0.  \\
\end{split}
\end{equation}
Multiplying decaying factor, we define,
\begin{equation} \label{def decaying}
\begin{split} 
\Omega^\alpha_+ := e^{-\gamma t}\mathcal{W}_+,\,\,\,\Omega^\alpha_- := e^{-\gamma t}\mathcal{W}_-,
\end{split}
\end{equation}
We get
\begin{equation} \label{10.32}
\begin{split}
& a_0(\partial_t\Omega^\alpha_+ + \gamma\Omega^\alpha_+) - \varepsilon\partial_i \left(a_{ij}\partial_j \Omega^\alpha_+ \right) = 0,  \quad \Omega^\alpha_+|_{z=h_0} = e^{-\gamma t} \left((\omega^\alpha_{v,h}+\omega^\alpha_{B,h})\circ\Phi^{-1}\circ X\right)(t,y,h_0(y)),
\end{split}
\end{equation}
and similarly,
\begin{equation} \label{10.33}
\begin{split}
& b_0(\partial_t\Omega^\alpha_- + \gamma\Omega^\alpha_-) - \varepsilon\partial_i \left(b_{ij}\partial_j \Omega^\alpha_- \right) = 0,  \quad \Omega^\alpha_-|_{z=h_0} = e^{-\gamma t} \left((\omega^\alpha_{v,h}-\omega^\alpha_{B,h})\circ\Phi^{-1}\circ X\right)(t,y,h_0(y)).
\end{split}
\end{equation}
We will use Theorem 10.6 in \cite{NMFR1}. To do this we should first show that $Y_1,Y_2$ satisfy similar property as Lemma 10.5 in \cite{NMFR1}.

\begin{lemma} \label{lemma 10.4}
Let us assume that for $T\in [0,T^\varepsilon],\,\,T^\varepsilon \leq 1$, there exist $M>0$ such that the following holds.
\begin{equation} \label{micro cond}
\sup_{[0,T]}\Lambda_{6,\infty}(t) + \int_0^T \left(\varepsilon\|\nabla\mathcal{V}^6\|^2 + \varepsilon\|\nabla\mathcal{B}^6\|^2 + \varepsilon\|\nabla S^v_n\|_4^2 + \varepsilon\|\nabla S^B_n\|_4^2 \right) \leq M.
\end{equation}
Under above assumption, for $t\in [0,T]$ and $i=1,2$ we have the following estimates.
\begin{equation} \label{10.34}
\begin{split}
\left|J_i(t,x)\right|_{W^{1,\infty}} +  \left|1/J_i(t,x)\right|_{W^{1,\infty}} &\leq \Lambda_0,  \\
\left|\nabla Y_i(t)\right|_{L^\infty} + \left|\partial_t\nabla Y_i(t)\right|_{L^\infty} &\leq \Lambda_0 e^{\Lambda(M)t},  \\
\left|\nabla Y_i(t)\right|_{W^{1,\infty}} + \left|\partial_t\nabla Y_i(t)\right|_{W^{1,\infty}} &\leq \Lambda(M) e^{\Lambda(M)t},  \\
\sqrt{\varepsilon}\left\|\nabla^2Y_i\right\|_{L^\infty} + \sqrt{\varepsilon}\left\|\partial_t\nabla^2Y_i\right\|_{L^\infty} &\leq \Lambda(M)(1+t)^2 e^{\Lambda(M)t} ,
\end{split}
\end{equation}
where 
\[
	J_{i}(t,x) := |\det \nabla Y_{i}(t,x)|,\quad i=1,2.
\]
\end{lemma}
\begin{proof}
1) First one comes from the fact that $J_i(t,x) = J_i(0,x)$, since $u,H$ are both divergence free.\\
2) Secondly,
$$
\partial_t DY_i = D(v \mp B)D\Phi^{-1}DY_i.
$$
Taking $L^\infty$ and using Gronwall's inequality,
$$
\left|\nabla Y_i(t,\cdot)\right|_{L^\infty} \leq \Lambda_0 e^{\Lambda(M)t}.
$$
Again, taking $L^\infty$ to above chain rule and using the result, we get
$$
\left|\partial_t\nabla Y_i(t,\cdot)\right|_{L^\infty} \leq \Lambda_0 e^{\Lambda(M)t}.
$$
These two inequalities give second result.\\
3) We take conormal derivative $Z$ to above chain rule,
$$
\partial_t Z(DY_i) = ZD(v \mp B)D\Phi^{-1}DY_i + D(v \mp B)ZD\Phi^{-1}DY_i + D(v \mp B)D\Phi^{-1}ZDY_i.
$$
Now we use above results 1) and 2), and Gronwall's inequality to get
$$
\left|\nabla Y_i(t,\cdot)\right|_{W^{1,\infty}} \leq \Lambda(M) e^{\Lambda(M)t}.
$$
Again using conormally differentiated chain rule, 
$$
\left|\partial_t\nabla Y_i(t,\cdot)\right|_{W^{1,\infty}} \leq \Lambda(M) e^{\Lambda(M)t}.
$$
4) We take $\sqrt{\varepsilon}\nabla$ to chain rule, to get
$$
\sqrt{\varepsilon}\partial_t D^2 Y_i = \sqrt{\varepsilon}D^2(v \mp B)D\Phi^{-1}DY_i + \sqrt{\varepsilon}D(v \mp B)D^2\Phi^{-1}DY_i + \sqrt{\varepsilon}D(v \mp B)D\Phi^{-1}D^2Y_i.
$$
We use Proposition \ref{proposition 9.9} and Gronwall's inequality, then
\begin{equation*}
\begin{split}
\sqrt{\varepsilon}\left\|\nabla^2Y_i\right\|_{L^\infty} &\leq \Lambda(M)(1+t)^2 e^{\Lambda(M)t} + \Lambda(M)\int_0^t\sqrt{\varepsilon}\left\|\nabla^2Y_i(s)\right\|_{L^\infty}, \\
\sqrt{\varepsilon}\left\|\nabla^2Y_i\right\|_{L^\infty} &\leq \Lambda(M)(1+t)^2 e^{\Lambda(M)t}.
\end{split}
\end{equation*}
Using this result, chain rule, and Gronwall's inequality, estimate for $\p_t\nabla Y_{i}\vert_{i=1,2}$ are also same.
$$
\sqrt{\varepsilon}\left\|\partial_t\nabla^2Y_i\right\|_{L^\infty} \leq \Lambda(M)(1+t)^2 e^{\Lambda(M)t}.
$$
\end{proof}

We are ready to apply Theorem 10.6 in \cite{NMFR1}.
\begin{theorem} \label{theorem 10.5}
There exist $\gamma_0$ which depends on $M$ defined by (\ref{micro cond}) such that for $\gamma \geq \gamma_0$, solution of (\ref{10.32}) and (\ref{10.33}) satisfy the following estimate, respectively.
\begin{equation} \label{10.38}
\left\|\Omega^{m-1}_{\pm}\right\|_{H^{\frac{1}{4}}(0,T;L^2)}^2 \leq \Lambda(M)\sqrt{\varepsilon}\int_0^T \left|(\Omega^{m-1}_{\pm})^b\right|_{L^2}^2.
\end{equation}
\end{theorem}
\begin{proof}
From Lemma \ref{lemma 10.4}, we can apply Theorem 10.6 in \cite{NMFR1} directly for both (\ref{10.32}) and (\ref{10.33}).
\end{proof}

\noindent Now we state estimate for homogeneous part.
\begin{proposition}  \label{proposition 10.6}
Under the same assumption in \ref{lemma 10.4}, we get, for $|\alpha| \leq m-1$,
\begin{equation} \label{10.39}
\left\|\omega^{\alpha}_{v,h}\right\|^2_{L^4(0,T;L^2)},\left\|\omega^{\alpha}_{B,h}\right\|^2_{L^4(0,T;L^2)} \leq \Lambda(M)\int_0^T \left( \|\mathcal{V}^m\|^2 + \|\mathcal{B}^m\|^2 + |h|_m^2 \right) + \frac{\varepsilon}{2}\int_0^T \left( \left\|\nabla\mathcal{V}^m\right\|^2 + \left\|\nabla\mathcal{B}^m\right\|^2 \right).
\end{equation}
\end{proposition}
\begin{proof}
We define,
\[
	\|f\|^2_{H^{\frac{1}{4}}_{tan}(0,T;L^2)} = \inf\{ \|Pf\|_{H^{\frac{1}{4}}(\mathbb{R},L^2(S))}, Pf = f\quad\text{on}\quad[0,T]\times S \}.
\]
By sobolev embedding, 
$$
\left\|\Omega^{\alpha}_{\pm}\right\|_{L^4(0,T;L^2)}^2 \leq C \left\|\Omega^{m-1}_{\pm}\right\|_{H^{\frac{1}{4}}(0,T;L^2)}^2 \leq \Lambda(M)\sqrt{\varepsilon}\int_0^T e^{-2\gamma t}\left((\omega^\alpha_{v,h}\pm\omega^\alpha_{B,h})\circ\Phi^{-1}\circ X\right)^2.
$$
For general sobolev embedding, $C$ may blow up as $T\rightarrow 0$. But here, $C$ is independent to $T$. This is because we use sobolev embedding on $\mathbb{R}_+$. 
We can change variable to function on $S$ using $\Phi$ and estimate $J_i$ of Lemma \ref{lemma 10.4}, to get
\begin{equation}
\begin{split}
\left\| \omega_{v,h}^\alpha \pm \omega_{B,h}^\alpha \right\|_{L^4(0,T;L^2)}^2 &\leq \Lambda(M)\sqrt{\varepsilon}\int_0^T \left\| (Z^\alpha\omega_v)^b \pm (Z^\alpha\omega_B)^b \right\|^2_{m-1}  \\
&\leq \Lambda(M)\sqrt{\varepsilon}\int_0^T \left( \left\| (Z^\alpha\omega_v)^b \right\|^2 + \left\| (Z^\alpha\omega_B)^b \right\|^2 \right).
\end{split}
\end{equation}
We apply (\ref{10.7}) for boundary integrals to finish the proof.
\end{proof}

\section{Uniform estimate}

We prove uniform energy estimate for $\varepsilon=\lambda$. For sufficiently large $m \geq 6$, let initial data satisfy
\begin{equation} \label{11.1}
\begin{split}
	\mathcal{I}_m(0) &:= \|v_0\|_{E^m} +  \|B_0\|_{E^m} + |h_0|_m + \sqrt{\varepsilon}|h_0|_{m+\frac{1}{2}} + \|v_0\|_{E^{2,\infty}} + \|B_0\|_{E^{2,\infty}} 	\\
	&\quad + \sqrt{\varepsilon}\|\partial_{zz}v_0\|_{L^\infty} + \sqrt{\varepsilon}\|\partial_{zz}B_0\|_{L^\infty} < \infty.	 	
\end{split}
\end{equation}

\noindent By the assumption of Theorem \ref{theorem 1.4}, for smoothed data $(v_{0}^{\varepsilon,\delta}, B_{0}^{\varepsilon,\delta}, h_{0}^{\varepsilon,\delta})$, we have local existence time $T^{\varepsilon,\delta}$. And $\forall T \leq T^{\varepsilon,\delta}$, we have additional regularity by parabolic estimate,
\begin{equation} \label{11.2}
\begin{split}
\mathcal{N}_m(T) &:= \sup_{[0,T]} \Big( \|v(t)\|_m^2 + \|B(t)\|_m^2 + |h(t)|_m^2 + \varepsilon |h(t)|_{m+\frac{1}{2}}^2 \\
&\quad + \|v(t)\|_{E^{2,\infty}}^2 + \|B(t)\|_{E^{2,\infty}}^2 + \varepsilon\|\partial_{zz}v(t)\|_{L^\infty}^2 + \varepsilon\|\partial_{zz}B(t)\|_{L^\infty}^2  \Big)  \\
&\quad + \left( \|\partial_z v\|_{L^4([0,T];H_{co}^{m-1})}^2 + \|\partial_z B\|_{L^4([0,T];H_{co}^{m-1})}^2 \right) \\
&\quad + \varepsilon\int_0^T \|\nabla v\|_m^2 + \varepsilon\int_0^T \|\nabla B\|_m^2 + \varepsilon\int_0^T \|\nabla\partial_z v\|_{m-2}^2 + \varepsilon\int_0^T \|\nabla\partial_z B\|_{m-2}^2 < \infty,
\end{split}
\end{equation}
where $(v,B,h)$ is solution for regularized initial data $(v_{0}^{\varepsilon,\delta}, B_{0}^{\varepsilon,\delta}, h_{0}^{\varepsilon,\delta})$. Without loss of generality, we can assume that (\ref{7.14}) holds by choosing sufficiently small $T^{\varepsilon,\delta} \ll 1$.
From section 7, section 8, and section 9, we know that it is equivalent to control 
\begin{equation} \label{11.3}
\begin{split}
\mathcal{E}_m(T) &:= \sup_{[0,T]} \Big( \|\mathcal{V}^m(t)\|^2 + \|\mathcal{B}^m(t)\|^2 + |h(t)|_m^2 + \varepsilon |h(t)|_{m+\frac{1}{2}}^2 \\
&\quad + \|S^v_n(t)\|_{E^{2,\infty}}^2 + \|S^B_n(t)\|_{E^{2,\infty}}^2 + \varepsilon\|\partial_z S^v_n(t)\|_{L^\infty}^2 + \varepsilon\|\partial_z S^B_n(t)\|_{L^\infty}^2  \Big)  \\
&\quad + \left( \|\omega_v\|_{L^4([0,T];H_{co}^{m-1})}^2 + \|\omega_B\|_{L^4([0,T];H_{co}^{m-1})}^2 \right) \\
&\quad + \varepsilon\int_0^T \|\nabla \mathcal{V}^m\|^2 + \varepsilon\int_0^T \|\nabla \mathcal{B}^m\|^2 + \varepsilon\int_0^T \|\nabla S^v_n\|_{m-2}^2 + \varepsilon\int_0^T \|\nabla S^B_n\|_{m-2}^2,
\end{split}
\end{equation}
instead of $\mathcal{N}_m(T)$. We know that
\begin{equation}
\begin{split}
	\mathcal{N}_{m}(T) \leq \Lambda(\frac{1}{c_0}, \mathcal{E}_{m}(T)),\quad 
	\mathcal{E}_{m}(T) \leq \Lambda(\frac{1}{c_0}, \mathcal{N}_{m}(T)),	
\end{split}
\end{equation}
where second inequality holds by product estimates. When two parameter $R$ and $c_0$ satisfy $\frac{1}{c_0} << R$, we define $T_*^{\varepsilon,\delta}$ as,
\begin{equation} \label{11.4}
T_*^{\varepsilon,\delta} = \sup\left\{ T\in[0,1] \bigg| \ \mathcal{E}_m(t)\leq R,\,\,|h(t)|_{2,\infty} \leq \frac{1}{c_0},\,\,\partial_z\varphi(t)\geq c_0,\,\,g-(\partial_z^\varphi q^E)|_{z=0} \geq \frac{c_0}{2},\,\,\forall t\in [0,T] \right\}.
\end{equation}

Now we combine our propositions and corollaries to get uniform energy estimate for both $\varepsilon$ and $\delta$. From Corollary \ref{corollary 9.1}, for $T\leq T^{\varepsilon,\delta}_{*}$,
$$
\Lambda_{6,\infty}(T) \leq \Lambda(R).
$$
And from Corollary \ref{9.58},
\[
	\int_{0}^{T} \sqrt{\varepsilon} \|\nabla^{2}v\|_{1,\infty} \leq \Lambda(R).
\]
Using (\ref{7.15}), (\ref{7.16}), Proposition (\ref{9.6}), Lemma (\ref{lemma 9.8}), Proposition \ref{proposition 10.3}, and Proposition \ref{proposition 10.6} to obtain
$$
\mathcal{E}_m(T) \leq \Lambda(\frac{1}{c_0},\mathcal{I}_m(0)) + \Lambda(R)\sqrt{T} + \Lambda(R)\int_0^T |(\partial_z\partial_t q^E)^b|_{L^\infty}.
$$
Using Proposition \ref{proposition 5.6} again,
\begin{equation} \label{11.5}
{E}_m(T) \leq \Lambda(\frac{1}{c_0},\mathcal{I}_m(0)) + \Lambda(R)\sqrt{T},
\end{equation}
which is independent to both $\varepsilon$ and $\delta$. Moreover, for $\forall t \leq T\leq T^{\varepsilon,\delta}_{*}$,
\begin{equation} \label{11.6}
\begin{split}
	|h(t)|_{2,\infty} &\leq |h(0)|_{2,\infty} + \Lambda(R)T, 	\\
	\partial_z\varphi(t) &\geq 1 - \int_0^t \|\partial_t\nabla\eta\|_{L^\infty} \geq 1 - \Lambda(R)T, 	\\
	g - (\partial_z^\varphi q^E)^b(t) &\geq g-(\partial_z^\varphi q^E)|_{z=0} - \Lambda(R)\int_0^t ( 1 + |(\partial_t\partial_z q^E)^b|_{L^\infty} ) 		\\
	&\geq g-(\partial_z^\varphi q^E)|_{z=0} - \Lambda(R)\sqrt{t}.
\end{split}	
\end{equation}


Note that above four inequalities in (\ref{11.5}) and (\ref{11.6}) are all $(\varepsilon, \delta)$-independent inequalities. Hence, we can choose proper $R = \Lambda(\mathcal{I}_m(0),|h|_{2,\infty})$ so that we can pick $T_*$ which is $(\varepsilon,\delta)$-independent and for all $T \leq \text{Min}(T_*,T_*^{\varepsilon,\delta})$, four inequalities in (\ref{11.4}) are satisfied. Now we can send regularizing parameter $\delta\rightarrow 0$, to get uniform time interval $T_*$ for initial data $\mathcal{I}_m(0)$. This is possible because for each $\delta$, $\mathcal{N}_m(T_*)$ is uniformly bounded in $\delta$ so we can use strong compactness.

\section{Uniqueness}
\subsection{Uniqueness of Theorem \ref{theorem 1.4}}
In the previous section, we proved existence of viscous system (\ref{1.4}). In this subsection we prove uniqueness of the system. We suppose two solutions $(v^1,B^1,\varphi^1,q^1),(v^2,B^2,\varphi^2,q^2)$ with same initial condition. We will perform $L^2$ energy estimate (zero order estimate) and use Gronwall's inequality to show that $L^2$ norm of differences are locally zero. We write
$$
\bar{v}^\varepsilon := v^\varepsilon_1 - v^\varepsilon_2,\,\,\,\bar{B}^\varepsilon := B^\varepsilon_1 - B^\varepsilon_2,\,\,\,\bar{h}^\varepsilon := h^\varepsilon_1 - h^\varepsilon_2,\,\,\,\bar{q}^\varepsilon := q^\varepsilon_1 - q^\varepsilon_2
$$
with initial condition $\bar{v}(0) = \bar{B}(0) = \bar{h}(0)=0$. Let both solutions satisfy on $[0,T^\varepsilon]$,
$$
\mathcal{N}^i_m(T^\varepsilon) \leq R,\,\,\,i=1,2.
$$
By divergence free condition, $\nabla^{\varphi_i}\cdot v_i^\varepsilon = 0,\,\,i=1,2$,
$$
\left(\partial_t + v_{y,i}^\varepsilon\cdot\nabla_y + V_{z,i}^\varepsilon\partial_z\right)v_{i}^\varepsilon + \nabla^{\varphi_i}q_i^\varepsilon - \varepsilon\triangle^{\varphi_i}v_i^\varepsilon = (B^\varepsilon_i\cdot\nabla^{\varphi_i})B_i.
$$
Then we have equation about $(\bar{v}^\varepsilon, \bar{B}^\varepsilon, \bar{h}^\varepsilon, \bar{q}^\varepsilon)$. First for Navier-Stokes,
\begin{equation} \label{12.1}
\left(\partial_t + v_{y,1}^\varepsilon\cdot\nabla_y + V_{z,1}^\varepsilon\partial_z\right)\bar{v}^\varepsilon + \nabla^{\varphi_1}\bar{q}^\varepsilon - \varepsilon\triangle^{\varphi_1}\bar{v}^\varepsilon - (B_1^\varepsilon \cdot P_1^* \nabla) \bar{B}^\varepsilon = F,
\end{equation}
where
\begin{equation*}
\begin{split}
F &:= (v_{y,2}^\varepsilon - v_{y,1}^\varepsilon)\cdot \nabla_y v_2^\varepsilon + (V_{z,2}^\varepsilon - V_{z,1}^\varepsilon)\partial_z v_2^\varepsilon - \left(\frac{1}{\partial_z \varphi_2^\varepsilon} - \frac{1}{\partial_z \varphi_1^\varepsilon}\right)\left(P_1^*\nabla q_2^\varepsilon\right) 	\\
&\quad + \frac{1}{\partial_z \varphi_2^\varepsilon}\left((P_2 - P_1)^*\nabla q_2^\varepsilon\right) + \varepsilon\left(\frac{1}{\partial_z \varphi_2^\varepsilon} - \frac{1}{\partial_z\varphi_1^\varepsilon}\right)\nabla\cdot(E_1\nabla v_2^\varepsilon) + \varepsilon\frac{1}{\partial_z\varphi_2^\varepsilon}\nabla\cdot\left((E_2-E_1)\nabla v_2^\varepsilon\right) 	\\
&\quad + (\bar{B}^\varepsilon\cdot P_2^*\nabla)B_2^\varepsilon + (B_1^\varepsilon\cdot(P_1^* - P_2^*)\nabla)B_2^\varepsilon.
\end{split}
\end{equation*}

For Faraday's law, similarly as above, we have
\begin{equation} \label{12.2}
\left(\partial_t + v_{y,1}^\varepsilon\cdot\nabla_y + V_{z,1}^\varepsilon\partial_z\right)\bar{B}^\varepsilon - \varepsilon\triangle^{\varphi_1}\bar{B}^\varepsilon - (B_1^\varepsilon \cdot P_1^* \nabla) \bar{v}^\varepsilon = E,
\end{equation}
where
\begin{equation*}
\begin{split}
E &:= (v_{y,2}^\varepsilon - v_{y,1}^\varepsilon)\cdot \nabla_y B_2^\varepsilon + (V_{z,2}^\varepsilon - V_{z,1}^\varepsilon)\partial_z B_2^\varepsilon + \varepsilon\left(\frac{1}{\partial_z \varphi_2^\varepsilon} - \frac{1}{\partial_z\varphi_1^\varepsilon}\right)\nabla\cdot(E_1\nabla B_2^\varepsilon) 	\\
&\quad + \varepsilon\frac{1}{\partial_z\varphi_2^\varepsilon}\nabla\cdot\left((E_2-E_1)\nabla B_2^\varepsilon\right) + (\bar{B}^\varepsilon\cdot P_2^*\nabla)v_2^\varepsilon + (B_1^\varepsilon\cdot(P_1^* - P_2^*)\nabla)v_2^\varepsilon.
\end{split}
\end{equation*}
For divergence-free condition,
\begin{equation} \label{12.3}
\nabla^{\varphi_1}\cdot \bar{v}^\varepsilon = - \left(\frac{1}{\partial_z \varphi_2^\varepsilon} - \frac{1}{\partial_z\varphi_1^\varepsilon}\right)\nabla\cdot\left(P_1 v_2^\varepsilon\right) - \frac{1}{\partial_z \varphi_2^\varepsilon} \nabla \cdot \left((P_2 - P_1) v_2^\varepsilon\right),\,\,\,\text{same\,\,for\,\,}B.
\end{equation}
For Kinematic boundary condition,
\begin{equation} \label{12.4}
\partial_t \bar{h}^\varepsilon - (v^\varepsilon)^b_{y,1}\cdot\nabla_{y} h + \left((v_{z,1}^\varepsilon)^b - (v_{z,2}^\varepsilon)^b \right) = -\left((v_{y,2}^\varepsilon)^b - (v_{y,1}^\varepsilon)^b \right)\cdot \nabla_{y} h_2^\varepsilon.
\end{equation}
Continuity of stress tensor condition becomes,
\begin{equation} \label{12.5}
\bar{q}^\varepsilon\mathbf{n}_1 - 2\varepsilon\left(S^{\varepsilon_1}\bar{v}^\varepsilon\right)\mathbf{n}_1 - g\bar{h}^\varepsilon = 2\varepsilon\left( \left( S^{\varphi_1} - S^{\varphi_2} \right) v_2^\varepsilon\right)\mathbf{n}_1 + 2\varepsilon\left( S^{\varphi_2} v_2^\varepsilon\right)\left( \mathbf{n}_1 - \mathbf{n}_2 \right).
\end{equation}
Using (\ref{12.1})-(\ref{12.5}), we get $L^2$- energy estimate. Details are nearly same as high order estimate which was shown throughout this paper. (Since initial condition is zero, no initial term appear on right hand side)
\begin{equation} \label{12.6}
\begin{split}
&\left\|\bar{v}^\varepsilon(t)\right\|_{L^2}^2 + \left\|\bar{B}^\varepsilon(t)\right\|_{L^2}^2 + \left|\bar{h}^\varepsilon(t)\right|_{L^2}^2 + \varepsilon\int_0^t \left( \left\|\nabla\bar{v}^\varepsilon\right\|_{L^2}^2 + \left\|\nabla\bar{B}^\varepsilon\right\|_{L^2}^2 \right) ds 	\\
&\leq \Lambda(R)\int_0^t \left( \left\|\bar{v}^\varepsilon(s)\right\|_{L^2}^2 + \left\|\bar{B}^\varepsilon(s)\right\|_{L^2}^2 + \left|\bar{h}^\varepsilon(s)\right|_{H^{\frac{1}{2}}}^2 + \|\nabla \bar{q}^\varepsilon\|_{L^(S)}\|\bar{v}^\varepsilon\|_{L^(S)} \right)ds.
\end{split}
\end{equation}

Using pressure estimate, we also get,
$$
\|\nabla \bar{q}^\varepsilon\|_{L^(S)} \leq \Lambda(R)\left( |\bar{h}^\varepsilon|_{H^{1/2}} + \|\bar{v}^\varepsilon\|_{H^1(S)} + \|\bar{B}^\varepsilon\|_{H^1(S)} \right).
$$
Now with help of Proposition \ref{proposition 3.9}, we get the result.
\begin{equation} \label{12.7}
\begin{split}
&\left\|\bar{v}^\varepsilon(t)\right\|_{L^2}^2 + \left\|\bar{B}^\varepsilon(t)\right\|_{L^2}^2 + \left|\bar{h}^\varepsilon(t)\right|_{L^2}^2 + \sqrt{\varepsilon}\left|\bar{h}^\varepsilon(t)\right|_{H^{1/2}}^2 + \varepsilon\int_0^t \left( \left\|\nabla\bar{v}^\varepsilon\right\|_{L^2}^2 + \left\|\nabla\bar{B}^\varepsilon\right\|_{L^2}^2 \right) ds 	\\
&\leq \Lambda(R)\int_0^t \left( \left\|\bar{v}^\varepsilon(s)\right\|_{L^2}^2 + \left\|\bar{B}^\varepsilon(s)\right\|_{L^2}^2 + \left|\bar{h}^\varepsilon(s)\right|_{H^{\frac{1}{2}}}^2 \right)ds.
\end{split}
\end{equation}
In above equations for $(\bar{v}^\varepsilon,\bar{B}^\varepsilon,\bar{h}^\varepsilon, \bar{q}^\varepsilon)$, right hand side does not have low order than $L^2$ energy. However we have already uniform bound of high-order energy, so we can collect high order terms into $\Lambda(R)$. 

\subsection{Uniqueness of Theorem \ref{theorem 1.5}.}
For two solutions $(v_1, B_1, h_1, q_1),(v_2, ,B_2, h_2, q_2)$ with same initial condition. Suppose,
\begin{equation} \label{12.8}
\sup_{[0,T]}\left( \|v_i\|_m + \left\|\partial_z v_i\right\|_{m-1} + \left\|\partial_z v_i\right\|_{1,\infty} + \|B_i\|_m + \left\|\partial_z B_i\right\|_{m-1} + \left\|\partial_z B_i\right\|_{1,\infty} + |h_i|_m \right) \leq R,\,\,\,\,i=1,2.
\end{equation}
Define $\bar{v} := v_1 - v_2,\,\,\,\bar{B} := B_1 - B_2,\,\,\,\bar{h} := h_1 - h_2,\,\,\,\bar{q} := q_1 - q_2$ and we write equation of $(\bar{v}, \bar{h}, \bar{q})$, as before. Euler equation becomes,
\begin{equation} \label{12.9}
\left(\partial_t + v_{y,1}\cdot\nabla_y + V_{z,1}\partial_z\right)\bar{v} + \nabla^{\varphi_1}\bar{q} - (B_1\cdot\nabla^{\varphi_1})\bar{B} = F',
\end{equation}
where
\begin{equation*}
\begin{split}
F' &:= (v_{y,2} - v_{y,1})\cdot \nabla_y v_2 + (V_{z,2} - V_{z,1})\partial_z v_2 - \left(\frac{1}{\partial_z \varphi_2} - \frac{1}{\partial_z \varphi_1}\right)\left(P_1^*\nabla q_2\right) + \frac{1}{\partial_z \varphi_2}\left((P_2 - P_1)^*\nabla q_2\right) 	\\
&\quad + (\bar{B}\cdot P_2^*\nabla)B_2 + (B_1\cdot(P_1^* - P_2^*)\nabla)B_2.
\end{split}
\end{equation*}

For Faraday's law, similarly,
\begin{equation} \label{12.10}
\left(\partial_t + v_{y,1}\cdot\nabla_y + V_{z,1}\partial_z\right)\bar{v} - (B_1\cdot\nabla^{\varphi_1})\bar{v} = E',
\end{equation}
where
$$
E' := (v_{y,2} - v_{y,1})\cdot \nabla_y v_2 + (V_{z,2} - V_{z,1})\partial_z v_2 + (\bar{B}\cdot P_2^*\nabla)v_2 + (B_1\cdot(P_1^* - P_2^*)\nabla)v_2.
$$
For divergence-free condition, (same for $B$)
\begin{equation} \label{12.11}
\nabla^{\varphi_1}\cdot \bar{v} = - \left(\frac{1}{\partial_z \varphi_2} - \frac{1}{\partial_z\varphi_1}\right)\nabla\cdot\left(P_1 v_2\right) - \frac{1}{\partial_z \varphi_2} \nabla \cdot \left((P_2 - P_1) v_2\right).
\end{equation}
For Kinematic boundary condition,
\begin{equation} \label{12.12}
\partial_t \bar{h} - v^b_{y,1}\cdot\nabla_{y} h + \left( v_{z,1}^b - v_{z,2}^b \right) = -\left( v_{y,2}^b - v_{y,1}^b \right)\cdot \nabla_{y} h_2.
\end{equation}
Continuity of stress tensor condition becomes simply,
\begin{equation} \label{12.13}
\bar{q}\mathbf{n}_1 = g\bar{h}.
\end{equation}
By performing basic $L^2$-estimate and pressure estimate,
\begin{equation} \label{12.14}
\left\|\bar{v}(t)\right\|_{L^2}^2 + \left\|\bar{B}(t)\right\|_{L^2}^2 + \left|\bar{h}(t)\right|_{L^2}^2 \leq \Lambda(R)\int_0^t \left( \left\|\bar{v}(s)\right\|_{H^1}^2 + \left\|\bar{B}(s)\right\|_{H^1}^2 + \left|\bar{h}(s)\right|_{H^{\frac{1}{2}}}^2 \right)ds.
\end{equation}
We should control $\|v\|_1$ on right hand side. But, since there are no dissipation on the left hand side, we cannot make it absorbed. Instead, we use vorticity. Let's define vorticity $\omega_v = \nabla^{\varphi} \times v$, $\omega_B = \nabla^{\varphi} \times B$ (which is equivalent to $\omega_v = (\nabla\times u)(t,\Phi)$ and $\omega_B = (\nabla\times H)(t,\Phi)$ ). We have (same for $B$),
$$
\omega_v\times\mathbf{n} = \frac{1}{2}\left( D^\varphi v\mathbf{n} - (D^\varphi v)^T \mathbf{n} \right) = \mathbf{S}^{\varphi} v\mathbf{n} - (D^\varphi v)^T \mathbf{n} = \frac{1}{2}\partial_{\mathbf{n}}u - g^{ij}\left(\partial_j v\cdot \mathbf{n}\right)\partial_{y^i}.
$$
Hence, it suffice to estimate $\omega_v$ instead of $\partial_z v$, 
\begin{equation} \label{12.15}
\left\|\partial_z v\right\|_{L^2} + \left\|\partial_z B\right\|_{L^2} \leq \Lambda(R)\left( \|\omega_v\|_{L^2} + \|\omega_B\|_{L^2} + \|v\|_1 + |h|_{\frac{1}{2}} \right).
\end{equation}
To estimate $\omega$ (both for $v,B$), we use vorticity equation
\begin{equation} \label{12.16}
\left( \partial_t^{\varphi_i} + v_i\cdot\nabla^{\varphi_i} \right) \omega_i = \left( \omega_i\cdot\nabla^{\varphi_i} \right) v_i .
\end{equation}
$L^2$ energy estimates of $\bar{\omega_v}$ and $\bar{\omega_B} $ are 
\begin{equation} \label{12.17}
\left\|\bar{\omega_v}(t)\right\|_{L^2}^2 + \left\|\bar{\omega_B}(t)\right\|_{L^2}^2 \leq \Lambda(R)\int_0^t \left( |\bar{h}(s)|_1^2 + \|\bar{v}(s)\|_{H^1(S)}^2 + \|\bar{\omega_v}(s)\|_{L^2}^2 + \|\bar{B}(s)\|_{H^1(S)}^2 + \|\bar{\omega_B}(s)\|_{L^2}^2\right) ds.
\end{equation}
So we finish the proof.

\section{Inviscid limit}
Proof for this section is exactly same as \cite{NMFR1}. For interior, it is clear that we can just add $B$-related terms those have same regularity as $v$. For boundary condition, since we have $B=0$ on the boundary, definition of weak solution makes sense also. At result, we get sequence (up to subsequence) $(v^\varepsilon(t),B^\varepsilon(t),h^\varepsilon(t))$, which converges to $(v,B,h)$ in weak $L^2(S)\times L^2(S) \times L^2(\mathbb{R}^2)$. By Proposition \ref{L2 est}, we have $L^2$ energy conservation, so this it converges in strong $L^2(S)\times L^2(S) \times L^2(\mathbb{R}^2)$. $L^\infty$ follows from $L^2$ convergence, uniform bounds of energy, and anisotropic embedding Proposition \ref{proposition 3.2}.

\section{Appendix}

In appendix, we show well-posedness of the system (\ref{1.2}). Full detail is given in \cite{DHL}, so we explain scheme of the proof, instead of full detail.

\subsection{Lagrangian coordinate and main statement}
 We use standard Lagrangian map $X(t,\cdot):\xi \rightarrow x $ and define  
\begin{equation} \label{lagrangian}
\begin{split}
v(t,\xi) &:= u(t,X(t,\xi)) = u(t,x),  \\
B(t,\xi) &:= H(t,X(t,\xi)) = H(t,x),  \\
q(t,\xi) &:= p(t,X(t,\xi)) = p(t,x).
\end{split}
\end{equation}

\noindent We should rewrite the system (\ref{1.2}) in terms of $(v,B,q)$ and $(t,\xi)$. Let $\Pi^x_{\xi}$ be corresponding transform from $\Omega(t), \ S_F(t)$ to $\Omega, \ S_F$. Then (\ref{1.2}) and (\ref{1.3}) can be rewritten in $(t,\xi)$ in $\Omega$.
\begin{equation} \label{app Larg system}
\begin{cases}
v_t - \nu\triangle_v v + \nabla_v q = B\cdot\nabla_v B\,\,\,\,\text{in}\,\,\,\,Q_T, \\
B_t - \lambda\triangle_v B = B\cdot\nabla_v v,\,\,\,\,\text{in}\,\,\,\,Q_T, \\
\nabla_v\cdot v = 0,\,\,\,\,\text{in}\,\,\,\,Q_T, \\
\nabla_v\cdot B = 0,\,\,\,\,\text{in}\,\,\,\,Q_T, \\
q - 2\nu\mathbf{D}_v(v)\mathbf{n}^{(v)}\cdot\mathbf{n}^{(v)} = gh\,\,\,\,\text{on}\,\,\,\,S_{F,T}, \\
B = 0,\quad \{S_{F} \cup \{\mathbb{R}^{3} \backslash \O \}\} \times [0,T),  	\\
v(0)=u_0,\,\,\,\,B(0)=H_0,\,\,\,\,\Omega\times\{t=0\},
\end{cases}
\end{equation}
with compatibility conditions
\begin{equation} \label{app Larg compatibility}
\begin{cases}
\nabla_v \cdot v_0 = 0, \quad\text{in} \ \Omega, \\
\nabla_v \cdot B_0 = 0, \quad\text{in} \ \Omega, \\
B_0 = 0,\quad\text{on} \quad \{S_{F} \cup \{\mathbb{R}^{3} \backslash \O \}\},  	\\
\mathbf{\Pi}^{(v)}\mathbf{D}_v(v_0)\mathbf{n}^{(v)} = 0,\quad\text{on} \ S_F, \\
\end{cases}
\end{equation}
where $Q_T := \Omega \times [0,T)$, $ S_{F,T} := S_F\times [0,T) $, and
\begin{equation} \label{app symbols def}
\begin{split}
f^{(v)} &= f^{(v)}(t,\xi) = \Pi_{\xi}^x f(t,x),  \\
\nabla_v &:= (\nabla_{v,1},\nabla_{v,2},\nabla_{v,3}) = \mathcal{G}\nabla = \mathcal{G}(\nabla_1,\nabla_2,\nabla_3),   \\
\mathcal{G} &:= \mathcal{G}^{(v)} = \left(\frac{\partial X_v}{\partial\xi}\right)^{-t} = \left( I + \int_0^t (Dv) d\tau \right)^{-t},  \\
\mathbf{D}_v(f) &:= \frac{1}{2}((\nabla_v f) + (\nabla_v f)^t).  \\
\end{split}
\end{equation}

We will solve the system (\ref{app Larg system}) with (\ref{app Larg compatibility}) in Sobolev-Slobodetskii, fractional sobolev space.

\begin{definition}
	By $W_2^l(\Omega)$, we define,	
	\begin{equation} \label{app slobodetskii}
	\|u\|^2_{W_2^l(\Omega)} := \sum_{0 \leq |\alpha| < l}\|D^\alpha u\|^2_{L^2(\Omega)} + \|u\|^2_{\dot{W}_2^l(\Omega)},
	\end{equation}
	where
	$$
	\|u\|^2_{\dot{W}_2^l(\Omega)} := 
	\begin{cases} \sum_{|\alpha|=l}\|D^\alpha u\|^2_{L^2(\Omega)}\,\,\,\,\text{if}\,\,l\in \mathbb{Z}, \\
	\sum_{|\alpha|=[l]} \int_{\Omega}\int_{\Omega} \frac{|D^\alpha u(x) - D^\alpha u(y)|^2}{|x-y|^{n+2\{l\}}}\,\,\,\,\text{if}\,\,l=[l] + \{l\}\notin \mathbb{Z},0<{l}<1 .
	\end{cases}
	$$
\end{definition}

\begin{definition}
	We also define the anisotropic space $W_2^{l,l/2}(Q_T)$ as 
	\begin{equation} \label{app aniso}
	\begin{split}
	\|u\|^2_{W_2^{l,l/2}(Q_T)} & :=  \|u\|^2_{W_2^{l,0}(Q_T)} + \|u\|^2_{W_2^{0,l/2}(Q_T)} \\
	& = \int_0^T \|u(t,\cdot)\|^2_{W_2^l(\Omega)} dt + \int_\Omega \|u(\cdot,x)\|^2_{W_2^{l/2}(0,T)} dx.
	\end{split}
	\end{equation}
\end{definition}

\begin{definition}
	By $H^{l,l/2}_\gamma(Q_T),\gamma \geq 0$, we define,
	\begin{equation} \label{app H space}
	\|u\|^2_{H^{l,l/2}_\gamma(Q_T)} := \|u\|^2_{H^{l,0}_\gamma(Q_T)} + \|u\|^2_{H^{0,l/2}_\gamma(Q_T)},
	\end{equation}
	where
	\begin{equation}
	\begin{split}
	\|u\|^2_{H^{l,0}_\gamma(Q_T)} &:= \int_0^T e^{-2\gamma t} \|u(t,\cdot)\|^2_{\dot{W}_2^l(\Omega)}dt,  \\
	\|u\|^2_{H^{0,l/2}_\gamma(Q_T)} &:= \gamma^l\int_0^T e^{-2\gamma t} \|u(t,\cdot)\|^2_{L^2(\Omega)}dt  \\
	&+ \int_0^T e^{-2\gamma t} dt \int_0^\infty \left\| \left( \frac{\partial}{\partial t}\right)^k u_0(t,\cdot) - \left( \frac{\partial}{\partial t}\right)^k u_0(t-\tau,\cdot) \right\|^2_{L^2(\Omega)}\frac{d\tau}{\tau^{1+l-2k}},
	\end{split}
	\end{equation}
	if $l/2$ is not an integer, $k=[l/2], \ u_0(t,x) = u(t,x)(t>0), \ u_0(t,x) = 0(t<0)$. If $l/2$ is an integer, then the double integral in the norm should be replaced by 
	$$
	\int_{-\infty}^T e^{-2\gamma t} \left\| \left( \frac{\partial}{\partial t} \right)^{l/2}u(t,\cdot) \right\|^2_{L^2(\Omega)} dt,
	$$
	and $\left( \frac{\partial}{\partial t}\right)^j u |_{t=0}(j=0,1,2,\cdots,l/2-1)$ should be satisfied.
\end{definition}

\begin{definition}
	We also introduce the space $H_\gamma^{\ell+\frac{1}{2},\frac{1}{2},\frac{\ell}{2}}(S_{F,T})$ to treat trace on the boundary $S_{F,T}$.
	\begin{equation} \label{app trace space}
	\begin{split}
	\|u\|^2_{H^{\ell+\frac{1}{2},\frac{1}{2},\frac{\ell}{2}}_\gamma(S_{F,T})} &:= \int_0^T e^{-2\gamma t} \left( \|u\|^2_{\dot{W}_2^{\ell+\frac{1}{2}}(S_{F})} + \gamma^{\ell} \|u\|^2_{\dot{W}_2^{1/2}(S_{F})} \right) dt   \\
	& + \int_0^T e^{-2\gamma t} dt \int_0^\infty \left\| \left( \frac{\partial}{\partial t}\right)^k u_0(t,\cdot) - \left( \frac{\partial}{\partial t}\right)^k u_0(t-\tau,\cdot) \right\|^2_{W_2^{1/2}(S_{F})}\frac{d\tau}{\tau^{1+\ell-2k}},
	\end{split}
	\end{equation}
	if $\ell/2$ is note an integer, $k=[\ell/2], \ u_0(t,x) = u(t,x)(t>0), \ u_0(t,x) = 0(t<0)$. If $\ell/2$ is an integer, then the double integral in the norm should be replaced by 
	$$
	\int_{-\infty}^T e^{-2\gamma t} \left\| \left( \frac{\partial}{\partial t} \right)^{\ell/2}u(t,\cdot) \right\|^2_{W_2^{l/2}(S_{F})} dt,
	$$
	and $\left( \frac{\partial}{\partial t}\right)^j u |_{t=0}(j=0,1,2,\cdots,\frac{\ell}{2}-1)$ should be satisfied. 
\end{definition}

\begin{definition}
	We define,
	\begin{equation} \label{app mixed space}
	(\|u\|_{Q_T}^{(l+2)})^{2} := (\|u\|_{Q_T}^{(l)})^{2} + \sum_{|s|=2} (\|D_x^s u\|_{Q_T}^{(l)})^{2} + \sum_{|s|=0}^1\|D_x^s u\|_{L^2(Q_T)}^{2} ,
	\end{equation}
	where
	$$
	(\|u\|_{Q_T}^{(l)})^{2} := \|u\|^2_{W_2^{l,l/2}(Q_T)} + T^{-l}\|u\|^2_{L^2(Q_T)}.
	$$
	And, since we deal nonlinear terms on right hand side of (\ref{app Larg system}), we need $L^\infty _T$ type norm. We define
	\begin{equation} \label{app infty space}
	\|u\|^2_{H^{l+2,l/2+1}(Q_T)} := \|u\|^{(l+2)2}_{Q_T} + \sup_{t<T}\|u(t)\|^2_{W_2^{l+1}(\Omega)}.
	\end{equation}
\end{definition}

We state a Lemma for above function spaces from Proposition 1 in \cite{PS2}.
\begin{lemma} \label{app interpolation lemma}
	For smooth $u(x)$ and $v(x)$, (defined in a domain $\Omega \subset \mathbb{R}^n$), they satisfy the following inequalities,
	\begin{equation}
	\begin{split}
	\|uv\|_{W_2^l(\Omega)} &\leq c \|u\|_{W_2^l(\Omega)} \|v\|_{W_2^s(\Omega)},\quad s>n/2, \quad l<n/2,  \\
	\|uv\|_{L^2(\Omega)} &\leq c \|u\|_{W_2^l(\Omega)} \|v\|_{W_2^{n/2-l}(\Omega)},\quad l<n/2, \\
	\|uv\|_{W_2^l(\Omega)} &\leq c \left( \|u\|_{W_2^l(\Omega)} \|v\|_{W_2^s(\Omega)} + \|v\|_{W_2^l(\Omega)} \|u\|_{W_2^s(\Omega)}\right),\quad l,s>n/2.
	\end{split}
	\end{equation}
\end{lemma}

Using above functional spaces, we claim the following well-posedness result.
\begin{theorem} \label{app main theorem}
	Let $l \in (\frac{1}{2},1)$, and initial conditions $h_0 \in W_2^{l+5/2}(S_F), \ u_0=v_0 \in W_2^{l+1}(\Omega)$, and $\ H_0=B_0 \in W_2^{l+1}(\Omega)$, with compatibility conditions
	\begin{equation*}
	\begin{cases}
	\nabla\cdot v_0 = \nabla\cdot B_0 = 0,\quad\text{in}\quad\Omega, \\
	\mathbf{D}(v_0)\mathbf{n}_0 - (\mathbf{D}(v_0)\mathbf{n}_0\cdot \mathbf{n}_0)\mathbf{n}_0 = 0,\,\,\,\,\text{on}\,\,\,\,S_{F}, \\
	H_0 = 0,\quad\text{on}\quad S_F. \\
	\end{cases}
	\end{equation*}
	Then there exist a unique solution $(v,B,q)$ to the system (\ref{app Larg system}) such that
	$$
	\|v\|_{H^{l+2,l/2+1}(Q_{T^*})} + \|B\|_{H^{l+2,l/2+1}(Q_{T^*})} + \|\nabla q\|_{W_2^{l,l/2}(Q_{T^*})} + \|q\|_{W_2^{l+1/2,l/2+1/4}(S_{F,{T^*}})}
	$$
	$$
	\leq C_0 \left( \|v_0\|_{W_2^{l+1}(\Omega)} + \|B_0\|_{W_2^{l+1}(\Omega)} + \|h_0\|_{W_2^{l+3/2}(\mathbb{R}^2)} \right),
	$$
	for some $T^* > 0$. Moreover for any $T>0$, we can choose sufficiently small $\varepsilon(T)$ such that if 
	\[
	\|v_0\|_{W_2^{l+1}(\Omega)} + \|B_0\|_{W_2^{l+1}(\Omega)} + \|h_0\|_{W_2^{l+3/2}(\mathbb{R}^2)}  \leq \varepsilon(T),
	\]
	then we have a unique solution for (\ref{app Larg system}) and (\ref{app Larg compatibility}) on the time interval $[0,T]$.
\end{theorem}

\subsection{Linear problem} We study linearized problem for two main PDE of (\ref{app Larg system}).
\subsubsection{Stokes problem}
From (\ref{app Larg system}), we start from linear Stokes problem. We use result by A.Tani \cite{AT} for the following Stokes system.
\begin{equation} \label{app stokes system}
\begin{cases}
v_t - \nu\triangle v + \nabla q = f, \quad \text{in}\quad Q_T, \\
\nabla\cdot v = \rho,\quad \text{in}\quad Q_T, \\
v(0) = v_0,\quad \text{in}\quad \Omega\times\{t=0\}, \\
2\nu[\mathbf{D}(v)\mathbf{n} - (\mathbf{D}(v)\mathbf{n}\cdot \mathbf{n})\mathbf{n}] = d,\quad \text{on}\quad S_{F,T}, \\
-q + 2\nu\mathbf{D}(v)\mathbf{n}\cdot \mathbf{n} = b,\quad \text{on}\quad S_{F,T}. \\
\end{cases}
\end{equation}

From Theorem 2.1 in \cite{AT}, we state
\begin{proposition} \label{app stokes with zero initial}
	(Linear Stokes problem with zero initial data)
	Let $l>\frac{1}{2}$, $0<T\leq\infty$, and $\gamma$ be sufficiently large so that $\gamma\geq\gamma_0\geq 1$. And we assume $S_F \in W_2^{l+3/2}$. Then, for
	$$
	(f,\rho,d,b)\in H_{\gamma}^{l,l/2}(Q_T) \times H_{\gamma}^{l+1,(l+1)/2}(Q_T) \times H_{\gamma}^{l+1/2,l/2+1/4}(S_{F,T}) \times H_{\gamma}^{l+1/2,1/2,l/2}(S_{F,T}),
	$$
	satisfying compatibility conditions $d\cdot n_0 = 0$, $(\rho,d)|_{t=0} = 0$, $\rho=\nabla\cdot R$, and $R\in H_{\gamma}^{l+1,1,l/2}(Q_T)$, there exist a unique solution $(v,q)\in H_{\gamma}^{l+2,l/2+1}(Q_T) \times H_{\gamma}^{l+1,1,l/2}(Q_T)$ to the problem (\ref{app stokes system}) with zero initial condition $v_0 = 0$. Moreover we have the following estimate,
	\begin{equation*}
	\begin{split}
	\|v\|_{H_{\gamma}^{l+2,l/2+1}(Q_T)} + \|q\|_{H_{\gamma}^{l+1,1,l/2}(Q_T)}  &\leq  C \{ \|f\|_{H_{\gamma}^{l,l/2}(Q_T)} + \|\rho\|_{H_{\gamma}^{l+1,(l+1)/2}(Q_T)} \\
	&+ \|R\|_{H_{\gamma}^{0,l/2+1}(Q_T)} + \|d\|_{H_{\gamma}^{l+1/2,l/2+1/4}(S_{F,T})} + \|b\|_{H_{\gamma}^{l+1/2,1/2,l/2}(S_{F,T})} \},
	\end{split}
	\end{equation*}
	where $\|q\|^2_{H_{\gamma}^{l+1,1,l/2}(Q_T)} := \|q\|^2_{H_{\gamma}^{l,l/2}(Q_T)} + \|\nabla q\|^2_{H_{\gamma}^{l,l/2}(Q_T)}$.
\end{proposition}

For general data, $v_0 \neq 0$, we first produce a function $U_0(t,x)$ which has same initial data data $U_0(0) = v_0$ and its space-time sobolev norm is bounded by $v_0$. We gain the following Lemma.

\begin{lemma} \label{app extention lemma}
	1) For $v\in W_2^{l+2,l/2+1}(Q_T)$ where $1/2<l<1$, there is extension $U \in W_2^{r,r/2}(Q_{\infty})$ such that
	$$
	\|U\|^{(l+2)}_{Q_\infty} \leq C \|v_0\|_{W_2^{l+1}(\Omega)},
	$$
	where $\|U\|^{(l+2)}$ is defined in definition (\ref{app mixed space}).  \\
	2) Let $w(0)=0$ and $w \in W_2^{l+2,l/2+1}(Q_T)$, then there exist a extension $w_{\text{ext}} \in W_2^{l+2,l/2+1}(Q_\infty)$ such that
	$$
	\|w_{\text{ext}}\|_{W_2^{l+2,l/2+1}(Q_\infty)} \leq C \|w\|_{W_2^{l+2,l/2+1}(Q_T)}.
	$$
\end{lemma}

Using above Lemma we gain linear result for general data.

\begin{proposition} \label{app prop general stokes}
	(Linear stokes problem with general initial data)
	Let $l \in (\frac{1}{2},1)$, $0<T<\infty$ and $S_F \in W_2^{l+3/2}$. And $(f,\rho,u_0,(b,d))$ in (\ref{app stokes system}) satisfy
	\begin{equation*}
	\begin{split}
	(f,\rho,u_0,(b,d)) &\in W_2^{l,l/2}(Q_T) \times W_2^{l+1,(l+1)/2}(Q_T) \times W_2^{l+1}(\Omega) \times W_2^{l+1/2,l/2+1/4}(S_{F,T}),  \\
	\rho &= \nabla\cdot R,\quad\text{where} \quad R \in L^2(Q_T)\quad \text{and} \quad R_t \in W_2^{0,l/2}(Q_T), \\
	\nabla\cdot u_0 &= \rho|_{t=0},\\
	d|_{t=0} &= 2\nu[\mathbf{D}(v)\mathbf{n}_0 - (\mathbf{D}(v)\mathbf{n}_0\cdot \mathbf{n}_0)\mathbf{n}_0]|_{S_F},  \\
	d\cdot \mathbf{n}_0 &= 0.
	\end{split}
	\end{equation*}
	Then system (\ref{app stokes system}) with general initial data $v_0$ has a solution 
	\[
	v \in W_2^{l+2,l/2+1}(Q_T) \cap C_T W_2^{l+1}(\Omega), \quad q \in W_2^{l,l/2}(Q_T), \quad \nabla q \in W_2^{l,l/2}(Q_T),\quad q \in W_2^{l+1/2,l/2+1/4}(S_{F,T}),
	\]
	with the estimate
	\begin{equation}
	\begin{split}
	& \|v\|_{H^{l+2,l/2+1}(Q_T)} + \|q\|^{(l)}_{Q_T} + \|\nabla q\|^{(l)}_{Q_T} + \|q\|_{W_2^{l+1/2,l/2+1/4}(S_{F,T})}  \\
	&\leq C(T) \{ \|f\|^{(l)}_{Q_T} + \|\rho\|_{W_2^{l+1,(l+1)/2}(S_{F,T})} + \|R\|_{W_2^{0,l/2+1}(Q_T)} + T^{-l/2}\|R_t\|_{L^2(Q_T)}  \\
	&+ \|(b,d)\|_{W_2^{l+1/2,l/2+1/4}(S_{F,T})} + T^{-l/2}\|b\|_{W_2^{l/2,0}(S_{F,T})} + \|u_0\|_{W_2^{l+1}(\Omega)} \},
	\end{split}
	\end{equation}
	where $C(T)$ is time dependent constant on $T$ non-decreasingly.
\end{proposition}

\begin{proof}
	Let us write $w = v - U$, where $U$ satisfies 
	\begin{equation} \label{app U ineq}
	\int_{0}^\infty \left( \|\partial_t U(\cdot,t)\|^2_{L^2(\Omega)} + \sum_{|s|=2} \|D^s U(\cdot,t)\|^2_{L^2(\Omega)}\right)\frac{1}{|t|^l} dt + \|U\|^2_{W_2^{l+2,l/2+1}(Q_\infty)} \leq C \|v_0\|^2_{W_2^{l+1}(\Omega)},
	\end{equation} 
	$U(0)=v_0$, and $U \in W^{l+2,\frac{l}{2}+1}(Q_{\infty})$. Then $(w,q)$ solves linear stokes problem with zero initial data,
	\begin{equation}\label{app w eq}
	\begin{cases}
	\frac{\partial w}{\partial t} - \nu\triangle w + \nabla q = f - \frac{\partial U}{\partial t} + \nu\triangle U := f^{\prime},\quad\text{in}\quad Q_T, \\
	\nabla\cdot w = \rho - \nabla\cdot U := \rho^{\prime},\quad \text{in}\quad Q_T, \\
	2\nu[\mathbf{D}(w)\mathbf{n} - (\mathbf{D}(w)\mathbf{n}\cdot\mathbf{n})\mathbf{n}] = d - 2\nu[\mathbf{D}(U)\mathbf{n} - (\mathbf{D}(U)\mathbf{n}\cdot\mathbf{n})\mathbf{n}] := d^{\prime},\quad \text{on}\quad S_{F,T}, \\
	q - 2\nu\mathbf{D}(w)\mathbf{n}\cdot\mathbf{n} := - b^{\prime},\quad \text{on}\quad S_{F,T} \\
	w(0) = 0,\quad\Omega\times\{t=0\}. \\
	\end{cases}
	\end{equation}
	For source terms $(f',\rho',d',b')$ in the right hand sides, we can estimate the followings.
	\begin{equation} \label{app f' rho'}
	\begin{split}
	\|f'\|^{(l)}_{Q_T} &\leq C( \|f\|^{(l)}_{Q_T} + \|u_0\|_{W_2^{l+1}(\Omega)} ),  \\
	\|\rho'\|_{W_2^{l+1,0}(Q_T)} &\leq C ( \|\rho\|_{W_2^{l+1,0}(Q_T)} + \|u_0\|_{W_2^{l+1}(\Omega)} ),	\\
	\|d'\|_{H_0^{l+1/2,l/2+1/4}(S_{F,T})} &\leq C \left( \|d\|_{W_2^{l+1/2,l/2+1/4}(S_{F,T})} + \|u_0\|_{W_2^{l+1}(\Omega)} \right),	\\
	\|b'\|_{H_0^{l+1/2,1/2,l/2}(S_{F,T})} &\leq C \left( \|b\|_{W_2^{l+1/2,l/2+1/4}(S_{F,T})} + T^{-l/2}\|b\|_{W_2^{1/2,0}(S_{F,T})} + \|u_0\|_{W_2^{l+1}(\Omega)} \right).
	\end{split}
	\end{equation}

	\noindent We apply the result of Proposition \ref{app stokes with zero initial} with (\ref{app f' rho'}) to get, 
	\begin{equation} \label{app w}
	\begin{split}
	&\|w\|_{H_0^{l+2,l/2+1}(Q_T)} + \|q\|_{H_0^{l+1,1,l/2}(Q_T)} \leq e^{\gamma T} \left( \|w\|_{H_{\gamma}^{l+2,l/2+1}(Q_T)} + \|\nabla q\|_{H_{\gamma}^{l,l/2}(Q_T)} \right)  \\
	&\quad \leq C e^{\gamma T} \{ \|f'\|_{H_{\gamma}^{l,l/2}(Q_T)} + \|\rho'\|_{H_{\gamma}^{l+1,(l+1)/2}(Q_T)} + \|R''\|_{H_{\gamma}^{0,l/2+1}(Q_T)} + \|d'\|_{H_{\gamma}^{l+1/2,l/2+1/4}(S_{F,T})} + \|b'\|_{H_{\gamma}^{l+1/2,1/2,l/2}(S_{F,T})} \}\\
	&\quad \leq C(T) \{ \|f\|^{(l)}_{Q_T} + \|\rho\|_{W_2^{l+1,(l+1)/2}(Q_T)} + T^{-l/2}\|R_t\|_{L^2(Q_T)} + \|R\|_{W_2^{0,l/2+1}(Q_T)}  \\
	&\quad \quad + \|(d,b)\|_{W_2^{l+1/2,l/2+1/4}(S_{F,T})} + T^{-l/2}\|b\|_{W_0^{1/2,0}(S_{F,T})} + \|v_0\|_{W_2^{l+2}(\Omega)} \},  \\
	\end{split}
	\end{equation}
	where $\gamma$ is a fixed constant, found in Proposition \ref{app stokes with zero initial} and $C(T)$ depends on time $T$, non-decreasingly on $T$, which means it does not blow up as $T\rightarrow 0$. Meanwhile, from boundary condition,
	\begin{equation} \label{app q}
	\|q\|_{W_2^{0,l/2+1/4}(S_{F,T})} \leq \|2\nu\mathbf{D}(v)n_0\cdot n_0\|_{W_2^{0,l/2+1/4}(S_{F,T})} + \|b\|_{W_2^{0,l/2+1/4}(S_{F,T})}.
	\end{equation}
	Lastly, we should estimate $L^\infty_T$ type estimate of $w,U$. Since $w$ has zero initial data, we use Lemma \ref{app extention lemma},
	\begin{equation} \label{app wC}
	\begin{split}
	\|w\|_{C_T W_2^{l+1}(\Omega)} &\leq \|w_{\text{ext}}\|_{C_\infty W_2^{l+1}(\Omega)} \leq C \left( \|w_{\text{ext}}\|_{L^2_\infty W_2^{l+2}(\Omega)} + \|\partial_t w_{\text{ext}}\|_{L^2_\infty W_2^{l}(\Omega)} \right)  \\
	&\leq C \|w_{\text{ext}}\|_{W_2^{l+1,l/2+1}(Q_\infty)} \leq C \|w\|_{W_2^{l+1,l/2+1}(Q_T)},
	\end{split}
	\end{equation}
	where $w_{\text{ext}} \in W_2^{l+2,l/2+1}$ is an extension of $w$ which satisfies (we can choose time-independent $C$)
	$$
	\|w_{\text{ext}}\|_{W_2^{l+1,l/2+1}(Q_\infty)} \leq C \|w\|_{W_2^{l+1,l/2+1}(Q_T)}.
	$$ 
	Similarly,
	\begin{equation} \label{app U}
	\begin{split}
	\|U\|_{C_T W_2^{l+1}(\Omega)} &\leq \|U_{\text{ext}}\|_{C_\infty W_2^{l+1}(\Omega)} \leq C \left( \|U_{\text{ext}}\|_{L^2_\infty W_2^{l+2}(\Omega)} + \|\partial_t U_{\text{ext}}\|_{L^2_\infty W_2^{l}(\Omega)} \right) \\
	&\leq C \|U_{\text{ext}}\|_{W_2^{l+1,l/2+1}(Q_\infty)} \leq C \|v_0\|_{W_2^{l+1}(\O)},
	\end{split}
	\end{equation}
	where $U_{\text{ext}} \in W_2^{l+2,l/2+1}$ is an extension of $U$ which satisfies
	\begin{equation*} \label{app U ext}
	\|U_{\text{ext}}\|_{W_2^{l+1,l/2+1}(Q_\infty)} \leq C \|v_0\|_{W_2^{l+1}(\Omega)},
	\end{equation*}
	from Lemma \ref{app extention lemma}. From (\ref{app infty space}),
	$$
	\|v\|^{(l+2)2}_{Q_T} + \sup_{t<T} \|v\|^2_{W_2^{l+1}(\Omega)} \leq C\left\{ \|U\|^{(l+2)2}_{Q_T} + \sup_{t<T}\|U\|^2_{W_2^{l+1}(\Omega)} + \|w\|^{(l+2)2}_{Q_T} + \sup_{t<T}\|w\|^2_{W_2^{l+1}(\Omega)} \right\},
	$$
	so putting Lemma \ref{app extention lemma}, (\ref{app w}), (\ref{app q}), (\ref{app wC}), and (\ref{app U}) together, we finish the proof.
\end{proof}

\subsubsection{Heat equation}
For linearized equatino for Faraday's law, we solve,
\begin{equation} \label{app heat system}
\begin{cases}
B_t - \lambda\triangle B = g,\quad\text{in}\quad Q_T, \\
B = 0,\quad\text{on}\quad \{S_F\cup \{\mathbb{R}^{3}\backslash\O \}\}\times [0,T), \\
B(0) = B_0 = H_0,\quad \text{in}\quad \Omega\times\{t=0\}. \\
\end{cases}
\end{equation}

\begin{proposition} \label{app prop general heat}
	For the system (\ref{app heat system}), we have a unique solution with the following estimates.
	\begin{equation*}
	\|B\|_{H^{l+2,l/2+1}(Q_T)} \leq C(T) \left( \|g\|^{(l)}_{Q_T} + \|B_0\|_{W_2^{l+1}(\Omega)} \right),
	\end{equation*}
	where $C(T)$ depends on time $T$ non-decreasingly.
\end{proposition}

\subsection{constant coefficient nonlinear problem}
We solve constant coefficient nonlinear problem of (\ref{app stokes system}) and (\ref{app heat system}).
\begin{equation} \label{app cst coeff nonlin}
\begin{cases}
v_t - \nu\triangle v + \nabla q = (B\cdot\nabla)B + f ,\quad\text{in}\quad Q_T, \\
\nabla\cdot v = \rho,\quad \text{in}\quad Q_T, \\
v(0) = v_0,\quad\text{in}\quad\Omega\times\{t=0\}, \\
2\nu[\mathbf{D}(v)\mathbf{n} - (\mathbf{D}(v)\mathbf{n}\cdot \mathbf{n})\mathbf{n}] = d,\quad \text{on}\quad S_{F,T}, \\
-q + 2\nu\mathbf{D}(v) \mathbf{n}\cdot \mathbf{n} = b \,\,\,\,\text{on}\,\,\,\,S_{F,T}, \\
B_t - \lambda\triangle B = (B\cdot\nabla)v + g,\quad\text{in}\quad Q_T ,\\
B = 0,\quad\text{on}\quad \{ S_{F} \cup \{\mathbb{R}^{3} \backslash \O \}\} \times [0,T), \\
B(0) = B_0 ,\quad \text{in}\quad \Omega\times\{t=0\}. \\
\end{cases}
\end{equation}

Main problem is to control nonlinear forcing terms $(B\cdot\nabla)B$ and $(B\cdot\nabla)v$ on the right hand sides. We first prove the following Lemma.
\begin{lemma} \label{app nonlinear est}
	When $l>\frac{1}{2}$, we have the following nonlinear estimate.
	$$
	\|F\nabla G\|^{(l)}_{Q_T} \leq C ( T + T^{(1-l)/2} + T^{1/2} ) \|F\|_{H^{l+2,l/2+1}(Q_T)} \|G\|_{H^{l+2,l/2+1}(Q_T)}.
	$$
\end{lemma}
\begin{proof}
	\begin{equation*}
	\begin{split}
	\|F\nabla G\|^{(l)2}_{Q_T} &= \|F\nabla G\|^2_{W_2^{l,l/2}(Q_T)} + T^{-l}\|F\nabla G\|^2_{L^2(Q_T)}  \\
	&= \|F\nabla G\|^2_{W_2^{l,0}(Q_T)} + \|F\nabla G\|^2_{W_2^{0,l/2}(Q_T)} + T^{-l}\|F\nabla G\|^2_{L^2(Q_T)}.
	\end{split}
	\end{equation*}
	Using Lemma (\ref{app interpolation lemma}),
	\begin{equation} \label{F delta G}
	\begin{split}
	\|F \nabla G\|^{(l)}_{Q_T} &\leq \|F \nabla G\|_{W_2^{l,0}(Q_T)} + \|F \nabla G\|_{W_2^{0,l/2}(Q_T)} + T^{-l/2}\|F \nabla G\|_{L^2(Q_T)}  \\
	&\leq C \left\{ T \|F\|_{C_T W_2^{l+1}(\Omega)}\|G\|_{C_T W_2^{l+1}(\Omega)} + \underbrace{ T^{-l/2}\|F\nabla G\|_{L^2(Q_T)} }_{(II)} \right\}  \\
	&\quad + \left( \underbrace{ \int_0^T \int_0^T \frac{\left\| (F \nabla G)(t)-(F \nabla G)(s) \right\|^2_{L^2(\Omega)}}{\left| t-s \right|^{1+l}} dtds }_{(I)} \right)^{1/2} .
	\end{split}
	\end{equation}
	We focus on the last term $(I)$. Let us write $t-s = h$. Domain can be divided symmetrically into two region $t>s,s>t$. By changing order of integral and Lemma (\ref{app interpolation lemma}),
	\begin{equation} \label{(I)}
	\begin{split}
	(I)&= \int_0^T \int_0^T \frac{\left\| (F \nabla G)(t)-(F \nabla G)(s) \right\|^2_{L^2(\Omega)}}{\left| t-s \right|^{1+l}} dtds  \\
	&\leq C \int_0^T \frac{dh}{h^{1+l}} \int_h^T \left\|F(s)\nabla (G(s+h)-G(s)) \right\|_{L^2(\Omega)}^2 ds   \\
	&\quad + C \int_0^T \frac{dh}{h^{1+l}} \int_h^T \left\|(F(s+h)-F(s))\nabla G(s)\right\|_{L^2(\Omega)}^2 ds \\
	&\leq C \|F\|^2_{C_T W_2^{l+1}(\Omega)}\int_0^T \frac{dh}{h^{1+l}} \int_h^T \left\|G(s+h) - G(s)) \right\|_{W_2^1(\Omega)}^2 ds  \\
	&\quad + C \|G\|_{C_T W_2^{l+1}(\Omega)}\int_0^T \frac{dh}{h^{1+l}} \int_h^T \left\|F(s+h)-F(s)\right\|_{W_2^1(\Omega)}^2 ds  \\
	&\leq CT \|F\|^2_{C_T W_2^{l+1}(\Omega)}\int_0^T \frac{dh}{h^{1+(l+1)}} \int_h^T \left\|G(s+h) - G(s)) \right\|_{W_2^1(\Omega)}^2 ds  \\
	&\quad+ CT \|G\|_{C_T W_2^{l+1}(\Omega)}\int_0^T \frac{dh}{h^{1+(l+1)}} \int_h^T \left\|F(s+h)-F(s)\right\|_{W_2^1(\Omega)}^2 ds  \\
	&\leq CT \|F\|^2_{C_T W_2^{l+1}(\Omega)} \|G\|^2_{H^{1+2(\frac{l+1}{2}),\frac{1}{2}+\frac{l+1}{2}}(Q_T)} + CT \|G\|^2_{C_T W_2^{l+1}(\Omega)} \|F\|^2_{H^{1+2(\frac{l+1}{2}),\frac{1}{2}+\frac{l+1}{2}}(Q_T)}  \\
	&\leq CT \|F\|^2_{H^{l+2,l/2+1}(Q_T)} \|G\|^2_{H^{l+2,l/2+1}(Q_T)}.
	\end{split}
	\end{equation}
	For $(II)$,
	\begin{equation} \label{(II)}
	\begin{split}
	(II) &= T^{-l/2}\|F\nabla G\|_{L^2(Q_T)} \\
	&\leq CT^{-l/2} \left( \int_0^T \|F\|^2_{W_2^{\frac{3}{2}-l}(\Omega)}\|\nabla G\|^2_{W_2^{l}(\Omega)} dt \right)^{1/2}  \\
	&\leq C T^{(1-l)/2} \|F\|_{C_T W_2^{l+1}(\Omega)} \|G\|_{C_T W_2^{l+1}(\Omega)}.
	\end{split}
	\end{equation}
	Hence, from (\ref{F delta G}), (\ref{(I)}), and (\ref{(II)}), we have the following estimate.
	\begin{equation*}
	\begin{split}
	\|F \nabla G\|^{(l)}_{Q_T} &\leq C \big\{ T \|F\|_{C_T W_2^{l+1}(\Omega)}\|G\|_{C_T W_2^{l+1}(\Omega)} + T^{(1-l)/2} \|F\|_{C_T W_2^{l+1}(\Omega)} \|G\|_{C_T W_2^{l+1}(\Omega)}  \\
	&\quad + T^{1/2} \|F\|_{H^{l+2,l/2+1}(Q_T)} \|G\|_{H^{l+2,l/2+1}(Q_T)} \big\}  \\
	&\leq C ( T + T^{(1-l)/2} + T^{1/2} ) \|F\|_{H^{l+2,l/2+1}(Q_T)} \|G\|_{H^{l+2,l/2+1}(Q_T)}.
	\end{split}
	\end{equation*}
\end{proof}

We use Proposition \ref{app prop general stokes}, \ref{app prop general heat}, and Lemma \ref{app nonlinear est} to solve system (\ref{app cst coeff nonlin}).

\begin{proposition} \label{prop cst coeff nonlin}
	Let $l\in(1/2,1)$, and $S_F \in W_2^{l+3/2}$. Assume that $(f,\rho,v_0,(b,d))$ in $(\ref{app cst coeff nonlin})$ satisfy
	\begin{equation*}
	\begin{split}
	(f,\rho,u_0,(b,d)) &\in W_2^{l,l/2}(Q_T) \times W_2^{l+1,(l+1)/2}(Q_T) \times W_2^{l+1}(\Omega) \times W_2^{l+1/2,l/2+1/4}(S_{F,T}),  \\
	\rho &= \nabla\cdot R, \quad R \in L^2(Q_T), \\
	R_t &\in W_2^{0,l/2}(Q_T).
	\end{split}
	\end{equation*}
	Also assume that the following compatibility conditions hold.
	$$
	\nabla\cdot v_0 = \rho|_{t=0},\quad d|_{t=0} = 2\nu[\mathbf{D}(v)\mathbf{n}_0 - (\mathbf{D}(v)\mathbf{n}_0\cdot \mathbf{n}_0)\mathbf{n}_0]|_{S_F},\quad d\cdot \mathbf{n}_0 = 0.
	$$
	Then system (\ref{app cst coeff nonlin}) has a solution 
	\begin{equation*}
	\begin{split}
	v &\in W_2^{l+2,l/2+1}(Q_T) \ \cap \ C_T W_2^{l+1}(\Omega),\quad q \in W_2^{l,l/2}(Q_T), \\
	\nabla q &\in W_2^{l,l/2}(Q_T), \quad q \in W_2^{l+1/2,l/2+1/4}(S_{F,T}).
	\end{split}
	\end{equation*}
\end{proposition}
\begin{proof}
	We construct iteration scheme as following.
	\begin{equation}
	\begin{cases}
	v^{(m+1)}_t - \nu\triangle v^{(m+1)} + \nabla q^{(m+1)} = (B^{(m)}\cdot\nabla)B^{(m)} + f ,\quad\text{in}\quad Q_T, \\
	\nabla\cdot v^{(m+1)} = \rho,\quad\text{in}\quad Q_T, \\
	v^{(m+1)}(0) = v_0,\quad \text{in}\quad \Omega\times\{t=0\}, \\
	2\nu[\mathbf{D}(v^{(m+1)})\mathbf{n} - (\mathbf{D}(v^{(m+1)})\mathbf{n}\cdot \mathbf{n})\mathbf{n}] = d,\quad\text{on}\quad S_{F,T}, \\
	-q + 2\nu\mathbf{D}(v^{(m+1)})\mathbf{n}\cdot \mathbf{n} = b,\quad\text{on}\quad S_{F,T}, \\
	B_t^{(m+1)} - \lambda\triangle B^{(m+1)} = (B^{(m)}\cdot\nabla)v^{(m)} + g ,\quad\text{in}\quad Q_T, \\
	B^{(m+1)} = 0, \quad\text{on}\quad S_{F,T}\cup\{\mathbb{R}^{3}\backslash Q_{T}\}, \\
	B^{(m+1)}(0) = B_0,\quad\text{in}\quad\Omega\times\{t=0\}, \\
	(v^{(0)},q^{(0)},B^{(0)}) = (0,0,0).
	\end{cases}
	\end{equation}
	
	From Proposition \ref{app prop general stokes} and \ref{app prop general heat} in section 3, we have a unique solution $(v^{(m+1)},q^{(m+1)},B^{(m+1)})$, for given data $(v^{(m)},B^{(m)})$,
	\begin{equation}
	\begin{split}
	v^{(m+1)} &\in W_2^{l+2,l/2+1}(Q_T)\cap C_T W_2^{l+1}(\Omega), \\
	B^{(m+1)} &\in W_2^{l+2,l/2+1}(Q_T)\cap C_T W_2^{l+1}(\Omega), \\
	q^{(m+1)} &\in W_2^{l,l/2}(Q_T),  \\
	\nabla q^{(m+1)} &\in W_2^{l,l/2}(Q_T),  \\
	q &\in W_2^{l,l/2}(S_{F,T}).
	\end{split}
	\end{equation}
	To get uniform bounds, we first define 
	\begin{equation} \label{mathcal A}
	\begin{split}
	\mathcal{A}^{(m+1)} &:= \|v^{(m+1)}\|_{H^{l+2,l/2+1}(Q_T)} + \|B^{(m+1)}\|_{H^{l+2,l/2+1}(Q_T)} + \|q^{(m+1)}\|^{(l)}_{Q_T} \\
	&\quad + \|\nabla q^{(m+1)}\|^{(l)}_{Q_T} + \|q^{(m+1)}\|_{W_2^{l+1/2,l/2+1/4}(S_{F,T})} .
	\end{split}
	\end{equation}
	Then we use estimates of Proposition \ref{app prop general stokes} and \ref{app prop general heat} to get,
	\begin{equation} \label{A est}
	\begin{split}
	\mathcal{A}^{(m+1)} &\leq C(T) \{ \|(B^{(m)}\cdot\nabla)B^{(m)}\|^{(l)}_{Q_T} + \|(B^{(m)}\cdot\nabla)v^{(m)}\|^{(l)}_{Q_T} + \|f\|^{(l)}_{Q_T} + \|g\|^{(l)}_{Q_T}  \\
	&\quad + \|\rho\|_{W_2^{l+1,(l+1)/2}(S_{F,T})} + \|R\|_{W_2^{0,l/2+1}(Q_T)} + T^{-l/2}\|R\|_{L^2(Q_T)}  \\
	&\quad + \|(b,d)\|_{W_2^{l+1/2,l/2+1/4}(S_{F,T})} + T^{-l/2}\|b\|_{W_2^{l/2,0}(S_{F,T})} + \|u_0\|_{W_2^{l+1}(\Omega)} + \|H_0\|_{W_2^{l+1}(\Omega)}\}.
	\end{split}
	\end{equation}
	
	We also define data part as
	\begin{equation} \label{mathcal D}
	\begin{split}
	\mathcal{D}(T) &:= \|f\|^{(l)}_{Q_T} + \|g\|^{(l)}_{Q_T} + \|\rho\|_{W_2^{l+1,(l+1)/2}(S_{F,T})} + \|R\|_{W_2^{0,l/2+1}(Q_T)}  \\
	&\quad + T^{-l/2}\|R\|_{L^2(Q_T)} + \|(b,d)\|_{W_2^{l+1/2,l/2+1/4}(S_{F,T})} \\
	&\quad + T^{-l/2}\|b\|_{W_2^{l/2,0}(S_{F,T})} + \|v_0\|_{W_2^{l+1}(\Omega)} + \|B_0\|_{W_2^{l+1}(\Omega)}.
	\end{split}
	\end{equation}
	Using this definition, (\ref{A est}) can be controlled by,
	\begin{equation} \label{mathcal A est}
	\begin{split}
	\mathcal{A}^{(m+1)} \leq  C(T)\left( \mathcal{D}(T) + ( T + T^{(1-l)/2} + T^{1/2} ) \mathcal{A}^{(m)2} \right).
	\end{split}
	\end{equation}
	Now we suffice to show uniform bound and contraction mapping to apply fixed point argument. \\
	$\mathbf{Uniform\,\,bound} \ $ By cauchy sequence argument, we can pick sufficiently small $T_0 > 0$ such that
	\begin{equation} \label{unif bdd 1}
	\mathcal{A}^{(m)} \leq 2(C(0)+1)\mathcal{D}(T_0) \doteq M_{T_0},\,\,\,\,\forall m \in \mathbb{N}.
	\end{equation}
	$\mathbf{Contraction\,\,mapping} \ $ Let us define difference,
	$$
	v^{(m+1)} - v^{(m)} := \mathcal{V}^{(m+1)},\quad B^{(m+1)} - B^{(m)} := \mathcal{B}^{(m+1)},\quad q^{(m+1)} - q^{(m)} := \mathcal{Q}^{(m+1)}.
	$$
	
	Again, using Proposition \ref{app prop general stokes}, \ref{app prop general heat}, and Lemma \ref{app nonlinear est},
	\begin{equation}
	\begin{split}
	\bar{\mathcal{A}}^{(m+1)} &:= \|\mathcal{V}^{(m+1)}\|_{H^{l+2,l/2+1}(Q_T)} + \|\mathcal{B}^{(m+1)}\|_{H^{l+2,l/2+1}(Q_T)}  \\
	&\quad + \|\mathcal{Q}^{(m+1)}\|^{(l)}_{Q_T} + \|\nabla \mathcal{Q}^{(m+1)}\|^{(l)}_{Q_T} + \|\mathcal{Q}^{(m+1)}\|_{W_2^{l+1/2,l/2+1/4}(S_{F,T})}   \\
	&\leq C(T) \{ \|(B^{(m)}\cdot\nabla)\mathcal{B}^{(m)}\|^{(l)}_{Q_T} + \|(\mathcal{B}^{(m)}\cdot\nabla)B^{(m-1)}\|^{(l)}_{Q_T}  \\
	&\quad + \|(B^{(m)}\cdot\nabla)\mathcal{V}^{(m)}\|^{(l)}_{Q_T} + \|(\mathcal{B}^{(m)}\cdot\nabla)v^{(m-1)}\|^{(l)}_{Q_T} \}  \\
	&\leq C(T)( T + T^{(1-l)/2} + T^{1/2} ) M_{T_0} \bar{\mathcal{A}}^{(m)}.
	\end{split}
	\end{equation}
	We can find sufficiently small $T_1 > 0$, (Without loss of generality, we pick this so that smaller than $T_0$), such that
	$$
	C(t)( t + t^{(1-l)/2} + t^{1/2} ) M_{T_0} < 1,\,\,\,\,\forall t < T_1,
	$$
	since $C(T)$ is time dependent on $T$ non-decreasingly, so that finite near $T=0$. Hence we have an unique solution via fixed point argument.
\end{proof}

\subsection{Proof of theorem \ref{app main theorem}}
In this subsection, we give finish the proof of Theorem \ref{app main theorem}. First, we solve (\ref{app Larg system}) and (\ref{app Larg compatibility}). And then, propagation of divergence free condition of magnetic field from initial data will be justified. 
\subsubsection{Fully nonlinear system}
We want to solve
\begin{equation} \label{app fully nonlin system}
\begin{cases}
v_t - \nu\triangle_v v + \nabla_v q = B\cdot\nabla_v B,\quad\text{in}\quad Q_T, \\
B_t - \lambda\triangle_v B = B\cdot\nabla_v v,\quad\text{in}\quad Q_T, \\
\nabla_v\cdot v = 0,\quad\text{in}\quad Q_T, \\
v(0) = v_0,\quad\Omega\times\{t=0\}, \\
B = 0 ,\quad \{S_{F}\cup \{\mathbb{R}^{3}\backslash \O \}\} \times [0,T),  \\
B(0) = B_0,\quad\Omega\times\{t=0\}, \\
q - 2\nu\mathbf{D}_v(v)\mathbf{n}^{(v)}\cdot\mathbf{n}^{(v)} = gh,\quad\text{on}\quad S_{F,T}, \\
2\nu\mathbf{D}(v)\mathbf{n}^{(v)} - 2\nu(\mathbf{D}_v(v)\mathbf{n}^{(v)}\cdot\mathbf{n}^{(v)})\mathbf{n}^{(v)} = 0,\quad\text{on}\quad S_{F,T}. \\
\end{cases}
\end{equation}

Note that this system does not contain divergence free condition of $B$, since (\ref{app heat system}) does not include any condition about divergence free of $B$. We first state a lemma which resembles Lemma (\ref{app nonlinear est}). 

\begin{lemma} \label{fully nonlin est}
	For $l>\frac{1}{2}$ we have the following nonlinear estimate.
	$$
	\|FG\|^{(l)}_{Q_T} \leq C(T) \|F\|_{H^{l+1,(l+1)/2}(Q_T)} \|G\|_{H^{l+1,(l+1)/2}(Q_T)},
	$$
	where $C(T)$ depends on $T$ non-decreasingly.
\end{lemma}

We alrady have a unique solution $(\mathbf{v},\mathbf{B},\mathbf{q})$ for $0 \leq t < T_1$ to system (\ref{app cst coeff nonlin}) with $(f,\rho,u_0,d,b) = (0,0,0,0,-gh)$. We find a solution of the form $ (v,B,q)=(\mathbf{v}+v^*,\mathbf{B}+B^*,\mathbf{q}+q^*)$. Then system (\ref{app fully nonlin system}) becomes,
\begin{equation} \label{app perturb}
\begin{cases}
v^*_t - \nu\triangle_v v^* + \nabla_v q^* = \nu(\triangle_v - \triangle)\mathbf{v} - (\nabla_v - \nabla)\mathbf{q} + \mathbf{B}\cdot(\nabla_v - \nabla)\mathbf{B}, \\
\,\,\,\,\,\,\,\,\,\,\,\,\,\,\,\,\,\,\,\,\,\,\,\,\,\,\,\,\,\,\,\,\,\,\,\,\,\,\,\,\,\,\,\,\,\,\,\,\,\,\,\,\,\,\,\,\,\,\,\, + \mathbf{B}\cdot\nabla_v B^* + B^*\cdot\nabla_v \mathbf{B} + B^*\cdot\nabla_v B^*,\quad\text{in}\quad Q_T, \\ 
B^*_t - \lambda\triangle_v B^* = \lambda(\triangle_v - \triangle)\mathbf{B} + \mathbf{B}\cdot(\nabla_v - \nabla)\mathbf{v}, \\
\,\,\,\,\,\,\,\,\,\,\,\,\,\,\,\,\,\,\,\,\,\,\,\,\,\,\,\,\,\,\,\,\,\,\,\,\,\,\,\,\,+ \mathbf{B}\cdot\nabla_v v^* + B^*\cdot\nabla_v \mathbf{v} + B^*\cdot\nabla_v v^*,\quad\text{in}\quad Q_T, \\ 
\nabla_v\cdot v^* = -\nabla_v\cdot \mathbf{v},\quad\text{in}\quad Q_T, \\
v^*(0) = v_0,\quad B^*(0) = B_0,\quad\text{in}\quad \Omega\times\{t=0\}, \\
-q^* + 2\nu\mathbf{D}_v(v^*)\mathbf{n}^{(v)}\cdot\mathbf{n}^{(v)} = -2\nu\mathbf{D}_v(\mathbf{v})\mathbf{n}^{(v)}\cdot\mathbf{n}^{(v)} + 2\nu\mathbf{D}_v(\mathbf{v})n_0\cdot n_0,\quad\text{on}\quad S_{F,T}, \\
2\nu[\mathbf{D}(v^*)\mathbf{n}^{(v)} - (\mathbf{D}_v(v^*)\mathbf{n}^{(v)}\cdot\mathbf{n}^{(v)})\mathbf{n}^{(v)}] = -2\nu[\mathbf{D}(\mathbf{v})\mathbf{n}^{(v)} - (\mathbf{D}_v(\mathbf{v})\mathbf{n}^{(v)}\cdot\mathbf{n}^{(v)})\mathbf{n}^{(v)}],\quad\text{on}\quad S_{F,T}. \\
\end{cases}
\end{equation}

We make the following iteration scheme to solve above (\ref{app perturb}). Note that $\nabla\cdot \mathbf{v} = 0$.

\begin{equation} \label{app perturb iteration}
\begin{cases}
v^{*(m+1)}_t - \nu\triangle v^{*(m+1)} + \nabla q^{*(m+1)} = \nu(\triangle_m - \triangle)v^{(m)} - (\nabla - \nabla_m)q^{(m)} + \mathbf{B}\cdot(\nabla_m - \nabla)\mathbf{B}, \\
\quad\quad\quad\quad\quad\quad\quad\quad\quad\quad\quad\quad\quad\quad\quad + \mathbf{B}\cdot\nabla_m B^{*(m)} + B^{*(m)}\cdot\nabla_m \mathbf{B} + B^{*(m)}\cdot\nabla_m B^{*(m)} \\
\quad\quad\quad\quad\quad\quad\quad\quad\quad\quad\quad\quad\quad\quad\quad := f^{(m)},\quad\text{in}\quad Q_T, \\ 
B^{*(m+1)}_t - \lambda\triangle B^{*(m+1)} = \lambda(\triangle_m - \triangle)B^{(m)} + \mathbf{B}\cdot(\nabla_m - \nabla)\mathbf{v}, \\
\quad\quad\quad\quad\quad\quad\quad\quad\quad\quad + \mathbf{B}\cdot\nabla_m v^{*(m)} + B^{*(m)}\cdot\nabla_m \mathbf{v} + B^{*(m)}\cdot\nabla_m v^{*(m)} \\
\quad\quad\quad\quad\quad\quad\quad\quad\quad\quad := g^{(m)},\quad\text{in}\quad Q_T, \\ 
\nabla \cdot v^{*(m+1)} = (\nabla-\nabla_m)\cdot v^{(m)} := \rho^{(m)},\quad\text{in}\quad Q_T, \\
v^{*(m+1)}(0) = u_0,\quad B^{*(m+1)}(0) = H_0,\quad\Omega\times\{t=0\}, \\
-q^{*(m+1)} + 2\nu\mathbf{D}(v^{*(m+1)})\mathbf{n}\cdot \mathbf{n} = 2\nu [\mathbf{D}(v^{(m)})\mathbf{n}\cdot \mathbf{n} - \mathbf{D}_m(v^{(m)})\mathbf{n}^{(m)}\cdot \mathbf{n}^{(m)}]   \\
\quad\quad\quad\quad\quad\quad\quad\quad\quad\quad\quad\quad\quad\quad := b^{(m)},\quad\text{on}\quad S_{F,T}, \\
2\nu[\mathbf{D}(v^{*(m+1)})\mathbf{n} - (\mathbf{D}_v(v^{*(m+1)})\mathbf{n}\cdot \mathbf{n})\mathbf{n}] = -2\nu \{ [\mathbf{D}_m(v^{(m)})\mathbf{n}^{(m)} - (\mathbf{D}_m(v^{(m)})\mathbf{n}^{(m)}\cdot \mathbf{n}^{(m)})\mathbf{n}^{(m)}], \\
\quad\quad\quad\quad\quad\quad\quad\quad\quad\quad\quad\quad\quad\quad\quad\quad\quad - [\mathbf{D}(v^{*(m)})\mathbf{n} - (\mathbf{D}(v^{*(m)})\mathbf{n} \cdot \mathbf{n})\mathbf{n}] \} \\
\quad\quad\quad\quad\quad\quad\quad\quad\quad\quad\quad\quad\quad\quad\quad\quad\quad := d^{(m)},\quad\text{on}\quad S_{F,T} ,\\
(v^{*(0)},B^{*(0)},q^{*(0)}) = (0,0,0),
\end{cases}
\end{equation}
where
\begin{equation}
\begin{split}
(v^{(m)},B^{(m)},q^{(m)}) &:= (\mathbf{v} + v^{*(m)},\mathbf{B} + B^{*(m)},\mathbf{q} + q^{*(m)}),  \\
\nabla_m &:= \nabla_{v^{(m)}},\quad \mathbf{D}_m := \mathbf{D}_{v^{(m)}},\quad n^{(m)} := n^{v^{(m)}} \\
\rho^{(m)} &:= \nabla\cdot R^{(m)},\quad R^{(m)} := (\mathbf{I} - \mathcal{G}^{(m)})v^{(m)},\quad \mathcal{G}^{(m)} := \mathcal{G}^{v^{(m)}} ,
\end{split}
\end{equation}
and 
\begin{equation} \label{app delta 12}
T^{1/2} \|v^{(m)}\|^{(l+2)}_{Q_{T}} \leq \delta_1,\,\,\,\,T^{1/2} \|B^{(m)}\|^{(l+2)}_{Q_{T}} \leq \delta_2.
\end{equation}
Using Proposition \ref{app prop general stokes} and \ref{app prop general heat},
\begin{equation}
\begin{split}
&\|v^{*(m+1)}\|_{H^{l+2,l/2+1}(Q_T)} + \|B^{*(m+1)}\|_{H^{l+2,l/2+1}(Q_T)}  \\
&\quad + \|q^{*(m+1)}\|^{(l)}_{Q_T} + \|\nabla q^{*(m+1)}\|^{(l)}_{Q_T} + \|q^{*(m+1)}\|_{W_2^{l+1/2,l/2+1/4}(S_{F,T})}  \\
&\leq C_*(T) \{ \|f^{(m)}\|^{(l)}_{Q_T} + \|g^{(m)}\|^{(l)}_{Q_T} + \|\rho^{(m)}\|_{W_2^{l+1,(l+1)/2}(S_{F,T})} + \|R^{(m)}\|_{W_2^{0,l/2+1}(Q_T)}  \\
&\quad + T^{-l/2}\|R_t^{(m)}\|_{L^2(Q_T)} + \|(b^{(m)},d^{(m)})\|_{W_2^{l+1/2,l/2+1/4}(S_{F,T})}  \\
&\quad + T^{-l/2}\|b^{(m)}\|_{W_2^{l/2,0}(S_{F,T})} + \|u_0\|_{W_2^{l+1}(\Omega)} + \|H_0\|_{W_2^{l+1}(\Omega)} \}.
\end{split}
\end{equation}

Similar as before, we claim uniform bound and contraction mapping property.  \\
$\mathbf{Uniform\,\,bound}$ 
Using Lemma \ref{app nonlinear est}, we can derive
\begin{equation} \label{app f_m}
\begin{split}
\|f^{(m)}\|^{(l)}_{Q_T} &\leq C(T,\delta_1,\delta_2) Z_m ( 1 + Z_m ),
\end{split}
\end{equation}
where
\begin{equation} \label{app Z_m}
\begin{split}
Z_m &:= \|v^{(m)}\|_{H^{l+2,l/2+1}(Q_T)} + \|B^{(m)}\|_{H^{l+2,l/2+1}(Q_T)} + \|q^{(m)}\|^{(l)}_{Q_T} \\
&\quad + \|\nabla q^{(m)}\|^{(l)}_{Q_T} + \|q^{(m)}\|_{W_2^{l+1/2,l/2+1/4}(S_{F,T})},
\end{split}
\end{equation}
and $C(T,\delta_1,\delta_2)$ is positive constant depending increasingly on both arguments with property that $ C(T,\delta_1,\delta_2)\rightarrow 0$ as $T\rightarrow 0$ and may vary line to line. Using exactly same argument, we have the same estimate for $g^{(m)}$, $\rho^{(m)}$, $ R^{(m)}$, $R_t^{(m)}$, $d^{(m)}$, and $b^{(m)}$, 
\begin{equation} \label{app others}
\begin{split}
\|g^{(m)}\|^{(l)}_{Q_T} &\leq C(T,\delta_1,\delta_2) Z_m ( 1 + Z_m ),  \\
\|\rho^{(m)}\|_{W_2^{l+1,(l+1)/2}(S_{F,T})} &\leq C(T,\delta_1,\delta_2) Z_m ( 1 + Z_m ),  \\
\|R^{(m)}\|_{W_2^{0,l/2+1}(Q_T)} &\leq C(T,\delta_1,\delta_2) Z_m ( 1 + Z_m ) \\
&\quad  + C(T,\delta_1,\delta_2) Z_m ( \|v^{(m)}\|_{W_2^{l+1,0}(Q_T)} + \|Dv^{(m)}\|_{W_2^{0,l/2}(Q_T)} ),  \\
T^{-l/2}\|R_t^{(m)}\|_{L^2(Q_T)} &\leq C(T,\delta_1,\delta_2)T^{(1-l)/2} Z_m ( 1 + Z_m ) \leq C(T,\delta_1,\delta_2) Z_m ( 1 + Z_m ),  \\
\|d^{(m)}\|_{W_2^{l+1/2,l/2+1/4}(S_{F,T})} &\leq C(T,\delta_1,\delta_2) Z_m ( 1 + Z_m ), \\
\|b^{(m)}\|_{W_2^{l+1/2,l/2+1/4}(S_{F,T})} + T^{-l/2}\|b^{(m)}\|_{W_2^{l/2,0}(S_{F,T})} &\leq C(T,\delta_1,\delta_2) Z_m ( 1 + Z_m ).
\end{split}
\end{equation}
Similar as (\ref{app Z_m}), we define,
\begin{equation} \label{app Z_m *}
\begin{split}
Z_m^* &:= \|v^{*(m)}\|_{H^{l+2,l/2+1}(Q_T)} + \|B^{*(m)}\|_{H^{l+2,l/2+1}(Q_T)} + \|q^{*(m)}\|^{(l)}_{Q_T} \\
&+ \|\nabla q^{*(m)}\|^{(l)}_{Q_T} + \|q^{*(m)}\|_{W_2^{l+1/2,l/2+1/4}(S_{F,T})}.
\end{split}
\end{equation}
Then, combining (\ref{app f_m}) and (\ref{app others}), 
$$
Z_{m+1}^* \leq h_0 + h_1 Z^*_m + h_2 Z_m^{*2},
$$
with positive constant $h_0,h_1,h_2$ with following properties.\\
1) $h_0 = h_0(T,\delta_1,\delta_2)$ is monotone increasing function with all its argument.\\
2) $h_{1,2} = h_{1,2}(T,\delta_1,\delta_2) \rightarrow 0$ as $T\rightarrow 0$.\\
From cauchy sequence argument, there exist $z^*$ such that if $Z_m^* \leq z^*$, then
$$
Z_{m+1}^* \leq h_0 + h_1 z^*_m + h_2 (z_m)^{*2} \leq z^*.
$$
Hence we have uniform bound,
\begin{equation} \label{app uniform bound 2}
Z_m^* \leq z^*,\,\,\,\,\forall m
\end{equation}
$\mathbf{Contraction\,\,mapping}$ We use the similar notation in section 4, to denote difference. 
$$
v^{*(m+1)} - v^{*(m)} \doteq \mathcal{V}^{*(m+1)},\,\,\,\,B^{*(m+1)} - B^{*(m)} \doteq \mathcal{B}^{*(m+1)},\,\,\,\,q^{*(m+1)} - q^{*(m)} \doteq \mathcal{Q}^{*(m+1)}.
$$
We make equation of $\mathcal{V}^{*(m)}$, $\mathcal{B}^{*(m)}$, and $\mathcal{Q}^{*(m)}$, and use Proposition \ref{app prop general stokes} and \ref{app prop general heat} with $v_0 = B_0 = 0$, to get
\begin{equation} \label{app Y ineq}
\begin{split}
Y_{m+1}^* &:= \|\mathcal{V}^{*(m+1)}\|_{H^{l+2,l/2+1}(Q_T)} + \|\mathcal{B}^{*(m+1)}\|_{H^{l+2,l/2+1}(Q_T)}  \\
&\quad + \|\mathcal{Q}^{*(m+1)}\|^{(l)}_{Q_T} + \|\nabla \mathcal{Q}^{*(m+1)}\|^{(l)}_{Q_T} + \|\mathcal{Q}^{*(m+1)}\|_{W_2^{l+1/2,l/2+1/4}(S_{F,T})}  \\
&\leq C_*(T) \{ \|f^{*(m)}\|^{(l)}_{Q_T} + \|g^{*(m)}\|^{(l)}_{Q_T} + \|\rho^{*(m)}\|_{W_2^{l+1,(l+1)/2}(S_{F,T})} + \|R^{*(m)}\|_{W_2^{0,l/2+1}(Q_T)} \\
&\quad + T^{-l/2}\|R_t^{*(m)}\|_{L^2(Q_T)} + \|(b^{*(m)},d^{*(m)})\|_{W_2^{l+1/2,l/2+1/4}(S_{F,T})} + T^{-l/2}\|b^{*(m)}\|_{W_2^{l/2,0}(S_{F,T})} \},
\end{split}
\end{equation}
where $C_{*}(T)$ is non-decreasingly time dependent constant, which means it does not blow up as $T\rightarrow 0$. SImilar as before, we can control every terms in right hand side, to get 
\begin{equation}
Y_{m+1}^* \leq \chi Y_m^*,
\end{equation}
where $\chi < 1$ if we pick a $\varepsilon$ and sufficiently small $T$ which is smaller than $T_1$ of (\ref{app delta 12}). Hence, by contraction mapping principle, we solve (\ref{app fully nonlin system}). So far, we proved Theorem \ref{app main theorem}, except $\nabla\cdot H = 0$. 

\subsubsection{Divergence free of $H$}
To show propagation of divergence free property, we appeal to maximum principle of convection-diffusion equation.
\begin{equation*}
\begin{cases}
H_t + (u\cdot\nabla)H - (H\cdot\nabla)u = \lambda\triangle H, \\
\nabla\cdot u = 0 ,  \\
\nabla\cdot H_0 = 0.
\end{cases}
\end{equation*}
Taking divergence to above equation and using notation  $\mathcal{H} := \nabla\cdot H $,
$$
\mathcal{H}_t + (u\cdot\nabla)\mathcal{H} + (\nabla u):(\nabla H)^t - (\nabla H):(\nabla u)^t - (H\cdot\nabla)(\nabla\cdot u) = \lambda\triangle \mathcal{H},
$$
where $A:B := \sum_{i,j} A_{ij} B_{ij}$. Hence,
$$
\mathcal{H}_t + (u\cdot\nabla)\mathcal{H} - \lambda\triangle \mathcal{H} = 0.
$$
Then by maximum principle of convection-diffusion equation, 
$$
\|(\nabla\cdot H) (t)\|_{L^\infty} = \|\mathcal{H}(t)\|_{L^\infty} \leq \|\mathcal{H}(0)\|_{L^\infty} = \|\nabla\cdot H_0\|_{L^\infty} = 0,
$$
during the time interval for the solution $H$.

\section*{Acknowledgments}
The author would like to thank Professor N. Masmoudi for helpful discussion and suggestion. The author  was  partially supported by NSF grant DMS-1211806. This research was performed when the author was at Courant Institute of Mathematical Science.

\vspace{1cm}
\indent\indent \author{$\textsc{Donghyun\,\,\,Lee}$}\\
\indent \address{$\textsc{University\,\,\,of\,\,\,Wisconsin-Madison,\,\,\,Van Vleck Hall, 480 Lincoln Drive, Madison, Wi  53706}$} \\
\indent \email{$\mathtt{dlee374@wisc.edu}$}

\end{document}